\documentclass[12pt]{article}
\usepackage{amssymb,amsmath,amscd}

\input epsf.tex

\title{Stability structures, motivic Donaldson-Thomas invariants and cluster transformations}

\author {Maxim Kontsevich, Yan Soibelman}

\begin{document}
\maketitle

\newcommand{\op}[1]{\operatorname{#1}}

\newcommand{\CC}{{\mathcal C}}

\renewcommand{\O}{{\mathcal O}}

\newcommand{\E}{{\mathcal E}}
\newcommand{\F}{{\mathcal F}}

\newcommand{\g}{{\mathfrak g}}
\newcommand{\h}{{\mathfrak h}}

\renewcommand{\k}{{\bf k}}
\newcommand{\kk}{{\overline{\bf k}}}

\newtheorem{defn}{Definition}
\newtheorem{thm}{Theorem}
\newtheorem{lmm}{Lemma}
\newtheorem{rmk}{Remark}
\newtheorem{prp}{Proposition}
\newtheorem{conj}{Conjecture}
\newtheorem{exa}{Example}
\newtheorem{cor}{Corollary}
\newtheorem{que}{Question}
\newtheorem{ack}{Acknowledgments}
\newcommand{\C}{{\bf C}}
\newcommand{\K}{{\bf k}}
\newcommand{\R}{{\bf R}}
\newcommand{\N}{{\bf N}}
\newcommand{\Z}{{\bf Z}}
\newcommand{\Q}{{\bf Q}}
\newcommand{\G}{\Gamma}
\newcommand{\A}{A_{\infty}}

\newcommand{\epi}{\twoheadrightarrow}
\newcommand{\mono}{\hookrightarrow}
\newcommand\ra{\rightarrow}
\newcommand\uhom{{\underline{\op{Hom}}}}
\renewcommand\O{{\cal O}}
\newcommand{\epp}{\varepsilon}

\tableofcontents

\section{Introduction}

\subsection{Counting problems for $3$-dimensional Calabi-Yau varieties}

Let $X$ be a compact complex $3$-dimensional K\"ahler manifold such that $c_1(T_X)=0\in Pic(X)$
(hence by Yau theorem $X$ admits a Calabi-Yau metric).
We can associate with $X$ several moduli spaces which have the virtual dimension zero:

a) moduli of holomorphic curves in $X$ with  fixed genus and degree;

b) moduli of holomorphic vector bundles on $X$ (or, more generally, of coherent sheaves) with a fixed Chern character;

c) moduli of special Lagrangian submanifolds\footnote{Recall that a Lagrangian submanifold $L\subset X$ is called {\it special} iff the restriction to $L$ of a holomorphic volume form on $X$ is a real volume form on $L$.} with a fixed homology class endowed with a $U(1)$ local system.

In order to have a well-defined virtual number of points of the moduli space one needs compactness and a perfect obstruction theory with virtual dimension zero (see \cite{BF}, \cite{T1}, \cite{T2}).\footnote{ The latter means that the deformation theory of a point is controlled by a differential-graded Lie algebra $\g$ such that $H^i(\g)=0$ for $i\ne 1,2$ and
$\dim H^1(\g)=\dim H^2(\g)$.}
The compactification is known in the case a). It is given by the moduli of stable maps. The corresponding virtual numbers are Gromov-Witten invariants (GW-invariants for short). Donaldson and Thomas in \cite{DoT},\cite{T2} addressed the cases b) and c). Analytical difficulties there are not completely resolved. The most understood example is the one
of torsion-free sheaves of rank one with the fixed Chern character of the form $(1,0,a,b)\in H^{ev}(X)$.
The corresponding virtual numbers are called Donaldson-Thomas invariants (DT-invariants for short).
One sees that the number of (discrete) parameters describing GW-invariants is equal to $1+\dim H^2(X)$ (genus and degree)
and coincides with the number of parameters describing DT-invariants. The conjecture from \cite{MaNOP} (proved in many cases) says that GW-invariants and DT-invariants  can be expressed one through another. The full putative virtual numbers in the case b) should depend on as twice as many parameters (i.e. $\dim  H^{ev}(X)$). By mirror symmetry one reduces the case c) to the case b) for the dual Calabi-Yau manifold. Unlike to GW-invariants and DT-invariants these virtual numbers should depend on some choices (the K\"ahler structure in the case b) and the complex structure in the case c), see \cite{T2}). In particular, in the case c), for a compact $3d$ Calabi-Yau manifold $X$ we should have an even function $$\Omega_{SLAG}: H^3(X,\Z)\setminus \{0\}\to \Q\,\,,$$ which depends  on the complex structure on $X$ in such a way that for any non-zero $\gamma \in H^3(X,\Z)$ the number $\Omega_{SLAG}(\gamma)$ is a constructible function
with respect to a real analytic stratification
of the moduli space of complex structures. Moreover this number
is integer for a generic complex structure. The invariant $\Omega_{SLAG}(\gamma)$ is the virtual number of special Lagrangian submanifolds $L\subset X$ in the  class $\gamma$ (or more generally, special Lagrangian submanifolds endowed with  
 local systems of arbitrary rank).

Our aim in this paper is to describe a framework for ``generalized Donaldson-Thomas invariants" and their wall-crossing formulas in the case
of non-commutative compact $3d$ Calabi-Yau varieties.
A choice of polarization
(``complexified K\"ahler structure")
will be encoded into a choice of ``stability condition" on ${\cal C}$. Then we define a generalized Donaldson-Thomas invariant $\Omega(\gamma)$  as  the
``number" of stable  objects in ${\cal C}$ with a fixed class $\gamma$ in the $K$-group.
Similar problem for abelian categories was addressed in the series of papers by Joyce  \cite{Jo1}\cite{Jo2}\cite{Jo3}
 and in the recent paper of Bridgeland and Toledano Laredo \cite{BrTL}. Our paper can be thought of as a generalization to the case of triangulated categories (the necessity of such a generalization is motivated by both mathematical and physical applications,
see e.g. \cite{T1}, \cite{Dou}). One of motivations for our counting formula was the microlocal formula by K.~Behrend (see
 \cite{B}) for the virtual number in the case of so called symmetric obstruction theory (see \cite{BF}), which is the
 case for objects in $3d$ Calabi-Yau categories. 
The above example b) corresponds to the bounded derived category $D^b(X)$ of coherent sheaves on $X$ (more precisely to its $\A$-enrichment). The example c) corresponds to the Fukaya category. In that case the wall-crossing formulas describe the behavior of
$\Omega_{SLAG}$. Even in the geometric situation our formalism extends beyond the case of smooth compact Calabi-Yau varieties.

\subsection{Non-commutative  varieties with polarization}

All $\A$-categories in this paper will be ind-constructible. This roughly means that their spaces of objects are countable inductive limits of constructible sets (for more details see  Section 3).
We define a {\it non-commutative proper algebraic variety over a base field $\k$} as an $Ext$-finite ind-constructible $\k$-linear triangulated $\A$-category $\CC$. For two objects
$E$ and $F$  we denote by $\op{Hom}^{\bullet}(E,F)$ the complex of morphisms and by
$\op{Ext}^{\bullet}(E,F)$ its cohomology.

Here are few examples of such categories.

\begin{exa}
a)  $\A$-version of $D^b(X)$, the bounded derived category of the category of coherent sheaves
on a smooth projective algebraic variety $X/\k$. In this case $D^b(X)$ coincides with the
triangulated category $Perf(X)$ of perfect complexes on $X$.

b) More generally, for a (not necessarily proper) smooth variety $X$
endowed with a closed proper  subset $X_0\subset X$,  the corresponding triangulated category is the full subcategory of $Perf(X)$ consisting
 of complexes of sheaves with cohomology supported on $X_0$.

c) Also for a (not necessarily proper) smooth variety $X$ we can consider the
 the full subcategory of $Perf(X)$ consisting
 of complexes of sheaves with compactly supported cohomology.
 
d) The $\A$-version of the category $Perf(X)$ of perfect complexes on a proper, not necessarily smooth scheme $X$ over $\k$.

e) If $A$ is an $\A$-algebra with finite-dimensional cohomology then ${\cal C}=Perf(A)$ is the category
of perfect $A$-modules.

f) If $\k$ is the field of characteristic zero and $A$ is finitely generated in the sense of \cite{ToVa}, (in particular it is homologically smooth, see \cite{KoSo3}) then
${\cal C}$ is the category of $A$-modules of finite dimension over $\k$.

g) If the category ${\cal C}$ is ind-constructible and $E\in Ob({\cal C})$ then left and right orthogonal
to the minimal triangulated subcategory generated by $E$ are also ind-constructible (since the conditions $\op{Ext}^{\bullet}(X,E)=0$
and $\op{Ext}^{\bullet}(E,X)=0$ are ``constructible").
\end{exa}

Let us make few comments on the list. Example a) is a particular case
of examples b),c),d). Using the results of \cite{BVdB}
we can reduce geometric examples b),d) to the algebraic example e),
 and also the example c) to the example f).
Let us discuss a typical (and most important) example e) at the
level of objects of the category. We claim that the set of isomorphism classes of  objects of $\CC$ 
 can be covered by an
inductive limit of constructible sets. First, replacing $A$ by its
minimal model we may assume that $A$ is finite-dimensional. Basic
examples of perfect $A$-modules are direct sums of shifts of $A$,
i.e. modules of the type $$M=A[n_1]\oplus A[n_2]\oplus\dots\oplus
A[n_r],\,\,r\geqslant 0 $$
 and their ``upper-triangular deformations" (a.k.a. twisted
complexes). The latter are described by solutions to the Maurer-Cartan
equations $$\sum_{1\le l\le r-1}m_l(\alpha,\dots,\alpha)=0\,\,$$
 where
$\alpha=(a_{ij})_{i<j}$ is an upper-triangular $r\times r$ matrix
with coefficients in $A$ and $\op{deg}a_{ij}=n_i-n_j+1$.

This gives a closed scheme. For a point $x$ of this scheme we have the corresponding $\A$-module $M_x$ over $A$.  In order to describe all other objects of $Perf(A)$ we need to take direct summands (up to homotopy) of such modules $M_x$. The set of such summands is parametrized by all
$\A$-morphisms of the non-unital algebra $P=\k\cdot p$ with the product given by $p^2=p$  (and with the trivial differential) to $\op{End}^\bullet(M_x)$.
Every such morphism is described by a collection of linear maps
$f_n: P^{\otimes n}\to \op{End}^{1-n}(M_x)$ satisfying a system of polynomial equations. Notice that if $n$ is large then $f_n=0$ because the corresponding negative graded components of $\op{End}^\bullet (M_x)$ are trivial.
Therefore we again have a finite system of polynomial equations.

For given $N=\sum_{1\le j\le r}(|n_j|+1)$ we obtain a scheme of finite type $Mod_N$ parametrizing some objects of $Perf(A)$. Clearly $Ob(Perf(A))=\cup_{N\geqslant 1}Mod_N$. Each isomorphism class of an object appears in the union for infinitely many values of $N$. In order to avoid the ``overcounting" we define a subscheme of finite type $Mod_N^0\subset Mod_N$ consisting of objects not isomorphic to objects from $Mod_{N^{\prime}}$ for $N^{\prime}<N$.
We conclude that objects of $Perf(A)$ form an ind-constructible set (more precisely, an ind-constructible stack). 
One can take care about morphisms in the category in a similar way. This explains the example e).

We define a {\it polarization} on a non-commutative proper algebraic variety over $\k$ (a version of Bridgeland stability condition, see \cite{Br1})
by the following data and axioms:

\begin{itemize}

\item{an ind-constructible homomorphism $\op{cl}: K_0(\CC)\to \Gamma$, where $\Gamma\simeq\Z^n$ is a free abelian group of finite rank endowed with a bilinear form\footnote{ In physics literature $\Gamma$ is called the charge lattice.} $\langle \bullet, \bullet\rangle:\Gamma\otimes \Gamma\to \Z$ such that for any two objects $E,F\in Ob(\CC)$ we have $$\langle \op{cl}(E), \op{cl}(F)\rangle=\chi(E,F):=\sum_i(-1)^i\dim \op{Ext}^i(E,F)\,\,,$$}
\item{an additive map $Z:\Gamma \to \C$, called the central charge,}
\item{a collection $\CC^{ss}$ of (isomorphism classes of) non-zero objects in $\CC$ called the
semistable ones, such that $Z(E)\ne 0$ for any $E\in
\CC^{ss}$, where we write $Z(E)$ for $Z(\op{cl}(E))$,}
\item{a choice $\op{Log} Z(E) \in \C$ of the logarithm of $Z(E)$ defined for any $E\in \CC^{ss}$.}
\end{itemize}

Making a connection with \cite{Br1} we say that the last three items define a {\it stability structure (or stability condition)} on the category ${\cal C}$.

For $E\in \CC^{ss}$ we denote by $\op{Arg}(E)\in \R$ the imaginary part of $\op{Log} Z(E)$.

 The above data satisfy the following axioms:
 \begin{itemize}
 \item{for all $E\in  \CC^{ss}$ and for all $ n\in \Z$ we have $E[n]\in \CC^{ss} $ and
$$\op{Arg} Z(E[n])=\op{Arg} Z (E)+\pi  n\,\,, $$}
 \item{ for all $E_1,E_2\in \CC^{ss}$ with $\op{Arg}(E_1)> \op{Arg} (E_2)$ we have $${\op{Ext}_{\cal
C}}^{\le 0}\,(E_1,E_2)=0\,\,,$$}
 \item{ for any object $E\in Ob(\CC)$ there exist $n\geqslant 0$ and a chain of morphisms
  $0=E_0\to E_1 \to \dots \to E_n=E$ (an analog of  filtration)
  such that the corresponding ``quotients'' $F_i:=Cone( E_{i-1}\to E_{i}),\,\,i=1,\dots,n$ are semistable and
   $\op{Arg}(F_1)>\op{Arg}(F_2)> \dots > \op{Arg}(F_{n})$,}
   \item{for each $\gamma \in \Gamma\setminus \{0\}$ the set of isomorphism classes of a
$\CC^{ss}_\gamma\subset Ob({\cal C})_{\gamma}$ consisting of semistable objects $E$ such that
   $\op{cl}(E)=\gamma$ and $\op{Arg}(E)$ is fixed, is a constructible set,}
   \item{(Support Property) Pick a norm $\parallel\cdot \parallel $ on $\Gamma \otimes
\R$, then there exists $C>0 $ such that
    for all $E\in \CC^{ss}$ one has
   $ \parallel E \parallel\le C|Z(E)|$.}
 \end{itemize}

In the above definition one can allow $\Gamma$ to have a torsion.
 In geometric examples
a), d) for $\k=\C$ one can take $\Gamma=K^0_{top}(X(\C))$
 where $K^0_{top}$ denotes the topological $K^0$-group. Similarly, in examples b),c)
  one should take the $K^0$-groups with appropriate supports.
 Another choice for $\Gamma$ is the image of the algebraic $K^0$-group under the Chern
character. Yet another
 choice is $\Gamma=K_0^{num}(\CC)$, which is the quotient of the group $K_0(\CC)$ by the
intersection of the left and right kernels of the Euler form $\chi(E,F)$. Finally one can pick a finite collection of ind-constructible functors $\Phi_i: \CC\to Perf(\k), 1\le i\le n$ and define
$$\op{cl}(E)=(\chi(\Phi_1(E)),\dots,\chi(\Phi_n(E)))\in \Z^n=:\Gamma\,\,,$$ where $\chi:K_0(Perf(\k))\to \Z$ is the isomorphism of groups given by the Euler characteristic.

\begin{rmk} The origin of the Support Property is geometric and can be explained in the case of the category of $A$-branes (the derived Fukaya category $D^b({\cal F}(X))$)
of a compact $3$-dimensional Calabi-Yau manifold $X$. Let us fix a Calabi-Yau metric $g_0$ on $X$. Asymptotically, in the large volume limit (as the rescaled symplectic form approaches infinity) it gives rise
to the stability condition on $D^b({\cal F}(X))$ such that
stable objects are special Lagrangian submanifolds,
and $|Z(L)|$ is the volume of $L$ with respect to $g_0$. Then for any harmonic form $\eta$ one has
$|\int_{L}\eta|\le C|Z(L)|$. It follows that the norm of the cohomology class
of $L$ is bounded (up to a scalar factor) by the norm of the linear functional $Z$.
\end{rmk}

The Support Property implies that the set $\{Z(E)\in \C\,| \,E\in \CC^{ss}\}$ is a discrete subset of $\C$ with at most polynomially
growing density at infinity. It also implies that the stability condition is locally finite in the sense of Bridgeland (see \cite{Br1}).
Any stability condition gives a bounded $t$-structure on $\CC$ with the corresponding heart consisting of
  semistable objects $E$ with $\op{Arg}(E) \in (0,\pi]$ and their extensions.

  \begin{rmk} The case of the classical Mumford notion of stability with respect to an ample
line bundle (and its refinement for coherent sheaves defined by Simpson) is not an example of the Bridgeland stability condition, it is rather
  a limiting degenerate case of it (see \cite{Ba}, \cite{Tod1} and Remark at the end of Section 2.1).

\end{rmk}

  For given $\CC$ and a homomorphism  $\op{cl}: K_0(\CC)\to \Gamma$ as above, let us denote by $Stab(\CC):=Stab(\CC,\op{cl})$ the set of
stability conditions $(Z,\CC^{ss}, (\op{Log} Z(E))_{E\in \CC^{ss}})$. Space $Stab(\CC)$
 can be endowed with a  Hausdorff topology, which we discuss in detail in Section 3.4.
Then we have an ind-constructible version of
the following fundamental result of Bridgeland (see \cite{Br1}).

 \begin{thm}  The forgetting map $Stab(\CC)\to \C^n\simeq { \op{Hom}}(\Gamma,\C)$ given by
  $( Z,\CC^{ss}, (\op{Log} Z(E))_{E\in \CC^{ss}})\mapsto Z$,
 is a local homeomorphism.
 \end{thm}

Hence, $Stab(\CC)$ is a complex manifold, not necessarily connected.
Under appropriate assumptions one can show also that the group of autoequivalences $\op{Aut}(\CC)$ acts properly and discontinuously on $Stab(\CC)$.
On the quotient orbifold $Stab(\CC)/ \op{Aut} (\CC)$ there is a natural non-holomorphic action of
the group $GL^+(2,\R)$ of orientation-preserving $\R$-linear automorphisms of $\R^2\simeq \C$.

\subsection{Donaldson-Thomas invariants for non-commutative $3d$ Calabi-Yau varieties}

Recall that a non-commutative Calabi-Yau variety of dimension $d$ (a.k.a Calabi-Yau category
of dimension $d$) is given
by an $Ext$-finite triangulated $\A$-category ${\cal C}$ which carries a functorial
non-degenerate pairing $$(\bullet, \bullet): \op{Hom}^\bullet_{\CC}(E,F)\otimes \op{Hom}^\bullet_{\CC}(F,E)\to \k[-d]$$
 (see e.g. \cite{KoSo3}, \cite{So1}, \cite{Laz}),
such that the polylinear forms $(m_n(f_0,\dots,f_n),f_{n+1})$ defined on
$\otimes_{0\le i\le n+1}\op{Hom}^\bullet_{\CC}(E_i,E_{i+1})$ by higher compositions $m_n$ are cyclically invariant.
We will discuss mainly the case $d=3$ and assume that our non-commutative $3d$ Calabi-Yau variety is ind-constructible and endowed with polarization.

Under these assumptions we define  {\it motivic Donaldson-Thomas invariants} which take values in certain Grothendieck groups of algebraic varieties (more details are given in Sections 4 and 6). Assuming some ``absence of poles" conjectures, which we discuss in detail in Section 7 one can pass to the ``quasi-classical limit" which corresponds to the taking of Euler characteristic of all relevant motives. In this way we obtain the putative {\it numerical DT-invariants} $\Omega(\gamma)\in \Q, \gamma\in \Gamma\setminus \{0\}$.
Morally, $\Omega(\gamma)$ counts semistable objects of $\CC$ with a given class $\gamma\in \Gamma\setminus \{0\}$.

There is a special case when our formulas can be compared with those from \cite{B} (see Section 7.1). Namely, let us define a {\it Schur object} $E\in Ob(\CC)$ as such that $$\op{Ext}^{<0}(E,E)=0, \op{Ext}^0(E,E)=\k\cdot Id_E\,\,.$$
 By the Calabi-Yau property in the dimension $d=3$ we know that the only possibly non-trivial groups 
$\op{Ext}^i(E,E), i=0,1,2,3$ are $$\op{Ext}^0(E,E)\simeq \op{Ext}^3(E,E)\simeq \k\,,\,\,\,\op{Ext}^1(E,E)\simeq (\op{Ext}^2(E,E))^{\ast}\,\,.$$
 In other words the ranks are $(1,a,a,1),\, a\in \Z_{\geqslant 0}$. Recall (see \cite{KoSo2}, \cite{KoSo3}) that the deformation theory of any object $E\in Ob(\CC)$ is controlled by a differential-graded Lie algebra (DGLA for short) $\g_E$ such that $H^i(\g_E)\simeq \op{Ext}^i(E,E), i\in \Z$. For a given Schur object $E$ instead of $\g_E$ we can use a DGLA
$\widehat{\g}_E=\tau_{\le 2}(\g_E)/\tau_{\le 0}(\g_E)$ where $\tau_{\le i}$ is the truncation functor. This makes sense since $\tau_{\le 0}(\g_E)$ is an ideal (in the homotopy sense) in $\g_E$. The modified deformation theory gives rise to a perfect obstruction theory in the sense of \cite{B}, \cite{BF}. The corresponding moduli space is the same as the original one, although controlling DGLAs are not quasi-isomorphic. The contribution of Schur objects to $\Omega(\gamma)$ can be identified with the Behrend microlocal formula for DT-invariants.
From this point of view objects of the category $\CC$ should be interpreted as critical points of the function (called the potential), which is obtained from the solution to the so-called classical master equation. The latter has a very natural interpretation in terms of the non-commutative formal symplectic dg-scheme defined by the $\A$-category $\CC$ endowed with a Calabi-Yau structure (see \cite{KoSo3}).

\subsection{Multiplicative wall-crossing formula}

The wall-crossing formulas for the numerical Donaldson-Thomas invariants do not depend on their origin and can be expressed in terms of graded Lie algebras. This is explained in Section 2. Our main application is the case of  $3d$ Calabi-Yau categories. Let us recall that
if ${\cal C}$ is an $Ext$-finite Calabi-Yau category of the odd dimension $d$ (e.g. $d=3$)
then the Euler form

$$\chi:K_0(\CC)\otimes K_0(\CC)\to \Z,\,\,\,\,\chi(E,F):=\sum_{n\in \Z} (-1)^n \dim \, \op{Ext}^n(E,F)$$
is skew-symmetric.
In this case we also assume that if $\CC$ is endowed with polarization, then a skew-symmetric bilinear form
$\langle\bullet,\bullet\rangle: \Gamma\otimes \Gamma\to \Z$ is given and satisfies
$$\langle \op{cl}(E),\op{cl}(F)\rangle=\chi(E,F)\,\,\,\,\,\forall E,F\in Ob(\CC)\,.$$

In general, having a free abelian group $\Gamma$ of finite rank endowed with an integer-valued skew-symmetric form $\langle\bullet,\bullet\rangle$, we
define a Lie algebra over $\Q$ $\g_{\Gamma}:=\g_{(\Gamma,\langle\bullet,\bullet\rangle)}$, with the basis $(e_\gamma)_{\gamma\in \Gamma}$ and the Lie bracket

 $$ [ e_{\gamma_1}, e_{\gamma_2}]=(-1)^{\langle \gamma_1,\gamma_2\rangle}\langle
\gamma_1,\gamma_2\rangle\,e_{\gamma_1+\gamma_2}\,\,.$$

This Lie algebra is isomorphic (non-canonically) to the Lie algebra of regular functions
on the algebraic Poisson torus ${ \op{Hom}}(\Gamma, {{\bf G}}_m)$ endowed with the natural translation-invariant Poisson bracket.\footnote{Later we will use the multiplication  as well.
It is given explicitly by
$e_{\gamma_1}e_{\gamma_2}=(-1)^{\langle \gamma_1,\gamma_2\rangle}\,e_{\gamma_1+\gamma_2}.$}

An additive map $Z:\Gamma\to \C$ is called {\it generic}
if  there are no two $\Q$-independent elements of the lattice $\Gamma$ which are mapped by $Z$
into the same straight line in $\R^2=\C$. The set of non-generic maps is a countable union of real hypersurfaces in $\C^n=\op{Hom}(\Gamma, \C)$. These hypersurfaces are called {\it walls}.

Let us choose such an additive map $Z$ and an arbitrary norm $\parallel \bullet \parallel$ on the real vector space $\Gamma_{\R}=\Gamma\otimes \R$. We will keep the same notation for the $\R$-linear extension
of $Z$ to $\Gamma_{\R}$.  Finally, assume that we are given an even map $\Omega:\Gamma\setminus \{0\}\to \Z$
supported on the set  of $\gamma\in \Gamma$ such that $\parallel \gamma \parallel\le C|Z(\gamma) |$ for some given constant $C>0$.

Let $(Z_t)_{t\in [0,1]}$ be a generic piece-wise smooth path in $\C^n=\op{Hom}(\Gamma, \C)$ such that $Z_0$ and $Z_1$ are generic. The wall-crossing formula calculates the function $\Omega_1$ corresponding to $Z_1$ in terms of the function $\Omega=\Omega_0$ corresponding to $Z_0$.  This is analogous to the analytic continuation of a holomorphic function expressed in terms of its Taylor coefficients. The continuation is unique if it exists, and is not changed under a small deformation of the path with the fixed endpoints.

Let us call {\it strict} a sector in $\R^2$ with the vertex at the origin $(0,0)$ which is less than $180^\circ$.
With a strict sector  $V\subset \R^2$ we associate a group element $A_V$ given by the infinite product
$$ A_V:=\prod_{\gamma\in Z^{-1}(V )\cap \Gamma}^{\longrightarrow}
    \exp\left(- \Omega(\gamma)\sum_{n=1}^\infty \frac{e_{n\gamma}}{n^2}\right).$$

The product takes value in a certain pro-nilpotent Lie group $G_V:=G_{V,\Gamma, \langle \bullet,\bullet,\rangle}$, which we will describe below.
The right arrow in the product sign means that the product is taken in the {\it clockwise } order on the set  of rays $\R_+\cdot Z(\gamma)\subset V\subset \C$. For the product in the anti-clockwise order
we will use the left arrow.

Let us describe the  Lie algebra $\g_V=Lie(G_V)$ of the pro-nilpotent Lie group $G_V$.
We denote by $C(V)$ a convex cone in $\Gamma_{\R}$ which is the convex hull of the set of points $v \in Z^{-1}(V)$ such that $\parallel v \parallel\le C|Z(v)|$. The Lie algebra   $Lie(G_V)$ is the
infinite product $\prod_{\gamma \in \Gamma \cap  C(V)} \Q\cdot e_\gamma$ equipped with the above Lie bracket.

Now we can  formulate the wall-crossing formula.
It says (roughly) that $A_V$ does not change as long as no lattice point $\gamma\in  \Gamma$
with $\Omega_t(\gamma)\ne 0$ crosses the boundary of the cone $Z^{-1}_t(V)$ (here $\Omega_t$ corresponds to
the point $t\in [0,1]$). By our assumptions, if $t=t_0$ corresponds to a non-generic central charge $Z_{t_0}$ then there exists a $2$-dimensional lattice $\Gamma_0\subset \Gamma$  such that $Z_{t_0}(\Gamma_0)$ belongs to a real line ${\R}e^{i\alpha}$ for some $\alpha\in [0,\pi]$.

The wall-crossing formula describes  the change of values $\Omega(\gamma)$ for $\gamma \in \Gamma_0$ and depends only on the restriction
 $\Omega_{|\Gamma_0}$ of $\Omega$ to the  lattice $\Gamma_0$. Values $\Omega(\gamma)$ for $\gamma\notin \Gamma_0$ 
do not change at $t=t_0$. 
Denote by $k\in \Z$ the value of the form $\langle\bullet,\bullet\rangle$ on a fixed basis $\gamma_1,\gamma_2$ of $\Gamma_0\simeq \Z^2$ such that $C(V)\cap \Gamma_0\subset \Z_{\geqslant 0}\cdot\gamma_1\oplus \Z_{\geqslant 0}\cdot\gamma_2$. We assume that $k\ne 0$, otherwise there will be no jump in values of $\Omega$ on $\Gamma_0$.
The group  elements which we are interested in
can be identified with products of the following 
automorphisms\footnote{Here we write an automorphism as acting
on elements of the algebra of functions. The corresponding automorphism on {\emph{points}} is given by the inverse formula.} of $\Q[[x,y]]$ preserving the symplectic form $k^{-1}(xy)^{-1} dx\wedge dy$:
$$
T_{a,b}^{(k)}:(x,y)\mapsto $$
$$\mapsto \left( x\cdot (1-(-1)^{kab} x^a y^b)^{-kb}, y\cdot(1-(-1)^{kab} x^a y^b)^{ka}\right),\medskip a,b\geqslant 0, a+b\geqslant 1\,.
$$
For $\gamma=a\gamma_1+b\gamma_2$ we have
$$T_{a,b}^{(k)}=\exp\left(-\sum_{n\geqslant 1}\frac{e_{n\gamma}}{n^2}\right)$$
in the above notation.
Any exact symplectomorphism $\phi$ of $\Q[[x,y]]$ can be decomposed uniquely into a clockwise and an anti-clockwise product:
$$\phi=\prod_{a,b}^{\longrightarrow} \left(T_{a,b}^{(k)}\right)^{c_{a,b}}=\prod_{a,b}^{\longleftarrow} \left(T_{a,b}^{(k)}\right)^{{d}_{a,b}}$$
with certain exponents $c_{a,b},d_{a,b}\in \Q$. These exponents should be interpreted as the 
limiting values
of the functions $\Omega_{t_0}^{\pm}=\lim_{t\to t_0\pm 0}\Omega_t$ restricted to $\Gamma_0$.
The passage from the clockwise order (when the slope $a/b \in [0,+\infty]\cap {\mathbb P}^1(\Q)$ decreases) to the anti-clockwise order (when the slope increases) gives the change of $\Omega_{|\Gamma_0}$ as we cross the wall.
It will be convenient to denote $T_{a,b}^{(1)}$ simply by $T_{a,b}$. The pro-nilpotent group generated by 
 transformations $T_{a,b}^{(k)}$ coincides with the one generated by transformations $T_{a,|k|b}$.

The compatibility of the wall-crossing formula with the {\it integrality} of the numbers $\Omega(\gamma)$ is not obvious but follows from:

\begin{conj}
If for $k>0$ one decomposes the product $T_{1,0}\cdot T_{0,k}$ in the opposite order:
$$T_{1,0}\cdot T_{0,k}= \prod_{a/b \text{ increases}}  (T_{a,kb})^{d(a,b,k)},$$
then $d(a,b,k)\in \Z$ for all $a,b,k$.
\end{conj}

An equivalent form of this conjecture says that
 if one decomposes $T_{1,0}^k \cdot T_{0,1}^k$ in the opposite order then all exponents will belong to $k\Z$.

Here are decompositions for $k=1,2$ 
$$T_{1,0} \cdot T_{0,1} =T_{0,1}\cdot T_{1,1}\cdot T_{1,0}\,\,,$$
$$T_{1,0}^{(2)} \cdot T_{0,1}^{(2)}= T_{0,1}^{(2)}\cdot T_{1,2}^{(2)}\cdot T_{2,3}^{(2)}\cdot\dots
    \cdot (T_{1,1}^{(2)})^{-2} \cdot \dots\cdot T_{3,2}^{(2)} \cdot T_{2,1}^{(2)}
\cdot T_{1,0}^{(2)}\,\,,$$
or equivalently 
$$T_{1,0} \cdot T_{0,2}= T_{0,2}\cdot T_{1,4}\cdot T_{2,6}\cdot\dots
    \cdot T_{1,2}^{-2} \cdot \dots\cdot T_{3,4} \cdot T_{2,2}
\cdot T_{1,0}\,\,.$$
Greg Moore and Frederik Denef pointed out that the factors in the last formula correspond to the BPS spectrum of $N=2, d=4$ super Yang-Mills model studied by Seiberg and Witten in \cite{SeibW}. A ``physical" explanation of our formulas in this context was given in \cite{GaMNeit}, see also our Section 2.7.

For $k\geqslant 3$ or $k\le -1$ the  decomposition of $T_{1,0}^{(k)}\cdot T_{0,1}^{(k)}$ is  not yet known completely.
Computer experiments give a conjectural formula for the diagonal term with the slope $a/b=1$.
The corresponding symplectomorphism is given by the map
$$(x,y)\mapsto (x \cdot F_k(xy)^{-k}, y \cdot F_k(xy)^{k}),$$
where the series $F_k=F_k(t)\in 1+ t\Z[[t]]$ is an algebraic hypergeometric series given for $k\geqslant 3$ by the formulas:
$$\sum_{n=0}^\infty \binom {(k-1)^2n +k-1}{n}\frac {t^n}{(k-2)n+1}=
        \exp\left(\sum_{n=1}^\infty \binom{(k-1)^2n}{n}
\frac{k}{(k-1)^2}\frac{t^n}{n}\right).$$

The function $F_k$ satisfies the equation $$F_k(t)\left(1-t F^{k-2}_k(t)\right)^k-1=0\,\,.$$

\begin{rmk}
The above example for $k=1$ is compatible with the expected behavior of Donaldson-Thomas invariants when we have two spherical
objects $E_1,E_2\in \CC$ (sphericity means that ${\op{Ext}}^{\bullet}(E_i,E_i)=H^{\bullet}(S^3)$) such that

$${\op{Ext}}^1 (E_2,E_1)=\k, \,\,\, {\op{Ext}}^n(E_2,E_1)=0 \text{ for } n\ne 1\,\,.$$
In this case on the one side of the wall we have two semistable objects $E_1,E_2$, and on the other side we have three
semistable objects $E_1,E_2, E_{12}$ where $E_{12}$ is the extension of $E_2$ by $E_1$. In the case of the derived of the Fukaya category  the objects $E_i, i=1,2$ can correspond to embedded special Lagrangian spheres intersecting transversally at one point. Then $E_{12}$ corresponds to their Lagrangian connected sum.

\end{rmk}

The automorphisms $T_{a,b}$ are  a special case of the more general ones.
Namely, we can consider the following rational automorphisms of $\g_{\Gamma}$ (considered as a Poisson algebra):

$$ T_{\gamma}: e_{\mu}\mapsto (1-e_{\gamma})^{\langle \gamma, \mu\rangle}e_{\mu},\,\,\,\,\, \gamma,\mu \in \Gamma\,\,.$$
The group element $A_V$ in the above notation has the form
$$A_V=\prod_{\gamma\in Z^{-1}(V )\cap \Gamma}^{\longrightarrow} T_\gamma^{\,\Omega(\gamma)}$$
and acts on a completion of $\g_\Gamma$. 
It is easy to quantize this Poisson algebra.
The corresponding algebra (quantum torus) is additively generated by quantum generators $\hat{e}_{\gamma}, \gamma\in \Gamma$ subject to the relations
$$\hat{e}_{\gamma}\hat{e}_{\mu}=q^{{1\over{2}}\langle \gamma,\mu\rangle}\hat{e}_{\gamma+\mu}\,\,,$$
where $q$ is a parameter (with the classical limit $q^{1\over{2}}\to -1 $). Then one has formulas similar to the above for the ``quantum"  analogs of automorphisms
$T_{\gamma}, \gamma \in \Gamma$ (see Sections 6.4 and 7.1).

 For general $k\ge 2$ the decomposition of the product $T_{1,0}\cdot T_{0,k}$ as in Conjecture 1, describes numerical DT-invariants
of the Calabi-Yau category associated with the Kronecker quiver consisting of two vertices and $k$ parallel arrows (see Section 8 for a general theory). 
 Recent paper \cite{Re1} gives an explicit formula for $d(a,b,k)$  in terms of the Euler characteristic of the framed moduli space of semistable representations of the quiver. Moreover, a weak form of the integrality Conjecture 1 is proved in \cite{Re1}.

\subsection{Some analogies and speculations}

The above formulas for symplectomorphisms are  partially motivated
by \cite{KoSo1}, Section 10, where similar formulas appeared in a different problem.
Both formulas involve  Hamiltonian vector fields associated with the dilogarithm function.
The problem discussed in \cite{KoSo1} was the reconstruction of the rigid analytic K3 surface from its skeleton,
which is a sphere $S^2$ equipped with an integral affine structure, singular at a finite set of points.
The group which is very similar to the pro-nilpotent group $G_V$ was introduced in the loc. cit.
where we assigned symplectomorphisms to edges of a certain tree in $S^2$.  That tree should be thought of as an analog of the walls in the space
of stability structures. Edges of the tree (we called them ``lines" in \cite{KoSo1}) correspond to pseudo-holomorphic
discs with the boundary on the Lagrangian toric fibers of the dual K3 surface. When we approach the ``large complex structure limit" cusp in the moduli space of K3 surfaces, the discs degenerate into  gradient lines of some smooth functions on $S^2$, thus defining edges of the tree.
Hence the reconstruction problem for K3 surfaces (and for higher-dimensional
Calabi-Yau manifolds, see \cite{GS1}, \cite{GS2}) is governed by the counting of rational curves in the mirror dual Calabi-Yau manifold.
This observation suggests by analogy the questions below.

1) First, let us note that we may assume that the bilinear form $\langle \bullet,\bullet \rangle$ is non-degenerate on $\Gamma$
by replacing $\Gamma$ by a ``larger" lattice  (e.g. by $\Gamma\oplus \Gamma^{\vee}$, where $\Gamma^{\vee}=\op{Hom}(\Gamma, \Z)$ is the dual lattice, see Section 2.6). Then the Lie algebra $\g_{\Gamma}$ will be realized as the Lie algebra of exact Hamiltonian vector fields on the algebraic symplectic torus $\op{Hom}(\Gamma, {\bf G}_m)$. The collection of formal symplectomorphisms $A_V$ defined above give rise to a rigid analytic space ${\cal X}^{an}$ over any non-archimedean field, similarly to \cite{KoSo1}. This space carries an analytic symplectic form and describes ``the behavior at infinity" of a (possibly non-algebraic) formal smooth symplectic scheme over $\Z$. String Theory suggests that there exists an actual complex symplectic manifold ${\cal M}$ (vector or hyper multiplet moduli space) admitting a (partial) compactification
$\overline{\cal M}$ and such that 
$${\cal X}^{an}(\C((t)))=\overline{\cal M}(\C[[t]])\setminus ({\cal M}(\C[[t]])\cup (\overline{\cal M}\setminus {\cal M})(\C[[t]]))\,\,,$$
i.e. it is the space of formal paths hitting the compactifying divisor but not belonging to it).
In the case of the Fukaya category of a complex $3d$ Calabi-Yau manifold $X$ the space ${\cal M}$ looks ``at infinity" as a deformation of a complex symplectic manifold ${\cal M}^{cl}$ where $\dim  {\cal M}=\dim  {\cal M}^{cl}=\dim  H^3(X)$. The latter is the total space of the bundle ${\cal M}^{cl}\to {\cal M}_X$, where
${\cal M}_X$ is the moduli space of complex structures on $X$. The fiber of the bundle is isomorphic to the space 
$$(H^{3,0}(X)\setminus \{0\})\times (H^3(X,\C)/H^{3,0}(X)\oplus H^{2,1}(X)\oplus H^3(X,\Z))$$
 parametrizing  pairs {\it (holomorphic volume element, point of the intermediate Jacobian)}.\footnote{In a very interesting paper \cite{GaMNeit} a construction of the hyperk\"ahler structure on ${\cal M}$ was suggested by means of our wall-crossing formulas.}
Furthermore, as we discuss in Section 7.2, we expect that there is a complex integrable system associated with an arbitrary homologically smooth $3d$ Calabi-Yau category and the fiber being the ``Deligne cohomology" of the category.

2) Is it true that the counting of the  invariants $\Omega(\gamma)$ for ${\cal C}$ is equivalent to the counting of (some) holomorphic discs ``near infinity" in $\overline{\cal M}$? Is it possible to construct an $\A$-category associated with those discs and to prove
that it is a $3$-dimensional Calabi-Yau category?

3) The study of the dependence of BPS states on a point of the moduli space
of vector and hyper multiplets given in \cite{DiM} and \cite{DM} suggests that ${\cal M}$ is hyperk\"ahler and the  invariants $\Omega(\gamma)$  for ${\cal C}$ (counting of objects of $\CC$) can be interpreted as
the counting of some ``quaternion curves" in ${\cal M}$. Punctures ``at infinity" of those curves can be interpreted as $4d$ black holes.
It would be nice
to think about the problem of counting
such maps as  a ``quaternionic analog" of the counting of rational Gromov-Witten invariants.
Hopefully (by the analogy with the ``Gromov-Witten story")  one can define an appropriate
$\A$-category (``quaternionic Fukaya category") and prove that it is a $3$-dimensional Calabi-Yau category.
This would relate our invariants $\Omega(\gamma)$ with ``quaternionic" Gromov-Witten invariants.

4) Geometry similar to the one discussed in this paper also appears in the theory of moduli spaces
of holomorphic abelian differentials (see e.g. \cite{Zo}). The moduli space of abelian differentials
is a complex manifold, divided by real ``walls" of codimension one into pieces glued from
convex cones. It also carries a natural non-holomorphic action of the group $GL^+(2,\R)$. There is an analog of the central charge $Z$ in the story. It is given
by the integral of an abelian differential over a path between marked points in a complex curve.
This makes plausible the idea that the  moduli space of abelian differentials associated with
a complex curve with marked points, is isomorphic
to the moduli space of stability structures on the (properly defined) Fukaya category of this curve.

5) We expect that our wall-crossing formulas are related to those in the Donaldson theory of $4d$ manifolds with $b_2^+=1$ (cf. e.g. recent paper \cite{NaYo}) as well as with Borcherds hyperbolic Kac-Moody algebras and multiplicative automorphic forms. The formulas from \cite{ChVer} also look very similar.

\subsection{About the content of the paper}

In Section 2 we work out in  detail the approach to the invariants $\Omega(\gamma)$ and the wall-crossing formula sketched in the Introduction in the framework of graded Lie algebras.
It is based on the notion of stability data which admits two equivalent descriptions: in terms of a collection of elements $a(\gamma)$ of a graded Lie algebra $\g=\oplus_{\gamma\in \Gamma}\g_{\gamma}$ and in terms of a collection of group
elements $A_V$ which satisfy
the ``Factorization Property". The latter says  that $A_{V_1}A_{V_2}=A_{V}$ for any strict sector $V$ and its decomposition into two sectors $V_1,V_2$ (in the clockwise order) by a ray emanating from the vertex.
We define the topology on the space of stability data. It immediately leads to the wall-crossing
formula. Then we discuss a special case when the lattice  carry an integer-valued skew-symmetric bilinear form.
The skew-symmetric form on the lattice $\Gamma$ gives rise
to a Poisson structure on the torus $\op{Hom}(\Gamma,{\bf G}_m)$ of its characters.
Then we introduce a {\it double symplectic} torus, which corresponds to  the lattice $\Gamma\oplus \Gamma^{\vee}$.
This allows us to construct  an embedding of the pro-nilpotent groups $G_V$ (see Section 1.4) into the group of
formal symplectomorphisms
of the double torus. We show how the ``numerical DT-invariants" $\Omega(\gamma)$ arise from a collection of elements $A_V$ which satisfy
the Factorization Property $A_{V_1}A_{V_2}=A_{V}$ for any strict sector $V$. We introduce the notion of the ``wall of second kind" such that
(in the categorical framework) crossing such a wall corresponds to a change of the $t$-structure.
Then the multiplicative wall-crossing formula is equivalent to the triviality of the monodromy of a ``non-linear connection" on the space of numerical stability data. Also we discuss the relationship with the works of Joyce, and Bridgeland and Toledano-Laredo by introducing (under certain conditions) a  connection with irregular singularities on $\C$. In Section 2.7. we explain how 
stability data arise from complex
integrable systems. We illustrate our consideration by an example
of Seiberg-Witten curve. Arising geometry is the same as in the
``string junction" interpretation of Seiberg-Witten model (see e.g.
\cite{MikNekSet}).
The last section is devoted to stability data on $\mathfrak{gl}(n,\Q)$. It is related to the study of vacua in $N=2$ supersymmetric Quantum Field Theories (see \cite{CeVa}).

Section 3 is devoted to some basics on  ind-constructible categories, including
the definition of the topology on the space of stability structures. Also we discuss the notion of the potential of an object of Calabi-Yau category and the categorical version of the wall-crossing formula.
The way it is formulated is close intuitively to the physics considerations: we look how the ``motive" of the moduli space of semistable objects changes when some of exact triangles become unstable.

Section 4 is devoted to motivic  functions and motivic Milnor fiber. We start by recalling basics on motivic functions and motivic integration, including their equivariant versions (motivic stack functions, see also \cite{Jo4}).  Then we discuss the notion of motivic Milnor fiber introduced by Denef and Loeser as well as its $l$-adic incarnation. Rough idea is to use the motivic Milnor fiber of the potential of the $3d$ Calabi-Yau category in order to define invariants of the ind-constructible set of semistable objects.
The technical question arises: there might be two quadrics with the same rank and determinant but different Chow motives. In order to resolve this difficulty we introduce certain equivalence relation on motivic functions, so that in the quotient such quadrics are the same. Also, we discuss an important integral identity which 
will play the key role in Section 6.

Section 5 is devoted to an additional structure, which we call orientation data. It is a super line bundle on the space of objects of our category. Roughly, it is a square root of the super line bundle of cohomology. Although the numerical DT-invariants do not depend on the orientation data, the motivic DT-invariants introduced in Section 6 depend on it in an essential way.

Section 6 is devoted to the definition of motivic DT-invariants. 
First we define the {\it motivic Hall algebra} of an ind-constructible triangulated $\A$-category and prove its associativity. It generalizes the derived Hall algebra introduced by To\"en in \cite{To}. We define the motivic version $A_V^{\op{Hall}}$ of the   element $A_V$ as an invertible element of the completed
motivic Hall algebra associated with the sector $V$. The elements $A_V^{\op{Hall}}$ satisfy
the Factorization Property.
Basic idea behind the
Factorization Property (and hence the multiplicative wall-crossing formula) is that the infinite product
in the latter corresponds to the integration over the space of {\it all} objects of the category $\CC_V$ generated
by extensions of semistable objects with the central charge in $V$. The latter can be easily controlled when we cross the wall.

Motivic DT-invariants appear as elements of a certain quantum torus with the coefficient ring given by the equivalence classes of motivic functions. Basic fact is the theorem which says that
in the case of $3d$ Calabi-Yau category there is a homomorphism of the motivic Hall algebra into the motivic quantum torus defined in terms of the motivic Milnor fiber of the potential. In many cases the images of the elements $A_V^{\op{Hall}}$ can be computed explicitly in terms of the motivic version of the quantum dilogarithm function.
The images of $A_V^{\op{Hall}}$ are denoted by $A_{V}^{mot}$. This collection (one element for every strict sector $V$)
is called the motivic DT-invariant. The collection of these elements satisfy the Factorization Property. Replacing motives  by their Serre polynomials,
 we obtain $q$-analogs
of Donaldson-Thomas invariants,  denoted by $A_{V,q}$. We discuss their properties as well as the ``quasi-classical limit"  $A_V$ as $q^{1/2}\to -1$. We formulate the conjectures about the existence of the limit (absence of poles conjecture) and integrality property of the limits (integrality conjecture). The latter are related to the Conjecture 1 from Section 1.4. These conjectures are discussed in detail in Section 7, where we present various arguments and computations in their favor. Presumably, the technique developed by D.~Joyce  can lead to the proof
of  our conjectures. The numerical DT-invariants $\Omega:\Gamma\setminus\{0\}\to \Z$ are defined as coefficients in the decomposition of symplectomorphism $A_V$ into the product of powers $T_{\gamma}^{\,\Omega(\gamma)}$ in the clockwise order.

In Section 8 we consider in detail the case of $3$-dimensional Calabi-Yau category endowed with a finite
collection of spherical generators satisfying some extra property (cluster collection). Such categories correspond to quivers with potentials (Theorem 9).
Applying general considerations from the previous sections we formulate some results about quivers
and mutations. They are almost obvious in the categorical framework, but seem to be new in the framework
of quivers. Finally we explain that cluster transformations appear naturally as birational
symplectomorphisms of the double torus in the case when crossing of the wall of second kind corresponds to
a mutation at a vertex of the quiver (equivalently, to a mutation at the corresponding
spherical object of the Calabi-Yau category).

Several parts of the theory presented here have to be developed in more detail. This concerns  ind-constructible categories
 and motivic stack functions. Also, we present only a sketch of the proof of the $l$-adic version of the main identity in Section 4.4, leaving aside few 
technical details (which are not difficult to restore), and   the definition of the orientation data
 for cluster categories in Section 8.2 is left as a conjecture (although there is no doubt that it should be true).

{\it Acknowledgments.}  We thank to Mina Aganagic, Roma Bezrukavnikov, Tom Bridgeland, Frederik Denef, Emanuel Diaconescu, Pierre Deligne,  Sasha Goncharov, Mark Gross, Dominic Joyce, Greg Moore, Andrew Neitzke, Nikita Nekrasov, Andrei Okounkov, Rahul Pandharipande, Markus Reineke, Balazs Szendr\"oi, Don Zagier  for useful discussions and correspondence.  Y.S. thanks to IHES and the University Paris-6 for excellent research conditions.
His work was partially supported by an NSF grant.

\section{Stability conditions for graded Lie algebras}

\subsection{Stability data}

Let us fix a free abelian group $\Gamma$ of finite rank,
and a graded Lie algebra $\g=\oplus_{\gamma\in \Gamma}\g_{\gamma}$ over $\Q$.\footnote{In examples $\g$ is a $R$-linear Lie algebra, where $R$ is a commutative unital $\Q$-algebra.}

\begin{defn} Stability data on $\g$ is a pair
$\sigma=(Z,a)$ such that:

1) $Z: \Gamma\to \R^2\simeq \C$ is a homomorphism of abelian groups called the central charge;

2) $a=(a(\gamma))_{\gamma\in \Gamma\setminus \{0\}}$ is a collection of elements $a(\gamma)\in \g_{\gamma}$,

satisfying the following

\vspace{2mm}

{\bf Support Property}:
\vspace{2mm}

Pick a  norm $\parallel\bullet\parallel$ on the vector space $\Gamma_{\R}=\Gamma\otimes_{\Z}\R$.
Then there exists $C>0$ such that for any $\gamma\in \op{Supp}a $ (i.e. $a(\gamma)\ne 0$) one has
$$\parallel\gamma\parallel\le C|Z(\gamma)|\,\,.$$

\end{defn}

Obviously the Support Property does not depend on the choice of the norm. We will denote the set of all stability
data on $\g$ by $Stab(\g)$. Later we will equip this set with a Hausdorff topology.

The Support Property is equivalent to the following condition (which we will also call the Support Property):

{\it There exists a  quadratic form $Q$ on $\Gamma_{\R}$ such that

1) $Q_{|\op{Ker}\, Z}<0$;

2) $\op{Supp} a\subset \{\gamma\in \Gamma\setminus \{0\}|\,\,Q(\gamma)\geqslant 0\}$,

where we use the same notation $Z$ for the natural extension of $Z$ to $\Gamma_{\R}$.}

Indeed, we may assume that the norm $\parallel\bullet\parallel$ is the Euclidean norm in a chosen basis and
take $Q(\gamma)=-\parallel\gamma\parallel^2+C_1|Z(\gamma)|^2$ for sufficiently large positive
constant $C_1$. Generically $Q$ has signature $(2,n-2)$, where $n=\op{rk}\Gamma$. In degenerate
cases $Q$ can have signature $(1,n-1)$ or $(0,n)$.

For a given quadratic form $Q$ on $\Gamma_{\R}$ we denote by $Stab_Q(\g)\subset Stab(\g)$ the set of stability data satisfying the above conditions 1) and 2). Obviously $Stab(\g)=\cup_QStab_Q(\g)$, where the union is taken over all quadratic forms $Q$.

\begin{rmk}
In the case of a $3$-dimensional Calabi-Yau manifold $X$ there is a natural candidate for the quadratic
form $Q$ of the signature $(2,n-2)$ needed to formulate the Support Property.
Namely, identifying $H^3(X,\R)$ with $H^{3,0}(X,\C)\oplus H^{2,1}(X,\C)$ we can equip $H^3(X,\R)$
with the complex structure. Furthermore, the natural symplectic form coming from the Hodge structure gives rise to a pseudo-hermitian
form on $H^3(X,\R)$ of the signature $(2,n-2)$, where $n=\dim _{\R}H^3(X,\R)$. One can ask whether 
 this form is positive on every special Lagrangian submanifold of $X$. If this is true, then
the Support Property gives rise to a bound on the support of the function $\Omega$ discussed in Section 1.4.
\end{rmk}

Support Property implies the following estimate for the number of points in the $\op{Supp}a$ with the central charge
inside of the disc of radius $R$:
$$\#\left(Z(\op{Supp}a)\cap \{z\in \C|\,\,|z|\le R\}\right)=O(R^n)\,\,,$$
where $R\to \infty$ and $n=\op{rk}\Gamma$, Therefore the set $Z(\op{Supp}a)$ is discrete in $\C$
and does not contain zero.

\begin{rmk} It seem reasonable to consider ``limiting cases'' of stability data when the Support Property is not satisfied. Then the numbers $\op{Re}Z$ and $\op{Im}Z$ are allowed to take values in arbitrary totally ordered fields, e.g. $\R((t))$ (here $t$ is a formal parameter such that $t>0$ and $t<x$ for any $x\in \R_{>0}$). Some of our considerations below make sense in this situation. In the framework of stability conditions on triangulated categories such structures appeared in \cite{Ba}, \cite{Tod1}.

\end{rmk}

\subsection{Reformulation of the stability data}

In what follows we will consider various cones in $\Gamma_{\R}$ and in $\R^2$ i.e. subsets, which are
closed under addition and multiplication by a positive real number. We assume that the vertex of the cone (i.e. the zero of the vector space)
does not belong to the cone. We will call a cone {\it strict}
if it is non-empty and does not contain a straight line. In particular, all strict cones on the plane (we will call them strict sectors) are sectors, which are
smaller than 180 degrees (not necessarily closed or open). We allow the sector to be degenerate (which means that it is a ray
with the vertex at the origin).
We orient the plane (and hence all sectors) in the clockwise direction. We write $l_1\le l_2$
if the rays $l_1,l_2$ bound a strict closed sector and $l_1$ precedes $l_2$ in the clockwise order (we allow
$l_1=l_2$).

Let us fix a quadratic form $Q$ on $\Gamma_\R$.
We are going to describe below another set of data and will show that it is naturally isomorphic to the
set $Stab_Q(\g)$. Let ${\cal S}$ be the set of strict sectors in $\R^2$
possibly degenerate (rays).

We denote by $\widehat{Stab}_Q(\g)$ the set of pairs $(Z, A)$ such that:

{\it a) $Z: \Gamma\to \R^2$ is an additive map such that
$Q_{|\op{Ker}Z}<0$;

b) $A=(A_V)_{V\in {\cal S}}$ is a collections of elements $A_V\in G_{V,Z,Q}$, where
$G_{V,Z,Q}$ is a pro-nilpotent
group with the pro-nilpotent graded Lie algebra
$$\g_{V,Z,Q}=\prod_{\gamma\in\Gamma\cap C(V,Z,Q)}\g_{\gamma}\,\,,$$
where $C(V,Z,Q)$ is the convex cone generated by the set
$$S(V,Z,Q)=\{x\in \Gamma_{\R}\setminus \{0\}|\, Z(x)\in V, Q(x)\geqslant 0\}\,\,.$$}

The above definition makes sense because the cone $C(V,Z,Q)$ is strict, as one can easily see
by elementary linear algebra.
Hence for a triangle $\Delta$ which is cut from $V$ by a straight line,  any $\gamma\in Z^{-1}(\Delta)$ can be represented as a sum of other elements
of $\Gamma\cap C(V,Z,Q)$ in finitely many ways.
Furthermore, the triangle $\Delta$ defines an ideal $J_{\Delta}\subset \g_{V,Z,Q}$ consisting
of elements $y=(y_{\gamma})\in \g_{V,Z,Q}$ such that for every component $y_{\gamma}$
the corresponding $\gamma$ does not belong to the convex hull
of $Z^{-1}(\Delta)$. Then the quotient $\g_{\Delta}:=\g_{V,Z,Q}/J_{\Delta}$ is a  nilpotent
Lie algebra, and $\g_{V,Z,Q}=\varprojlim_{\Delta\subset V}g_{\Delta}$.

Let $G_{\Delta}=\exp(\g_{\Delta})$ be the nilpotent  group corresponding to the Lie algebra $\g_{\Delta}$. 
Then $G_{V,Z,Q}=\varprojlim_{\Delta}G_{\Delta}$ is a pro-nilpotent group.
If $V=V_1\sqcup V_2$ (in the clockwise order) then there are natural embeddings
$G_{V_i,Z,Q}\to G_{V,Z,Q}, i=1,2$.

We impose the following axiom on the set of pairs $(Z,A)$:

{\bf Factorization Property:}

{\it The element $A_V$ is given by the product
$A_V=A_{V_1}A_{V_2}$ where the equality is understood in $G_{V,Z,Q}$.}

We remark
that if $Q_1\le Q$ and both forms $Q,Q_1$ are negative on $\op{Ker}Z $ 
then $G_{V,Z,Q_1}\subset G_{V,Z,Q}$ for any $V\in {\cal S}$.
We say that the $(Z,A)\in \widehat{Stab}_Q(\g) $ and $(Z^{\prime},A^{\prime})\in \widehat{Stab}_{Q^{\prime}}(\g)$ are equivalent if
$Z=Z^{\prime}:=Z$ and there exists $Q_0$ such that $Q\le Q_0, Q^{\prime}\le Q_0$, ${Q_0}_{|\op{Ker}Z}<0$ and 
moreover for any $V\in {\cal S}$ we have $A_V=A^{\prime}_V$  as elements of the group
$G_{V,Z,Q_0}$.

\begin{thm} 1) For a fixed $Q$ there is a natural bijection between sets $\widehat{Stab}_Q(\g)$
and ${Stab}_Q(\g)$.

2) Any two elements of $\widehat{Stab}_Q(\g)$ and $\widehat{Stab}_{Q^{\prime}}(\g)$ are equivalent if and only if they define  the same element in $Stab(\g)$.
\end{thm}

{\it Proof.} Suppose that we are given a pair $(Z,A)\in \widehat{Stab}_Q(\g)$. In order to construct
the corresponding element in $Stab_Q(\g)$ we take the same $Z$ as a homomorphism $\Gamma\to \R^2$. 
What is left is to construct a collection
$a({\gamma})\in \g_{\gamma}$. We define it such as follows.

a) If $Z(\gamma)=0$ then we set $a({\gamma})=0$.

b) Suppose $Z(\gamma)\ne 0$. Let us consider the ray $l=\R_{>0}Z(\gamma)$.
Then we have an element $\log(A_l)\in \g_{l,Z,Q}\subset \prod_{\gamma\in \Gamma}\g_{\gamma}$.
We denote by $a({\gamma})$ the component of $\log(A_l)$ which belongs to $\g_{\gamma}$.
This assignment gives rise to stability data $(Z,a)\in {Stab}_Q(\g)$. In order to show that it is injective, we observe
that the Factorization Property implies that $A_V=\prod_{l\subset V}^{\longrightarrow}A_l$, where
the product is taken in the clockwise order over the set of all rays $l$ which belong to $V$.
Indeed, let us consider the image of $A_V$ in $G_{\Delta}$. Then only finitely many rays contribute
to the product $\prod_{l\subset V}^{\longrightarrow}A_l$, and the product formula follows from the Factorization Property.
Since $G_{V,Z,Q}=\varprojlim_{\Delta}G_{\Delta}$ the desired equality holds.

Conversely, if we have stability data $(Z,a)\in{Stab}_Q(\g)$  , then we construct a pair $(Z,A)$ taking the same $Z$
and $Q$, and for any ray $l$ we set
$$A_l=\exp\left(\sum_{\gamma\in \Gamma\cap C(l,Z,Q)}a({\gamma})\right).$$
Notice that $A_l=1$  if there are no elements $\gamma$ such that $Z(\gamma)\in l$. We define $A_V$ for any $V\in a$ using
the Factorization Property, i.e. $A_V=\prod_{l\subset V}^{\longrightarrow}A_l.$ This proves part 1) of the theorem. Part 2) follows immediately from definitions.
The theorem is proved. $\blacksquare$

\begin{rmk} We will use the same name ``stability data" for either of the set of data
which appear in the above theorem and will denote either set by $Stab(\g)$.

\end{rmk}

\begin{rmk} Let $\R^2\setminus \{(0,0)\}=\sqcup_{1\le i\le n}V_i$, where $V_i, 1\le i\le n$ are strict (semiclosed) sectors. Then  the stability data with
a given central charge $Z$ are uniquely determined by an arbitrary collection of elements
$A_{V_i}\in G_{V_i,Z,Q}$ for some quadratic form $Q$. 
\end{rmk}

There exists a generalization of stability data suitable for motivic Hall algebras. Namely, let us assume that
the Lie algebra $\g$ carries an automorphism $\eta$ such that $\eta(\g_{\gamma})=\g_{-\gamma}$
for any $\gamma\in \Gamma$.

\begin{defn} Symmetric stability data for $(\g,\eta)$ is a pair
$(Z,\widehat{a})$ where $Z:\Gamma\to \C$ is an additive map
and $\widehat{a}$ is a map $(\gamma,\varphi)\mapsto\widehat{a}(\gamma,\varphi)\in \g_{\gamma}$ where
$\varphi \in \R$,
$\gamma\in \Gamma$ is such that $Z(\gamma)\in \R_{>0}e^{i\varphi}$ and
$$\widehat{a}(\gamma,\varphi+\pi)=\eta(\widehat{a}(\gamma,\varphi))\,\,.$$

\end{defn}

All the considerations about stability data admit a straightforward generalization to the symmetric case.
We will use them without further comments.

\begin{rmk} Let $H_{\Gamma}$ be a $\Gamma$-graded unital associative algebra considered as a graded Lie algebra.
Then the pro-nilpotent groups $G_{V,Z,Q}$ discussed above are the groups of invertible elements of the form $f=1+\dots$
in the pro-nilpotent associative algebras which are completions of $H_{\Gamma}$.

\end{rmk}

\begin{rmk} Decomposition $\g=\oplus_{\gamma\in \Gamma}\g_{\gamma}$ and the Lie algebras $\g_{V,Z,Q}$ are similar to the root decomposition and nilpotent subalgebras in Kac-Moody Lie algebras.
The involution $\gamma\mapsto -\gamma$ is similar to the ``Cartan involution".
 These analogies deserve further study, since Donaldson-Thomas invariants (more precisely, counting functions for BPS states) appear in physics as
a kind of character formulas (see e.g. \cite{DiM}, formula (2.7)). In particular our multiplicative wall-crossing formulas in the case of wall of second kind should be related
to automorphic forms of Borcherds (see \cite{Bor}). The motivic Hall algebra defined below in Section 6 could be thought  of as the motivic version of the algebra of BPS states (see \cite{HM}).

\end{rmk}

\subsection{Topology and the wall-crossing formula}

Here we are going to introduce a Hausdorff topology on the set of stability data in such a way
that the forgetting map 
$$Stab(\g)\to \op{Hom}(\Gamma,\C)\simeq \C^n,\,\, (Z,a)\mapsto Z$$
 will be a local homeomorphism.
In particular $Stab(\g)$ carries a structure of a complex manifold (in general with an uncountable number of components, each of which is paracompact). In order to define the topology we define
the notion of a continuous family of points in $Stab(\g)$.

Let $X$ be a topological space, $x_0\in X$ be a point, and $(Z_x,a_x)\in Stab(\g)$ be a family parametrized by $X$.

\begin{defn} We say that the family is continuous at $x_0$ if the following conditions are satisfied:

a) The map $X\to \op{Hom}(\Gamma,\C), x\mapsto Z_x$ is continuous at $x=x_0$.

b) Let us choose a quadratic form $Q_0$ such that
$(Z_{x_0},a_{x_0})\in Stab_{Q_0}(\g)$. Then there exists an open
neighborhood $U_0$ of $x_0$ such that
$(Z_{x},a_{x})\in Stab_{Q_0}(\g)$ for all $x\in U_0$.

c) For any closed strict sector $V$ such that $Z(\op{Supp}a_{x_0})\cap \partial V=\emptyset$ the map
$$x\mapsto \log\,A_{V,x,Q_x}\in \g_{V,Z_x,Q_x}\subset \prod_{\gamma\in \Gamma}\g_{\gamma}$$
 is continuous
at $x=x_0$. Here we endow the vector space $\prod_{\gamma\in \Gamma}\g_{\gamma}$ with the product topology of discrete sets, and $A_{V,x,Q_x}$ is the group element associated with $(Z_x,a_x)$, sector $V$ and a quadratic form $Q_x$ such that 
$(Z_x,a_x)\in Stab_{Q_x}(\g)$.

\end{defn}

\begin{rmk} Part c) of the Definition 3 means that for any $\gamma\in \Gamma\setminus\{0\}$ the
component of $\log\,A_{V,x,Q_x}$ belonging to $\g_{\gamma}$ is locally constant as a function of $x$ in a neighborhood of $x_0$.

\end{rmk}
The element $\log(A_{V,x,Q_x})\in \prod_{\gamma\in \Gamma}\g_{\gamma}$ does not depend on $Q_x$, e.g. we can take $Q_x:=Q_0$ for $x$ close to $x_0$. 
The continuity means informally that for
any closed triangle $\Delta\subset \R^2$ with one vertex at the origin, the projection of
$\log\,A_{V,Z_x,Q_x}$ into the vector space $\oplus_{\gamma\in \Delta}\g_{\gamma}$ does not depend
on $x\in X$ as long as there is no element $\gamma\in \op{Supp}a_x$ such that $Z(\gamma)$ crosses the boundary
$\partial\Delta$.

It is easy to see that the above definition gives rise to a topology on $Stab(\g)$.
\begin{prp} This topology is Hausdorff.

\end{prp}
{\it Proof.} Let $(Z,a)$ and $(Z^{\prime},a^{\prime})$ be two limits of a sequence $(Z_n,a_n)$ as $n\to \infty$. We have to prove that $(Z,a)=(Z^{\prime},a^{\prime})$.
It is clear that $Z=Z^{\prime}$ since $\op{Hom}(\Gamma,\C)$ is Hausdorff. Let us now choose quadratic forms
$Q$ and $Q^{\prime}$ which are compatible with $a$ and $a^{\prime}$ respectively in the sense of Definition 3. Then there exists a quadratic form ${Q}_0$ such that ${Q}_0$ is negative on $\op{Ker}Z=\op{Ker}Z^{\prime}$ and also $Q\le {Q}_0, Q^{\prime}\le {Q}_0$. Then for all sufficiently large $n$ the form ${Q}_0$ is compatible with $a_n$.

For a generic sector $V\subset \R^2$ its boundary rays do not intersect $Z(\Gamma)$. By part c) of the Definition 3 we have: $A_{V,Z, {Q}_0}=A_{V,Z^{\prime}, {Q}_0}^{\prime}$ since the product
$\prod_{\gamma \in \Gamma}\g_{\gamma}$ is Hausdorff. Since any ray in $\R^2$ with the vertex at the origin can be obtained as an intersection of  generic sectors then we conclude that $a=a^{\prime}$. The Proposition is proved. $\blacksquare$

Let us fix an element $Z_0\in \op{Hom}(\Gamma,\C)$ and a quadratic form $Q_0$ compatible with $Z_0$ (i.e. negative on its kernel).
We denote by $U_{Q_0,Z_0}$ the connected component containing $Z_0$ in the domain $\{Z\in \op{Hom}(\Gamma,\C)|\,\, (Q_0)_{|\op{Ker} Z}<0\}$.
In what follows we will frequently use the following elementary observation.

\begin{prp} If $Q$ is a quadratic form on a finite-dimensional vector space $\Gamma_{\R}$ and $Z:\Gamma_{\R}\to \C$ is an $\R$-linear map such that $Q_{|\op{Ker}Z}<0$ then the intersection
$\{x\in \Gamma_{\R}|Q(x)\geqslant 0\}\cap Z^{-1}(l)$ is a convex cone (possibly empty) for any ray $l\subset \C$ with the vertex at the origin.

\end{prp}

Let $\gamma_1,\gamma_2\in \Gamma\setminus \{0\}$ be two $\Q$-linearly independent elements such that
$Q_0(\gamma_i)\geqslant 0, Q_0(\gamma_1+\gamma_2)\geqslant 0, i=1,2$.
We introduce the set
$${\cal W}^{Q_0}_{\gamma_1,\gamma_2}=\{Z\in U_{Q_0,Z_0}|\,\,\R_{>0}Z(\gamma_1)=\R_{>0}Z(\gamma_2)\}\,\,.$$
In this way we obtain a countable collection of hypersurfaces
${\cal W}^{Q_0}_{\gamma_1,\gamma_2}\subset U_{Q_0,Z_0}$ called the {\it walls corresponding to $Q_0,\gamma_1,\gamma_2$}. We denote their union by ${\cal W}_1:={\cal W}_1^{Q_0}$ and sometimes call it
{\it the wall of first kind} (physicists call it the wall of marginal stability).

Let us consider a continuous path $Z_t, 0\le t\le 1$ in $U_{Q_0,Z_0}$ which intersects each of these walls for finitely many values of $t\in [0,1]$. Suppose that we have a continuous lifting path $(Z_t,a_t)$ of this path such that $Q_0$ is compatible with each $a_t$ for all $0\le t\le 1$. Then for any
$\gamma\in \Gamma\setminus \{0\}$ such that $Q_0(\gamma)\geqslant 0$ the element $a_t(\gamma)$ does not change as long as $t$ satisfies the condition $$Z_t(\gamma)\notin \cup_{\gamma_1,\gamma_2\in \Gamma\setminus \{0\},\, \gamma_1+\gamma_2=\gamma}{\cal W}_{\gamma_1,\gamma_2}^{Q_0}\,\,.$$ If this condition is not satisfied we say that $t$ is a discontinuity point for $\gamma$.
For a given $\gamma$ there are finitely many discontinuity points.

Notice that for each $t\in [0,1]$ there exist limits 
$$a^{\pm}_t(\gamma)=\lim_{\varepsilon\to 0, \,\varepsilon>0}a_{t\pm \varepsilon}(\gamma)$$ (for $t=0$ or $t=1$ only one of the limits is well-defined).
Then the continuity of the lifted path $(Z_t,a_t)$ is equivalent to the following {\it wall-crossing formula}
which holds for {\it any $t\in [0,1]$} and arbitrary $\gamma\in \Gamma\setminus\{0\}$:
$$\prod_{\mu\in \Gamma^{prim},\,Z_t(\mu)\in l_{{\gamma},t}}^{\longrightarrow}\exp\left(\sum_{n\geqslant 1}a_t^{-}(n\mu)\right)=$$
$${}=\exp\left(\sum_{\mu\in \Gamma^{prim},\,Z_t(\mu)\in l_{{\gamma},t}, \,n\geqslant 1}a_t(n\mu)\right)=\prod_{\mu\in \Gamma^{prim},\,Z_t(\mu)\in l_{{\gamma},t}}^{\longrightarrow}\exp\left(\sum_{n\geqslant 1}a^{+}_t(n\mu)\right),$$
where $l_{\gamma,t}=\R_{>0}Z_t(\gamma)$, and $\Gamma^{prim}\subset \Gamma$ is the set of primitive vectors.
The first and the last products are taken in the clockwise order of $\op{Arg}(Z_{t-\varepsilon})$ and
$\op{Arg}(Z_{t+\varepsilon})$ respectively, where $\varepsilon>0$ is sufficiently small.
Moreover, for each $\gamma$ we have $a_t^{-}(\gamma)=a_t^{+}(\gamma)=a_t(\gamma)$ unless there exist non-zero $\gamma_1,\gamma_2$ such that $\gamma=\gamma_1+\gamma_2$ and $Z_t\in {\cal W}_{\gamma_1,\gamma_2}^{Q_0}$.

\begin{rmk} Informally speaking, the wall-crossing formula says that for a very small sector $V$ containing the ray $l_{{\gamma},t}$ the corresponding element $A_V$, considered as a function of time, is locally constant in a neighborhood of $t$.

\end{rmk}

For each $\gamma\in \Gamma\setminus \{0\}$ the wall-crossing formula allows us to calculate $a_{1}(\gamma)$ is terms of $a_{0}(\gamma^{\prime})$ for a finite collection of elements $\gamma^{\prime}\in \Gamma\setminus \{0\}$. Morally it is an inductive procedure on the ordered set of discontinuity points $t_i\in [0,1]$. The only thing we need to check is that for each $\gamma\in \Gamma\setminus \{0\}$ the computation involves finitely many elements of $\Gamma$. For that we need some preparation. First we introduce a partial order on the set $S_{Q_0}=(\Gamma\cap Q^{-1}_0(\R_{\geqslant 0}))\times [0,1]$ generated by the following relations:

a) $(\gamma,t)\geqslant (\gamma,t^{\prime})$ if $t\geqslant t^{\prime}$;

b) if $\gamma=\sum_{1\le i\le m}\gamma_i, Q_0(\gamma_i)\geqslant 0, Z_t(\gamma_i)\in \R_{>0}Z_t(\gamma), 1\le i\le m, m\geqslant 2$, where not all $\gamma_i$ belong to
$\Q\cdot\gamma$, then $(\gamma,t)\geqslant (\gamma_i,t)$ for all $1\le i\le m$.

\begin{lmm} For any $(\gamma,t)\in S_{Q_0}$ the set $(\gamma^{\prime},t^{\prime})\in S_{Q_0}$
such that $(\gamma^{\prime},t^{\prime})\le (\gamma,t)$ is a finite union of  sets of the form
$\{\gamma_{\alpha}\}\times [0,t_{\alpha}]$.

\end{lmm}

The Lemma immediately implies the desired result.

\begin{cor} The element $a_t({\gamma})$ is a finite Lie expression of the elements $a_{0}(\gamma_{\alpha})$.

\end{cor}

{\it Proof of the Lemma}. Let us assume the contrary. Then we have an infinite sequence $t_1>t_2>t_3>\dots$
such that 
$$(\gamma_1,t_1)\geqslant _{a)}(\gamma_2,t_2)\geqslant_{b)}(\gamma_3,t_3)\geqslant_{a)}(\gamma_4,t_4)\geqslant_{b)}\dots\,\,,$$
where the subscript $a)$ or $b)$ denotes the two different possibilities for the partial order defined above.
Let $t_{\infty}=\lim_{n\to \infty}t_n$. It is easy to see that there exists a Euclidean norm $\parallel\bullet\parallel$ on $\Gamma_{\R}$ such that for any $v_1,v_2\in \Gamma_{\R}$ satisfying the properties $Q_0(v_i)\geqslant 0, Z_{t_{\infty}}(v_1)\in \R_{>0}Z_{t_{\infty}}(v_2)$ we have the inequalities
$\parallel v_i\parallel <\parallel v_1+v_2\parallel$ for $i=1,2$.

%\vspace{3mm}
%FIGURE 4 (Vectors and their sum in the cone)

%\centerline{\epsfbox{vectors.eps}}

%\vspace{3mm}

Moreover the same property holds if we replace the map $Z_{t_{\infty}}$  by an additive map $Z$ which is close to it. Then we conclude that $$\parallel \gamma_{2n}\parallel>\parallel \gamma_{2n+1}\parallel=
\parallel \gamma_{2n+2}\parallel$$ for all sufficiently large $n$. This contradicts to the fact the lattice $\Gamma$ is discrete in $\Gamma_{\R}$. The Lemma is proved. $\blacksquare$

The previous discussion allows us to lift a generic path $Z_t, 0\le t\le 1$ as above to
a unique continuous path $(Z_t,a_t)\in Stab(\g), 0\le t\le 1$ which starts at a given point $(Z_0,a_0)\in Stab(\g)$.
In other words, we have the notion of a {\it parallel transport} along a generic  path.
This observation is a part of the following more general statement.

\begin{thm} For given quadratic form  $Q_0$ and $(Z_0,a_0)\in Stab_{Q_0}(\g)$ there exists a unique continuous map $\phi: U_{Z_0,Q_0}\to Stab(\g)$ such that it is a section of the natural projection $Stab(\g)\to U_{Z_0, Q_0}$, and  $\phi(Z_0)=(Z_0,a_0)$.

\end{thm}

{\it Proof.} We have already proved the existence of a lifted path $(Z_t,a_t)$ for a generic path $Z_t$ provided the beginning point $Z_0$ is fixed.
What is left to prove that the endpoint $(a_1,Z_1)$  does not depend on a choice of the generic path $Z_t$. We are going to sketch the proof leaving the details for the reader.

Let us consider an infinitesimally small loop around the intersection point $Z$ of two or more walls. We would like to prove that the monodromy of the parallel transport along the loop is trivial.
There are two possibilities:

a) there are two different sublattices $\Gamma_1,\Gamma_2\subset \Gamma$ of ranks $\geqslant 2$ such that $Z(\Gamma_i), i=1,2$ belong to two different lines in the plane $\R^2$;

b) there exists a sublattice $\Gamma_3\subset \Gamma$ such that $\op{rk}\Gamma_3\geqslant 3$ and $Z(\Gamma_3)$ belongs to a line in $\R^2$.

In the case a) the corresponding Lie subalgebras of the completion of $\g_{\Gamma}$ are graded by non-intersecting subsets of $\Gamma$. Hence the corresponding wall-crossing transformations commute.

In the case b) let us choose a decomposition $\R^2\setminus \{(0,0)\}=\sqcup_{1\le i\le 4}V_i$, where $V_i, 1\le i\le 4$ are strict sectors such that $\R\cdot Z(\Gamma_3)\subset V_1\sqcup V_3\sqcup \{(0,0)\}$. When we move around the infinitesimally small loop the element $a(\gamma)$ can change only for $\gamma\in \Gamma_3$. Hence we can replace $\Gamma$ by $\Gamma_3$ in all computations. The wall-crossing formula implies that the elements $A_{V_i}, 1\le i\le 4$ do not change along the loop (moreover, by our assumption we have $A_{V_2}=A_{V_4}=1$). By Remark 7
from Section 2.2 we conclude that the stability data with the central charge $Z^{\prime}$ which is close to $Z$ are uniquely determined by $Z^{\prime}$ and the collection of elements $A_{V_i}, 1\le i\le 4$.
Hence the monodromy around the loop is trivial.

Finally one has to check the the {\it global} monodromy around a loop in $U_{Z_0,Q_0}$ is trivial. It follows from the fact that the fundamental group $\pi_1(U_{Z_0,Q_0})$ is generated by the loop $Z\mapsto Ze^{2\pi it}, t\in [0,1]$. But the monodromy around this  loop is trivial for generic $Z$, because the loop does not intersect the walls. $\blacksquare$.

We can write the wall-crossing formula in the way similar to the one from the Introduction.
In the case of generic path we have at a discontinuity point $t_0\in [0,1]$ a two-dimensional lattice $\Gamma_0\simeq \Z^2$ which is projected by $Z_{t_0}$ into a real line in $\R^2$.
We choose an isomorphism $\Gamma_0\simeq \Z^2$ in such a way that $Q_0^{-1}(\R_{\geqslant 0})\cap (\Gamma_0\setminus \{0\})$ is contained in $\Z_{>0}^2\cup \Z^2_{<0}$. Also we assume that the orientation of $\Gamma_0\otimes \R$ defined by $Z_t$ agrees with the one on $\Z^2_{>0}$ for $t=t_0-\varepsilon$ and is opposite to it for $t=t_0+\varepsilon$, where $\varepsilon>0$ is sufficiently small.

Then if
$\gamma=(m,n)\in \Z^2_{>0}$, and $a^{\pm}_{t_0}(\gamma):=a^{\pm}(m,n)$, we can write the
the wall-crossing formula in the following way:
$$\prod_{(m,n)=1}^{\longrightarrow}\exp\left(\sum_{k\geqslant 1}a^{-}(km,kn)\right)=\prod_{(m,n)=1}^{\longleftarrow}\exp\left(\sum_{k\geqslant 1}a^{+}(km,kn)\right),$$
where in the LHS we take the product over all coprime $m,n$ in the increasing order of $m/n\in \Q$,
while in the RHS we take the product over all coprime $m,n$ in the decreasing order.
Both products are equal to
$\exp(\sum_{m,n> 0}a_{t_0}(m,n))$.

\subsection{Crossing the wall of  second kind}

Here we will interpret the parallel transport in a different way, introducing a wall of another kind.
We use the notation from the previous section. In particular, we fix the quadratic form $Q_0$ and the connected component $U$ of the set $\{Z\in \op{Hom}(\Gamma,\C)|\,\, (Q_{0})_{|\op{Ker} Z}<0\}$.

For a given primitive $\gamma\in \Gamma\setminus \{0\}$ we introduce the set ${\cal W}_{\gamma}^{Q_0}=\{Z\in U|\,Z(\gamma)\in \R_{>0}\}$. It is a hypersurface in $U$. We call it a {\it wall of second kind} associated with $\gamma$. We call the union $\cup_{\gamma}{\cal W}_{\gamma}^{Q_0}$ the {\it wall of second kind} and denote it by
${\cal W}_2$.

\begin{defn} We say that a path $\sigma=(Z_t)_{0\le t\le 1}\subset U$ is short if the convex cone
$C_{\sigma}$ which is the convex hull of $\left(\cup_{0\le t\le 1}Z_t^{-1}(\R_{>0})\right)\cap \{Q_0\geqslant 0\}$
is strict.

\end{defn}

With a short path we  associate a pro-nilpotent group $G_{C_{\sigma}}$ with the Lie algebra
$\g_{C_{\sigma}}=\prod_{\gamma\in C_{\sigma}\cap \Gamma}\g_{\gamma}$.

The following result is obvious.
\begin{prp} For a generic short path $\sigma=(Z_t)_{0\le t\le 1}$ there exists no more than countable
set $t_i\in [0,1]$ and corresponding primitive $\gamma_i\in \Gamma\setminus \{0\}$ such that $Z_{t_i}\in {\cal W}_{\gamma_i}^{Q_0}$. For each $i$ we have: $\op{rk} Z_{t_i}^{-1}(\R)\cap \Gamma=1$.

\end{prp}

Let us recall the continuous lifting map $\phi: U\to Stab(\g)$ from the previous section. In the notation of the previous Proposition we define for any $t_i$ a group element
$$A_{t_i}=\exp\left(\varepsilon_i\sum_{n\geqslant 1}a_{t_i}({n\gamma_i})\right)\in G_{C_{\sigma}}\,\,,$$
where $\varepsilon_i=\pm 1$ depending on the direction in which the path $Z_t(\gamma_i)$ crosses $\R_{>0}$ for $t$ sufficiently close to $t_i$.

\begin{thm} For any short loop the monodromy $\prod_{t_i}^{\longrightarrow}A_{t_i}$ is equal to the identity
(here the product is taken in the increasing order of the elements $t_i$).

\end{thm}

{\it Proof.} Here we also present a sketch of the proof. Similarly to the proof of the Theorem 3 we consider the case of infinitesimally small loop $\sigma$
around a point $Z$ such that $\op{rk} \Gamma_2=2$ where $\Gamma_2:=Z^{-1}(\R)\cap \Gamma$ (i.e. $Z$ is a point where two, and hence infinitely many, walls of second kind intersect). Since $\sigma$ is infinitesimally small we can replace $\Gamma$ by $\Gamma_2$. Then we have the space $\op{Hom}(\Gamma_2,\C)\simeq \R^4$ which contains a countable collection of walls consisting of those $Z: \Gamma_2\to \C$ for which there exists $\gamma\in \Gamma_2\setminus \{0\}$ such that $Q_0(\gamma)\geqslant 0$ and $Z(\gamma)\in \R$. All the hypersurfaces contain $\R^2= \op{Hom}(\Gamma_2,\R)\subset \op{Hom}(\Gamma_2,\C)$. Factorizing by this subspace $\R^2$ we obtain a collection of lines with rational slopes in the union of two opposite strict sectors  $S\cup (-S)\subset \R^2=\op{Hom}(\Gamma_2,i\R)$.

%\vspace{3mm}

%FIGURE 5 (Two symmetric sectors with lines)

%\centerline{\epsfbox{sectors.eps}}

%\vspace{3mm}

We have to prove that the product over a loop surrounding $(0,0)$ is the identity element. But it is easy to see that the product over the rays belonging to each of the sectors is equal to the left (resp. right) hand side of the wall-crossing formula. $\blacksquare$

Let us now introduce a set
${\cal X}_1\subset \Gamma\times U$ which consists of pairs $(\gamma,Z)$ such that $\gamma\in\Gamma\setminus\{0\}$ is a non-zero  element, $Q_0(\gamma)\geqslant 0, Z(\gamma)\in \R_{>0}$ and $Z^{-1}(\R_{>0})\cap \Gamma=(\Q_{>0}\cdot \gamma)\cap \Gamma$.

\begin{prp} The set of continuous sections $\psi: U\to Stab(\g)$ such that $\psi(Z)$ is compatible with $Q_0$ for any $Z\in U$ is in one-to-one correspondence with functions
$\tilde{a}: {\cal X}_1\to \g$ such that $\tilde{a}(\gamma,Z)\in \g_{\gamma}$ satisfying the property that for any small loop $\sigma$ the monodromy defined in the previous theorem is equal to the identity.

\end{prp}

{\it Proof.} The bijection is given by the formula $\tilde{a}(\gamma, Z)=a_{\psi(Z)}(\gamma)$. By the previous theorem the corresponding monodromy is trivial. Conversely, the triviality of the monodromy is equivalent to the wall-crossing formula in the special case when a $2$-dimensional sublattice of $\Gamma$ is mapped by $Z$ into the line $\R\subset \C$. The general case of an arbitrary line can be reduced to this one by a rotation $Z\mapsto Ze^{2\pi it}$ (it does not change the values $a(\gamma)$ because we do not cross the wall of first kind). $\blacksquare$

Let us also introduce a set ${\cal X}_2\subset \Gamma\times U$ which consists  of such pairs $(\gamma,Z)$ that $Q_0(\gamma)\geqslant 0, Z(\gamma)>0$ and there are no non-zero $\Q$-independent elements $\gamma_1,\gamma_2\in \Gamma$ with the property $\gamma=\gamma_1+\gamma_2, Q_0(\gamma_i)\geqslant 0, Z(\gamma_i)>0, i=1,2$. Since ${\cal X}_2$ is a locally-closed hypersurface in $U\times \op{Hom}(\Gamma,\C)$ it has finitely many connected components. Obviously, we have ${\cal X}_1\subset {\cal X}_2$.

It follows from the wall-crossing formula that for a continuous section $\psi:U\to Stab(\g)$ the restriction of the function $a$ to ${\cal X}_2$ is locally-constant and uniquely determines the section $\psi$. Therefore, the values of the restriction $a_{|\pi_0({\cal X}_2)}$ provides a countable coordinate system (satisfying a countable system of equations) on the set of continuous sections
$\{\psi: U\to Stab(\g)\}$ as above. It can be compared with another countable coordinate system (with no constraints) given the value $\psi(Z_0)$ for $Z_0\in U$. The latter coordinate system is not very convenient since one has to choose a generic $Z_0$.

\subsection{Invariants $\Omega(\gamma)$ and the dilogarithm}

Let $\Gamma$ be a free abelian group of finite rank $n$ as before, endowed
with a skew-symmetric integer-valued  bilinear form $\langle \bullet, \bullet\rangle : \Gamma \times \Gamma\to \Z$. Recall the Lie algebra
$\g_{\Gamma}=\g_{\Gamma,\langle \bullet,\bullet \rangle}=\oplus_{\gamma\in \Gamma}\Q\cdot e_{\gamma}$ with the Lie bracket
$$[e_{\gamma_1},e_{\gamma_2}]=(-1)^{\langle \gamma_1, \gamma_2\rangle}\langle \gamma_1, \gamma_2\rangle e_{\gamma_1+\gamma_2}\,\,.$$
Let us introduce a commutative associative product on $\g_{\Gamma}$ by the formula
$$e_{\gamma_1}e_{\gamma_2}=(-1)^{\langle \gamma_1, \gamma_2\rangle}e_{\gamma_1+\gamma_2}\,\,.$$
We denote by ${\mathbb T}_{\Gamma}:={\mathbb T}_{\Gamma,\langle \bullet, \bullet\rangle}$ the spectrum of this commutative algebra. It is easy to see that ${\mathbb T}_{\Gamma}$ is a torsor over the algebraic torus $\op{Hom}(\Gamma, {\bf G}_m)$. Moreover ${\mathbb T}_{\Gamma}$ is an algebraic Poisson manifold with the Poisson bracket $$\{a,b\}:=[a,b]\,.$$ The Poisson structure on ${\mathbb T}_{\Gamma}$ is invariant with respect to the action of $\op{Hom}(\Gamma, {\bf G}_m)$.

We can specify the results of the previous sections to the Lie algebra $\g_{\Gamma}$.
For stability data $(Z,a)$ we can write uniquely (by the M\"obius inversion formula) 
$$a(\gamma)=-\sum_{n\geqslant 1, {1\over{n}}\gamma\in \Gamma\setminus \{0\}}{\Omega(\gamma/n)\over{n^2}}e_{\gamma}\,\,,$$
where $\Omega:\Gamma\setminus \{0\}\to \Q$ is a function. Then we have
$$\exp\left(\sum_{n\geqslant 1}a(n\gamma)\right)=\exp\left(-\sum_{n\geqslant 1} \Omega(n\gamma)\sum_{k\geqslant 1}{e_{kn\gamma}\over{k^2}}\right):=
\exp\left(-\sum_{n\geqslant 1}\Omega(n\gamma)\op{Li}_2(e_{n\gamma})\right),$$
where $\op{Li}_2(t)=\sum_{k\geqslant 1}{t^k\over{k^2}}$ is the dilogarithm function.

The Lie algebra $\g_{\Gamma}$ acts on ${\mathbb T}_{\Gamma}$ by Hamiltonian vector fields. Let us denote by $T_{\gamma}$ the formal Poisson automorphism
$$T_{\gamma}=\exp(\{-\op{Li}_2(e_{\gamma}),\bullet\})\,\,,\,\,\, T_\gamma(e_\mu)=(1-e_\gamma)^{\langle \gamma,\mu \rangle} e_\mu$$
considered as an automorphism of algebra of functions.

More precisely for any strict convex cone $C\subset \Gamma_{\R}$ containing $\gamma$ the element $T_{\gamma}$ acts on the formal scheme $Spf(\prod_{\mu\in \Gamma \cap C}\Q e_{\mu})$. Moreover $T_{\gamma}$ is the Taylor expansion of a birational automorphism of ${\mathbb T}_{\Gamma}$.

Finally, in the case when $\Gamma$ comes from a $3d$ Calabi-Yau category the numbers $\Omega(\gamma)$ are (conjecturally) integers for $(\gamma,Z)\in {\cal X}_2$ in notation of Section 2.4. They provide generalization of DT-invariants (BPS degeneracies in physics language).

\subsection{Symplectic double torus}

If the skew-symmetric bilinear form on $\Gamma$ is degenerate, then the
action of $\g_{\Gamma}$ on
${\mathbb T}_{\Gamma}$ is not exact. In order to remedy the problem we can embed
$(\Gamma, \langle\bullet,\bullet\rangle)$ into a larger symplectic lattice. A possible choice is $\Gamma\oplus \Gamma^{\vee}$, where $\Gamma^{\vee}=\op{Hom}(\Gamma,\Z)$. The corresponding non-degenerate bilinear form is 
$$\langle(\gamma_1,\nu_1),(\gamma_2,\nu_2)\rangle=\langle \gamma_1,\gamma_2\rangle+\nu_2(\gamma_1)-\nu_1(\gamma_2)\,\,.$$

Let us choose a basis $e_i, 1\le i\le n=\op{rk}\Gamma$ of $\Gamma$. It gives rise to the coordinates $y_i, 1\le i\le n$ on ${\mathbb T}_{\Gamma}$.
The Poisson structure on ${\mathbb T}_{\Gamma}$ can be written as
$$\{y_i,y_j\}=b_{ij}y_iy_j\,\,,$$
where $b_{ij}=\langle e_i,e_j \rangle$.

Let us also introduce additional coordinates $x_j, 1\le j\le \op{rk}\Gamma^{\vee}$ in such a way that
$(y_i,x_j), 1\le i,j\le n$ will be the coordinates on the double torus $D({\mathbb T}_{\Gamma})$ with the Poisson brackets
$$\{x_i,x_j\}=0,\{y_i,x_j\}=\delta_{ij}y_i x_j\,\,.$$
There is a projection
$$\pi:D({\mathbb T}_{\Gamma})\to {\mathbb T}_{\Gamma}, \pi((y_i)_{1\le i\le n},(x_j)_{1\le j\le n})=(y_i)_{1\le i\le n}\,\,.$$
Notice that  $\pi$ is a Poisson morphism of the symplectic manifold
$D({\mathbb T}_{\Gamma})$ onto the Poisson manifold ${\mathbb T}_{\Gamma}$.

Let $C\subset \Gamma_{\R}$ be a closed convex strict cone.
Let us choose
a closed convex strict cone $C_1\subset (\Gamma\oplus\Gamma^{\vee})\otimes {\R}$ which contains $C\oplus\{0\}$.
With the cone $C_1$ we associate the Poisson algebra $\Q[[C_1]]$ consisting of series
$\sum_{\gamma,\delta\in C_1\cap (\Gamma\oplus\Gamma^{\vee})}c_{\gamma,\delta}y^{\gamma}x^{\delta}$. The pro-nilpotent group $G_C=\exp(\prod_{\gamma\in C\cap \Gamma}\g_{\gamma})$ acts by Poisson automorphisms
of $\Q[[C_1]]$.

Let us consider a closed algebraic submanifold $N\subset D({\mathbb T}_{\Gamma})$ defined by the equations
$$y_i\prod_jx_j^{b_{ij}}=-1, \,\,\,\,\,\,1\le i\le n\,\,.$$

\begin{lmm} The image of the group $G_C$ preserves the corresponding completion of $N$.

\end{lmm}

{\it Proof.} It suffices to check that the image of the Lie algebra $\g_{C}$ preserves the equations
of $N$. Notice that this image belongs to the Lie algebra of Hamiltonian vector fields
on $D({\mathbb T}_{\Gamma})$ generated by $\{y^{\gamma},\bullet\}$, where $\gamma\in \Gamma$ and
$y^{\gamma}=y_1^{\gamma_1}\dots y_n^{\gamma_n}$. Taking logarithms we see that
$$\{\log(y^{\gamma}),\log(y_i)+\sum_j b_{ij}\log(x_j)\}=\sum_j\gamma_j b_{ji}+\sum_j b_{ij}\gamma_j=0\,\,.$$
This concludes the proof. $\blacksquare$

\begin{rmk} It is clear that the action of the image of $G_{C}$ also commutes with the map $\pi$.
Moreover the image of $G_C$ in the group of exact symplectomorphisms of the completion of ${\mathbb T}_{\Gamma}$ corresponding to $C$ can be  characterized by the property that it preserves the completion of
$N$ and commutes with $\pi$.

\end{rmk}

Let us finally make a remark about a possible non-archimedean geometry interpretation of our construction.
Let us choose a complete non-archimedean field $K$ with the residue field of characteristic zero. Extending scalars we can think of the algebraic variety $D({\mathbb T}_{\Gamma})$ as of variety over $K$. We denote by
$D({\mathbb T}_{\Gamma})^{an}$ the corresponding non-archimedean $K$-analytic space in the sense of Berkovich (see \cite{KoSo1} for the explanation of the relevance of Berkovich approach to the large complex structure limit of Calabi-Yau varieties). Then the group $G_C$ acts on the analytic subset of $D({\mathbb T}_{\Gamma})^{an}$ given by inequalities $\{|e_{\gamma}|<1, \gamma\in C\setminus \{0\}\}$. Here we interpret
$e_{\gamma}$ as a Laurent monomial on $D({\mathbb T}_{\Gamma})$.

The symplectic double torus together with submanifold $N$ will be used again only in Section 8.

\subsection{Complex integrable systems and stability data}

In this section we explain how complex integrable systems (with some additional structures)
give rise to stability data in the graded Lie algebra $\mathfrak{g}_\Gamma$ associated with a symplectic lattice. In particular, Seiberg-Witten differential can be interpreted as the central charge
for a complex integrable system, while the BPS degeneracies are interpreted via our ``numerical" Donaldson-Thomas invariants as the number of certain gradient trees on the base of a complex integrable system.

Recall that a {\it complex integrable system} is a holomorphic map $\pi:X\to B$ where $(X,\omega^{2,0}_X)$ is a holomorphic symplectic manifold, $\dim X=2\dim B$,  and the generic fiber of $\pi$ is a Lagrangian submanifold, which is a polarized abelian variety. We assume (in order to simplify the exposition) that the polarization is
 principal. The fibration $\pi$ is non-singular outside of a closed subvariety $B^{sing}\subset B$ of codimension at least one.
It follows that on the open subset $B^{sm}:=B\setminus B^{sing}$ we have a local system ${\bf \Gamma}$ of symplectic lattices with the fiber over $b\in B^{sm}$ equal to $\Gamma_b:=H_1(X_b,\Z), X_b=\pi^{-1}(b)$ (the symplectic structure on $\Gamma_b$ is given by the  polarization).

Furthermore, the set $B^{sm}$ is locally (near each point $b\in B^{sm}$) embedded as a holomorphic Lagrangian subvariety into an affine symplectic space parallel to $H_1(X_b,\C)$. Namely, let us choose a symplectic basis $\gamma_i\in \Gamma_b, 1\le i\le 2n$. Then we have a collection of holomorphic closed
$1$-forms $\alpha_i=\int_{\gamma_i}\omega^{2,0}_X, 1\le i\le 2n$ in a neighborhood of $b$. There exists (well-defined locally up to an additive constant) holomorphic functions $z_i, 1\le i\le 2n$ such that $\alpha_i=dz_i, 1\le i\le 2n$. They define an embedding of a neighborhood of $b$ into $\C^{2n}$. The collection of $1$-forms $\alpha_i$ gives rise to an element $\delta\in H^1(B^{sm}, {\bf \Gamma}^{\vee}\otimes \C)$. {\it We  assume that $\delta=0$}.
This assumption is equivalent to an existence of a section
$Z\in \Gamma(B^{sm},{\bf \Gamma}\otimes {\cal O}_{B^{sm}})$ such that $\alpha_i=Z(\gamma_i), 1\le i\le 2n$.
\begin{defn} We call $Z$ the central charge of the integrable system.

\end{defn}

Hence, for every point $b\in B^{sm}$ we have a symplectic lattice $\Gamma_b$ endowed with an additive map
$Z_b:\Gamma_b\to \C$. Our goal will be to define a continuous family of stability data on graded Lie algebras
  $\mathfrak{g}_{\Gamma_b}$ with central charges $Z_b$.

First, we show an example of section $Z$.
\begin{exa} (Seiberg-Witten curve)

Let $B=\C$ be a complex line endowed with a complex coordinate $u$. We denote by $X^0=T^{\ast}\left(\C\setminus\{0\}\right)$ the cotangent bundle to the punctured line. We endow it with the coordinates $(x,y), y\ne 0$ and the symplectic form 
$$\omega^{2,0}=dx\wedge{dy\over{y}}\,\,.$$ There is a projection $\pi^0:X^0\to B$ given by 
$$\pi(x,y)={1\over{2}}(x^2-y-{c\over{y}})\,\,,$$
 where $c$ is a fixed constant. Fibers of $\pi^0$ are punctured elliptic curves 
$$y+{c\over{y}}=x^2-2u\,\,.$$ We denote by $X$ the compactification of $X^0$ obtained by the compactifications of the fibers. We denote by $\pi:X\to B$ the corresponding projection. Then $Z_u\in H^1(\pi^{-1}(u),\C)$ is represented by a meromorphic $1$-form $\lambda_{SW}={xdy\over{y}}$ (Seiberg-Witten form). The form $\lambda_{SW}$ has zero residues, hence it defines an element of $H^1(\pi^{-1}(u),\C)$ for any $u\in B^{sm}$, where $B^{sm}=B\setminus\{b_{-},b_{+}\}$ consists of points where the fiber of $\pi$ is a
non-degenerate elliptic curve.

\end{exa}

The dense open set $B^{sm}\subset B$ carries a K\"ahler form 
$$\omega_{B}^{1,1}=\op{Im}\left(\sum_{1\le i\le n}\alpha_i\wedge\overline{\alpha_{n+i}}\right)\,\,.$$ We denote by $g_{B}$ the corresponding K\"ahler metric.

For any $t\in \C^{\ast}$ we define an integral affine structure on $C^{\infty}$-manifold $B^{sm}$ given by a collection of closed $1$-forms $\op{Re}(t\alpha_i), 1\le i\le 2n$. For any simply-connected open subset $U\subset B^{sm}$ and a covariantly constant section $\gamma\in \Gamma(B^{sm},{\bf \Gamma})$ we have a closed $1$-form
$$\alpha_{\gamma,t}=\op{Re}\left(t\int_{\gamma}\omega_X^{2,0}\right)=d\op{Re}(tZ(\gamma))\,\,,$$
and the corresponding gradient vector field $v_{\gamma,t}=g_B^{-1}(\alpha_{\gamma,t})$.
 Notice that this vector field is a constant field with integral direction in the integral affine structure associated with closed $1$-forms $\op{Im}(t\alpha_i), 1\le i\le 2n$.

Similarly to \cite{KoSo1} we can construct  infinite oriented trees lying in $B$ such that its external vertices belong to $B^{sing}$, and edges are {\it positively oriented } trajectories of vector fields $v_{\gamma,t}$. All internal vertices have valency at least $3$, and  every such vertex should be thought of as a splitting point: a trajectory of the vector field $v_{\gamma,t}$ is split at a vertex into several trajectories of vector fields $v_{\gamma_1,t},\dots,v_{\gamma_k,t}$ such that $\gamma=\gamma_1+\dots+\gamma_k$.

The restriction of the function $Z$ to a tree gives rise to a $\C$-valued function such that on the trajectory of vector field $v_{\gamma,t}$ it is equal to the restriction of $Z(\gamma)$ to this trajectory. We assume that this function approach to zero as long as we approach an external vertex of the tree (which belongs to $B^{sing}$).
It is easy to see that  $tZ(\gamma)$ is a positive number at any other point of the tree (hence it defines a length function).
We expect that for any point $b\in B^{sm}$ and $\gamma\in \Gamma_b$ there exist finitely many such trees which pass the point $b$ in the direction of $\gamma$ (we can think of $b$ as a root of the tree, hence we can say above that we consider oriented trees such that all external vertices except of the root belong to $B^{sing}$). Here we choose an affine structure with $t\in \R_{>0}(Z(\gamma)_b^{-1})$.
Probably the number of such trees for fixed $b,\gamma$ is finite, since their lengths
should be bounded.\footnote{In \cite{KoSo1} we modified the gradient fields near $B^{sing}$ in order to guarantee the convergence of infinite products in the adic topology. It seems that we were too cautious, and the convergence holds without any modification.}

For a fixed $t\in \C^{\ast}$ the union $W_t$ of all trees as above is in fact a countable union of real hypersurfaces in $B^{sm}$. They are analogs of the walls of second kind. The set $W_t$ depends on $\op{Arg}t$ only. The union
$\cup_{\theta\in [0,2\pi i)}W_{te^{i\theta}}$ swap the whole space $B^{sm}$. Let us denote by $W^{(1)}$ the union over all $t\in \C^{\ast}/\R_{>0}$ of the sets of internal vertices of all trees in $W^{(1)}$ (splitting points of the gradient trajectories). This is an analog of the wall of first kind.

In \cite{KoSo1} we suggested a procedure of assigning rational multiplicities to edges of trees (see also \cite{GS1},\cite{GS2}).
 This leads to the following picture. Consider the total space $tot({\bf \Gamma})$ of the local system ${\bf \Gamma}$. It follows from above assumptions and considerations that we have a locally constant function $\Omega:tot({\bf \Gamma})\to \Q$ which jumps at the subset consisting of the lifts of the wall $W^{(1)}$ to $tot({\bf \Gamma})$.  Then for a fixed $b\in
B^{sm}$ the pair $(Z,\Omega)$ defines stability data on the graded Lie algebra $\g_{\Gamma_b}$ of the group of formal symplectomorphisms of the symplectic torus ${\mathbb T}_{\Gamma_b}$. In this way we obtain a local embedding $B^{sm}\hookrightarrow Stab(\g_{\Gamma_b})$.

In the above example of Seiberg-Witten curve, the wall $W^{(1)}$ is an oval-shaped curve which contains two singular points 
$b_{\pm}\in B^{sing}$. 
A typical $W_t$ looks such as follows.

\vspace{3mm}

\centerline{\epsfbox{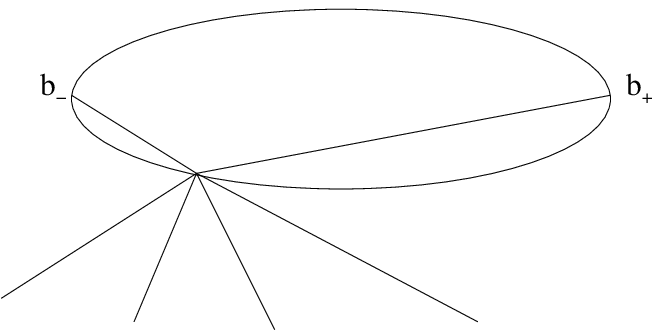}}

\vspace{3mm}

\vspace{3mm}

The wall-crossing formula coincides with the one for $T_{a,b}^{(2)}$ (see Introduction).

\begin{rmk} 1) We expect that the above considerations are valid for a large class of complex integrable systems, e.g. Hitchin system.

2) In the case when we have a $3d$ complex compact Calabi-Yau manifold $X$, the moduli space ${\cal M}_X$ of complex structures on $X$ is locally embedded into the projective space ${\bf P}(H^3(X,\C))$ as a base of a Lagrangian cone ${\cal L}_X\subset H^3(X,\C)$. It carries a K\"ahler metric (Weil-Petersson metric). We can repeat the above considerations given for integrable systems, replacing the gradient flows by the attractor flow (see e.g. \cite{DM}). The above case of integrable systems is obtained in the limit,  when the cone becomes ``very sharp".

\end{rmk}

\subsection{Relation with the works of Joyce, and of Bridgeland and Toledano-Laredo}

Let $\g, \Gamma$ be as in Section 2.1. We assume that the ground field is $\C$.
Suppose that
$C\subset \Gamma_{\R}$ is a strict convex cone. We are interested in such stability data $(Z,a)$ that
$\op{Supp}a\subset C\cup (-C)$.
We define $D$ as an open subset of $\op{Hom}(\Gamma, \C)$ which consists
of additive maps such that $C\cap \Gamma$ is mapped into the upper-half plane ${\cal{H}}_+=\{
z\in \C\,|\,\op{Im}(z)>0\}$.
We interpret $D$ as an open subset of $Stab(\g)$.
Every $\alpha\in C\cap {\Gamma}$ gives rise to an invertible
function (coordinate) $z_{\alpha}\in {\cal O}(D)^{\times}$ such that $z_{\alpha}(Z,a)=Z(\alpha)$.

Recall the pro-nilpotent Lie algebra $\g_C=\prod_{\gamma \in C\cap \Gamma}\g_{\gamma}$ and the corresponding pro-nilpotent group $G_C$.

In the paper \cite{Jo1} by D.~Joyce  the following system of differential
equations for a collection of holomorphic functions
$(f_{\alpha})_{\alpha \in C\cap {\Gamma}}, f_{\alpha}\in {\cal O}(D)\otimes \g_{\alpha}$ was considered:
$$\forall \alpha \in C\cap \Gamma\,\,\,\,\,\,\,df_{\alpha}=-{1\over{2}}\sum_{\beta+\gamma=\alpha}[f_{\beta},f_{\gamma}]\,\,d\log{z_{\beta}
\over{z_{\gamma}}}\,\,.$$

It follows that if $(f_{\alpha})$ satisfies the above system of equations then
the differential $1$-form 
$$\omega=\sum_{\alpha}f_{\alpha}d\log\,z_{\alpha}\in \Omega^1(D)\widehat{\otimes} {\g}_C:=\prod_{\alpha \in C\cap \Gamma}(\Omega^1(\g)\otimes \g_{\alpha})$$
 gives rise to the flat connection,
since $$d\omega+{1\over{2}}[\omega,\omega]=0\,\,.$$
 Moreover, setting $F=\sum_{\alpha}f_{\alpha}$ we observe that
$$dF+[\omega,F]=0\,\,,$$
i.e. $F$ is a flat section of this connection in the adjoint representation.

One can check by induction that there exists a unique solution to the above system of differential
equations (modulo constants for each function $f_{\alpha}$). This means that the set
of solutions is isomorphic to ${\g}_C$ (non-canonically).

For any $n\geqslant 0$ and pairwise different numbers $x_i\in \C\setminus \{0,1\},\,i=1,\dots,n$ we introduce the following function (multilogarithm) which is holomorphic when all $x_i$ lie outside of the interval $[0,1]$:
$$L_n(x_1,\dots,x_n):=v.p.\int_{0<t_1<t_2<\dots<t_n<1}\prod_{1\le i\le n}{dt_i\over{t_i-x_i}}\,\,,$$
where $v.p.$ means ``principal value".
Then $L_0=1,L_1(x)=\log(1-1/x)$ where we define the branch of the logarithm by taking the cut along the ray $(-\infty,0)$.

For a given collection $(f_{\alpha})$ as above, Joyce defined a collection of  functions ${\cal E}_{\alpha}$
on $D$ with values in the completed universal enveloping algebra $U(\g_C)$:
$${\cal E}_{\alpha}=\sum_{n\geqslant 1}\sum_{\alpha_1+\dots+\alpha_n=\alpha}f_{\alpha_1}\dots
f_{\alpha_n}I_n(z_{\alpha_1},z_{\alpha_2},\dots,z_{\alpha_n})\,\,,$$
where for  $z_1,\dots,z_n\in \C$ such that $0<\op{Im}z_1<\op{Im}z_2<\dots<\op{Im}z_n$ we set\footnote{This formula was proposed in \cite{BrTL}
 as the inversion of the Joyce formula which expressed $f_{\alpha}$'s in terms of ${\cal E}_{\alpha}$'s.}

$$\begin{array}{l}
I_n(z_1,z_2,\dots,z_n):=\\
{}=2\pi i(-1)^{n-1} L_{n-1}\left({z_1\over{z_1+\dots+z_n}},{z_1+z_2\over{z_1+\dots+z_n}},\dots,
{z_1+z_2+\dots+z_{n-1}\over{z_1+\dots+z_n}}\right)\,\,.
\end{array}$$

One can show that in fact ${\cal E}_{\alpha}\in \g_C$, and it is a locally constant along strata of the stratification defined by 
the walls $z_{\beta}/z_{\gamma}\in \R$ where $\alpha=\beta+\gamma$ with $ \beta,\gamma\in C\cap \Gamma$ and $\beta$ is not parallel to $\gamma$.

For a solution $(f_{\alpha})$ of the above system of differential equations we define
a differential $1$-form on $D\times \C^{\ast}$ such that
$$\widehat{\omega}:=\sum_{\alpha}f_{\alpha} e^{vz_{\alpha}}d\log(vz_{\alpha})\,\,,$$
where $v$ is the standard coordinate on $\C^{\ast}$.

Then one checks that 
$$d\widehat{\omega}+{1\over{2}}[\widehat{\omega},\widehat{\omega}]=0\,\,.$$
Let $M_{(f_{\alpha})}\in {G}_C$ be the monodromy of the corresponding flat connection computed
along a closed loop in the complex $v$-plane, which starts at $+i\infty$ and goes in the anti-clockwise direction around $v=0$. The flatness implies that the monodromy does not depend on the point of $D$.

On the other hand let us consider the element $N_{(f_{\alpha})}\in {G}_C$ defined as
$$\prod_{l\subset {\cal H}_+}^{\longrightarrow}\exp\left(\sum_{\alpha\in C\cap {\Gamma},z_{\alpha}\in l}{\cal E}_{\alpha}\right),$$
where the product is taken over all rays $l\subset {\cal H}_+$ with the vertex at the origin.

\begin{conj} We have $M_{(f_{\alpha})}=N_{(f_{\alpha})}$.

\end{conj}

The conjecture implies that the elements $a_{\alpha}:={\cal E}_{\alpha}$ satisfy the wall-crossing formula. The element $M_{(f_{\alpha})}$ is equal (in our notation) to the element
$A_V$, where $V$ is a strict sector in ${\cal H}_+$ containing $Z(C)$.
We will discuss below a sequence of identities which imply the conjecture. But we need to introduce certain functions first.

Let  $\varphi:(0,1)\to \C\setminus\{0\}$ be the infinite contour which starts and ends at $+i\infty$, goes in the anti-clockwise direction and surrounds the point $0\in \C$. With the contour $\varphi$ we associate the following function on  
$\left({\cal H}_+\right)^n,\,n\geqslant 1$:
$$K_n(z_1,\dots,z_n):=\int_{0<t_1<\dots<t_n<1}\exp\left(\sum_{1\le i\le n}\varphi(t_i)z_i\right)\prod_{1\le i\le n}\varphi^{\prime}(t_i)/\varphi(t_i)dt_i\,\,.$$
Notice that this function can be written as a Chen iterated integral

$$K_n(z_1,\dots,z_n)=\int_{\varphi}\omega_1\circ\omega_2\circ\dots\circ \omega_n\,\,,$$
where $\omega_i=e^{v z_i}dv/v, 1\le i\le n,\,\, v\in \C\setminus \{0\}$.

Let us fix $n\geqslant 1$  and a collection of complex numbers $z_i\in {\cal H}_+, \,1\le i\le n$. We call a sequence $0=i_0<i_1<\dots<i_{k-1}<i_k$ {\it admissible} if
$$\op{Arg}(z_1+\dots+z_{i_1})\geqslant  \op{Arg}(z_{i_1+1}+\dots+z_{i_2})\geqslant \dots\geqslant \op{Arg}(z_{i_{k-1}+1}+\dots+z_{i_k})\,\,.$$
For a fixed admissible sequence we have a partition $k=l_1+\dots+l_m$ where $l_1,l_2,\dots,l_m$ are the numbers of consecutive equalities in the above sequence of inequalities for the arguments.
Let $\Omega_{k,l_1,\dots,l_m}(z_1,\dots,z_n)$ be the set of all admissible sequences $0=i_0<i_1<\dots<i_{k-1}<i_k$
with the given partition $k=l_1+\dots+l_m$.
Under these assumptions and notation one can see that the previous Conjecture 2 is equivalent to

\begin{conj} We have
$$K_n(z_1,\dots,z_n)=\sum_{\Omega_{k,l_1,\dots,l_m}(z_1,\dots,z_n)}
\prod_{1\le j\le m}{1\over{l_j!}}{I}_{i_1}(z_1,\dots,z_{i_1})\cdot{I}_{i_2-i_1}(z_{i_1+1},\dots,z_{i_2})\cdot\dots$$
$$\cdot{I}_{i_k-i_{k-1}}(z_{i_{k-1}+1},\dots,z_{n})\,\,.$$

\end{conj}
Indeed, for $z_i=z_{\alpha_i},\,i=1,\dots,n$ the l.h.s. of the formula is the contribution of  the term 
$f_{\alpha_1}\dots f_{\alpha_n}$ in the expansion of $M_{(f_{\alpha})}$. Similarly, the r.h.s. is the contribution of the same
term in $N_{(f_{\alpha})}$.

Here we give a proof of the above conjecture in the special case:

\begin{prp} If $0<\op{Arg}\,z_1<\dots<\op{Arg}\,z_n<\pi$ then
$$K_n(z_1,\dots,z_n)=I_n(z_1,z_2,\dots,z_n)\,\,.$$

\end{prp}

{\it Proof.} For $n=1$ both sides are equal to $2\pi i$. For $n\geqslant 2$ we proceed by induction.
First one checks directly that
$$dK_n(z_1,\dots,z_n)=-\sum_{i=1}^{n-1}d\log\left(\frac{z_{i+1}}{z_i}\right)\, K_{n-1}(z_1,\dots, z_i+z_{i+1},\dots,z_n)\,\,.$$
The same formula holds if we replace $K_n,K_{n-1}$ by $I_n,I_{n-1}$ respectively. Thus we see by induction that
$K_n-I_n=const_n$. We want to prove that $const_n=0$. In order to do that we take $z_j=z_j(\varepsilon), 1\le j\le n$,
such as follows: 
$$z_1(\varepsilon)=i+{1\over {\varepsilon}},\, z_n(\varepsilon)=i-{1\over {\varepsilon}},\,z_k(\varepsilon)=i-k, \,
2\le k\le n-1\,\,.$$ Here $i=\sqrt{-1}$. Then 
$$0<\op{Arg}\,z_1(\varepsilon)<\op{Arg}\,z_2(\varepsilon)<\dots<\op{Arg}\,z_n(\varepsilon)<\pi$$ and
$|\sum_{1\le j\le k}z_j(\varepsilon)|\to \infty$ as $\varepsilon\to 0$ for $k=1,\dots,n-1$, and moreover $|\sum_{1\le j\le n}z_j(\varepsilon)|$ is a constant function of $\epsilon$. Therefore, $$I_n(z_1(\varepsilon),\dots, z_n(\varepsilon))\to 0$$ as $\varepsilon\to 0$,
since all the arguments of the function $L_{n-1}$ in the definition of $I_n$ approach infinity.

Hence in order to finish the proof it suffices to show that 
$$K_n(z_1(\varepsilon),\dots, z_n(\varepsilon))\to 0$$ as $\varepsilon\to 0$.
Here is the sketch of the proof.\footnote{We thank Andrei Okounkov for the idea of the proof.}
Notice that $$\int_{+i\infty}^{v_2}e^{v_1z_1(\varepsilon)}dv_1/v_1={1\over{z_1(\varepsilon)}}e^{v_1z_1(\varepsilon)}/v_2+
r_1(\varepsilon)\,\,,$$ where $r_1(\varepsilon)\to 0$ as $\varepsilon\to 0$. Repeating we obtain that
$$K_n(z_1(\varepsilon),\dots, z_n(\varepsilon))={1\over{z_1(\varepsilon)}}
{1\over{z_1(\varepsilon)+z_2(\varepsilon)}}\dots{1\over{z_1(\varepsilon)+z_2(\varepsilon)+\dots+z_{n-1}(\varepsilon)}}\times$$
$${}\times\int_{\varphi}e^{v_n(z_1(\varepsilon)+\dots+z_n(\varepsilon))}dv_n/v_n+r_n(\varepsilon)\,\,,$$
where the integral is taken over the contour $\varphi$ described before, and $r_n(\varepsilon)\to 0$ as $\varepsilon\to 0$.
It follows from our choice of numbers $z_j(\varepsilon), 1\le j\le n$ that
$K_n(z_1(\varepsilon),\dots, z_n(\varepsilon))\to 0$ as $\varepsilon\to 0$. $\blacksquare$

%FIGURE 2 (Integration contour)

%\centerline{\epsfbox{contour-bis.eps}}

One can hope that the technique developed in \cite{Jo1} helps in proving the general case.

A  relationship  between Joyce formulas and iterated integrals is discussed in \cite{BrTL} in a slightly different form. In that paper
the elements $N_{(f_{\alpha})}$ are interpreted as Stokes multipliers for a different system of differential
equations on $\C$ (with coordinate $t$) with values in the Lie algebra which is an extension of $\g_C$ by the abelian Lie algebra $\op{Hom}(\Gamma,\C)$ (an analog of the Cartan subalgebra).  It has irregular singularity at the origin given by $Z\over{t^2}$, where $Z$ is the central charge of the stability structure. In fact the connection from \cite{BrTL} reduces to our connection after
the change of variables $v=1/t$ and the conjugation by $\exp(-vZ)$.

\subsection{Stability data on ${\mathfrak{gl}}(n)$}

Let $\g={\mathfrak{gl}}(n,\Q)$ be the Lie algebra of the general linear group. We consider it as a $\Gamma$-graded Lie algebra $\g=\oplus_{\gamma \in \Gamma}\g_{\gamma}$,  where $$\Gamma=\{(k_1,\dots,k_n)|\,\,k_i\in \Z, \sum_{1\le i\le n}k_i=0\}$$ is the root lattice. We endow $\g$ with the Cartan involution $\eta$. Algebra $\g$ has 
 the standard basis $E_{ij}\in \g_{\gamma_{ij}}$  consisting of matrices with the single non-zero entry at the place $(i,j)$ equal to $1$. Then $\eta(E_{ij})=-E_{ji}$.
In what follows we are going to consider symmetric (with respect to $\eta$) stability data on $\g$.

We notice that  $$\op{Hom}(\Gamma,\C)\simeq \C^n/\C\cdot(1,\dots,1)\,\,.$$ We define a subspace $\op{Hom}^\circ(\Gamma,\C)\subset \op{Hom}(\Gamma,\C)$ consisting (up to a shift by the multiples of the vector $(1,\dots,1)$) of vectors $(z_1,\dots,z_n)$ such that $z_i\ne z_j$ if $i\ne j$. Similarly we define a subspace $\op{Hom}^{\circ\circ}(\Gamma,\C)\subset \op{Hom}(\Gamma,\C)$ consisting (up to the same shift) of such $(z_1,\dots,z_n)$ that there is no $z_i,z_j,z_k$ belonging to the same real line as long as $i\ne j\ne k$. Obviously there is an inclusion $\op{Hom}^{\circ\circ}(\Gamma,\C)\subset \op{Hom}^\circ(\Gamma,\C)$.

For $Z\in \op{Hom}(\Gamma,\C)$ we have $Z(\gamma_{ij})=z_i-z_j$. If $Z\in \op{Hom}^{\circ\circ}(\Gamma,\C)$ then symmetric stability data with such $Z$ is the same as a skew-symmetric matrix $(a_{ij})$ with rational entries determined from the equality $a(\gamma_{ij})=a_{ij}E_{ij}$. Every continuous path in $\op{Hom}^\circ(\Gamma,\C)$ admits a unique lifting to $Stab(\g)$ as long as we fix the lifting of the initial point. The matrix $(a_{ij})$ changes when we cross walls in  $\op{Hom}^\circ(\Gamma,\C)\setminus \op{Hom}^{\circ\circ}(\Gamma,\C)$. A typical wall-crossing corresponds to the case when in the above notation the point $z_j$ crosses a straight segment joining $z_i$ and $z_k, i\ne j\ne k$. In this case the only change in the matrix $(a_{ij})$ is of the form:
$$a_{ik}\mapsto a_{ik}+a_{ij}a_{jk}\,\,.$$
This follows from the multiplicative wall-crossing formula which is of the form:
$$\begin{array}{l}\exp(a_{ij}E_{ij})\exp(a_{ik}E_{ik})\exp(a_{jk}E_{jk})=\\
\qquad{}=\exp(a_{jk}E_{jk})\exp((a_{ik}+a_{ij}a_{jk})E_{ik})\exp(a_{ij}E_{ij})\,\,.
\end{array}$$
Same wall-crossing formulas appeared in \cite{CeVa} in the study of the change of the number of solitons in $N=2$ two-dimensional supersymmetric QFT. In \cite{CeVa} the numbers $a_{ij}$ were integers, and the wall-crossing preserved integrality. In our considerations, for any $Z\in \op{Hom}^{\circ\circ}(\Gamma,\C)$ the fundamental group
$\pi_1(\op{Hom}^\circ(\Gamma,\C),Z)$ acts on the space of skew-symmetric matrices by polynomial transformations with integer coefficients. It can be identified with the well-known actions of the pure braid group on the space of upper-triangular matrices in the theory of Gabrielov bases of isolated singularities and in the theory of triangulated categories endowed with exceptional collections. Furthermore, the matrices $\exp(a_{ij}E_{ij})=1+a_{ij}E_{ij}$ can be interpreted as Stokes matrices of a certain connection in a neighborhood of $0\in \C$, which has irregular singularities ($tt^{\ast}$-connection from \cite{CeVa}, see also \cite{GaMNeit}). This observation should be compared with the results about the irregular connection from the previous section.

\section{Ind-constructible categories and stability structures}

\subsection{Ind-constructible categories}

Here we introduce an ind-constructible version of the notion of a (triangulated)  $\A$-category.
Let $\k$ be a field, $\overline{\k}$ be its algebraic closure. By a variety over $\k$ (not necessarily irreducible)
we mean a reduced separated scheme of finite type over $\k$. 
Recall  the following definition.

\begin{defn} Let $S$ be a variety over $\k$. A subset $X\subset S(\kk)$ is called constructible over $\k$ 
if it belongs to the Boolean algebra generated by $\kk$-points of open (equivalently closed) subschemes of $S$.

\end{defn}

Equivalently, a constructible set is the union of a finite collection  of $\kk$-points of
disjoint locally closed subvarieties $(S_i\subset S)_{i\in I}$.

For any field extension $\k\subset \k^{\prime}\subset \kk$ we define the set of $\k^{\prime}$-points $X(\k')$ of the constructible set $X$
as $(X\cap S(\k'))\subset S(\kk)$. In particular, $X(\kk)=X$.

We define the category ${\cal CON}_{\k}$ of constructible sets over $\k$ as a category with objects
$(X,S)$, where $X$ and $S$ as above. The set of morphisms $\op{Hom}_{{\cal CON}_{\k}}((X_1,S_1),(X_2,S_2))$ is defined to be the set of maps  
 $f:X_1\to X_2$ such that there exists a decomposition of $X_1$ into the finite disjoint union of $\kk$-points of 
varieties $(S_i\subset S_1)_{i\in I}$ such that
 the restriction of $f$ to each $S_i(\kk)$ is a morphism of schemes $S_i\to S_2$.
We see that there is a natural faithful functor from ${\cal CON}_{\k}$ to the category of sets equipped with
the action of $\op{Aut}(\overline{\k}/\k)$.

\begin{defn} An ind-constructible set over $\k$ is given by a chain of embeddings of constructible sets
$X:=(X_1\to X_2\to X_3\to\dots)$ over $\k$. A morphism of ind-constructible sets  is defined as
a map $g:\cup_iX_i(\overline{\k})\to \cup_iY_i(\overline{\k})$ such that for any $i$ there exists $n(i)$ such  that
$g_{|X_i(\overline{\k})}:X_i(\overline{\k})\to Y_{n(i)}(\overline{\k})$ comes from a constructible map.

\end{defn}

Ind-constructible sets form a full subcategory ${\cal IC}_{\k}$ of the category of ind-objects
in ${\cal CON}_{\k}$.

\begin{rmk} Equivalently, we can consider a countable collection $Z_i=X_i\setminus X_{i-1}$ of non-intersecting
constructible sets. Then a morphism $\sqcup_{i\in I}Z_i\to \sqcup_{j\in J}Z_j^{\prime}$ is given by a collection of constructible
maps $f_i:Z_i\to \sqcup_{i\in J_i}Z_j^{\prime}$, where each $J_i$ is a finite set.

\end{rmk}

The category of constructible (or ind-constructible) sets has fibered products.
 There is a notion
of constructible (or ind-constructible) vector bundle (i.e. the one with the fibers which are affine spaces of various finite dimensions).

\begin{defn} An  ind-constructible  $\A$-category  over $\k$ 
is defined by the following data:

1) An ind-constructible set ${\cal M}=Ob(\CC)=\sqcup_{i\in I}X_i$ over $\k$, called the set of objects.

2) A collection of ind-constructible vector bundles $${\cal HOM}^n\to {\cal M}\times {\cal M}, \,n\in \Z$$
called the bundles of morphisms of degree $n$.  The restriction ${\cal HOM}^n\to X_i\times X_j$ is  a
 finite-dimensional constructible vector bundle for any $n\in \Z, i,j\in I$,
and the restriction
${\cal HOM}^n\to X_i\times X_j$ is a zero bundle for $n\le C(i,j)$, where $C(i,j)$
is some constant.

3) For any $n\geqslant 1, l_1,\dots,l_n\in \Z$,  ind-constructible  morphisms  of ind-constructible bundles
$$m_n:p_{1,2}^{\ast}{\cal HOM}^{l_1}\otimes\dots\otimes p_{n,n+1}^{\ast}{\cal HOM}^{l_n}\to p_{1,n+1}^{\ast}{\cal HOM}^{l_1+\dots+l_n+2-n}\,\,,$$
where $p_{i,i+1}$ denote natural projections
of ${\cal M}^{n+1}$ to ${\cal M}^2$. These morphisms are called higher composition maps.

\end{defn}
The above data are required to satisfy the following axioms A1)-A3):

{\it A1) Higher associativity property for $m_n, n\geqslant 1 $ in the sense of $\A$-categories.} We leave for the reader to write down the corresponding well-known identities (see \cite{Kel},\cite{KoSo3}).

This axiom implies that we have a small $\kk$-linear non-unital $A_\infty$-category $\CC(\kk)$ with the set of objects ${\cal M}(\kk)$ and morphisms
 ${\cal HOM}^\bullet(\kk)$. 

{\it A2) (weak unit) There is a constructible section $s$ of the ind-constructible bundle
${\cal HOM}^0_{|Diag}\to {\cal M}$ such that the image of $s$ belongs to
the kernel of $m_1$ and gives rise to the identity morphisms in $\Z$-graded $\kk$-linear category $H^{\bullet}(\CC(\kk))$.}

Alternatively, instead of A2) one can use the axiomatics of $A_\infty$-categories with strict units (see \cite{KoSo3}, \cite{LyuMan}).

An ind-constructible $\A$-category $\CC$ gives rise to a collection of ind-constructible bundles over $Ob(\CC)\times Ob(\CC)$ given by  $${\cal EXT}^i:=H^{i}({\cal HOM}^\bullet),\, i\in \Z$$ with the fiber over a pair of objects $(E,F)$ equal to 
$$\op{Ext}^i(E,F):=H^{i}({\cal HOM}_{E,F}^\bullet)\,\,.$$ The cohomology groups are taken with respect to the differential $m_1$.

 {\it A3) (local regularity) There exists a countable collection of schemes $(S_i)_{i\in I}$ of finite type over $\k$,  
 a collection of algebraic $\k$-vector
bundles $\op{HOM}^n_{i}, n\in \Z$ over $S_i\times S_i$ for all $i$,
and  ind-constructible identifications 
$$\sqcup_i S_i(\kk)\simeq {\cal M}, \op{HOM}^n_{i}\simeq {\cal HOM}_{|S_i\times S_i}^n, n\in \Z$$
 such that all higher compositions $m_n,n\geqslant 2$,
considered for objects from $S_i$ for any given $i\in I$, become morphisms of algebraic vector bundles.}

We will often call ind-constructible $\A$-categories simply by ind-constructible categories. 
The basic example of an ind-constructible category is the category $Perf(A)$ of perfect $A$-modules where $A$ an $A_\infty$-algebra over $\k$ with finite-dimensional cohomology
(see the discussion after the Example 1 in Section 1.2 of Introduction).

We define a {\it functor} between two ind-constructible categories mimicking the usual definition of an $\A$-functor.

A functor $\Phi:{\cal C}_1\to {\cal C}_2$ is called an {\it equivalence} if $\Phi$ is a full embedding,  i.e. it
induces an isomorphism
$$\op{Ext}^\bullet(E,F)\simeq \op{Ext}^\bullet(\Phi(E),\Phi(F))\,\,\forall E,F\in Ob(\CC_1)(\kk)$$
 and moreover, there exists an   ind-constructible {\it over }$\kk$ map $$s:Ob(\CC_2)(\kk)\to Ob(\CC_1)(\kk)$$ such
that for any object $E\in Ob(\CC_2)(\kk)$ we have $E\simeq \Phi(s(E))$.

 Using the notions of a functor and of an  equivalence we can define the
  property of an ind-constructible weakly unital $\A$-category $\cal C$ to be triangulated.
  For example, the property to have exact triangles can be formulated as follows.
  Consider a finite $\A$-category ${\cal C}_3$ consisting of $3$ objects $E_1,E_2,E_3$ with non-trivial morphism spaces 
  $$\op{Hom}^0(E_i,E_i)=\k\cdot \op{id}_{E_i},\,\op{Hom}^0(E_1,E_2)\simeq \op{Hom}^0(E_2,E_3)\simeq \op{Hom}^1(E_3,E_1)\simeq \k$$
  equivalent to the full subcategory of the category of representations of the quiver $A_2$ 
  consisting of modules of dimensions $(0,1),(1,1),(1,0)$.
  Let ${\cal C}_2\subset {\cal C}_3$ be the full subcategory consisting of first two objects.
  
  It is easy to see directly from the definitions that for any ind-constructible 
  category $\cal C$ there are natural ind-constructible categories $\op{Fun}({\cal C}_i,{\cal C}),\,i=2,3$ whose objects over $\overline{\k}$ are the usual $\A$-categories of functors
   from ${\cal C}_i(\overline{\k})$ to ${\cal C}(\overline{\k})$ as defined e.g. in \cite{KoSo2},\cite{KoSo3} and \cite{Kel}. There exists a natural restriction functor
   $$r_{32}:\op{Fun}({\cal C}_3,{\cal C}) \to\op{Fun}({\cal C}_2,{\cal C})\,.$$
   Similarly to the setting of usual $\A$-categories,  the ind-constructible version of the axiom of exact triangles says that $r_{32}$
   is an equivalence. In the same manner one can define other properties of triangulated $\A$-categories (i.e.  the existence of shift functors, finite sums,  see \cite{KoSo3}, \cite{So1}) in the ind-constructible setting.

In Sections 5,6 we will use a simplified notation $Cone(f)$ for a cone of morphism $f$ in $\CC(\kk)$ ``pretending'' that cones are 
functorial. The precise prescription is to 
 take an object in $\op{Fun}({\cal C}_2,{\cal C})$ corresponding to $f$, find an isomorphic object in $\op{Fun}({\cal C}_3,{\cal C})$, and then take
the image in $\CC(\kk)$ of the object $E_3$. All this can be properly formulated using the language of constructible stacks, see 3.2 and 4.2.

Let us call an ind-constructible $\A$-category {\it minimal on the diagonal} if the restriction
of $m_1$ to the diagonal $\Delta\subset {\cal M}\times {\cal M}$ is trivial.
One can show that any ind-constructible $\A$-category is equivalent to a one which is minimal on the diagonal.

\begin{rmk} Typically in practice one has a decomposition $Ob(\CC)=\sqcup_{i\in I}X_i$ where $X_i$ are schemes, not just constructible sets. Moreover, for any $E\in X_i(\kk)$ there is a natural map $T_E X_i\to \op{Ext}^1(E,E)$. The reason for this is the fact that the deformation theory of the object $E$ should be controlled by the DGLA $\op{Ext}^{\bullet}(E,E)$. We did not include the above property into the list of axioms since it does not play any role in our constructions.
\end{rmk}

\subsection{Stack of objects}

In this section we assume that the ground field $\k$ is perfect, i.e. the $\kk$ is a Galois field over $\k$.
Our goal in this section is to explain how to associate with an 
ind-constructible $\A$-category $\CC$ over $\k$ an ordinary $\k$-linear $\A$-category $\CC(\k)$, in such a way
 that ind-constructible equivalences will induce the usual equivalences.
For any field extension $\k'\supset \k$ (e.g. for $\k'=\k$) one can define $\A$-category $\CC^{naive}(\k')$
to be the small $\k'$-linear category with the set of objects given by $(Ob(\CC))(\k')$ and obvious morphisms and compositions.
This is not a satisfactory notion because
 in the definition of the equivalence we demand only the surjectivity on isomorphism classes of objects  over $\kk$.
The naive category $\CC^{naive}(\k)$ will be a full subcategory of the ``correct'' category $\CC(\k)$.
 One should read carefully brackets, as in our notation
 $$Ob(\CC(\k))\ne (Ob(\CC))(\k)=:Ob(\CC^{naive}(\k))\,,$$
contrary to the case of $\kk$ where we have 
 $$Ob(\CC(\kk))=(Ob(\CC))(\kk)\,.$$
We will see also the ``set of isomorphism classes of objects'' in $\CC$ should be better understood as an ind-constructible stack\footnote{Even in the case when $\k=\kk$ it is important to keep track on automorphisms groups of objects (and not only on the set of isomorphism classes), e.g. for the correct definition of the motivic Hall algebra in 6.1.}.

Let $\k'\subset \kk$ be a finite Galois extension of $\k$ and consider an element $E'\in (Ob(\CC))(\k')\subset Ob(\CC(\kk))$ such that $\sigma(E)$ is isomorphic to 
$E$ for all $\sigma \in \op{Gal}(\k'/\k)$. We would like to define the descent data for such $E'$, which should be data necessary to define
an object
in (not yet defined) $\k$-linear category $\CC(\k)$ which becomes isomorphic to $E'$ after the extension of scalars from $\k$ to $\k'$.

First, for a finite non-empty collection of objects $(E_i)_{i\in I}$ of any $\A$-category $\CC'$ linear over a field ${\k}'$ (not necessarily a perfect one)
we define an {\it identification data} for objects of this collection to be 
an $\A$-functor $\Phi$ from the $\A$-category $\CC_{I,{\k}'}$ describing $I$ copies of the same object: 
$$Ob(\CC_I)=I, \op{Hom}^\bullet(i,j) =\op{Hom}^0(i,j)\simeq {\k}'  $$
 to $\CC'$. In plain terms, to give such a functor is to give a closed morphism of degree 0 for any pair of objects $E_i,E_j$ (representing  the identity
 $\op{id}_{E_i}$ in $H^\bullet(\CC(\kk))$ for $i=j$), a homotopy for any triple of objects, homotopy between homotopies
 for any quadruple of objects, etc. Thus, we in a sense identify all the objects of the collection $(E_i)_{i\in I}$ and hence can treat it is a
 new object (canonically isomorphic to all $(E_i)_{i\in I}$), without choosing any specific element $i\in I$.

Returning to the case of $E'\in (Ob(\CC))(\k')$, we define the descent data
 as the identification of the collection of objects $(\sigma(E'))_{\sigma\in \op{Gal}(\k'/\k)}$ of the category $\CC^{naive}(\k')$
equivariant with respect to the action of $\op{Gal}(\k'/\k)$ acting both on the collection and on the coefficients in the identification.
 
We define the set $Ob(\CC(\k))$ of objects of $\CC(\k)$ to be the inductive limit over finite Galois extensions $\k'/\k$ of descent data as above.
Also one can define morphisms and higher compositions. 
We leave the following Proposition without a proof.
\begin{prp} There is a natural structure of a $\k$-linear $\A$-category on $\CC(\k)$ containing $\CC^{naive}(\k)$ as a full subcategory.
 Any equivalence $\Phi:\CC_1\to \CC_2$ in ind-constructible sense induces an equivalence $\CC_1(\k)\to \CC_2(\k)$.
 If $\CC$ is triangulated in ind-constructible sense then $\CC(\k)$ is also triangulated.
\end{prp}

If $E$ is an object of $\CC(\k)$ then any other object $E'$ of $\CC(\k)$ which is isomorphic to $E$ after the extension of scalars to $\CC(\kk)$
is in fact isomorphic to $E$ in $\CC(\k)$ (in other words, there are no non-trivial twisted forms).
The reason is that (as follows directly from definitions) the set of such ``$\k$-forms" of $E$
is classified by $H^1(\op{Gal}(\overline{\k}/\k),{\cal G}_E)$, where ${\cal G}_E$ is a simplicial group
associated with the $\A$-algebra $\op{End}^\bullet_{\CC(\k)}(E,E)$.
There is a spectral sequence
which converges to this set and has the second term $E_2=(E_2^{pq})$  given by
$$\begin{array}{l}
H^1(\op{Gal}(\overline{\k}/\k),(\op{Ext}_{\CC({\kk})}^0(E,E))^{\times})\,\,,\\ H^2(\op{Gal}(\overline{\k}/\k),\op{Ext}_{\CC({\kk})}^{-1}(E,E))\,\,,\\
H^3(\op{Gal}(\overline{\k}/\k),\op{Ext}_{\CC({\kk})}^{-2}(E,E)),\dots\,\,.
\end{array}$$
We observe that all Galois cohomology groups with coefficients in $\op{Ext}_{\CC({\kk})}^{<0}(E,E)$ are trivial
(since  $\op{Ext}_{\CC({\kk})}^{i}(E,E), i<0$  are just sums of copies of the additive group ${\bf G}_a(\kk)$). Also the set 
$H^1(\op{Gal}(\overline{\k}/\k),(\op{Ext}_{\CC({\kk})}^0(E,E))^{\times})$ is the one-element set, because for any finite-dimensional algebra
 $A$ over $\k$ we have $H^1(\op{Gal}(\overline{\k}/\k),A^\times)=0$ (a version of Hilbert 90 theorem, see also section 2.1 in \cite{Jo4}). 
One can deduce from the above spectral sequence an important corollary: the set of isomorphism classes of objects of $\CC(\k)$ is in a natural bijection 
 with the set of isomorphism classes of the usual descent data in  category $H^{\bullet}({\CC(\kk)})$ endowed with the strict action of $\op{Gal}(\kk/\k)$.

Finally, we will explain how to associate an ind-constructible stack to an ind-constructible category $\CC$ over $\k$.
First of all, we can always assume that  $\CC$ satisfies the following axiom

{\it A4) There exists a decomposition ${\cal M}=\sqcup_{i\in I}X_i$ into the countable disjoint union of constructible sets over $\k$ such that
  any  two isomorphic objects of $H^{\bullet}(\CC(\kk))$  belong to the same part $X_i(\kk)$ for some $i\in I$.}

Indeed, if we choose any decomposition ${\cal M}=\sqcup_{i\in I}X'_i$ into disjoint union of constructible sets over $\k$ and identify
$I$ with the set of natural number $\{1,2,\dots\}$, then we can shrink $X'_i(\kk)$ to the subset consisting of objects which are not isomorphic to objects
 from $\cup_{j<i}X'_j(\kk)$.

By axioms A3),A4) we may assume that $Ob(\CC)(\kk)$ is decomposed into the union of $\kk$-points of schemes $(S_i)_{i\in I}$ (as in Axiom A3)) such that any isomorphism class
 of $\CC(\kk)$ belongs only to one of the sets $S_i(\kk)$.  Let us call it a locally regular subdivision of $Ob(\CC)$. 
Moreover, we can assume
that all $S_i$ are smooth and equidimensional, and such that for any $i\in I$ there exists $\delta(i)\in \Z_{\geqslant 0}$
such that for  any $E\in S_i(\kk)$ the subset
 of objects in $S_i(\kk)$ isomorphic to to $E$ has dimension $\delta(i)$. This can be achieved by subdividing each $S_i$ into smaller pieces, 
 and by removing some unnecessary pieces consisting of objects which belong to other pieces.
Then taking a generic slice of codimension $\delta(i)$ (and thus shrinking $\CC$ to an equivalent subcategory),
 and taking further subdivisions,
 one may assume that we have a locally regular subdivision of $Ob(\CC)$ such that any isomorphism class of objects in $S_i(\kk)$ is finite.
Moreover, we may assume that the cardinality $c_i$ of all isomorphism classes in $S_i(\kk)$ depend only on $i$, and also the dimension $d_i$
 of the algebra $\op{Ext}^0(E,E)$ for $E\in S_i(\kk)$ also depends only on $i$.

For any given $i\in I$ let us consider the  constructible set $Z_i$ over $\k$ parametrizing isomorphism classes of objects in $S_i(\kk)$. There is a natural 
constructible (over $\k$) bundle of finite-dimensional unital associative algebras $\cal A$, with the fiber ${\cal A}_x$ over any full collection $x=(E_1,\dots, E_{c_i})$
 (up to permutation) of different isomorphic objects equal to 
$$\oplus_{1\le j_1,j_2\le c_i} \op{Ext}^0(E_{j_1},E_{j_2})\,.$$
The above algebra is Morita equivalent to $\op{Ext}^0(E_j,E_j)$ for every $j\le c_i$, and in fact is isomorphic to the matrix algebra 
 $${\cal A}_x\simeq \op{Mat}(c_i\times c_i, \op{Ext}^0(E_j,E_j))\,\,\,\forall j\le c_i\,.$$
Informally speaking, the ``stack'' of objects from $S_i$
 is the stack of projective modules $M$ over algebra ${\cal A}_x$ for some $x\in Z_i(\kk)$ which are isomorphic after Morita equivalence to
 a free module of rank one over $\op{Ext}^0(E_j,E_j)$ where $E_j$ is some representative of the equivalence class $x$, i.e. $M$ is isomorphic
to the standard module $(\op{Ext}^0(E_j,E_j))^{\oplus c_i}$ over the matrix algebra for every $j\le c_i$. We see
 that $M$ has dimension 
$$N_i:=c_i d_i$$ over $\kk$. This leads to the following construction. Define a constructible set $Y_i$ over $\k$ to be the
 set of pairs $(x,f)$ where $x\in Z_i(\kk)$ is a point and $f$ is a homomorphism of $A_x$ to the algebra of matrices $\op{Mat}(N_i\times N_i,\kk)$
 such that the resulting structure of ${\cal A}_x$-module on $\kk^{N_i}$ belongs to the isomorphism class of projective ${\cal A}_x$-modules discussed above.
 The group $GL(N_i,\kk)$ acts naturally on $Y_i$ by changing the basis in the standard coordinate space $\kk^{N_i}$.
 The quotient set is naturally identified with $Z_i(\kk)$, and the stabilizer of every point is isomorphic to  $\op{Ext}^0(E_j,E_j)^\times$ in the above notation.

The essential element of the presented construction is that everything is equivariant with respect to the action of $\op{Gal}(\kk/\k)$.
Hence, we come to the conclusion that one associates (making many choices) with an ind-constructible category
$\CC$ over $\k$ a countable collection of varieties $(Y_i)_{i\in I}$ (we can assume that $Y_i$ are not just constructible sets but varieties
after making further subdivisions) endowed with algebraic actions of affine algebraic groups $GL(N_i)$ such that the groupoid of isomorphism classes
 of $\CC(\kk)$ is naturally equivalent to the groupoid of the disjoint union of sets $Y_i(\kk)$ endowed with $GL(N_i,\kk)$-actions.
If we replace $\CC$ by an equivalent ind-constructible category, or make different choices in the construction, we obtain an equivalent in an obvious sense 
``ind-constructible stack''.  We will discuss ind-constructible stacks later, in Section 4.2.

Moreover, using the fact that the first Galois cohomology with coefficients in $GL(N_i)$ vanish, one can see that the same is true for $\CC(\k)$ (and replace $Y_i(\kk)$, $GL(N_i,\kk)$-actions by $Y_i(\k)$ and $GL(N_i,\k)$-action $\forall i\in I$).
In general, for any field $\k'$, $\,\k\subset \k'\subset \kk$ one can define the descent data for $\k'$ and a $\k'$-linear $\A$-category $\CC(\k')$ (which is triangulated 
if $\CC$ is triangulated in the ind-constructible sense). The groupoid of isomorphism classes of objects of $H^\bullet(\CC(\k'))$ is 
equivalent to the groupoid of the disjoint union of sets $Y_i(\k')$ endowed with $GL(N_i,\k')$-actions.
More generally, one can define the category $\CC(\k')$ for {\it any} field extension $\k'\supset \k$, not necessarily an algebraic one.
In the case $\k'=\kk$ we get a non-fatal crash of notations, because the $\A$-category $\CC(\kk)$ in last sense is equivalent to the previously
 defined $\CC(\kk)$.

In what follows, we will assume for convenience that $Ob(\CC)$ for an ind-constructible category $\CC$ is described by schemes $Y_i$ with $GL(N_i)$-actions.
In particular, for any extension $\k'\supset \k$ we will have a bijection
 $$Iso(\CC(\k'))\simeq \sqcup_{i \in I} Y_i(\k')/GL(N_i,\k')$$
between the set of isomorphism classes in $\CC(\k')$ and the set of orbits.

\begin{rmk} In fact, objects of an $\A$-category form not a stack but a higher stack, i.e. one should speak about isomorphisms between isomorphisms etc.
Passing to the level of ordinary stacks we make a truncation. Presumably, for a proper treatment of ind-constructible categories
 and problems like non-functoriality of cones, one should introduce higher constructible stacks. Looking on the guiding example of 
 identification data for finite non-empty collections, one can guess  an appropriate notion of a higher constructible stack.
 Namely, it should be a simplicial constructible set $X_\bullet$ which satisfies a constructible version of the Kan property
(i.e. there exists a constructible lifting from horns to simplices) and such that 

1) for any $k\geqslant 2$ the constructible map $(\partial_0,\dots,\partial_k):X_k\to (X_{k-1})^{k+1}$ is a constructible vector bundle over its image 
(i.e. there exists a constructible identification of non-empty fibers of this map with vector spaces).

2)  $\exists \,k_0$ such that $\forall k\geqslant k_0$ the above map is an inclusion.

The reason for the first property is that in the case of identification on each step (except first two) we have to solve linear equations.
The second property comes from the property ${\cal HOM}^n_{|X_i\times X_i}=0$ for $n\ll 0$ in our axiomatics of ind-constructible categories.

\end{rmk}

\subsection{Ind-constructible Calabi-Yau categories and potentials}

Let $\k$ be a field of characteristic zero.
Recall that a  Calabi-Yau category of dimension $d$ is a weakly unital $\k$-linear triangulated $\A$-category
${\CC}$ (see \cite{KoSo2}, \cite{KoSo3}, \cite{So1}), such that for any two objects $E,F$ the $\Z$-graded vector space
$\op{Hom}^\bullet(E,F)=\oplus_{n\in \Z}\op{Hom}^n(E,F)$ is finite-dimensional
(hence the space $\op{Ext}^\bullet(E,F)$ is also finite-dimensional)
and moreover:

1) We are given a  non-degenerate pairing 
$$(\bullet, \bullet): \op{Hom}^{\bullet}(E,F)\otimes \op{Hom}^{\bullet}(F,E)\to \k[-d]\,\,,$$
which is symmetric with respect to interchange of objects $E$ and $F$;

2) For any $N\geqslant 2$ and a sequence of objects $E_1,E_2,\dots,E_N$ we are given
a polylinear $\Z/N\Z$-invariant map
$$W_N: \otimes_{1\le i\le N}\left(\op{Hom}^{\bullet}(E_i,E_{i+1})[1]\right)\to \k[3-d]\,\,,$$
where $[1]$ means the shift in the category
of $\Z$-graded vector spaces, and we set $E_{N+1}=E_1$;

3) We have:
$$W_N(a_1,\dots,a_{N})=(m_{N-1}(a_1,\dots,a_{N-1}),a_N)\,\,,$$
where $m_n: \otimes_{1\le i\le n}\op{Hom}^{\bullet}(E_i,E_{i+1})\to \op{Hom}^{\bullet}(E_1,E_{n+1})[2-n]$ are higher
composition maps.

The collection $(W_N)_{N\geqslant 2}$ is called the {\it potential of $\CC$}.
If $d=3$ then for any object $E\in Ob(\CC)$ we define a formal series $W_{E}^{tot}$ at $0\in \op{Hom}^{\bullet}(E,E)[1]$
by the formula:
$$W_E^{tot}(\alpha)=\sum_{n\geqslant 2}{W_n(\alpha,\dots, \alpha)\over{n}}\,\,.$$
We call $W_E^{tot}$ the {\it total (or full) potential of the object $E$}. We call the {\it potential of $E$}
the restriction of $W_E^{tot}$ to the subspace $\op{Hom}^1(E,E)$. We will denote it by $W_E$. 

The notion of a Calabi-Yau category admits a natural generalization to the ind-constructible case (the pairing is required to be a morphism of constructible  vector bundles). It follows from the Axiom A3)  that there exists a decomposition of $Ob(\CC)\simeq \sqcup S_i$ into the disjoint union of schemes such that all Taylor components $W_N$ of the potential are symmetrizations of regular sections of cyclic powers of algebraic vector bundles on
schemes $S_i$. Therefore we can treat the family of potentials $W_{\CC}=(W_E)_{E\in Ob(\CC)}$ as a function, which is regular with respect to the variable $E$ and formal in the direction $\alpha\in \op{Hom}^1(E,E)$
(or $\alpha\in \op{Ext}^1(E,E)$ if our category is minimal on the diagonal). 
Also the  potential $W_E$ considered as a function of $E\in S_j$ becomes a section of the pro-algebraic vector bundle $\prod_{n\geqslant 2}Sym^n({\cal HOM}^1_{|Diag(S_j)\subset S_j\times S_j})^\star$, where $Diag$ denotes here the diagonal embedding.

\begin{prp} In the case of $3d$ Calabi-Yau category $\CC$ consisting of one object $E$ the potential $W_E$ admits (after a formal change of coordinates) a splitting:
$$W_E=W_E^{min}\oplus Q_E\oplus N_E\,\,,$$
where $W_E^{min}$ is the potential of the minimal model $\CC^{min}$ (i.e. it is a formal series on $\op{Ext}^1(E,E)$),
the quadratic form $Q_E$ is defined on the vector space $\op{Hom}^1(E,E)/\op{Ker}(m_1:\op{Hom}^1(E,E)\to \op{Hom}^2(E,E))$ by the formula
$Q_E(\alpha,\alpha)={m_2(\alpha,\alpha)\over{2}}$, and $N_E$ is the zero function on the image of the map $m_1:\op{Hom}^0(E,E)\to \op{Hom}^1(E,E)$.
In the above splitting formula we use the notation $(f\oplus g)(x,y)=f(x)+g(y)$ for the direct sum of formal functions $f$ and $g$.

\end{prp}

The above Proposition follows from the minimal model theorem for Calabi-Yau algebras (i.e. Calabi-Yau categories with only one object). In its formulation below we are going to use the language of formal non-commutative geometry from \cite{KoSo3}. We assume that the ground field has characteristic zero.

\begin{thm} a) Let $(X,x_0,\omega,d_X)$ be a $\Z/2\Z$-graded non-commutative formal pointed manifold $(X,x_0)$ endowed with an odd symplectic form $\omega$ and homological vector field $d_X$ which preserves $\omega$ and vanishes at $x_0$. Then it is isomorphic to the product
$$(X^{\prime}, x_0^{\prime},\omega^{\prime}, d_{X^{\prime}})\times (X^{\prime\prime}, x_0^{\prime\prime},\omega^{\prime\prime}, d_{X^{\prime\prime}})\,\,,$$
 where
$(X^{\prime}, x_0^{\prime},\omega^{\prime}, d_{X^{\prime}})$ is minimal in the sense that
$(Lie_{d_{X^{\prime}}})_{|T_{x_0^{\prime}X^{\prime}}}=0$ (i.e. $d_{X^{\prime}}$ vanishes quadratically at $x_0^{\prime}$),
and the second factor satisfies the following property:
there exists a finite-dimensional super vector space $V$ endowed with an even non-degenerate quadratic form $Q_V$ such that $(X^{\prime\prime}, x_0^{\prime\prime},\omega^{\prime\prime}, d_{X^{\prime\prime}})$ is isomorphic to the non-commutative formal pointed manifold associated with $V\oplus \Pi V^{\ast}$ (here $\Pi$ is the change of parity functor) endowed with a constant symplectic form $\omega_V$ coming from the natural pairing between $V$ and
$\Pi V^{\ast}$, and homological vector field $d_V$ is the Hamiltonian vector field associated with the pull-back of $Q_V$ under the natural projection $V\oplus \Pi V^{\ast}\to V$.

b) In the $\Z$-graded case when $X$ corresponds to a $3d$ Calabi-Yau algebra (i.e. $\omega$ has degree $-1$) a similar statement holds. In this case $V$ is $\Z$-graded vector space, $Q_V$ has degree $0$, and the tangent space $T_{x_0^{\prime\prime}}X^{\prime\prime}$ isomorphic to $V\oplus V^{\ast}[-1]$.
\end{thm}

{\it Proof.} One can prove  part a) similarly to the usual minimal model theorem for $\A$-algebras or $L_{\infty}$-algebras (it is induction by the order of the Taylor expansion, see e.g. \cite{KoSo2}). Part b) is a ${\bf G}_m$-equivariant version of part a). $\blacksquare$

The Proposition follows from part b) of the Theorem, since we have a decomposition
$\op{Hom}^1(E,E)\simeq \op{Ext}^1(E,E)\oplus V^0\oplus (V^1)^{\ast}$ where $V^i,i\in \Z$ are the graded components of $V$. The restriction of $W^{tot}_E$ to $\op{Hom}^1(E,E)$ is the direct sum of $W_E^{min}$, the restriction of $Q_V$ to $V^0$ (we identify $Q_V$ with $Q_E$) and the zero function on $(V^1)^{\ast}$.

\begin{cor} The minimal model potential $W_E^{min}$ does not depend on a choice of minimal model for $\op{End}^\bullet(E)$,
if considered up to a formal non-linear automorphism of the bundle ${\cal HOM}^1$ restricted to the diagonal
$Diag(S_j)\subset S_j\times S_j$.

\end{cor}
{\it Proof.} Change of the minimal model is a $\Z$-graded change of coordinates. It preserves the topological ideal generated by all coordinates of non-zero degrees. $\blacksquare$

We remark that there is a notion of Calabi-Yau category valid over a field $\k$ of arbitrary characteristic.
In the case of a category with one object $E$ let us denote by $A$ the $\A$-algebra $\op{Hom}^{\bullet}(E,E)$. We assume that $\op{Ext}^\bullet(E,E)=H^\bullet(A)$ is 
finite-dimensional. Then a Calabi-Yau structure of dimension $d$ on $A$ is given by a functional $Tr$ of degree $-d$ on the cyclic homology $HC_{\bullet}(A)$ such that the induced functional on $H^{\bullet}(A)/[H^{\bullet}(A),H^{\bullet}(A)]$ gives rise to a non-degenerate bilinear form $(a,b)\mapsto Tr(ab)$, where $a,b\in H^{\bullet}(A)$.

In the case of positive characteristic the notion of the potential does not exist in the conventional sense. This can be seen in the example $A={\bf F}_3\langle \xi\rangle/(\xi^4), \op{deg}\xi=+1$. The potential should have the form $W(\xi)=\xi^3/3+\dots$ which does not make sense over the field ${\bf F}_3$.

In general it seems that although the potential does not exist, its differential is well-defined as a closed $1$-form.

\begin{rmk} In the case of characteristic zero the cyclic homology $HC_{\bullet}(A)$ can be identified with the cohomology of the complex $\oplus_{n\geqslant 1}Cycl^n(A[1])$ of cyclically invariant tensors (see \cite{KoSo3}). Therefore the potential $W$ becomes a functional of degree $3-d$ on the latter complex, vanishing on the image of the differential. Hence it defines a class $[W]$ in $(HC_{\bullet}(A))^{\ast}$. The latter space is a $\k[[u]]$-module, where $u$ is a variable, $\op{deg} u=+2$ (see loc. cit.). The class $[W]$ is related to the functional $Tr$ discussed above by the formula $[W]=uTr$.
In the case of a Calabi-Yau algebra of dimension $d=2k+1$ it is natural to introduce a cyclic functional $W_k$ with the corresponding class $[W_k]=u^kTr$. It can be thought of as a higher-dimensional Chern-Simons action. In particular, it defines a formal power series $W_k^0$ of degree zero such that it vanishes with the first $k$ derivatives on the formal scheme of solutions to the Maurer-Cartan equation.

\end{rmk}

\subsection{Topology on the space of stability structures}

Let $\CC$ be an ind-constructible weakly unital $\A$-category over a field $\k$ of arbitrary characteristic.
Let $\op{cl}: Ob(\CC)\to \Gamma\simeq \Z^n$ be a map of ind-constructible sets (where $\Gamma$ is considered as a countable set of points) such that the induced map $ Ob(\CC)(\overline{\k})\to \Gamma$ factorizes through a group
homomorphism $\op{cl}_{\overline{\k}}: K_0(\CC(\overline{\k}))\to \Gamma$. It is easy to see that for any field extension
$\k^{\prime}\supset \k$ we obtain a homomorphism $\op{cl}_{\overline{\k}^{\prime}}: K_0(\CC(\overline{\k}^{\prime}))\to \Gamma$.

In the case when $\CC$ is a Calabi-Yau category we require that $\Gamma$ is endowed with an integer-valued bilinear
form $\langle \bullet,\bullet\rangle$ and the homomorphism $\op{cl}_{\overline{\k}}$ is compatible with $\langle \bullet,\bullet\rangle$ and the Euler form on $K_0(\CC(\overline{\k}))$.

For ind-constructible
triangulated $\A$-categories the notion of stability structure admits the following version.

\begin{defn} A constructible stability structure on $(\CC,\op{cl})$ is given by the following data (cf. Introduction, Section 1.2):

\begin{itemize}
\item{ an ind-constructible subset 
$$\CC^{ss}\subset Ob(\CC)$$  consisting of objects called semistable, and satisfying the condition that with each object it contains all isomorphic objects,}

\item{an additive map $Z:\Gamma\to \C$ called central charge, such that $Z(E):=Z(\op{cl}(E))\ne 0$ if $E\in \CC^{ss}$,}

\item{a choice of the branch of logarithm $\op{Log}Z(E)\in \C$ for any $E\in \CC^{ss}$ which is constructible as a function of $E$.}

\end{itemize}

These data are required to satisfy the corresponding axioms from Section 1.2 for the category
$\CC(\overline{\k})$.

In particular
\begin{itemize}
\item{the set of $E\in \CC^{ss}(\overline{\k})\subset Ob(\CC)(\overline{\k})$ with the fixed $\op{cl}(E)\in \Gamma\setminus \{0\}$ and fixed $\op{Log}Z(E)$ is a constructible set.}
\end{itemize}
\end{defn}

Before we proceed with the topology let us make a comparison with the ``Lie-algebraic" story of Section 2.
First, we observe that the set $\CC^{ss}$ can be thought of as an analog of the collection of elements $(a(\gamma))_{\gamma \in \Gamma\setminus \{0\}}$ from Section 2.1.
Then we give the following definition of another data and axioms which is equivalent to the one given above and can be thought of as an analog of the collection of the group elements $A_V$.

\begin{defn} A constructible stability structure on $(\CC,\op{cl})$ is given by the following data:

\begin{itemize}

\item{an additive map $Z:\Gamma\to \C$,}

\item{for any bounded connected set $I\subset \R$ an ind-constructible subset
${\cal P}(I)\subset Ob(\CC)(\kk)$ which contains which every object all isomorphic objects.}

\end{itemize}

These data are required to satisfy the following axioms\footnote{One should read expressions $Ob(\CC),{\cal P}(I)$ etc. in this list of axioms
 as sets of $\kk$-points.}:

\begin{itemize}
\item{the zero object of the category $\CC(\kk)$ belongs to all ${\cal P}(I)$,}
\item{ $\cup_{n\in \Z_{>0}}{\cal P}([-n,n])=Ob(\CC)(\kk)$,}
\item{if $I_1<I_2$ in the sense that every element of $I_1$ is strictly less than any
element of $I_2$ then for any $E_1\in {\cal P}(I_1)$ and $E_2\in {\cal P}(I_2)$ one has
$\op{Ext}^{\le 0}(E_2,E_1)=0$,}
\item{${\cal P}(I+1)={\cal P}(I)[1]$ where $[1]$ is the shift functor in $\CC(\kk)$,}
\item{(Extension Property) If $I=I_1\sqcup I_2$ and $I_1<I_2$ in the above sense then the ind-constructible set ${\cal P}(I)$ is isomorphic to the ind-constructible subset of such objects $E\in Ob(\CC)(\kk)$ which are extensions $E_2\to E\to E_1$ with $E_m\in {\cal P}(I_m), m=1,2$,}
\item{if $I$ is an interval of the length strictly less than one, $E\in {\cal P}(I), E\ne 0$,
then  $Z(\op{cl}(E))$ belongs to the strict sector 
$$V_I=\{z=re^{\pi i\varphi}\in \C^{\ast}|r>0,\varphi \in I\}\,\,,$$}
\item{there is a non-degenerate quadratic form $Q$ on $\Gamma_{\R}$ such that $Q_{|\op{Ker}Z}<0$, and for any interval $I$ of the length strictly less than $1$ the set
$$\{\op{cl}(E)\in \Gamma|E\in {\cal P}(I)\}\subset \Gamma$$
 belongs to the convex cone $C(V_I,Z,Q)$ defined in Section 2.2,}
\item{let $I$ be an interval of the length strictly less than $1$, and $\gamma\in \Gamma$. Then the set $\{E\in {\cal P}(I)| \op{cl}(E)=\gamma\}$ is constructible.}

\end{itemize}

\end{defn}

The equivalence of Definitions 9,10 can be proved similarly to the proof of Theorem 2.

With this equivalent description of a constructible stability condition we observe that the collection of sets ${\cal P}(I)$ considered for all intervals $I$ with the length less than $1$ are analogous to the collection of elements $A_V$ where $V=V_I$ (see Section 2.2) and the Extension Property is analogous to the Factorization Property.

One has the following result.

\begin{prp} For any constructible stability structure on $\CC$ and any field extension $\k\subset \k^{\prime}$ the   category $\CC(\k^{\prime})$ carries a locally finite stability structure in the sense of \cite{Br1} with the central charge given by $Z\circ \op{cl}_{{\k}^{\prime}}$ and the collection of additive subcategories ${\cal P}(I)({\k}^{\prime})$,
where $I$ runs through the set of bounded connected subsets of $\R$ as above.

\end{prp}

{\it Proof.} The proof is straightforward. Local finiteness in the sense of \cite{Br1} follows from our (stronger) assumption on the quadratic form $Q$. $\blacksquare$

Let us denote by $Stab(\CC,\op{cl})$ the set of constructible stability structures on $\CC$ with a fixed class map $\op{cl}$.
Our goal is to introduce a topology on $Stab(\CC,\op{cl})$. 

Let $\Delta\subset \C$ be a triangle with one vertex at the origin. We choose a branch  of the function $z\mapsto \op{Log} z$ for $z\in \Delta$. We denote the corresponding argument function by $\op{Arg}(z)$.
We denote by ${\CC}_{{\Delta},\op{Log}}$
an $\A$-subcategory of $\CC$ generated by the zero object ${\bf 0}$, semistable objects $E$ with $Z(E)\in \Delta, \op{Arg}(E)\in \op{Arg}(\Delta)$ as well as extensions $J$
of such objects satisfying the condition $Z(J)\in \Delta$. We allow the case $\Delta=V$ where $V$ is a sector, in which case we will use the notation $\CC_{V,\op{Log}}$.
It is easy to see that
$\CC_{{\Delta},\op{Log}}$ is an ind-constructible  category.
Notice that in the language of ind-constructible sets ${\cal P}(I)$ we have $Ob(\CC_{{V_I},\op{Log}})={\cal P}(I)$ for some choice of the branch $\op{Log}$. The condition of genericity of the sector $V_I$ corresponding to a closed interval $I=[a,b]$ of the length less than $1$ corresponds
to the following genericity condition of the set ${\cal P}(I)$: both ${\cal P}(\{a\})$ and
${\cal P}(\{b\})$ are zero categories (equivalently ${\cal P}([a,b])={\cal P}((a,b))$).

Let us fix a non-degenerate quadratic form $Q$ on $\Gamma_{\R}$ such that $Q_{|\op{Ker} Z}<0$ and $Q(\op{cl}(E))>0$ for any $E\in \CC^{ss}$.
 We introduce the topology on $Stab(\CC,\op{cl})$ in the following way. Let us consider a family
$\sigma_x=(Z_x,\CC^{ss}_x,\dots), x\in X$ of stability structures in a neighborhood of $x_0\in X$.
Then for every point $x$, a generic closed interval $I=[a,b]$ of the length less than $1$  we have the full category ${\cal P}(I)_x\subset Ob(\CC)$.
 For a given $\gamma\in \Gamma\setminus\{0\}$ we denote by ${\cal P}(I)_{x,\gamma}$ the constructible subset of objects $E\in {\cal P}(I)_x $ such that $\op{cl}(E)=\gamma$.

We say that a  family $\sigma_x=(Z_x,\CC^{ss}_x,\dots), x\in X$ of stability structures parametrized by a topological space $X$
is {\it continuous at a given stability structure $\sigma_{x_0}=(Z_0,\CC^{ss}_0,\dots)$ } if:

a) The map $x\to Z_x$ is continuous at $x=x_0$.

b) There exists a neighborhood $U$ of $x_0$ such that for any $E\in \CC^{ss}_x, x\in U$ we have $Q(\op{cl}(E))\geqslant 0$.

c) For any generic closed interval $I$ of the length strictly less than $1$ the
constructible set ${\cal P}(I)_{x,\gamma}$ is locally constant near $x_0$ (cf. Definition 3c)).

In this way we obtain a Hausdorff topology on $Stab(\CC,\op{cl})$.
We can define a parallel transport along a path $\sigma_t$ in the space $Stab(\CC,\op{cl})$ similarly to
the case of  stability structures in graded Lie algebras discussed in Section 2. Each time when we cross the wall
of first kind we use the above property c) in order to ``recalculate" the set
of semistable objects. In order to do this we use the following property:
$E\in Ob(\CC_{{\Delta},\op{Log}})$ is semistable iff there is no extension $E_2\to E\to E_1$ where $E_i,i=1,2$ are
non-zero objects
of $\CC_{{\Delta},\op{Log}}$ such that $\op{Arg}(E_2)>\op{Arg}(E_1)$. These considerations also ensure that the Theorem 1 from Introduction holds
(i.e. the natural projection of the space of stability conditions to the space of central charges is a local homeomorphism).

\section{Motivic functions and motivic Milnor fiber}

\subsection{Recollection on motivic functions}

Recall (see \cite{DeLo1}) that for any constructible set $X$ over $\k$ one can define an abelian group $Mot(X)$ of motivic functions as the group generated by symbols $[\pi:S\to X]:=[S\to X]$ where $\pi$ is a morphism of constructible sets, subject to the relations $$[(S_1\sqcup S_2)\to X]=[S_1\to X]+[S_2\to X]\,\,.$$
For any constructible morphism $f:X\to Y$ we have two homomorphisms of groups: 

1) $f_{!}: Mot(X)\to Mot(Y)$, such that
$[\pi:S\to X]\mapsto [f\circ \pi:S\to X]$;

2) $f^{\ast}: Mot(Y)\to Mot(X)$, such that
$[S^{\prime}\to Y]\mapsto [S^{\prime}\times_YX\to X]$.

Moreover, $Mot(X)$ is a commutative ring via the fiber product operation.
We denote by ${\mathbb L}\in Mot(Spec(\k))$ the element $[{\bf A}^1_{\k}]:=
[{\bf A}^1_\k \to  Spec(\k)]$. It is customary to add its formal inverse
${\mathbb L}^{-1}$ to the ring $Mot(Spec(\k))$ (or more generally to the ring $Mot(X)$ which is a $Mot(Spec(\k))$-algebra).

There are several ``realizations" of the theory of motivic functions which we are going to recall below.

(i) There is a homomorphism of rings 
$$\chi: Mot(X)\to Constr(X,\Z)\,\,,$$ where
$ Constr(X,\Z)$ is the ring of integer-valued constructible functions on $X$ endowed with the pointwise multiplication. More precisely, the element $[\pi:Y\to X]$ is mapped into $\chi(\pi)$, where
$\chi(\pi)(x)=\chi(\pi^{-1}(x))$, which is the Euler characteristic of the fiber $\pi^{-1}(x)$.

(ii) Let now $X$ be a scheme of finite type over a field $\k$, and $l\ne char\,\k$ be a prime number. There is a homomorphism of rings
$$Mot(X)\to K_0(D^b_{\op{constr}}(X,\Q_l))\,\,,$$
 where $D^b_{\op{constr}}(X,\Q_l)$ is the bounded derived category of \'etale $l$-adic  sheaves on  $X$ with constructible cohomology, such that
$$[\pi:S\to X]\mapsto \pi_{!}(\Q_l)\,\,,$$
which is the direct image of the constant sheaf $\Q_l$. Notice that $D^b_{\op{constr}}(X,\Q_l)$ is a tensor category, hence Grothendieck group $K_0$ is naturally a ring. The homomorphisms $f_{!}$ and $f^{\ast}$ discussed above correspond to the functors $f_{!}$ (direct image with compact support) and $f^{\ast}$ (pullback), which we will denote by the same symbols. We will also use the notation
$\int_X\phi:=f_{!}\,(\phi)$ for the canonical map $f:X\to Spec(\k)$.

(iii) In the special case $X=Spec(\k)$ the above homomorphism becomes a map
$$[S]\mapsto \sum_i(-1)^i[H^i_{c}(S\times_{Spec (\k)} Spec (\overline{\k}),\Q_l)]\in 
K_0(\op{Gal}(\overline{\k}/\k)-mod_{\Q_l})\,\,,$$
where $\op{Gal}(\overline{\k}/\k)-mod_{\Q_l}$ is the tensor category of finite-dimensional continuous
 $l$-adic representations of the Galois group $\op{Gal}(\overline{\k}/\k)$, and we take the \'etale cohomology of $S$ with compact support. 

(iv) If $\k={\bf F}_q$ is a finite field then for any $n\geqslant 1$ we have a homomorphism
$$Mot(X)\to \Z^{X({\bf F}_{q^n})}$$
 given by $$[\pi:Y\to X]\mapsto (x\mapsto \#\{y\in X({\bf F}_{q^n})\,|\,\pi(y)=x\})\,\, 
.$$
Operations $f^!,f_\ast$ correspond to  pullbacks and pushforwards of functions
on finite sets.

(v) If $\k\subset \C$ then the category of $l$-adic constructible sheaves on a scheme of finite type $X$ can be replaced in the above considerations by the Saito's category of mixed Hodge modules (see \cite{Sa}). 

(vi) In the case $X=Spec(\k)$ one has two additional homomorphisms:

a) Serre polynomial $$Mot(Spec(\k))\to \Z[q^{1/2}]$$ such that
$$[Y]\mapsto \sum_i(-1)^i\sum_{w\in {\Z_{\geqslant 0}}}\dim H^{i,w}_c(Y)q^{w/2}\,\,,$$
where
$H^{i,w}_c(Y)$ is the weight $w$ component in the $i$-th Weil cohomology group with compact support.

b) If $char\,\k=0$ then we have the Hodge polynomial $$Mot(Spec(\k))\to \Z[z_1,z_2]$$ such that
$$[Y]\mapsto \sum_{i\geqslant 0}(-1)^i\sum_{p,\,q\geqslant 0}\dim \op{Gr}^p_F(\op{Gr}^W_{p+q} H^i_{DR,c}(Y))z_1^{p}z_2^{q}\,\,,$$
where $\op{Gr}^W_{\bullet} $ and $\op{Gr}_F^{\bullet}$ denote the  graded components with respect to the weight and Hodge filtrations, and $H^i_{DR,c}$ denotes the de Rham cohomology with compact support.

Clearly the Hodge polynomial determines the Serre polynomial via the homomorphism $\Z[z_1,z_2]\to \Z[q^{1/2}]$ such that $z_i\mapsto q^{1/2}, i=1,2$.

\subsection{Motivic functions in the equivariant setting}

Here we give a short exposition of the generalization of the theory of motivic functions
in the equivariant setting (essentially due to  Joyce  \cite{Jo4},  here we use slightly different terms).

Let $X$ be a constructible set over a field $\k$ and $G$ be an affine algebraic group acting on $X$, in the sense that $G(\kk)$ acts on $X(\kk)$ and there
exists a $G$-variety $S$ over $\k$ with a constructible equivariant identification $X(\kk)\simeq S(\kk)$. 

We define the group $Mot^G(X)$ of $G$-equivariant motivic functions as abelian group generated by $G$-equivariant constructible maps $[Y\to X]$
 modulo the relations \begin{itemize}
\item $[(Y_1\sqcup Y_2)\to X]=[Y_1\to X]+[Y_2\to X]$,
\item $[Y_2\to X]=[(Y_1\times {\bf A}_\k^d)\to X]$ if $Y_2\to Y_1$ is a $G$-equivariant constructible vector 
bundle of rank $d$.
\end{itemize}
This group form a commutative ring via the fiber product, and a morphism of constructible sets with group actions induces a pullback homomorphism
 of corresponding rings. There is no natural operation of a pushforward for equivariant motivic functions, and for that one has to enlarge
   ring of functions.

  Consider the following 2-category of constructible stacks. First, its objects are pairs
$(X,G)$ as above\footnote{Strictly speaking,we should denote such stacks as triples $(X,G,\alpha)$ where $\alpha$ is the action of $G$ on $X$.}, and the objects of the category of 1-morphisms $\op{Hom}((X_1,G_1),(X_2,G_2))$ are pairs $(Z,f)$ where
$Z$ is a $G_1\times G_2$-constructible set such that $\{e\}\times G_2$ acts freely on $Z$ in such a way that we have the induced $G_1$-equivariant isomorphism
$Z/G_2\simeq X_1$, and $f:Z\to X_2$ is a $G_1\times G_2$-equivariant map ($G_1$ acts trivially on $X_2$). An element of
$\op{Hom}((X_1,G_1),(X_2,G_2))$ defines a map of sets $X_1(\kk)/G_1(\kk)\to X_2(\kk)/G_2(\kk)$.
Furthermore,  elements of $\op{Hom}((X_1,G_1),(X_2,G_2))$ form naturally objects of a groupoid, so we obtain a $2$-category ${\cal SCON}_{\k}$ of {\it constructible stacks over $\k$}.   The 2-category of constructible stacks carries a direct sum operation induced by the disjoint union of stacks $$(X_1,G_1)\sqcup (X_2,G_2)=
((X_1\times G_2\sqcup X_2\times G_1),G_1\times G_2)\,,$$ as well as the product induced by the Cartesian product
$$(X_1,G_1)\times (X_2,G_2)=(X_1\times X_2,G_1\times G_2)\,.$$

 The abelian group of stack motivic function $Mot_{st}((X,G))$ is generated by the group of isomorphism classes of 1-morphisms of stacks $[(Y,H)\to (X,G)]$ with the fixed target $(X,G)$, subject to the relations \begin{itemize}
\item
$[((Y_1,G_1)\sqcup (Y_2,G_2))\to (X,G)]=[(Y_1,G_1)\to (X,G)]+[(Y_2,G_2)\to (X,G)]$
\item $[(Y_2,G_1)\to (X,G)]=[(Y_1\times {\bf A}_\k^d,G_1)\to (X,G)]$ if $Y_2\to Y_1$ is a $G_1$-equivariant constructible vector 
bundle of rank $d$.
\end{itemize}
The ring $Mot^G(X)$ maps to $Mot_{st}((X,G))$. 
 Notice that every isomorphism class
$[(Y,H)\to (X,G)]$ corresponds to an ordinary morphism of constructible sets acted by algebraic groups. Indeed, in the notation of the definition of $1$-morphism of 
stacks with $(X_1,G_1):=(Y,H)$ and $(X_2,G_2):=(X,G)$, we can replace the source
  $(Y,H)$ by an equivalent stack $(Z,G_1\times G_2)$ and get an ordinary morphism
   $(Z,G_1\times G_2)\to (X,G)$ of constructible sets acted by algebraic groups.
    One can define the pullback, the pushforward and the product of elements of $Mot_{st}((X,G))$.

Finally, for a constructible stack ${\cal S}=(X,G)$ we define its class
 in the ring $ K_0(Var_{\k})[[{\mathbb L}]^{-1}, ([GL(n)]^{-1})_{n\geqslant 1}]$ as
$$[{\cal S}]={[(X\times GL(n))/G]\over {[GL(n)]}},$$
where we have chosen an embedding $G\to GL(n)$ for some $n\geqslant 1$, and $(X\times GL(n))/G$ is the ordinary quotient by the diagonal free action (thus in the RHS we have the quotient of motives of ordinary varieties). The result does not depend on the choice of embedding (see \cite{BDh}, Lemma 2.3). Then we define the integral
$\int_{\cal S}: Mot_{st}({\cal S})\to K_0(Var_{\k})[[{\mathbb L}]^{-1}, ([GL(n)]^{-1})_{n\geqslant 1}]$ as
$\int_{\cal S}[{\cal S}^{\prime}\to {\cal S}]=[{\cal S}^{\prime}]$.

If $\k$ is finite, one can associate with every constructible stack ${\cal S}=(X,G)$ a finite set ${\cal S}(\k)$, 
the set of orbits of $GL(n,\k)$ acting on $((X\times GL(n))/G)(\k)$.  There is a homomorphism of the algebra of stack motivic functions 
 to the algebra of $\Q$-valued functions
 on ${\cal S}(\k)$. The identity function 
$$\mathbf{1}_{(X,G)}:=[(X,G)\to (X,G)]$$
represented by the identity map,  when interpreted as a measure (for pushforwards) maps to  
the ``stack counting measure''  on ${\cal S}(\k)$  which is equal to $\#(GL(n,\k))^{-1}$ 
times the 
 direct image of the ordinary counting measure\footnote{For every finite set $S$ the counting measure on $S$ has weight 1 for each element $s\in S$.} on $((X\times GL(n))/G)(\k)$.  Its density with respect to the ordinary counting measure on 
${\cal S}(\k)$ 
 is given by
the inverse 
to the order of the stabilizer.

Our construction in Section 3.2 can be  rephrased as a construction of an ind-constructible stack of objects. Hence we can speak 
about stack motivic  functions on an ind-constructible category.

\subsection{Motivic Milnor fiber}

Let $M$ be a complex manifold, $x_0\in M$. Recall, that for a germ $f$ of an analytic function at $x_0$ such that $f(x_0)=0$ one can define its {\it Milnor fiber } $MF_{x_0}(f)$, which is a locally trivial $C^{\infty}$-bundle over $S^1$  of manifolds with the boundary (defined only up to a diffeomorphism):
$$ \{z\in M|\,\,dist(z,x_0)\le\varepsilon_1, |f(z)|=\varepsilon_2\}\to S^1=\R/2\pi\Z\,\,,$$
where $z\mapsto \op{Arg}\,f(z)$.
Here $dist$ is any smooth metric on $M$ near $x_0$, and there exists a constant $C=C(f,dist)$ and a positive integer $N=N(f)$ such that for $0<\varepsilon_1\le C$ and $0<\varepsilon_2<\varepsilon_1^N$ the
$C^{\infty}$ type of the bundle is the same for all $\varepsilon_1,\varepsilon_2,dist$.

In particular, taking the cohomology of the fibers we obtain a well-defined local system on $S^1$.

There are several algebro-geometric versions of this construction (theories of nearby cycles). They produce analogs of local systems on $S^1$, for example $l$-adic representations of the group $\op{Gal}(\k((t))^{sep}/\k((t)))$ where $l\ne char\,\k$.

There is a convenient model of the Milnor fiber in non-archimedean geometry. In order to describe it we note that the field $K=\k((t))$ is a non-archimedean field endowed with the (standard) valuation, and with the norm given by $|a|=c^{val(a)}$ for a given constant $c \in (0,1)$. Let 
$$f\in \k[[x_1,\dots,x_n]], f(0)=0$$ be a formal series considered as an element of $K[[x_1,\dots,x_n]]$. Clearly it is convergent in the non-archimedean sense in the domain $U\subset ({\bf A}^n_{\k})^{an}$ defined by inequalities $|x_i|<1, 1\le i\le n$.
The non-archimedean analog of the Milnor fiber is given (at the level of points) by the fibration 
$$\{x=(x_1,\dots,x_n)\in U|\max_i|x_i|\le\varepsilon_1, 0<|f(x)|\le \varepsilon_2\}\to $$ 
$$\to \{w\in({\bf A}^1_{\k})^{an}|\,\,0<|w|\le \varepsilon_2\}\,\,,$$ where $\varepsilon_1,\varepsilon_2$ are positive numbers as above (cf. \cite{NicS}).

Ideally, we would like to have the following picture. Let $V\to X$ be a vector bundle over a scheme of finite type $X/\k$, and 
$$f\in \prod_{n\geqslant 1}\Gamma(X,Sym^n(V^{\ast}))$$ be a function on the formal completion of zero section of $\op{tot}(V)$ vanishing on $X$.
 We would like to associate with such data a motivic Milnor fiber $$MF(f)\in Mot(X\times {\bf G}_m)\,\,.$$ 
 Here the factor ${\bf G}_m$ replaces the circle $S^1$ in the analytic picture.
 Moreover, the motivic function $MF(f)$ should be ``unramified" in ${\bf G}_m$-direction (i.e. it should correspond to a
${\bf G}_m$-invariant stratification of $X\times {\bf G}_m$). In the case $\k=\C$ , 
assuming that $f$ is convergent near zero section,
the value $MF_x(f)$ at a point $(x,\epsilon)\in X(\C)\times \C^\ast$  should be thought of as a representative of the alternating sum $$\sum_i(-1)^i[H^i(f^{-1}(\varepsilon)\cap B_{0,x})]$$ where $\varepsilon\ne 0$ is a sufficiently small complex number and $B_{0,x}$ is a small open ball around $0$ in the fiber $V_x$. Notice that here we use the usual cohomology and not the one
with compact support.

Also, we can consider the case when $X$ is a constructible set and $V\to X$ is a constructible vector bundle.
We say that $f\in \prod_{n\geqslant 1}\Gamma(X,Sym^n(V^{\ast}))$ is constructible if $(X,V,f)$
is constructibly isomorphic to an algebraic family of formal functions over a scheme of finite type as above.

This goal was achieved by Denef and Loeser (see \cite{DeLo1}) in the case $char\,\k=0$ by using motivic integration and resolution of singularities. In this case the group $Mot(X\times {\bf G}_m)$ is replaced by $Mot^{\mu}(X)$, where
$\mu=\varprojlim_{n}\mu_n$ and $\mu$ acts trivially on $X$.
Here  $\mu_n$ is the group of $n$-th roots of $1$ in $\k$. We will always assume that $\mu$-action  is ``good" in the sense that $\mu$ acts via a  finite quotient $\mu_n$ and every its orbit is contained in an affine open subscheme.
Notice that there is a homomorphism of groups $$Mot^{\mu}(X)\to Mot(X\times {\bf G}_m),\,\,\,
[\pi:Y\to X]\mapsto [\pi_1:(Y\times {\bf G}_m)/\mu_n\to X\times {\bf G}_m]\,\,,$$
 where $\mu$ acts on $Y$ via its quotient $\mu_n$ and $\pi_1(y,t)=\pi(y)t^n$.

As we work with constructible sets, it is sufficient to define the motivic Milnor fiber not for a family, but for an individual
 formal germ of a function. Let $M$ be a smooth formal scheme over $\k$ with closed point $x_0$ and $f$ be a function on $M$
vanishing at $x_0$ (e.g. $M$ could be the formal completion at $0$ of a fiber of vector bundle $V\to X$ in the above notation). We assume that $f$ is not identically equal to zero near $x_0$, otherwise the Milnor fiber would be empty. 
Let us choose a simple normal crossing resolution of singularities $\pi:M'\to M$ of the hypersurface in $M$ given by 
the equation $f=0$ with exceptional divisors $D_j, j\in J$. 
The explicit formula for the motivic Milnor fiber from \cite{DeLo1}   looks such as follows\footnote{In \cite{DeLo1} it was assumed that
$f$ is a regular function on a smooth scheme, but the formula and all the arguments work in the formal setting as well.}.

$$MF_{x_0}(f)=\sum_{I\subset J, I\ne \emptyset}(1-{\mathbb L})^{\#I-1}[\widetilde{{D_I}^0} \cap \pi^{-1}(x_0)]\in Mot^\mu
(Spec(\k))\,\,,$$
where $D_I=\cap_{j\in I}D_j$, $D_I^0$ is the complement in $D_I$ to the union of all other exceptional divisors, 
and $\widetilde{{D_I}^0}\to D_I^0$ is a certain Galois cover with the Galois group  $\mu_{m_I}$, where $m_I$ is the g.c.d.
of the multiplicities of all divisors $D_i,i\in I$ (see \cite{DeLo1} for the details). Informally speaking, the fiber of the cover
 $\widetilde{{D_I}^0}\to D_I^0$ is the set of connected components of a non-zero level set of function $f\circ \pi$ near a point of  $D_I^0$.

The space
$Mot^{\mu}(X)$ carries an associative product introduced by Looijenga (see \cite{Loo}) which is different from the one defined above. It is essential for the motivic Thom-Sebastiani theorem which will be discussed later.
Let us sketch  a construction of this product.

First, let us introduce the commutative ring $Mot(X\times {\bf A}^1_{\k})_{conv}$ which coincides as an abelian group with $Mot(X\times {\bf A}^1_{\k})$ but carries the ``convolution product"
$$[f_1: S_1\to X\times {\bf A}^1_{\k}]\circ[f_2: S_2\to X\times {\bf A}^1_{\k}]=[f_1\oplus f_2:S_1\times_X S_2\to X\times {\bf A}^1_{\k}],$$
where $(f_1\oplus f_2)(s_1,s_2)=(pr_X(f_1(s_1)),pr_{{\bf A}^1_{\k}}(f_1(s_1)) +pr_{{\bf A}^1_{\k}}(f_2(s_2)))$.
The ring $Mot(X\times {\bf A}^1_{\k})_{conv}$ contains the ideal 
$${\bf I}:=pr_X^{\ast}(Mot(X))\,\,.$$
By definition we have an epimorphism of abelian groups $Mot(X)\to {\bf I}$. Let 
$$i:X\to X\times {\bf A}^1_{\k}, x\mapsto (x,0),\,\,\,j:X\times ({\bf A}^1_{\k}\setminus \{0\})\to X\times {\bf A}^1_{\k}$$ be natural embeddings. They give rise to an isomorphism of abelian groups 
$$i^{\ast}\oplus j^{\ast}:Mot(X\times {\bf A}^1_{\k})\simeq Mot(X)\oplus Mot( X\times ({\bf A}^1_{\k}\setminus \{0\}))\,\,.$$ Since $i^{\ast}\circ pr_X^{\ast}=id_{Mot(X)}$ we see that the restriction of
$pr_X^{\ast}$ to ${\bf I}$ gives an isomorphism of abelian groups 
$${\bf I}\simeq Mot(X)\,\,,$$ and $j^{\ast}$ induces the isomorphism of groups 
$$Mot(X\times {\bf A}^1_{\k})/{\bf I}\simeq Mot( X\times ({\bf A}^1_{\k}\setminus \{0\}))\,\,.$$ Using the latter isomorphism we transfer the convolution product and endow $Mot( X\times ({\bf A}^1_{\k}\setminus \{0\}))$ with an associative product which we will call {\it exotic}.

Recall that we have a homomorphism of groups 
$$Mot^{\mu}(X)\to Mot(X\times {\bf G}_m)=Mot( X\times ({\bf A}^1_{\k}\setminus \{0\}))\,\,.$$ One can check that the image of $Mot^{\mu}(X)$ is closed under the exotic product.
Intuitively, the image consists of isotrivial families of varieties over $X\times ({\bf A}^1_{\k}\setminus\{0\})$
equipped with a flat connection which has  finite (i.e. belonging to some $\mu_n$) monodromy. The complicated formula from \cite{Loo} coincides with the induced product on $Mot^{\mu}(X)$.
In what follows we will use the notation
$${\cal M}^{\mu}(X):=(Mot^{\mu}(X), \mathrm{\, exotic\,\,\,product})\,\,. $$

Let $V\to X,\,V'\to Y$ be two constructible vector bundles endowed with constructible families $f,g$ of formal power series.
We denote by $f\oplus g$ the sum of pullbacks of $f$ and $g$ to the constructible vector
bundle 
$$pr_X^*V\oplus pr_Y^*V'\to X\times Y\,\,.$$
Then we have the following motivic version of Thom-Sebastiani theorem.

\begin{thm}(\cite{DeLo2}) One has
$$(1-MF(f\oplus g))=pr_X^*(1-MF(f))\cdot
pr_Y^*(1-MF(g))\in {\cal M}^\mu (X\times Y)\,\,.$$
\end{thm}

One can make similar constructions in the equivariant setting. 
Let  $X/\k$ be a constructible set endowed with the good action
of an  affine algebraic group $G$. We  endow $X$ also with the trivial $\mu$-action. Then, similarly to the above, we can equip $Mot^{G\times \mu}(X)$ with the exotic product (by considering $G$-equivariant families over $X$ in the previous considerations). We will denote the resulting ring by ${\cal M}^{G,\mu}(X)$.  Using the 
canonical resolution of singularities (see e.g. \cite{BiMi}) one can define the equivariant motivic Milnor fiber in 
the case of equivariant families of functions, and state the corresponding version of  Thom-Sebastiani theorem.

In the case of arbitrary $\k$ there is an $l$-adic version of the above results. More precisely, the theory of Milnor fiber is replaced by the theory of nearby cycles (see \cite{SGA}), with the convolution defined by Laumon \cite{Lau}. The Thom-Sebastiani theorem was proved in this case by Pierre Deligne and probably by Lei Fu (both unpublished).

There is an analog of the Hodge polynomial in the story (see \cite{DeLo1}, 3.1.3). Let us assume that $\k=\C$ for simplicity. Then we have a homomorphism of rings 
$${\cal M}^{\mu}(Spec(\k))\to \left\{\sum_{\alpha,\beta\in \Q,\alpha+\beta\in \Z}c_{\alpha,\beta}z_1^{\alpha}z_2^{\beta}\left|\,c_{\alpha,\beta}\in \Z\right.\right\}\,\,.$$
 Namely, for a smooth projective $\mu_n$-scheme $Y$ we set
$$[Y]\mapsto \sum_{p,q\geqslant 0, p,q\in \Z}(-1)^{p+q}\dim H^{p,q,0}(Y)z_1^pz_2^q+$$
$${}+\sum_{p,q\geqslant 0, p,q\in \Z}\sum_{1\le i\le n-1}(-1)^{p+q}\dim H^{p,q,i}(Y)z_1^{p+i/n}z_2^{q+1-i/n}\,\,,$$
where $H^{p,q,i}(Y)$ is the subspace of the cohomology $H^{p,q}(Y)$, where an element $\xi\in \mu_n$ acts by multiplication by $\xi^i$. The appearance of rational exponents was first time observed in the Hodge spectrum of a  complex isolated
singularity.
Taking $z_1=z_2=q^{1/2}$ we obtain the corresponding Serre polynomial.

\subsection{An integral identity}

In this section we are going to discuss the identity which will be crucial in the proof of the main theorem of Section 6.

Let $\k$ be a field of characteristic zero, and $V_1,V_2,V_3$ be finite-dimensional $\k$-vector spaces.

\begin{conj} Let $W$ be a formal series on the product $V_1\times V_2\times V_3$ of three vector spaces,
depending in a constructible way on finitely many extra parameters, such that $W(0,0,0)=0$ and $W$ has degree zero
with respect to the diagonal action of the multiplicative group ${\bf G}_m$ with the
weights $(1,-1,0)$. We denote by $\widehat{W}$ the ${\bf G}_m$-equivariant extension of $W$
to the formal neighborhood $\widehat{V_1}$ of $V_1\times \{0\}\times \{0\}\subset V_1\times V_2\times V_3$.
Then we have the following formula (where we denote the direct image by the integral):
$$\int_{v_1\in V_1}(1-MF_{(v_1,0,0)}(\widehat{W}))={\mathbb L}^{\dim  V_1}(1-MF_{(0,0,0)}(W_{|(0,0,V_3)}))\,\,,$$
where in the RHS we consider the motivic Milnor fiber  at $(0,0,0)$ of the restriction of  $W$  to the subspace $(0,0,V_3)$.
\end{conj}

Using the obvious equality $\int_{V_1} 1={\mathbb L}^{\dim  V_1}$ we can rewrite the identity as 
$$\int_{v_1\in V_1}MF_{(v_1,0,0)}(\widehat{W})={\mathbb L}^{\dim  V_1}\cdot MF_{(0,0,0)}(W_{|(0,0,V_3)})\,\,.$$

Let us discuss the $l$-adic version of the Conjecture. For simplicity we assume that the vector spaces do not depend on extra parameters.
Then we have a morphism of formal schemes $\pi: \widehat{V_1}\to Spf(\k[[w]])$ such that $w\mapsto \widehat{W}$ as well as an embedding $i_{V_1}:V_1\to \widehat{W}^{-1}(0)$. For any morphism 
$\pi: X\to Spf(\k[[w]])$ of formal schemes to we denote by $\R\psi_{\pi}$ the functor of nearby cycles.  It acts from the bounded derived category of $l$-adic constructible sheaves on $X$ to the bounded derived category of $l$-adic constructible sheaves on $X_0=\pi^{-1}(0)$ endowed with the action of the inertia group (hence  they can be informally thought of as $l$-adic constructible sheaves on $X_0\times_{Spec(\k)}\k((w))$).

\begin{prp} The complex $\R\Gamma_c(i^{\ast}_{V_1}\R\psi_{\pi}(\Q_l))$ is isomorphic
(as a complex of  $\op{Gal}(\overline{\k((w))}/\k((w)))$-modules) to the complex
$\R\Gamma_c(V_1,\Q_l)\otimes j^{\ast}\R\psi_{\widehat{\pi}}(\Q_l)$, where $j:Spec(\k)\to V_3, j(0)=0$ is the natural embedding and $\widehat{\pi}$ is the morphism of the formal completion of $0\in \{0\}\times\{0\}\times V_3$ to $Spf(\k[[w]])$ given by the restriction $W_{|\,\{0\}\times\{0\}\times V_3}$.

\end{prp}

{\it Proof.} We will give a sketch of the proof based on the non-archimedean model for the Milnor fiber described in Section 4.3. 

Let us consider the $\k((t))$-analytic space (in the sense of Berkovich) associated with
the scheme $(V_1\times V_2\times V_3)\times_{Spec(\k)}Spec(\k((t)))$. Let us choose sufficiently small positive numbers $\varepsilon_1, \varepsilon_2,\varepsilon_3,\varepsilon_4 $ (we will specify them later) and define an analytic subspace $Y=Y_{\varepsilon_1, \varepsilon_2,\varepsilon_3,\varepsilon_4}$ by the following inequalities: 
$$|v_1|\le 1+\varepsilon_1, \,\,|v_2|,|v_3|\le \varepsilon_2, \,\,\varepsilon_4\le |W(v_1,v_2,v_3)|\le \varepsilon_3\,\,.$$
Notice that the series $W(v_1,v_2,v_3)$ is convergent on $Y$ because of homogeneity property.
We introduce  another analytic space
$Y^{\prime}=Y_{\varepsilon_1, \varepsilon_2,\varepsilon_3,\varepsilon_4 }^{\prime}\subset Y_{\varepsilon_1, \varepsilon_2,\varepsilon_3,\varepsilon_4 } $ by changing the inequality $|v_1|\le 1+\varepsilon_1$ to the equality
$|v_1|=1+\varepsilon_1$ (all other inequalities remain unchanged). There is a natural projection $pr_{Y\to A}$ (resp. its restriction $pr_{Y'\to A}$) of $Y$
(resp. of $Y^{\prime}$) to the annulus
$A=\{w|\,\varepsilon_4\le |w|\le \varepsilon_3\}$. Let us now consider the complex
$$Cone((pr_{Y\to A})_\ast(\Q_l)\to  (pr_{Y'\to A})_\ast(\Q_l))[-1]\,.$$ 
This is a {\it lisse} sheaf on the annulus, i.e. a continuous $l$-adic representation of the fundamental group of the $\k((t))$-analytic space $A$. There is a tautological embedding $\k((w))\to {\cal O}^{an}(A)$. It induces the homomorphism of profinite groups $\pi_1(A)\to
\op{Gal}(\overline{\k((w))}/\k((w)))$. Then one can show that the complex 
$$Cone((pr_{Y\to A})_\ast(\Q_l)\to  (pr_{Y'\to A})_\ast(\Q_l))[-1]\simeq$$
$$\simeq (pr_{Y\to A})_\ast(Cone(\Q_l\to (i_{Y'\hookrightarrow Y})_*\Q_l))$$  on $A$ is quasi-isomorphic to the pull-back of the complex $\R\Gamma_c(i^{\ast}_{V_1}\R\psi_{\pi}(\Q_l))$ of $\op{Gal}(\overline{\k((w))}/\k((w)))$-modules.

We decompose the space $Y$ into a disjoint union $Y_0\sqcup Y_1$ where for $Y_0$ we have $|v_1||v_2|=0$ while for $Y_1$ we have $|v_1||v_2|\ne 0$.
Similarly, we have a decomposition $Y^{\prime}=Y_0^{\prime}\sqcup Y_1^{\prime}$. We claim
 that the complex 
$$(pr_{Y_0\to A})_\ast(Cone(\Q_l\to (i_{Y'_0\hookrightarrow Y_0})_*\Q_l))$$
is quasi-isomorphic to the pull-back of the complex $\R\Gamma_c(V_1,\Q_l)\otimes j^{\ast}\R\psi_{\widehat{\pi}}(\Q_l)$
of $\op{Gal}(\overline{\k((w))}/\k((w)))$-modules. Notice that $W_{|Y_0}$ depends on $v_3$ only.
Furthermore, $Y_0$ and $Y_0^{\prime}$ can be decomposed as the products
$$Y_0=Y_3\times Z_0,\,\,Y_0^\prime=Y_3\times Z_0^{\prime}\,\,,$$
where
$$\begin{array}{l}
Z_0:=\{(v_1,v_2)\in V_1^{an}\times V_2^{an}\,|\,\,v_1=0\, \,\mbox{or}\, \,v_2=0,\, |v_1|\le 1+\varepsilon_1, \,|v_2|\le\varepsilon_2\}\,\,,\\
Z_0^{\prime}:=\{(v_1,v_2)\in V_1^{an}\times V_2^{an}\,|\,\,v_2=0,\,|v_1|= 1+\varepsilon_1\}\,\,,\\
Y_3:=\{v_3\in V_3^{an}\,|\,\,|v_3|\le \varepsilon_2,\,\,\varepsilon_4\le |W(0,0,v_3)|\le \varepsilon_3\}\,\,.
\end{array}$$
Here we denote by $V_i^{an}, i=1,2,3$ the $\k((t))$-analytic space associated with the scheme $V_i\times_{Spec(\k)}Spec(\k((t)), i=1,2,3$.
Notice that $Z_0$ is the bouquet of two (non-archimedean) balls. It follows that  the inclusion of the ball 
$$Z_0^{\prime\prime}=\{(v_1,0)\,|\,|v_1|\le 1+\varepsilon_1\}$$ into $Z_0$ induces isomorphisms of the Berkovich \'etale cohomology groups of the $l$-adic sheaves on the analytic spaces. Therefore the cohomology groups of the pair $(Z_0,Z_0^{\prime})$ coincide with the cohomology groups of the pair $(Z_0^{\prime\prime},Z_0^{\prime})$.
The latter  are
equal to $\R\Gamma_c(V_1^{an},\Q_l)$ (which corresponds to the image of ${\mathbb L}^{\dim V_1}$ in $K_0(\op{Gal}(\overline{\k((t))}/\k((t)))-mod)$).

We have an obvious morphism of complexes of sheaves on $Y$:
$$ f:Cone(\Q_l\to (i_{Y'\hookrightarrow Y})_*\Q_l) \to 
Cone ((i_{Y_0\hookrightarrow Y})_*\Q_l\to (i_{Y'_0\hookrightarrow Y})_*\Q_l)\,\,.$$
In order to prove the Proposition we have to prove that the $(pr_{Y\to A})_*f$ is a quasi-isomorphism, i.e.
 $(pr_{Y\to A})_*(Cone(f))$ is zero. The compactness of spaces $Y,Y',Y_0,Y_0'$ implies that 
$$(pr_{Y\to A})_*(Cone(f))\simeq (pr_{Y_1\to A})_!(Cone(\Q_l\to (i_{Y_1'\to Y_1})_\ast \Q_l ))\,\,. $$

The (partially defined) actions of the group ${\bf G}_m$ on $Y_1$ and $Y_1^{\prime}$ are free, and the value of $W$ does not change under the action. 
More precisely, one can define easily analytic ``spaces of orbits'' $\widetilde{Y}_1\supset \widetilde{Y}_1^{\prime}$  
of ${\bf G}_m$ acting on $Y_1$ and $Y_1^{\prime}$ respectively. The projections
$$Y_1\to \widetilde{Y}_1,\,\,\, Y_1'\to \widetilde{Y}_1^{\prime}$$
are proper maps, and the map $W$ factors through them.
  Hence it is enough to check that
  $$(pr_{Y_1\to \widetilde{Y}_1})_*(Cone(\Q_l\to (i_{Y_1'\hookrightarrow Y_1})_*\Q_l))\simeq 0\,\,.    $$
This follows from the fact that every orbit in $Y_1$ is a closed annulus, its intersection with $Y_1^{\prime}$ is a circle,
 and the inclusion of a circle into an annulus induces an isomorphism of \'etale cohomology groups.

 This concludes the sketch of the proof. $\blacksquare$

\begin{rmk} 1) In the proof we did not specify the values of $\varepsilon_i, i=1,2,3,4$. We can take $\varepsilon_4=O(\varepsilon_3)$ (e.g. take $\varepsilon_4=\frac{1}{2}\varepsilon_3$), $\varepsilon_2=o(1)$ and $\varepsilon_3=O(\varepsilon_2^N),\varepsilon_1=O(\varepsilon_3^M)$ for some integers $N,M>0$.

2) In the proof we used the comparison of the cohomology of the sheaf of nearby cycles with the \'etale cohomology of subvarieties of $\k((t))$-analytic spaces (see \cite{NicS}).
\end{rmk}

We strongly believe that the analog of the Proposition holds at the level of motivic rings in the case $char\,\k=0$.

\subsection{Equivalence relation on motivic functions}

We start with a motivation for this section. There are many examples of pairs of constructible sets (or even varieties) $X_1,X_2$ over a field $\k$ such that their classes $[X_1]$ and $[X_2]$ in $Mot(Spec(\k))$ are different
(or at least not obviously coincide), but $X_1$ and $X_2$ coincide in each realization described in Section 4.1 (i)-(vi) (e.g.  when $X_1,X_2$ are isogeneous abelian varieties). In particular we will be interested in the case when
$X_l,l=1,2$ are affine quadrics given by the equations $\sum_{1\le i\le n}a_{i,l}x_i^2=1$,
such that they have equal determinants: $\prod_{1\le i\le n}a_{i,1}=\prod_{1\le i\le n}a_{i,2}\in \k^{\times}$.

Here we propose a modification of the notion of motivic function which is a version of the Grothendieck's approach to the theory of pure motives with numerical equivalence. Let us explain it in the case of $X=Spec(\k)$, where $\k$ is a field of characteristic zero (which we will assume throughout this section).

We start with the symmetric monoidal $\Q$-linear category ${\cal M}^{eff}(\k)$. Its objects are smooth projective varieties over $\k$ and
$$\op{Hom}_{{\cal M}^{eff}(\k)}(Y_1,Y_2)=\Q\otimes \op{Im}(Z^{\dim Y_2}(Y_1\times Y_2)\to H_{DR}^{2\dim Y_2}(Y_1\times Y_2))\,\,,$$
where $Z^n(X)$ denotes as usual the space of algebraic cycles in $X$ of codimension $n$, and we take the image of the natural map into the algebraic de Rham cohomology. Then $\op{Hom}_{{\cal M}^{eff}(\k)}(Y_1,Y_2)$ is a finite-dimensional $\Q$-vector space.
Instead of de Rham cohomology we can use Betti cohomology (for an embedding $\k\hookrightarrow \C$) or $l$-adic cohomology. Comparison theorems imply that the image of the group of cycles in the cohomology does not depend on a cohomology theory.

Composition of morphisms is given by the usual composition of correspondences, and the tensor product is given by the Cartesian product of varieties. We extend the category ${\cal M}^{eff}(\k)$ by adding formally finite sums (then it becomes an additive category), and finally taking the Karoubian envelope. The $K_0$-ring of the resulting category
contains the element ${\mathbb L}=[{\bf P}_{\k}^1]-[Spec(\k)]$. Adding formally the inverse ${\mathbb L}^{-1}$ we obtain the ring which we denote by $Mot_{coh}(Spec(\k))$. It is an easy corollary of Bittner theorem (see \cite{Bit}) that the natural map $Mot(Spec(\k))\to Mot_{coh}(Spec(\k))$ which assigns to a smooth projective variety its class in
$Mot_{coh}(Spec(\k))$ is a homomorphism of rings.

The above considerations can be generalized to the case of motives over constructible sets.

\begin{defn} Let $X$ be a constructible set over a field $\k, char\,\k=0$.
A constructible family of smooth projective varieties over $X$ is represented by a pair consisting of a smooth projective morphism
$h:Y\to X_0$ of schemes of finite type over $\k$ and a constructible isomorphism $j:X_0\stackrel{\mathrm constr}{\simeq} X$.
Two such representations
$$h:Y\to X_0,\,\,j:X_0\stackrel{\mathrm constr}{\simeq}X,\,\,\,\;\,h':Y'\to X_0',\,\,j':X_0'\stackrel{\mathrm constr}{\simeq}X$$
are identified if  we are given constructible isomorphisms $f:Y \stackrel{\mathrm constr}{\simeq}Y^{\prime}$, 
$g: X_0\stackrel{\mathrm constr}{\simeq}X_0^{\prime}$ such that $h^{\prime}\circ f=g\circ h,\,j'\circ g=j$, and for any point
  $x\in X(\overline{\k})$ the induced constructible isomorphism between smooth projective varieties $(j\circ h)^{-1}(x)$ and
$(j'\circ h')^{-1}(x)$ is an isomorphism of schemes.
\end{defn}
For a constructible family of smooth projective varieties over $X$ and a point $x\in X$ there is a well-defined smooth projective variety $Y_x$ over the residue field $\k(x)$ called {\it the fiber over $x$}.  Moreover, one can define the notion of constructible family of algebraic cycles of the fixed dimension. We say that such a family is homologically equivalent to zero if for any $x\in X$ the corresponding cycles in $Y_x$ map to zero in $H_{DR}^{\bullet}(Y_x)$. Also, having two constructible families of smooth projective varieties over $X$ one easily defines their product, which is again a  constructible family of smooth projective varieties over $X$. All that allows us to generalize our constructions from the case $X=Spec(\k)$ to the general case. In this way we obtain the ring $Mot_{coh}(X)$ as well as the natural homomorphism of rings $Mot(X)\to Mot_{coh}(X)$.

\begin{defn} We say that two elements of $Mot(X)$ are (cohomologically) equivalent if their images in $Mot_{coh}(X)$ coincide.

\end{defn}

The set of equivalence classes (in fact the ring) will be denoted by $\overline{Mot}(X)$. It is isomorphic to the image of $Mot(X)$ in $Mot_{coh}(X)$. In particular, the above-mentioned quadrics define the same element in $\overline{Mot}(Spec(\k))$.

Let now $X$ be a constructible set over $\k$, endowed with an action of an affine algebraic group $G$. We define an equivalence relation on $Mot^G(X)$ in the following way. First we choose an embedding $G\hookrightarrow GL(N)$.
We say that $f,g\in Mot^G(X)$ are equivalent if their pull-backs to $Mot((X\times GL(N))/G)$ have the same image in
$\overline{Mot}(X\times GL(N))/G)$. Using the fact that all $GL(N)$-torsors over a constructible set are trivial it is easy to show that this equivalence relation does not depend on the embedding $G\hookrightarrow GL(N)$.
The ring of equivalence classes is denoted by $\overline{Mot}^G(X)$.\footnote{In this way we have circumvented the problem of defining the category of $G$-equivariant motivic sheaves.} Similarly one defines the ring $\overline{Mot}^{G\times \mu}(X)$, where $\mu=\varprojlim_n\mu_n$. The exotic product descends to $\overline{Mot}^{G\times \mu}(X)$. Hence we obtain the ring $\overline{{\cal M}}^{G,\mu}(X)$ of equivalence classes as well as homomorphism of rings ${\cal M}^{G,\mu}(X)\to \overline{{\cal M}}^{G,\mu}(X)$.

\subsection{Numerical realization of motivic functions}

This section is not used in further consideration, its goal is only to show that  the abstractly defined ring 
$\overline{{\cal M}}^{G,\mu}(X)$ can be realized as certain ring of functions with numerical values. 

Let $Z$ be a scheme of finite type over a finite field $\k\simeq{\bf F}_q$ endowed with an action of the group $\mu_n$ of roots of $1$ such that $n<p=char \, {\bf F}_q$. We choose a prime $l\ne p$. Let us define $Y$ as the quotient
$(Z\times ({\bf A}^1_{\k}\setminus \{0\}))/\mu_n$ with respect to the diagonal action of $\mu_n$. Then we have a morphism 
$$\pi: Y\to {\bf A}^1_{\k}\setminus \{0\},\,\,\,(z,t)\mapsto t^n \,\,\forall (z,t)\in (Z\times {\bf A}^1_{\k}\setminus \{0\})(\overline{\bf F}_q)\,.$$
 Let $j:{\bf A}^1_{\k}\setminus \{0\}\to {\bf A}^1_{\k}$ be the natural embedding. We define the number
$$N_Z=\op{Tr}_{Fr}({\cal F}(j_{\ast}\pi_{!}({\Q}_{l,Y}))_{|s=1})\in \overline{\Q_{l}}\,\,,$$
which is the trace of the Frobenius $Fr$ of the fiber at $s=1$ of the Fourier transform of the
$l$-adic sheaf $j_{\ast}\pi_{!}({\Q}_{l,Y})$, where ${\Q}_{l,Y}$ denotes the constant sheaf $\Q_l$ on $Y$.
In fact the number $N_Z$ can be considered as an element of the  cyclotomic field $\Q(\eta_p)$, where
$\eta_p$ is a primitive $p$-th root of $1$:
$$1+\eta_p+\dots+\eta_p^{p-1}=0\,\,.$$

We have a canonical non-trivial character $\chi:{\bf F}_q\to \Q(\eta_p)^{\ast}$ given by the composition
of the trace $\op{Tr}_{{\bf F}_q\to{\bf F}_p}:{\bf F}_q\to {\bf F}_p$ with the additive character
$$\chi_p:{\bf F}_p\to \Q(\eta_p)^{\ast}, \,\,\,m\op{mod} p\mapsto \eta_p^m\,.$$
Then
$$N_Z=\sum_{s\in ({\bf A}^1_{\k}\setminus \{0\})({\bf F}_q)}\#(\pi^{-1}(s)({\bf F}_q))\,\chi(s)\,\,.$$
Notice that the last formula makes sense for constructible $Z$ as well.

Let now $X$ be a constructible set over a field $\k, char\,\k=0$, endowed with an action of an affine algebraic group $G$.

\begin{defn} We call a model for $(X,G)$ the following choices:
\begin{itemize}
 \item a finitely generated subring $R\subset \k$,
\item a scheme ${\cal X}\to Spec(R)$ of finite type,
\item an affine group scheme ${\cal G}\to Spec(R)$ together with an embedding ${\cal G}\hookrightarrow GL(N)_R$,
\item an action of ${\cal G}$ on ${\cal X}$,
\item a constructible identification over $\k$ of ${\cal X}\times_{Spec(R)}Spec(\k)$ with $X$, as well as an
isomorphism of groups ${\cal G}\times_{ Spec(R)}Spec(\k)\simeq G$ over $k$, compatible with the actions. 
\end{itemize}
\end{defn}
Such a model always exists, and models form a filterted system.
With a given model for $(X,G)$ we associate a commutative unital ring
$$K({\cal X})=\varinjlim_{\{\mathrm{open}\,\,\,U\subset ({\cal X}\times GL(n))/{\cal G}\}}\prod_{\{\mathrm{closed}\,\,\,x\in U\}}\Q(\eta_{char\, \k(x)})\,\,,$$
where $\k(x)$ is the residue field of $x$ (which is a finite field).

Suppose that we are given a model for $(X,G)$ and
let $f\in {\cal M}^{G,\mu}(X)$. As always we assume that the $\mu$-action on $X$ is good and factors through
the action of some $\mu_n$.

\begin{defn}

A model for $f$ compatible with the model $(R,{\cal X},{\cal G}) $ for $(X,G)$ consists of the following data:
\begin{itemize}
\item a finite set $J$, numbers $N_j,d_j, n_j\in \Z$, where $j\in J$, such that all numbers $N_j$ are positive and invertible in
 ${\cal O}({\cal X})$,
\item ${\cal G} \times \mu_{N_j}$-equivariant morphisms of constructible sets
${\cal Y}_j\to {\cal X}$ given for each $j\in J$, where ${\cal Y}_j\to {\cal X}$ are ${\cal G}\times \mu_{N_j}$-schemes of finite type, and we endow ${\cal X}$ with the trivial action of the group $\mu_{N_j}$,
\end{itemize}
These data are required to satisfy the condition that 
$$f=\sum_{j\in J}\,n_j\,[{\cal Y}_j\times_{Spec(R)}Spec(\k)\to X]\cdot{\mathbb L}^{d_j}\,.$$
\end{defn}
Models for $f$ always exist. Moreover, for any finite collection $(f_i)_{I\in I}$ of
elements of ${\cal M}^{G,\mu}(X)$ there exists a  model for $(X,G)$ with compatible models for $(f_i)_{I\in I}$ .

Having a model for $f$ we can associate with it an element $f_{num}\in K({\cal X})$ in the following way.
Let $x\in ({\cal X}\times GL(N)_R)/{\cal G}$ be a closed point. We can apply the considerations of the beginning of this section to the fiber $Z_{j,x}$ over the point $x$, of the map $({\cal Y}_j\times GL(N)_R)/{\cal G}\to ({\cal X}\times GL(N)_R)/{\cal G}$, where ${\cal Y}_j$ is the scheme from the definition of the model for $f$. Then for each $j$ we obtain an 
element $N_{Z_{j,x}}\in \Q(\eta_p), p=char\,\k(x)$. Finally, we set
$$f_{num}(x):=\sum_{j\in J}n_j\,q_x^{d_j}N_{Z_{j,x}}\,\,,$$
where $q_x:=\#\k(x)$. Hence we realize $f$ as a function with values in numbers.

\section{Orientation data on odd Calabi-Yau categories}

Considerations of this section are reminiscent of those in Quantum Field Theory when one tries to define determinants for the Gaussian integral in a free theory.

\subsection{Remarks on the motivic Milnor fiber of a quadratic form}

Although  the theory of motivic Milnor fiber was defined over a field of characteristic zero, an
essential part of considerations below has meaning over an arbitrary field $\k, char\,\k\ne 2$ if we replace the notion of motivic Milnor fiber by its $l$-adic version.

Let now $V$ be a $\k$-vector space endowed with a non-degenerate quadratic form $Q$. We define an element
$${I}(Q)=(1-MF_0(Q)){\mathbb L}^{-{1\over{2}}\dim V}\in {\cal M}^{\mu}(Spec(\k))[{\mathbb L}^{\pm 1/2}]\,\,,$$
 where
${\mathbb L}^{1/2}$ is a formal symbol which satisfies the relation $({\mathbb L}^{1/2})^2={\mathbb L}$, and $Q$ is interpreted as a function on $V$.
Then the motivic Thom-Sebastiani theorem implies
that $${I}(Q_1\oplus Q_2)={I}(Q_1){I}(Q_2)\,\,.$$
 Also we have ${I}(Q)=1$, if $Q$ is a split form: $Q=\sum_{1\le i\le n}x_iy_i$
for $V=\k^{2n}$. Therefore, 
 we have a homomorphism of groups $$I:Witt(\k)\to ({\cal M}^{\mu}(Spec(\k))[{\mathbb L}^{\pm 1/2}])^{\times}\,\,,$$
 where $Witt(\k)$ is the Witt group of the field $\k$. We can think of it as a multiplicative character. Let us denote by $J_2(\k):=J_2(Spec(\k))$ the quotient of the group $\Z\times \k^{\times}/(\k^{\times})^2$
by the subgroup generated by the element $(2,-1)$. There is an obvious homomorphism $Witt(\k)\to J_2(\k)$ given for a quadratic form $Q$ by $$[Q]\mapsto (\op{rk}Q, \op{det}(Q)\op{mod}(\k^{\times})^2)\,\,.$$
 Notice that all ``motivic realizations" of ${I}(Q)$ in the sense of Section 4.1 depend only on the image of $[Q]$ in $J_2(\k)$. This is similar to the classical formula (for $\k=\R$ and a positive definite form $Q$)
$$\int_V\exp(-Q(x))dx=(2\pi)^{-{1\over{2}}\dim V}(\op{det}(Q))^{-1/2}$$ in the sense that the answer depends on $\dim V$ and
$\op{det}(Q)$ only.
In particular, the homomorphism of rings $${\cal M}^{\mu}(Spec(\k))\to K_0(D^b_{\op{constr}}({\bf G}_m, \Q_l))$$ (see Section 4.1, (ii), (iii) and Section 4.2) induces (by combining with the above character) a homomorphism of abelian groups
$$Witt(\k)\to (K_0(D^b_{\op{constr}}({\bf G}_m, \Q_l))[{\mathbb L}^{\pm 1/2}])^{\times}\,\,.$$ It is easy to see that it factors through the homomorphism
$Witt(\k)\to J_2(\k)$. For example the element of $K_0$-group corresponding to the pair $(n,a), n\in 2\Z+1, a\in \k^{\times}$ is represented by ${\mathbb L}^{-1/2}[F]$ where $F$ is a local system on ${\bf G}_m$ associated with the double cover of ${\bf G}_m$ given by $y\mapsto y^2a, y\in {\bf G}_m(\kk)$.

\begin{que}
The above considerations give rise to the following question. Let us consider the family of quadratic forms
$Q_{a_1,a_2}(x,y)=a_1x_1^2+a_2x_2^2$ where $a_1,a_2\in \k^{\times}$. Is it true that
$[Q_{a_1,a_2}]=[Q_{a_1^{\prime},a_2^{\prime}}]$ in $K_0(Var_{\k})$ as long as $a_1a_2=a_1^{\prime}a_2^{\prime}$?

\end{que}

We expect that the answer to the Question is negative, and this is the main reason for introducing the equivalence relation for the motivic functions in Section 4.5.

The above considerations can be generalized to arbitrary constructible (or ind-constructible) sets. Namely, let $X$ be a constructible set over $\k$.  We define the group $J_2(X)$ as the quotient of the group 
$$Constr(X,\Z)\times \bigl(Constr(X,{\bf G}_m)/Constr(X,{\bf G}_m)^2\bigr)$$
 by the subgroup consisting of the elements $(2f,(-1)^f), f\in Constr(X,\Z)$, where we denote
by $Constr(X,Y)$  the set of constructible maps $X\to Y$. To a constructible vector bundle $V\to X$ endowed with a non-degenerate quadratic form $Q=(Q_x)_{x\in X}$ we associate the element 
$${I}(V,Q):=(1-MF_0(Q)){\mathbb L}^{-{1\over{2}}\dim V}\in ({\cal M}^{\mu}(X)[{\mathbb L}^{\pm 1/2}])^{\times}\,\,.$$ Here we treat each $Q_x$ as a formal power series on the fiber $V_x$ and $\dim V\in Constr(X,\Z)$. One can show that this correspondence gives rise to a homomorphism of groups 
$$\overline{I}: J_2(X)\to (\overline{{\cal M}}^{\mu}(X)[{\mathbb L}^{\pm 1/2}])^{\times}\,\,.$$ 
This fact has a simple ``numerical'' counterpart: for the case $\k\simeq {\mathbf F}_q$, 
two affine quadrics given by equations $Q_1(x)=0, Q_2(x)=0$ have the same number of points if $Q_1,Q_2$ are two 
non-degenerate quadratic forms
 of equal rank and determinant.

Let us consider a symmetric monoidal category $sPic_2(X)$ consisting of constructible super line bundles $L\to X$ endowed with an isomorphism $L^{\otimes 2}\simeq {\bf 1}_X$, where ${\bf 1}_X$ is a trivial even line bundle on $X$. It is easy to see that the group $J_2(X)$ is the group of isomorphism classes of objects of
$sPic_2(X)$.
If $V\to X$ is a constructible super vector bundle, $V=V^{even}\oplus \Pi V^{odd}$ then there is a well-defined super line bundle (called super determinant bundle) $\op{sdet}(V)\to X$ given by 
$$\op{sdet}(V):=\Pi^{\dim V^{even}-\dim V^{odd}}\left(\wedge^{top}V^{even}\otimes (\wedge^{top}V^{odd})^{\ast}\right)\,\,,$$ where
$\Pi$ is the parity change functor.

Recall canonical isomorphisms:
\begin{enumerate}
\item $\op{sdet}(V^{\ast})\simeq \op{sdet}(\Pi(V))\simeq (\op{sdet}(V))^{\ast}$,
\item $\op{sdet}(\oplus_{i\in I}V_i)\simeq \otimes_{i\in I}\op{sdet}(V_i)$,
\item if $V$ carries an odd differential $d$ then $\op{sdet}(V)\simeq \op{sdet}(H^{\bullet}(V,d))$,
\item for an exact triangle of complexes $(V_1^{\bullet},d_1)\to (V_2^{\bullet},d_2)\to (V_3^{\bullet},d_3)$ there is a canonical isomorphism
$$\op{sdet}(V_2^{\bullet})\simeq \op{sdet}(V_1^{\bullet})\otimes \op{sdet}(V_3^{\bullet})\,\,.$$
\end{enumerate}

In case if $V$ carries a non-degenerate quadratic form $Q=(Q_x)_{x\in X}$ we have a canonical isomorphism
$\op{sdet}(V)\simeq (\op{sdet}(V))^{\ast}$. Therefore in this case we have a well-defined object of $sPic_2(X)$. Its class in the group $J_2(X)$ is represented by the pair 
$$(\dim V,\,\,\{\op{det}(Q_x)\}_{x\in X}\op{ mod}(Constr(X,{\bf G}_m))^2)\,\,.$$

The above considerations can be generalized to the case when
$X$ is acted by an algebraic group $G$. Then one replaces the category $sPic_2(X)$ by the category $sPic_2(X,G)$ of $G$-equivariant constructible super line bundles $L$ endowed with a $G$-equivariant isomorphism $L^{\otimes 2}\simeq {\bf 1}_X$. The group of isomorphism classes of  $sPic_2(X,G)$ we denote by $J_2(X,G)$. In what follows we will often omit the word ``equivariant" in the considerations involving the category $sPic_2(X,G)$.

\begin{rmk}
Let us make an additional assumption that $\sqrt{-1}\in \k$.
In this case the quadratic form $x^2+y^2=(x+\sqrt{-1}y)(x-\sqrt{-1}y)$ is split. Then $MF_0(x^2+y^2)=1-{\mathbb L}$ and we can consider the element ${\mathbb L}^{1/2}:=1-MF_0(x^2)$ which enjoys the property
$({\mathbb L}^{1/2})^2={\mathbb L}$. Furthermore, the group $J_2(X)$ can be canonically identified with the product
$$Constr(X,\Z/2\Z)\times Constr(X,{\bf G}_m)/(Constr(X,{\bf G}_m))^2\,\,.$$ Therefore the isomorphism classes of objects of $sPic_2(X)$ can be identified with pairs {\it (constructible $\mu_2$-function, constructible $\mu_2$-torsor)}.

\end{rmk}

\subsection{Orientation data}

Let $\CC$ be an ind-constructible $\k$-linear $3$-dimensional Calabi-Yau category.\footnote{There is a notion of $\Z/2\Z$-graded odd or even Calabi-Yau category, see \cite{KoSo3}. Some considerations of this section can be generalized to $\Z/2\Z$-graded case.}
Then we have a natural ind-constructible super line bundle ${\cal D}$ over $Ob(\CC)$ with the fiber over $E$ given by
${\cal D}_E=\op{sdet}(\op{Ext}^{\bullet}(E,E))$. It follows that on the ind-constructible stack of exact triangles $E_1\to E_2\to E_3$ we have an isomorphism of the pull-backed line bundles which fiberwise reads as
$${\cal D}_{E_2}\otimes {\cal D}_{E_1}^{-1}\otimes {\cal D}_{E_3}^{-1}\simeq (\op{sdet}(\op{Ext}^{\bullet}(E_1,E_3)))^{\otimes 2}\,\,.$$
Let us explain this isomorphism. The multiplicativity of superdeterminants on exact triangles gives rise to a canonical isomorphism
$$\op{sdet}(\op{Ext}^{\bullet}(E_2,E_2))\simeq \op{sdet}(\op{Ext}^{\bullet}(E_1,E_1))\otimes \op{sdet}(\op{Ext}^{\bullet}(E_1,E_3))\otimes $$
$$\otimes\op{sdet}(\op{Ext}^{\bullet}(E_3,E_1))\otimes \op{sdet}(\op{Ext}^{\bullet}(E_3,E_3))\,\,.$$
By the Calabi-Yau property we have
$$\op{sdet}(\op{Ext}^{\bullet}(E_3,E_1))\simeq \op{sdet}(\Pi(\op{Ext}^{\bullet}(E_3,E_1)))^{\ast}\simeq \op{sdet}(\op{Ext}^{\bullet}(E_1,E_3))$$
which implies the desired formula.
When $Ob(\CC)=\sqcup_{i\in I}Y_i$ is a decomposition into the union of $GL(N_i)$-invariant constructible sets as at the end of 3.2, then the restriction ${\cal D}_{|Y_i}$ is a $GL(N_i)$-equivariant super line bundle and the above isomorphisms are also equivariant.

\begin{defn} Orientation data on $\CC$ consists of a choice of an ind-constructible super line bundle $\sqrt{{\cal D}}$ on $Ob(\CC)$ such that its restriction to each $Y_i, i\in I$ is $GL(N_i)$-equivariant, endowed on each $X_i$ with $GL(N_i)$-equivariant isomorphisms $(\sqrt{{\cal D}})^{\otimes 2}\simeq {\cal D}$ and such that for the natural pull-backs to the ind-constructible stack of exact triangles $E_1\to E_2\to E_3$ we are given equivariant isomorphisms:
$$\sqrt{{\cal D}}_{E_2}\otimes (\sqrt{{\cal D}}_{E_1})^{-1}\otimes (\sqrt{{\cal D}}_{E_3})^{-1}\simeq \op{sdet}(\op{Ext}^{\bullet}(E_1,E_3))$$
such that the induced equivariant isomorphism
$${\cal D}_{E_2}\otimes {\cal D}_{E_1}^{-1}\otimes {\cal D}_{E_3}^{-1}\simeq (\op{sdet}(\op{Ext}^{\bullet}(E_1,E_3)))^{\otimes 2}$$
coincides  with the one which we have a priori.

\end{defn}

We define the group $J_2(\CC):=\prod_{i\in I} J_2(Y_i,GL(N_i))$. We have a canonical equivariant super line bundle ${\cal D}_{\le 1}$ whose fiber at $E\in Ob(\CC)$ is
$${\cal D}_{\le 1,E}:=\op{sdet}(\tau_{\le 1}(\op{Ext}^{\bullet}(E,E)))\,\,,$$ where $\tau_{\le i}, i\in \Z$ denotes the standard truncation functor.
It is easy to see that we have an equivariant isomorphism ${\cal D}_{\le 1}^{\otimes 2}\simeq {\cal D}$. Then on the space of exact triangles $E_1\to E_2\to E_3$ we have an equivariant isomorphism of super line bundles fiberwise given by
$$({\cal D}_{{\le 1},E_2}\otimes {\cal D}_{{\le 1},E_1}^{-1}\otimes {\cal D}_{{\le 1},E_3}^{-1})^{\otimes 2}\simeq (\op{sdet}(\op{Ext}^{\bullet}(E_1,E_3)))^{\otimes 2}.$$

Let now $Fun({\CC}_3, \CC)$ be the ind-constructible category of $\A$-functors from the the category $\CC_3$ considered in 3.1. Its objects can be thought of as exact triangles $$E_1\to E_2\to E_3=Cone(E_1\to E_2)$$ in $\CC$. There are three functors $Funct({\CC}_3, \CC)\to \CC$ which associate to an exact triangle
$E_1\to E_2\to E_3$  the objects $E_1,E_2,E_3$ respectively. These functors induce three homomorphisms 
$$\phi_i: J_2(\CC)\to J_2(Funct(\CC_3, \CC)), i=1,2,3\,\,.$$
 The super line bundle $L$ with the fiber
$$L_{E_1\to E_2\to E_3}=({\cal D}_{{\le 1},E_2}\otimes {\cal D}_{{\le 1},E_1}^{-1}\otimes {\cal D}_{{\le 1},E_3}^{-1})\otimes (\op{sdet}(\op{Ext}^{\bullet}(E_1,E_3)))^{-1}$$ defines an element $l\in J_2(Funct({\cal A}_2, \CC))$, since
$L^{\otimes 2}\simeq {\bf 1}_{Funct(\CC_3, \CC)}$. Then a choice of orientation data on $\CC$ is equivalent to a choice of $h\in J_2(\CC)$ such that $-\phi_1(h)+\phi_2(h)-\phi_3(h)=l$. Indeed a choice of orientation data gives rise to a super line bundle $\sqrt{\cal D}$ such that $\sqrt{\cal D}^{\otimes 2}\simeq {\cal D}_{\le 1}^{\otimes 2}$. Therefore the quotient $h=\sqrt{\cal D}\otimes {\cal D}_{\le 1}^{-1}$ defines an element in $J_2(\CC)$, and the condition for the tensor squares of the super line bundles $\sqrt{\cal D}_{E_i}, i=1,2,3$ on the space of exact triangles is equivalent to the equation $-\phi_1(h)+\phi_2(h)-\phi_3(h)=l$.

\begin{rmk} All the above considerations admit a straightforward generalization to the case of Calabi-Yau category of arbitrary odd dimension $d$. In the case $d=1\,(\op{mod} 4)$ we have canonical orientation data given by
$$\sqrt{\cal D}_{\cal E}:=\op{sdet}(\tau_{\le {d-1\over{2}}}(\op{Ext}^{\bullet}(E,E)))\,\,.$$
 This is due to the observation that in the explicit description of the analog of the obstruction element $l$ defined above in terms of a super vector bundle endowed with a symmetric bilinear form, the super vector bundle turns out to be purely odd, hence the bilinear form is split. It follows that the obstruction element is trivial. In the case $d=3$ considered in this paper the obstruction does not have to vanish.

\end{rmk}

\subsection{Orientation data from a splitting of bifunctors}

Let $\cal C$ be a triangulated ind-constructive category over a field $\k$.  We will assume that all functors, bifunctors
etc. respect this structure. In this section we are going to discuss a special framework in which orientation data is easy to construct.

Let $F:\CC\times \CC^{op}\to Perf(Spec(\k))$ be a biadditive bifunctor and $d$ be an integer. We define the {\it dual bifunctor of degree $d$} as a bifunctor $F^{\vee}=F^{\vee,d}$ given by
$$F^{\vee}(E_2,E_1):=F(E_1,E_2)^{\ast}[-d]\,\,.$$
Clearly $F\mapsto F^{\vee}$ is an involution.
\begin{defn} A self-duality structure on $F$ of degree $d$ is an isomorphism $F\to F^{\vee}$ of bifunctors such that for any two objects $E_1,E_2$ the induced non-degenerate pairing
$$F(E_1,E_2)\otimes F(E_1,E_2)\to \k[-d]$$ is symmetric on the level of cohomology 
$H^\bullet(F(E_1,E_2))$. If $F$ is endowed with a self-duality structure of degree $d$ then we call it self-dual.

\end{defn}

For a Calabi-Yau category of dimension $d$  the bifunctor $(E_1,E_2)\mapsto \op{Hom}^{\bullet}(E_1,E_2)$ is self-dual.

For any self-dual bifunctor $F$ of odd degree $d$ we can repeat considerations of Section 5.2. Namely, we define an ind-constructible super line bundle ${\cal D}^F$ with the fiber  ${\cal D}^F_E:=\op{sdet}(F(E,E))$. Then for any exact triangle $E_1\to E_2\to E_3$ we have a canonical isomorphism

$$\sqrt{{\cal D}}_{E_2}^F\otimes (\sqrt{{\cal D}}_{E_1}^F)^{-1}\otimes (\sqrt{{\cal D}}_{E_3}^F)^{-1}\simeq \op{sdet}(F(E_1,E_3))^{\otimes 2}.$$

Then one can ask the same question: {\it is there an ind-constructible super line bundle $\sqrt{{\cal D}^F}$ which is compatible with the above isomorphism in the sense of Definition 15?}

The answer is positive for any bifunctor of the form $F\simeq H\oplus H^{\vee}$, with the obvious self-duality structure. In this case we set $$\sqrt{{\cal D}^F}_E:=\op{sdet}(H(E,E))\,\,.$$ More generally we can use an ${\bf A}^1$-homotopy in  this special case. More precisely, suppose we are given a bifunctor\footnote{In fact we would like to say that $G$ is ``ind-constructible'' in some sense. A sufficient, but not necessary condition would be
the existence of an ind-constructible functor $G'$ from $\CC\times \CC^{op}$ to $Perf({\bf P}^1_{\k})$ such that $G$ is isomorphic to the composition of
 $G'$ and the restriction functor $Perf({\bf P}^1_{\k})\to Perf({\bf A}^1_{\k})$.} $G: \CC\times \CC^{op}\to Perf({\bf A}^1_{\k})$. It can be thought of as a family $G_t: \CC\times \CC^{op}\to Perf(Spec(\k))$ of bifunctors, parametrized by $t\in {\bf A}^1(\k)$, namely $G_t=i_t^{\ast}\circ G$, where $i_t: Spec(\k)\to {\bf A}^1$ is the embedding corresponding to $t$. Since the category $Perf({\bf A}^1_{\k})$ has an obvious duality functor (taking dual to a complex of vector bundles) then the definition of self-duality structure extends naturally to families. Suppose that we have a family of self-dual bifunctors $G_t, t\in {\bf A}^1(\k)$ such that $G_0\simeq \op{Hom}^{\bullet}(\bullet,\bullet)$ and $G_1\simeq H\oplus H^{\vee}$ for some bifunctor $H$, and the isomorphisms preserve the self-duality structures. Then we have a canonical orientation data on $\CC$, since any super line bundle over ${\bf A}^1_{\k}$ is trivial and all fibers are canonically isomorphic.

\section{Motivic Donaldson-Thomas invariants}

\subsection{Motivic Hall algebra and stability data}

In this section the field $\k$ can have arbitrary characteristic.

Let $\CC$ be an
ind-constructible triangulated  $\A$-category over a field $\k$. We are going to describe a motivic generalization
of the derived Hall algebras from \cite{To}. 

As usual, we have  a constructible countable decomposition  $Ob(\CC)=\sqcup_{i\in I}Y_i$ with group $GL(N_i)$ acting on $Y_i$.
Let us consider a $Mot(Spec(\k))$-module $\oplus_i Mot_{st}(Y_i,GL(N_i))$ (see section 4.2) and extend it by adding negative powers ${\mathbb L}^n, n<0$ of the motive of the affine line  ${\mathbb L}$. We denote the resulting $Mot(Spec(\k))$-module by $H(\CC)$.  We understand elements of 
$H(\CC)$ as {\it measures} (and not as functions), because in the definition of the product we will use the pushforward maps.

We would like to make $H(\CC)$ into an associative algebra, called the {\it motivic Hall algebra}. We need some preparations for that. First we observe that if $[\pi_i: Z_i\to Ob(\CC)], i=1,2$ are two elements\footnote{Here we consider for simplicity the case when the groups acting on $Z_1,Z_2$ are trivial, the generalization to the case of non-trivial groups is straightforward.}
of $H(\CC)$ then one has a constructible set
$tot((\pi_1\times \pi_2)^{\ast}({\cal EXT}^1))$ which is the total space of the pull-back
of the ind-constructible bundle ${\cal EXT}^1$ over $Ob(\CC)\times Ob(\CC)$.
Then the map ${{\cal C}one}$ (see Section 3.1) after the shift $[1]$  maps the total space to $Ob(\CC)$.

For any $N\in \Z$ we introduce the ``truncated" Euler characteristic 
$$(E,F)_{\le N}:=\sum_{i\le N}(-1)^i\dim \op{Ext}^i(E,F)\,\,.$$
In the future we will use the notation $(E,F)_i$ for $\dim \op{Ext}^i(E,F)$, hence
$(E,F)_{\le N}=\sum_{i\le N}(-1)^i (E,F)_i$.

With the pair $[\pi_i: Y_i\to X_i], i=1,2$ as above we can associate a collection of
constructible sets
$$W_n  =\big\{(y_1,y_2,\alpha)\,|\,y_i\in Y_i, \alpha\in \op{Ext}^1(\pi_2(y_2),\pi_1(y_1))\,,\,
 (\pi_2(y_2),\pi_1(y_1))_{\le 0}=n\big\}\,\,, $$
 where $n\in \Z$ is arbitrary.
Clearly $$[tot((\pi_1\times \pi_2)^{\ast}({\cal EXT}^1))\to Ob(\CC)]=\sum_{n\in \Z}[W_n\to Ob(\CC)]\,\,.$$
We define the product
$$[Y_1\to Ob(\CC)]\cdot [Y_2\to Ob(\CC)]=\sum_{n\in \Z}[W_n\to Ob(\CC)]{\mathbb L}^{-n}\,\,,$$
where the map $W_n\to Ob(\CC)$ is given by the formula
$$(y_1,y_2,\alpha)\mapsto Cone(\alpha: \pi_2(y_2)[-1]\to \pi_1(y_1))\,\,.$$

\begin{prp} The above formula makes $H(\CC)$ into an associative algebra.

\end{prp}
{\it Proof.} We are going to prove the result for the ``delta functions" 
$$\nu_E=[pt\to Ob(\CC)], pt\mapsto E\,\,,$$ where $E$ is an object of $\CC(\kk)$. The case of equivariant families is similar. In other words, we would like to prove that
$$(\nu_{E_1}\cdot\nu_{E_2})\cdot \nu_{E_3}=\nu_{E_1}\cdot(\nu_{E_2}\cdot \nu_{E_3})\,\,.$$ 
 Replacing the category by its minimal model we may replace in all considerations $\op{Hom}^{\bullet}$ by $\op{Ext}^{\bullet}$. Let us also remark that an element $\alpha\in \op{Ext}^1(E,F)$ defines an extension $E_{\alpha}$ which we can interpret as a deformation of the object $E\oplus F$ (the trivial extension). Therefore for any object $G$ the group $\op{Ext}^{\bullet}(G,E_{\alpha})$ is equal to the cohomology of the complex $(\op{Ext}^{\bullet}(G,E\oplus F),d_{\alpha})$, where $d_{\alpha}$ is the operator of multiplication (up to a sign) by $\alpha$.

Notice that
$$\nu_{E_1}\cdot\nu_{E_2}={\mathbb L}^{-(E_2,E_1)_{\le 0}}[\op{Ext}^1(E_2,E_1)\to Ob(\CC)]:={\mathbb L}^{-(E_2,E_1)_{\le 0}}\int_{\alpha\in \op{Ext}^1(E_2,E_1)}\nu_{E_{\alpha}}\,\,,$$
where $E_{\alpha}$ is the object corresponding to the extension $\alpha$, i.e. 
$$E_{\alpha}=Cone(\alpha: E_2[-1]\to E_1)\,\,.$$
It follows that
$$(\nu_{E_1}\cdot\nu_{E_2})\cdot \nu_{E_3}={\mathbb L}^{-(E_2,E_1)_{\le 0}}\int_{
\alpha\in \op{Ext}^1(E_2,E_1),\,\beta\in \op{Ext}^1(E_3,E_{\alpha}) }{\mathbb L}^{-(E_3,E_{\alpha})_{\le 0}}\nu_{E_{\beta}}\,\,.$$

We observe that $$(E_3,E_{\alpha})_{\le 0}=(E_3,E_2\oplus E_1)_{\le 0}-l_{\alpha}=(E_3,E_2)_{\le 0}+ (E_3,E_1)_{\le 0}-l_{\alpha}\,\,,$$
 where the ``error term" $l_{\alpha}\geqslant 0$ can be computed in terms of the linear map $d_{\alpha}$.
Therefore one can write
$$(\nu_{E_1}\cdot\nu_{E_2})\cdot \nu_{E_3}={\mathbb L}^{-(E_2,E_1)_{\le 0}-(E_3,E_1)_{\le 0}-(E_3,E_2)_{\le 0}}\int_{\alpha\in \op{Ext}^1(E_2,E_1),\beta\in \op{Ext}^1(E_3,E_{\alpha})}{\mathbb L}^{l_{\alpha}}\nu_{E_{\beta}}\,\,.$$

One can write a similar expression for $\nu_{E_1}\cdot(\nu_{E_2}\cdot \nu_{E_3})$. In this case the ``error term"  will be denoted by $r_{\alpha}$ instead of $l_{\alpha}$.

Notice that the differential $$d_{\alpha}: \op{Ext}^0(E_3,E_2)\oplus \op{Ext}^0(E_3,E_1)\to \op{Ext}^1(E_3,E_1)\oplus \op{Ext}^1(E_3,E_2)$$ satisfies the property that the only non-trivial component is the map
$\alpha_R: \op{Ext}^0(E_3,E_2)\to \op{Ext}^1(E_3,E_1)$. Here we  denote by $\alpha_R$ the linear operator of multiplication by $\alpha\in \op{Ext}^1(E_2,E_1)$ from the right. We will use the same convention for the linear operator $\alpha_R: \op{Ext}^1(E_3,E_2)\to \op{Ext}^2(E_3,E_1)$.
Hence we see that
$$\dim \op{Ext}^1(E_3,E_{\alpha})= \dim \op{Ker}\left(\alpha_R: \op{Ext}^1(E_3,E_2)\to \op{Ext}^2(E_3,E_1)\right)+$$
$$+\dim \op{Coker}\left(\alpha_R: \op{Ext}^0(E_3,E_2)\to \op{Ext}^1(E_3,E_1)\right)\,\,.$$
Let us now consider the constructible set
$$X_{1,2,3}=\{(\alpha,\gamma, \delta)\in \op{Ext}^1(E_2,E_1)\oplus \op{Ext}^1(E_3,E_2)\oplus \op{Ext}^1(E_3,E_1)|\alpha\circ \gamma=0\}\,.$$

Notice that a triple $(\alpha,\gamma, \delta)\in X_{1,2,3}$ defines the deformation of the object $E_1\oplus E_2\oplus E_3$ preserving the  filtration
$$E_1\subset E_1\oplus E_2\subset E_1\oplus E_2\oplus E_3\,\,.$$ More precisely, the triple gives rise to a twisted complex, which is defined by the corresponding to $(\alpha,\gamma, \delta)$ solution to the Maurer-Cartan equation (strictly upper-triangular matrix acting on $E_1\oplus E_2\oplus E_3$). The latter observation means that there is an ind-constructible map $X_{1,2,3}\to Ob(\CC)$ which assigns to a point $(\alpha,\gamma, \delta)$ the corresponding twisted complex.

Let us now fix $\alpha\in \op{Ext}^1(E_2,E_1)$ and consider the ind-constructible subset
$X_{1,2,3}^{\alpha}\subset X_{1,2,3}$ which consists of the triples with fixed $\alpha$.
There is a natural projection $(\alpha,\gamma, \delta)\mapsto (\gamma,\delta)$,
which gives rise to the constructible map $X_{1,2,3}^{\alpha}\to \op{Ext}^1(E_3,E_{\alpha})$. 
This is a constructible affine bundle with the fibers isomorphic to $\op{Im}(\alpha_R: \op{Ext}^0(E_3,E_2)\to \op{Ext}^1(E_3,E_1))$. Also, one can see directly  that the dimension of the latter space is 
$$l_{\alpha}=(E_3,E_2\oplus E_1)_{\le 0}- (E_3,E_\alpha)_{\le 0}\,\,. $$ 
Hence we have the following identity in $H(\CC)$:
$$\int_{  \alpha\in \op{Ext}^1(E_2,E_1),\, \beta\in \op{Ext}^1(E_3,E_{\alpha})}{\mathbb L}^{l_{\alpha}}\nu_{E_{\beta}}=[X_{1,2,3}\to Ob(\CC)]\,\,.$$
Therefore, 
$$(\nu_{E_1}\cdot\nu_{E_2})\cdot \nu_{E_3}={\mathbb L}^{-(E_2,E_1)_{\le 0}-(E_3,E_1)_{\le 0}-(E_3,E_2)_{\le 0}}\cdot[X_{1,2,3}\to Ob(\CC)]\,\,.$$
Similar considerations show that
$$\nu_{E_1}\cdot(\nu_{E_2}\cdot \nu_{E_3})=  {\mathbb L}^{-(E_2,E_1)_{\le 0}-(E_3,E_1)_{\le 0}-(E_3,E_2)_{\le 0}}\cdot[X_{1,2,3}\to Ob(\CC)]\,\,.$$
This proves the associativity of the product in $H(\CC)$.
$\blacksquare$

For a constructible stability structure on $\CC$ with an ind-constructible class map $\op{cl}:K_0(\CC)\to \Gamma$, a central charge $Z:\Gamma\to \C$, a strict sector $V\subset \R^2$ and a branch $\op{Log}$ of the logarithm function on $V$ we have the category $\CC_V:=\CC_{V,\op{Log}}$ defined in Section 3.4. Hence we have
the {\it completion} 
$$\widehat{H}(\CC_V):=\prod_{\gamma\in (\Gamma\cap C(V,Z,Q))\cup\{0\}}H(\CC_V\cap \op{cl}^{-1}(\gamma))\,\,.$$

Then we have an invertible element $A_V^{\op{Hall}}\in \widehat{H}(\CC_V)$ such that
$$A_V^{\op{Hall}}:=1+\dots=\sum_{i\in I}  \mathbf{1}_{(Ob(\CC_V)\cap Y_i,GL(N_i))}\,\,,$$
where $\mathbf{1}_{\cal S}$ is the identity function (see 4.2) but interpreted as a counting measure\footnote{The same is true if one uses the language of higher stacks because for any $E\in Ob(\CC(V))(\kk)$ one has $\op{Ext}^{<0}(E,E)=0$.}.
In short, element $A_V^{\op{Hall}}$ is given by the counting measure restricted to
$\CC_V$. The summand $1$ comes from zero object.

\begin{prp} Elements $A_V^{\op{Hall}}$ satisfy the Factorization Property:
$$A_{V}^{\op{Hall}}=A_{V_1}^{\op{Hall}}\cdot A_{V_2}^{\op{Hall}}$$
for a strict sector $V=V_1\sqcup V_2$ (decomposition in the clockwise order).
\end{prp}

{\it Proof.} The proof follows from the following observations:

1) For any $E_i\in Ob(\CC_{V_i}(\kk)), i=1,2$ one has $(E_2,E_1)_{\le 0}=\dim\op{Ext}^0(E_2,E_1)$ because $\op{Ext}^i(E_2,E_1)=0$ for $i< 0$.

2) The set $\{[E]\in Iso(\CC_V(\kk))\}$ is in one-to-one correspondence with the set of isomorphisms classes of the triples $(E_1,E_2,\alpha)$ such that $E_i\in Ob(\CC_{V_i}(\kk)), i=1,2$ and $\alpha\in \op{Ext}^1(E_2,E_1)$ (the map between the sets assigns to the triple the extension $E_{\alpha}$).

3) The automorphism group of the triple $(E_1,E_2,\alpha)$ is the stabilizer of $\alpha$ for the natural action of the group $\op{Aut}(E_2)\times \op{Aut}(E_1)$ on the vector space $\op{Ext}^1(E_2,E_1)$.

4) There is an exact sequence of groups
   $$1\to \op{Ext}^0(E_2,E_1)\to \op{Aut}(E_\alpha)\to\op{Aut}( E_1,E_2,\alpha) \to 1   $$

In order to apply these observations one uses the fact that an object $E\in \CC_V(\kk)$ contains a unique subobject
$E_1\in Ob(\CC_{V_1}(\kk))$ such that the quotient object $E_2$ belongs to $\CC_{V_2}(\kk)$, and
then the factor ${\mathbb L}^{-(E_2,E_1)_{\le 0}}$ cancels the ratio between the stabilizer of $\alpha$ and the automorphism group of 
the extension $E_\alpha$. $\blacksquare$

\begin{cor} Let us endow $H(\CC)$ with an automorphism $\eta$ given by the shift functor $[1]$. Then the collection $(A_V^{\op{Hall}})$ gives rise to a symmetric
stability data on $H(\CC)$ considered as a graded Lie algebra (see  Definition 2 and Remark 8 in Section 2.2). Moreover we obtain a local homeomorphism $Stab(\CC,\op{cl})\to Stab(H(\CC))$.

\end{cor}

The above considerations can be illustrated in the case of finite fields. Namely, let us assume that $\CC$ is a triangulated category over a finite field ${\bf F}_q$. We define the Hall algebra $H(\CC)$ as an associative unital algebra over $\Q$, which is a $\Q$-vector space spanned by isomorphism classes $[E]$ of objects $E\in Ob(\CC)$.
The multiplication is given by the formula
$$[E]\cdot[F]=q^{-(F,E)_{\le 0}}\sum_{\alpha\in \op{Ext}^1(F,E)}[E_{\alpha}]\,\,,$$
where $E_{\alpha}$ is the extension corresponding to $\alpha\in \op{Ext}^1(F,E)$.

We define a stability condition on $\CC$ in the same way as in the Introduction (or Section 3.4) without imposing any constructibility condition (since we do not assume that our category is ind-constructible). Inside of the set $Stab(\CC)$ of stability conditions on $\CC$ we consider a subset $Stab^0(\CC)$ consisting of such stability conditions that the set $\{E\in \CC_{\gamma}^{ss}|\op{Arg}(E)=\varphi\}$ is finite for any $\gamma\in \Gamma, \varphi\in [0,2\pi)$. This property is analogous to  the one in the ind-constructible setting which says
 that $\CC_{\gamma}^{ss}$ is a constructible set.
Then for any strict sector $V$ and a choice of the branch $\op{Log}$ we have an element
$A_V^{\op{Hall}}\in \widehat{H}(\CC_V)$ given by
$$A_V^{\op{Hall}}=1+\dots=\sum_{[E]\in Iso(\CC_V)}{[E]\over {\#\op{Aut}(E)}}\,\,.$$
Similarly to the motivic case the collection of elements $(A_V^{\op{Hall}})$ satisfies the Factorization property. Hence it defines a stability data on the space $H(\CC)$ considered as a graded Lie algebra.

The relationship of our version of Hall algebra to the To\"en derived Hall algebra from \cite{To} is described in the following proposition.

\begin{prp} There is a homomorphism of rings $H(\CC)\to H_{To}(\CC)$,
where $H_{To}(\CC)$ is the derived Hall algebra over ${\bf F}_q$ defined by To\"en in \cite{To} (see also \cite{XiXu}),
such that
$$[\pi: Y\to Ob(\CC)]\mapsto \sum_{y\in Y({\bf F}_q)}[\pi(y)]\#\op{Aut}(y)({\bf F}_q)\,q^{(y,y)_{<0}}\,\,.$$
Moreover, for any strict sector $V$ the above homomorphism admits a natural
extension to the completed Hall algebras such that the element $A_V^{\op{Hall}}\in \widehat{H}(\CC_V)$ is mapped
to the element of the completed Hall algebra $\widehat{H}_{To}(\CC_V)$ given by $\sum_{[x]\in Iso(\CC_V)}[x]$.

\end{prp}

{\it Proof.} Straightforward.  $\blacksquare$

In fact, in the To\"en version of the Hall algebra  the factorization property 
$$A_V^{\op{Hall}}=A_{V_1}^{\op{Hall}}A_{V_2}^{\op{Hall}}$$ is essentially trivial. The reason is that the structure constants
in 
$\widehat{H}_{To}(\CC_V)$
for the elements of the basis corresponding to objects in a heart of a $t$-structure are the usual one, i.e.
 they count the number of 2-step filtrations of a given object with given isomorphism classes of the associate graded factors.
The factorization property means that any object in $\CC_V$ has a unique subobject in $\CC_{V_1}$ with quotient in $\CC_{V_2}$.

\begin{rmk}
One can try to go even further in an attempt to ``categorify" the motivic Hall algebra.
 Here one has to assume that objects of $\CC$ form not just an ind-constructible stack, but 
 a higher stack of locally finite type in the sense of To\"en and Vezzosi (see \cite{ToVe}).
The corresponding category will be the monoidal category of motivic sheaves on $Ob(\CC)$.
The motivic Hall algebra is the $K_0$-ring of this category. In the case of the non-commutative
variety endowed with polarization one can define (for any strict sector $V$) the subcategory
${\cal F}_V$ of ``motivic sheaves with central charges in $V$".
Nevertheless, the Factorization Property could fail since
the object $A_V^{\op{Hall}}$ can be non-isomorphic to  the object $A_{V_1}^{\op{Hall}}\otimes A_{V_2}^{\op{Hall}}$
(but their images in $K_0$ coincide).
\end{rmk}

Finally, we explain how to rephrase the factorization property in terms of $t$-structures,
without the use of stability conditions. 
Here we understand a $t$-structure $\alpha$
on a small triangulated category $\cal C$
as a pair of strictly full subcategories (i.e. a pair of sets of equivalence classes of objects)
$${\cal C}_{-,\alpha},\,\,{\cal C}_{+,\alpha}\subseteq {\cal C}$$
such that for any $E_-\in {\cal C}_{-,\alpha},\,\,E_+\in {\cal C}_{+,\alpha}$ we have
 $Ext^{\le 0}(E_-,E_+)=0$, and any object $E\in {\cal C}$ can be represented (uniquely) as an extension
 $$\tau_{-,\alpha}(E)\to E\to \tau_{+,\alpha}(E),\,\,\,\tau_{\pm,\alpha}(E)\in {\cal C}_{\pm,\alpha}\,\,.$$
 
 Any stability condition on $\cal C$ defines two $t$-structures $\alpha_l,\alpha_r$
  such that ${\cal C}_{-,\alpha_l}$ (resp. ${\cal C}_{-,\alpha_r}$) consists of extensions of semistable objects $E$ with $\op{Arg}(E)\geqslant 0$ (resp. with $\op{Arg}(E)>0$).
  These two $t$-structures do not change under the action of the group
  $$\left\{\left( 
  \begin{array}
  {cc} a_{11} & 0\\ a_{21} & a_{22}
  \end{array}\right) \bigl|\,\, a_{11},a_{22}>0\bigr.\right\}\subset \widetilde{GL_+}(2,\mathbf{R})$$
  of transformations preserving the upper half-plane.
  In particularly, we see that while a connected component
  in the $Stab({\cal C})$ is a real $2n$-dimensional manifold for $n:=\op{rank}(\Gamma)$,
   the set of corresponding $t$-structures is at most $(n-1)$-dimensional.

 Introduce an order on the set of $t$-structures by
 $$\alpha_1\le \alpha_2\,\,\Longleftrightarrow \,\,{\cal C}_{-,\alpha_1}\subseteq {\cal C}_{-,\alpha_2}\,\,\Longleftrightarrow\,\, {\cal C}_{+,\alpha_1}\supseteq {\cal C}_{+,\alpha_2}\,\,.$$
 The shift functor acts on $t$-structures, and 
 $\alpha[1]\le \alpha$ for any $t$-structure $\alpha$.

 Let now $\cal C$ be an ind-constructible category endowed with an 
 ind-constructible homomorphism $cl:K_0({\cal C}(\overline{\bf k}))\to \Gamma$
  and $\alpha_1,\alpha_2$
 are two ind-constructible $t$-structures. We say
  $$\alpha_1\le_{constr}\alpha_2$$
  iff 
  \begin{itemize}
  \item $\alpha_1\le \alpha_2\le \alpha_1[-1]\,\,,$
  \item $\forall\gamma\in \Gamma\,\,\,\,
    {\cal C}_{+,\alpha_1}\cap{\cal C}_{-,\alpha_2}\cap cl^{-1}(\gamma)$ is constructible,
    \item the cone generated by $\{\gamma\in \Gamma\,| \,{\cal C}_{+,\alpha_1}\cap{\cal C}_{-,\alpha_2}\cap cl^{-1}(\gamma)\ne 0\}$ is strict.
    \end{itemize}
    If $\alpha_1\le_{constr}\alpha_2$ then we define
   an element $A_{\alpha_1,\alpha_2}$ of an appropriately completed Hall
    algebra as the sum of the ``counting measure'' over the objects in
    ${\cal C}_{+,\alpha_1}\cap{\cal C}_{-,\alpha_2}$.
    Obviously, elements $A_V$ (for an open, or a closed, or a semi-open strict sector $V$)
    are of the form $A_{\alpha_1,\alpha_2}$ for  appropriate $t$-structures
     $\alpha_1,\alpha_2$. The factorization property
      generalizes to
      $$A_{\alpha_1,\alpha_3}=A_{\alpha_1,\alpha_2}\cdot A_{\alpha_2,\alpha_3}$$
      if $\alpha_1 \le_{constr}\alpha_2,\,\,\alpha_2\le_{constr}\alpha_3,\,\,\alpha_1\le_{constr}\alpha_3$.
      
      Notice that in the case of stability conditions the element $A_V$ 
    is preserved under the action of a subgroup of $\widetilde{GL_+}(2,\mathbf{R})$ conjugated to the group
     of positive diagonal matrices. This action on $Stab({\cal C})/Aut({\cal C})$
     has a good chance to be ergodic,
     as indicates a similar example with the moduli spaces of 
    curves with abelian differentials (see a review \cite{Zo}).

\subsection{Motivic weights and stability data on motivic quantum tori}

Let $\CC$ be a $3$-dimensional ind-constructible   Calabi-Yau category over a field of characteristic zero (see Section 3.3). 
In this section we are going to define motivic Donaldson-Thomas invariants associated with a constructible stability condition and an orientation data on $\CC$.

{\it Step 1}. 

Let us define the ring
$${D}^{\mu}={\cal M}^{\mu}(Spec(\k))[{\mathbb L}^{-1}, {\mathbb L}^{1/2},  ([GL(n)]^{-1})_{n\geqslant 1}]\,\,, $$
where  the ring ${\cal M}^{\mu}(Spec(\k))$ was defined in 4.3, and ${\mathbb L}=[{\bf A}^1_{\k}]$ is the motive of the affine line. The element
${\mathbb L}^{1/2}$ is a formal symbol satisfying the equation $({\mathbb L}^{1/2})^2={\mathbb L}$. Instead of inverting motives $$[GL(n)]=({\mathbb L}^n-1)({\mathbb L}^n-{\mathbb L})\dots({\mathbb L}^n-{\mathbb L}^{n-1})$$ of all general linear groups we can invert motives of all projective spaces
$$[{\bf P}^n]={{\mathbb L}^{n+1}-1\over {\mathbb L}-1}\,\,.$$  
We also will consider the ring $\overline{D^{\mu}}$ of equivalence classes of functions from $D^{\mu}$ by the equivalence relation defined in Section 4.5. 
 The ring $\overline{D^{\mu}}$ will play the role of the universal coefficient ring where motivic Donaldson-Thomas invariants take value.

{\it Step 2}. 

We define an algebra $\overline{{\cal M}^{\mu}}(Ob(\CC))$ associated with $\CC$ which will contain certain canonical element called the motivic weight.
First, we define 
$${\cal M}^{\mu}(Ob(\CC)):=\prod_i{\cal M}^{GL(N_i),\mu}(Y_i)[{\mathbb L}^{-1}, {\mathbb L}^{1/2}]\,\,,$$
where $(Y_i, GL(N_i))_{i\in I}$ is a decomposition of the stack of objects of $\CC$ as at the end of 3.2.
Algebra $\overline{{\cal M}^{\mu}}(Ob(\CC))$ is obtained from it by passing to the equivalence classes 
in the sense of Section 4.5. 

For any $GL(N_i)$-invariant constructible set $Z\subset Y_i$ for some $i\in I$, we have a ${\cal M}^\mu(Spec(\k))[{\mathbb L}^{-1}, {\mathbb L}^{1/2}]$-linear map
$$\int_Z: {\cal M}^{\mu}(Ob(\CC))\to {D}^{\mu} $$
which is the $\mu$-equivariant version of integral over stack $(Z,GL(N_i))$ (see 4.2) of the restriction  to $Z$.
 Explicitly, if $f_{|Y_i}$ is represented by a $\mu\times GL(N_i)$-equivariant map $X\to Y_i$ then
$$\int_Z f=[X\times_{Y_i} Z]/[GL(N_i)] \in D^\mu $$
where $[X\times_{Y_i} Z]$ is interpreted a constructible set with $\mu$-action.
By additivity we extend the integral to the case when $Z$ is a finite union of $GL(N_i)$-invariant constructible set $Z_i\subset Y_i$ for different $i\in I$.

{\it Step 3}.

Now we are going to define the motivic weight. 
Recall that for any $E\in Ob(\CC)(\kk)$ we have defined the potential $W_E^{min}$ which is a formal power series in $\alpha\in \op{Ext}^1(E,E)$ which starts with cubic terms. 
We denote by 
$$MF(E):=MF_0(W_E^{min})$$ the motivic Milnor fiber of $W_E^{min}$ at $0\in \op{Ext}^1(E,E)$. Then the assignment $E\mapsto MF(E)$ can be interpreted as the value of some function $MF\in {\cal M}^{\mu}(Ob(\CC))$.

Let us choose an orientation data $\sqrt{\cal D}$ for $\CC$. Recall that in Section 5.2 we defined the element $h\in J_2(\CC)$ represented by the equivariant super line bundle $\sqrt{\cal D}\otimes {\cal D}_{\le 1}^{-1}$ with trivialized tensor square. For a representative of $h$ given by a pair $(V,Q)$ we have $\overline{I}(h)=(1-MF_0(Q)){\mathbb L}^{-{1\over{2}}\op{rk} Q}$. Let us choose such a representative.

\begin{defn} The motivic weight $w\in {\cal M}^{\mu}(Ob(\CC))$ is the function defined on objects by the formula
$$w(E)={\mathbb L}^{{1\over{2}}\sum_{i\le 1}(-1)^i \dim \op{Ext}^i(E,E)}(1-MF(E))(1-MF_0(Q_E)){\mathbb L}^{{-{1\over{2}}\op{rk}Q_E}}\,\,.$$

\end{defn}

The image $\overline{w}\in \overline{{\cal M}^{\mu}}(Ob(\CC))$ does not depend on the choice of a representative of $h$ and is equal to
$$\overline{w}(E)={\mathbb L}^{{1\over{2}}\sum_{i\le 1}(-1)^i \dim \op{Ext}^i(E,E)}(1-MF(E))\overline{I}(h(E))\,\,,$$
where $h(E)$ is the value of the obstruction $h$ at the point $E$.

{\it Step 4}.

Let us now fix the following data:
\begin{itemize}
\item{a triple $(\Gamma, \langle\bullet,\bullet \rangle, Q)$ consisting of a free abelian group $\Gamma$
of finite rank endowed with a skew-symmetric bilinear form $\langle\bullet,\bullet \rangle: \Gamma\otimes \Gamma\to \Z$, and a quadratic form $Q$ on $\Gamma_{\R}=\Gamma\otimes \R$;}

\item{an ind-constructible, $\op{Gal}({\overline{\k}}/\k)$-equivariant homomorphism $$\op{cl}_{\kk}: K_0(\CC(\kk))\to \Gamma$$ compatible with the Euler form and the skew-symmetric bilinear form;}

\item{a constructible stability structure $\sigma\in Stab(\CC,\op{cl})$ compatible with the quadratic form $Q$ in the sense that $Q_{|\op{Ker}Z}<0$ and $Q(\op{cl}_{\kk}(E))\geqslant 0$ for $E\in \CC^{ss}(\overline{\k})$.}
\end{itemize}

In the next section we are going to define a homomorphism from the motivic Hall algebra to the associative unital algebra called motivic quantum torus. The latter is defined in the following way.

For any commutative unital ring $C$ which contains an invertible symbol ${\mathbb L}^{1/2}$ we introduce a $C$-linear associative algebra 
$${\cal R}_{\Gamma,C}:=\oplus_{\gamma\in \Gamma}C\cdot \hat{e}_{\gamma}$$
 where the generators $\hat{e}_{\gamma}, \gamma\in \Gamma$ satisfy the relations 
$$\hat{e}_{\gamma_1}\hat{e}_{\gamma_2}={\mathbb L}^{{1\over{2}}\langle \gamma_1,\gamma_2\rangle}\hat{e}_{\gamma_1+\gamma_2}, \,\,\,\hat{e}_0=1\,\,.$$
 We will call it the  {\it quantum torus} associated with $\Gamma$ and $C$. 
 
For any strict sector $V\subset \R^2$ we define
$${\cal R}_{V,C}:=\prod_{\gamma \in \Gamma\cap \,C_0(V,Z,Q)}C\cdot \hat{e}_{\gamma}$$
and call it the quantum torus associated with $V$. Here we introduce a notation which will be used later:
$$C_0(V,Z,Q):=C(V,Z,Q)\cup \{0\}$$ 
where the cone $C(V,Z,Q)$ was defined in 2.2.
Algebra  ${\cal R}_{V,C}$ is the natural completion of the subalgebra  ${\cal R}_{V,C}\cap {\cal R}_{\Gamma,C}\subset{\cal R}_{\Gamma,C}$.

Let us choose as $C$ the ring $\overline{D^{\mu}}$. We denote ${\cal R}_{\Gamma}:={\cal R}_{\Gamma, \overline{D^{\mu}}}$
the corresponding quantum torus and call it the {\it motivic quantum torus} associated with $\Gamma$. Similarly, we have motivic quantum tori $R_V$ associated with strict sectors $V$.

{\it Step 5}.

We define an element $A_V^{mot}\in {\cal R}_V:={\cal R}_{V,\overline{D^{\mu}}} $ in the following way. First, we fix a branch of the function $\op{Log}z$, where $ z\in V$ (the result will not depend on the choice of the branch). Recall the category $\CC_{V,\op{Log}}\subset \CC$ (see Section 3.4).
It follows from our assumptions that for any $\gamma\in \Gamma$ the set $\CC_{V,\gamma}=\{E\in Ob(\CC_{V,\op{Log}})|\op{cl}(E)=\gamma\}$ is constructible.

Finally, we define the desired element
$$A_V^{mot}=\sum_{\gamma}\left(\int_{\CC_{V,\gamma}}w \right)\cdot \hat{e}_{\gamma}\,\,.$$
The element $A_V^{mot}$ in fact depends only on $\overline{w}$.

Informally, one can write
$$A_V^{mot}=\sum_{E\in Iso(\CC_{V,\op{Log}})}{w(E)\over {[\op{Aut}(E)]}}\hat{e}_{\op{cl}(E)}=1+\dots\in {\cal R}_V,$$
where $Iso(\CC_{V,\op{Log}})$ denotes the set of isomorphism classes of objects of the category $\CC_{V,\op{Log}}$.

\begin{thm} Assuming  the integral identity, the collection of elements $(A_V^{mot})$ satisfies the Factorization Property:
if a strict sector $V$ is decomposed into a disjoint union $V=V_1\sqcup V_2$ (in the clockwise order) then
$$A_V^{mot}=A_{V_1}^{mot}A_{V_2}^{mot}\,\,. $$
Moreover we have a local homeomorphism $Stab(\CC)\to Stab\left({\cal R}_{\Gamma,\overline{D^{\mu}}}\right)$.

\end{thm}

This theorem  follows immediately from the statement of Proposition 11 (see 6.1) about the elements $A_V^{\op{Hall}}$, and 
  the Theorem 8 from the next section.

\subsection{From motivic Hall algebra to motivic quantum torus}

Assume  that $\CC$ is an ind-constructible  $3d$ Calabi-Yau category endowed with polarization and orientation data $\sqrt{\cal D}$. The Hall algebra of $\CC$ is graded by the corresponding lattice $\Gamma$: $H(\CC)=\oplus_{\gamma\in \Gamma}H(\CC)_{\gamma}$. Main result of this section is the following theorem.

\begin{thm} The map $\Phi: H(\CC)\to {\cal R}_{\Gamma}$ given by the formula
$$\Phi(\nu)=(\nu,w)\hat{e}_{\gamma}, \,\,\,\nu\in H(\CC)_{\gamma}$$
is a homomorphism of $\Gamma$-graded $\Q$-algebras. Here $w$ is the motivic weight and $(\bullet,\bullet)$ is the pairing between motivic measures and functions.

\end{thm}

In other words, the homomorphism $H(\CC)\to {\cal R}_{\Gamma}$ can be written as
$$ [\pi: Y\to Ob(\CC)]\mapsto $$
$$\mapsto \int_{Y}(1-MF(\pi(y)))\,(1-MF_0(Q_{\pi(y)}))\,{\mathbb L}^{-{1\over{2}}\op{rk}Q_{\pi(y)}}\,{\mathbb L}^{{1\over{2}}(\pi(y),\pi(y))_{\le 1}}\,\hat{e}_{\op{cl}(\pi(y))}\,\,,$$
where $\int_Y$ is understood as the direct image functor (see Section 4.2).

The natural extension of the above homomorphism to the completion of $\widehat{H}(\CC_V)$  maps the element $A_V^{\op{Hall}}$ to the
element $A_V^{mot}$ defined in Section 6.2.

{\it Proof.} For simplicity we will present the proof of the Theorem for 
$$\nu_E:=[\delta_E:pt\to Ob(\CC)]\,\,,$$
 where $\delta_E(pt)=E\in Ob(\CC(\k))$ is the ``delta-function". The general proof for equivariant constructible families is similar.
We will also assume that our category is minimal on the diagonal.

The proof will consists of several steps.

{\it Step 1.}

We have:

$$\nu_{E_1}\cdot\nu_{E_2}={\mathbb L}^{-(E_2,E_1)_{\le 0}}[\pi_{21}:\op{Ext}^1(E_2,E_1)\to Ob(\CC)]\,\,,$$
the map $\pi_{21}$ is the restriction of the cone map\footnote{Recall that we pretend that such a map exists. In fact, it is defined only as a 1-morphism of stacks.} 
$$tot(\op{Ker}(m_1:{\cal HOM}^0\to{\cal HOM}^1))\to Ob(\CC)$$  to the fiber over the point $(E_2[-1],E_1)$.
Under this map the element $\alpha\in \op{Ext}^1(E_2,E_1)$ is mapped to the object 
$$E_{\alpha}=Cone(\alpha: E_2[-1]\to E_1)\,\,.$$

Let us denote by $\gamma_i$ the class $\op{cl}(E_i)\in \Gamma, i=1,2$. Then we have:
$$\Phi(\nu_{E_i})={\mathbb L}^{{1\over{2}}(E_i,E_i)_{\le 1}}(1-MF_0(W_{E_i}^{min}))\,\overline{I}(h(E_i))\,\hat{e}_{\gamma_i},\,\,\, i=1,2\,\,,$$
where $h(E_i)$ is the value at $E_i$ of the element $h\in J_2(\CC)$ (i.e. the image of the restriction map to $J_2(\CC)$) given by the super line bundle $\sqrt{{\cal D}}\otimes {\cal D}_{\le 1}^{-1}$ with trivialized tensor square.

We have:

$$\Phi(\nu_{E_1})\Phi(\nu_{E_2})={\mathbb L}^{{1\over{2}}((E_1,E_1)_{\le 1}+(E_2,E_2)_{\le 1})}\times$$
$$\times\,{\mathbb L}^{{1\over{2}}((E_1,E_2)_{\le 1}-(E_2,E_1)_{\le 1})}(1-MF_0(W_{E_1}^{min}\oplus W_{E_2}^{min}))\overline{I}(h(E_1))\overline{I}(h(E_2))\,\hat{e}_{\gamma_1+\gamma_2}\,\,.$$
In order to obtain this formula  we used the Calabi-Yau property, which implies that
$$\begin{array}{rcl}
\langle \gamma_1,\gamma_2\rangle&=&\sum_{j\in \Z}(-1)^j(E_1,E_2)_j=(E_1,E_2)_{\le 1}+(E_1,E_2)_{\geqslant 2}=\\
&=&(E_1,E_2)_{\le 1}-(E_2,E_1)_{\le 1}\,\,,
\end{array}$$
where we employ the notation $(E_1,E_2)_{\geqslant m}=\sum_{j\geqslant m}(-1)^j(E_1,E_2)_j$.
Also we used the motivic Thom-Sebastiani theorem for the Milnor fibers and the product formula for the basis elements in the motivic quantum torus ${\cal R}_{\Gamma}$.

On the other hand, we can apply $\Phi$ to the product $\nu_{E_1}\cdot \nu_{E_2}$ and obtain:
$$\begin{array}{l}
\Phi(\nu_{E_1}\cdot \nu_{E_2})=\\
={\mathbb L}^{-(E_2,E_1)_{\le 0}}\int_{\alpha\in
\op{Ext}^1(E_2,E_1)}{\mathbb L}^{{1\over{2}}(E_{\alpha},E_{\alpha})_{\le
1}}(1-MF_0(W_{E_{\alpha}}^{min}))\overline{I}(h(E_{\alpha}))\hat{e}_{\gamma_1+\gamma_2}\,\,.
\end{array}$$

Using the identity $(E_2,E_1)_{\le 0}=(E_2,E_1)_{\le 1}+(E_2,E_1)_1$ (and also recall that $(E_2,E_1)_1=\dim \op{Ext}^1(E_2,E_1)$) and observing that
$$(E_1\oplus E_2,E_1\oplus E_2)_{\le 1}=(E_1,E_1)_{\le 1}+(E_2,E_2)_{\le 1}+(E_1,E_2)_{\le 1}+(E_2,E_1)_{\le 1}$$
we arrive to the following equality which is equivalent to 
$$\Phi(\nu_{E_1}\cdot \nu_{E_2})=\Phi(\nu_{E_1})\Phi(\nu_{E_2})$$ and hence implies the Theorem:

$${\mathbb L}^{(E_2,E_1)_1}(1-MF_0(W_{E_1}^{min}\oplus W_{E_2}^{min}))\overline{I}(h(E_1))\overline{I}(h(E_2))=$$
$$=\int_{\alpha\in \op{Ext}^1(E_2,E_1)}{\mathbb L}^{{1\over{2}}((E_{\alpha},E_{\alpha})_{\le
1}-(E_1\oplus E_2,E_1\oplus E_2)_{\le 1})}(1-MF_0(W_{E_{\alpha}}^{min}))\overline{I}(h(E_{\alpha}))\,\,.$$

{\it Step 2.}

Now we would like to express the difference

$$(E_{\alpha},E_{\alpha})_{\le 1}-(E_1\oplus E_2,E_1\oplus E_2)_{\le 1}$$ as the rank of a certain linear operator.
Recall that the object $E_{\alpha}$ can be thought of as a deformation of the object $E_0:=E_1\oplus E_2$.
Therefore, there is a spectral sequence which starts at $\op{Ext}^{\bullet}(E_1\oplus E_2,E_1\oplus E_2)$ and converges to
$\op{Ext}^{\bullet}(E_{\alpha},E_{\alpha})$. Using the $\A$-structure one can make it very explicit. Namely, let us denote by $d_{\alpha}:\op{Ext}^{\bullet}(E_1\oplus E_2,E_1\oplus E_2)\to \op{Ext}^{\bullet}(E_1\oplus E_2,E_1\oplus E_2) $ the differential of degree $+1$ given by the formula
$$d_{\alpha}=m_2(\alpha,\bullet)+m_2(\bullet,\alpha)+m_3(\alpha,\bullet,\alpha)\,\,.$$
Then the graded vector space $\op{Ext}^{\bullet}(E_{\alpha},E_{\alpha})$ is isomorphic to the cohomology of $d_{\alpha}$
(cf. e.g. \cite{KoSo3}, Remark 10.1.5).

It is clear that for any cohomological complex $(C^{\bullet},d)$ of finite-dimensional vector spaces we have the equality
$$\sum_{i\le 1}(-1)^i\dim H^i(C)-\sum_{i\le 1}(-1)^i\dim C^i
=\op{rk}d^{(1)},$$
where $d^{(1)}: C^1\to C^2$ is the component of $d$. Applying this observation to our complex we obtain that
$$(E_{\alpha},E_{\alpha})_{\le 1}-(E_1\oplus E_2,E_1\oplus E_2)_{\le 1}=\op{rk}d^{(1)}_{\alpha}.$$

{\it Step 3.}

Let us introduce a $\kk$-vector space
$${\cal M}_{E_1,E_2}=\op{Ext}^1(E_1\oplus E_2, E_1\oplus E_2)=$$
$$=\op{Ext}^1(E_1,E_1)\oplus \op{Ext}^1(E_2,E_1)\oplus \op{Ext}^1(E_1,E_2)\oplus \op{Ext}^1(E_2,E_2)\,\,.$$ It can be interpreted as the tangent space to the moduli space of formal deformations of the object $E_1\oplus E_2$.
We choose coordinates $(x,\alpha,\beta,y)$ on this space
in such a way that $x$ denotes the coordinates on $\op{Ext}^1(E_1,E_1)$, $\alpha$ denotes the coordinates on $\op{Ext}^1(E_2,E_1)$, $\beta$ denotes the coordinates on $\op{Ext}^1(E_1,E_2)$ and $y$ denotes the coordinates on $\op{Ext}^1(E_2,E_2)$. Then the point $(0,\alpha,0,0)$ corresponds (by abuse of notation) to the isomorphism class  $\alpha\in \op{Ext}^1(E_2,E_1)$
of an exact triangle $E_1\to E_{\alpha}\to E_2$. Later we are going to use the integral identity from Section 4.4 applying it to the formal neighborhood of the subspace consisting of the points $(0,\alpha,0,0)$. In order to do that we will relate the potential of the object $E_{\alpha}$ with a certain formal function on ${\cal M}_{E_1,E_2}$.

We may assume that the full subcategory $\CC(E_1,E_2)$ consisting of the pair of objects $E_1,E_2$ is minimal.
As in the case of one object the  potential of $\CC$ induces a formal power series $W_{E_1,E_2}=W(x,\alpha,\beta, y)$
on ${\cal M}_{E_1,E_2}$. It is defined as the abelianization of a series $\sum_{n\geqslant 3} W_n/n$ in cyclic paths in the quiver $Q_{E_1,E_2}$
with the vertices $E_1$ and $E_2$ and $(E_i,E_j)_1$ edges between vertices $E_i$ and $E_j$ for $i,j\in \{1,2\}$.
Since any cyclic path has the same number of edges in the direction $E_1\to E_2$ as in the direction $E_2\to E_1$ we conclude that the potential $W_{E_1,E_2}$ is ${\bf G}_m$-invariants with respect to the ${\bf G}_m$-action on the graded vector space ${\cal M}_{E_1,E_2}$ with the weights $wt\, x=wt\, y=0$ and $wt\, \alpha=-wt\, \beta=1$. The potential $W_{E_1,E_2}$ is obtained from the potential $W_{E_1\oplus E_2}$ by a formal change of variables.

It follows from ${\bf G}_m$-invariance of $W_{E_1,E_2}$ that it belongs to $\kk[\alpha][[x,\beta,\gamma]]$. Therefore it defines a function on the formal neighborhood of the affine subspace
$\{(0,\alpha,0,0)\}\subset {\cal M}_{E_1,E_2}$. In particular, for any $\alpha\in \op{Ext}^1(E_2,E_1)$ we obtain a formal power series $W_{E_1,E_2,\alpha}$ on ${\cal M}_{E_1,E_2}$ which is the Taylor expansion of $W_{E_1,E_2}$ at the point
$(0,\alpha,0,0)$. Similarly to the Proposition 7 from Section 3.3 the series $W_{E_1,E_2,\alpha}$ becomes (after a formal change of coordinates) a direct sum $W_{E_{\alpha}}^{min}\oplus Q_{E_{\alpha}}\oplus N_{E_{\alpha}}$, where $Q_{E_{\alpha}}$ is a non-degenerate quadratic form, $N_{E_{\alpha}}$ is the zero function on a vector subspace, and $W_{E_{\alpha}}^{min}$ does not contain terms of degree less than $3$ in its Taylor expansion. By the motivic Thom-Sebastiani theorem we have 
$$(1-MF_0(W_{E_1,E_2,\alpha}))=(1-MF(E_{\alpha}))(1-MF_0(Q_{E_{\alpha}}))\,\,.$$ 
Let us consider the quadratic form
$((W_{E_1,E_2})^{\prime\prime})_{|(0,\alpha,0,0)}$
on ${\cal M}_{E_1,E_2}$, where $(W_{E_1,E_2})^{{\prime\prime}}$ denotes the second derivative of the  potential 
with respect to the  affine coordinates. It follows from the above discussion that this quadratic form is equal
to the direct sum of  $Q_{E_{\alpha}}$ and the zero quadratic form.

It is easy to check that $
\left(((W_{E_1,E_2})^{\prime\prime})_{|(0,\alpha,0,0)}\right)(v)$ is equal to
$(d_{\alpha}^{(1)}v,v)$ for any $ v\in {\cal M}_{E_1,E_2}$. Hence $Q_{E_{\alpha}}$ can be identified with the quadratic form on $\op{I
m}(d_{\alpha}^{(1)})$ given by $(u, (d_{\alpha}^{(1)})^{-1}u)$.

{\it Step 4.}

Recall (see Section 5.2) that for any exact triangle $E_1\to E_{\alpha}\to E_2$ we have a super line bundle $L$ with a canonically trivialized square:

$$L_{E_1\to E_{\alpha}\to E_2}=({\cal D}_{{\le 1},E_{\alpha}}\otimes {\cal D}_{{\le 1},E_1}^{-1}\otimes {\cal D}_{{\le 1},E_2}^{-1})\otimes (\op{sdet}(\op{Ext}^{\bullet}(E_1,E_2)))^{-1}.$$
For a split triangle $E_{\alpha}\simeq E_1\oplus E_2$ (i.e. $\alpha=0$) this line bundle is canonically trivialized since
\begin{itemize}
\item by definition, for $E_{\alpha}\simeq E_1\oplus E_2$ we have
$${\cal D}_{\le 1,E_{\alpha}}\simeq {\cal D}_{\le 1,E_{1}}\otimes {\cal D}_{\le 1,E_{2}}\otimes $$
$$\otimes\op{sdet}(\op{Ext}^{\le 1}(E_1,E_2))\otimes \op{sdet}(\op{Ext}^{\le 1}(E_2,E_1))\,\,,$$
\item by the Calabi-Yau property we have 
$$\op{sdet}(\op{Ext}^{\le 1}(E_2,E_1))\simeq \op{sdet}(\op{Ext}^{\geqslant 2}(E_1,E_2))\,\,.$$
\end{itemize}
Therefore, for any exact triangle $E_1\to E_{\alpha}\to E_2$ we have an isomorphism
$$L_{E_1\to E_{\alpha}\to E_2}\simeq {\cal D}_{\le 1,E_{\alpha}}\otimes {\cal D}_{\le 1,E_{1}\oplus E_2}^{-1}\,\,.$$
On the other hand, considerations similar to those in Step 2 give rise to a canonical isomorphism 
$${\cal D}_{\le 1,E_{\alpha}}\otimes {\cal D}_{\le 1,E_{1}\oplus E_2}^{-1}\simeq \op{sdet}(\op{Im}(d_{\alpha}^{(1)}))\,\,.$$
One can see that the trivialization $(\op{sdet}(\op{Im}(d_{\alpha}^{(1)})))^{\otimes 2}\simeq {\bf 1}$ comes exactly 
from the non-degenerate quadratic form $Q_{E_{\alpha}}$.
Therefore, for an arbitrary exact triangle $E_1\to E_{\alpha}\to E_2$ we have an isomorphism of super lines compatible
 with the trivializations of squares: 
$$L_{E_1\to E_{\alpha}\to E_2}\simeq \op{sdet}(\op{Im}(d_{\alpha}^{(1)}))\,\,.$$
This implies that
$$\overline{I}(Q_{E_{\alpha}})=\overline{I}(l(E_1\to E_{\alpha}\to E_2))\,\,,$$
where $l\in J_2(Funct({\CC}_3, \CC))$ was defined in Section 5.2.

{\it Step 5.}

Let us apply the integral identity from Section 4.4 to the potential $W_{E_1,E_2}$. 
We put
$$V_1:=\op{Ext}^1(E_2,E_1),\,V_2:=\op{Ext}^1(E_1,E_2),\,V_3:=\op{Ext}^1(E_1,E_1)\oplus\op{Ext}^1(E_2,E_2)\,.  $$
We  have:
$$\int_{\alpha\in \op{Ext}^1(E_2,E_1)}(1-MF_{(0,\alpha,0,0)}(W_{E_1,E_2}))=$$ $$=
{\mathbb L}^{(E_2,E_1)_1}(1-MF_0((W_{E_1,E_2})_{|\op{Ext}^1(E_1,E_1)\oplus \op{Ext}^1(E_2,E_2)}))\,\,.$$

On the other hand the LHS of the integral identity is equal to
$$\int_{\alpha\in \op{Ext}^1(E_2,E_1)}(1-MF_0(Q_{E_{\alpha}}))(1-MF(W^{min}_{E_{\alpha}})=$$
$$=\int_{\alpha\in \op{Ext}^1(E_2,E_1)}{\mathbb L}^{{1\over{2}}\op{rk}Q_{E_{\alpha}}}\overline{I}(Q_{E_{\alpha}})(1-MF(W^{min}_{E_{\alpha}}))\,\,.$$

Recall that $\op{rk} Q_{E_{\alpha}}=\op{rk}(d_{\alpha}^{(1)})=(E_{\alpha},E_{\alpha})_{\le 1}-(E_1\oplus E_2,E_1\oplus E_2)_{\le 1}$ by Steps 2 and 3. Then the integral identity becomes the following equality:
$$\int_{\alpha\in \op{Ext}^1(E_2,E_1)}{\mathbb L}^{{1\over{2}}((E_{\alpha},E_{\alpha})_{\le 1}-(E_1\oplus E_2,E_1\oplus E_2)_{\le 1}) }\overline{I}(Q_{E_{\alpha}})(1-MF(W^{min}_{E_{\alpha}}))=$$
$$={\mathbb L}^{(E_2,E_1)_1}(1-MF_0(W^{min}_{E_1}\oplus W^{min}_{E_2}))\,\,.$$

Comparing this formula with the one we wanted to prove on Step 1 we see that they coincide if
$$\overline{I}(Q_{E_{\alpha}})={{\overline{I}(h(E_{\alpha}))}\over{{\overline{I}(h(E_{1}))}{\overline{I}(h(E_{2}))}}}\,\,.$$
Now using Step 4 we observe that this cocycle condition is equivalent to the main property of the orientation data on exact triangles.
This concludes the proof of the Theorem. $\blacksquare$

\begin{defn} Let $\CC$ be an ind-constructible $3$-dimensional Calabi-Yau category endowed with polarization, $\sigma\in Stab(\CC,\op{cl})$. We call the collection of elements $(A_V^{mot}\in {\cal R}_V)$ of the completed motivic quantum tori $({\cal R}_V)$ (for all strict sectors $V\subset \R^2$) the  motivic Donaldson-Thomas invariant of $\CC$.

\end{defn}

Let us consider the following unital $\Q$-subalgebra of $\Q(q^{1/2})$:
$$D_q:=\Z[q^{1/2},q^{-1/2},\left((q^n-1)^{-1}\right)_{n\geqslant 1}]\,.$$
 There is a homomorphism of rings $\overline{D^{\mu}}\to D_q$ given by the twisted Serre polynomial. Namely, it maps ${\mathbb L}^{1/2}\mapsto q^{1/2}$, and on
${\cal M}^{\mu}$ it is the composition of the Serre polynomial with the involution $q^{1/2}\mapsto -q^{1/2}$. 
Therefore, we have a homomorphism of algebras ${\cal R}_{\Gamma}\to {\cal R}_{\Gamma,q}$, where
${\cal R}_{\Gamma,q}$ is a $D_q$-algebra
generated by $\hat{e}_{\gamma}, \gamma\in \Gamma$, subject to the relations
$$\hat{e}_{\gamma}\hat{e}_{\mu}=q^{{1\over{2}}{\langle\gamma,\mu\rangle}}\hat{e}_{\gamma+\mu}\,,\,\,\,\hat{e}_0=1\,.$$
Similarly to the motivic case, we have the algebra ${\cal R}_{V,q}$ associated with any strict sector $V$.

The motivic DT-invariants give rise to stability data on the graded Lie algebra associated with
${\cal R}_{\Gamma,q}$. We will denote by $A_{V,q}\in {\cal R}_{V,q}$ the element corresponding to
$A_V^{mot}$.

\subsection{Examples}

1) Assume that a $3$-dimensional Calabi-Yau category $\CC$ is generated by one spherical object $E$ defined over $\k$. Therefore $R:=\op{Ext}^{\bullet}(E,E)\simeq H^{\bullet}(S^3,\k)$. In this case we take $\Gamma=K_0(\CC(\kk))\simeq \Z\cdot \op{cl}_{\kk}(E)$, and the skew-symmetric form on $\Gamma$ is trivial. In order to choose an orientation data, let us fix a basis $r_0=1,r_3$ in the algebra $R$ (the subscript indicates the degree). Consider $R$ as a bimodule over itself and denote this bimodule by $M$. The corresponding to $1,r_3$ bimodule basis will be denoted by $1_M,r_{3,M}$. Then we have a family $M_t, 0\le t\le 1$ of $R$-bimodule structures on $M$ such that 
$$1_{M}\cdot_t r_3=r_3\cdot_t 1_{M}=(1-t)r_{3,M}\,.$$
 Hence $M_0=M$ and $M_1\simeq N\oplus N^{\vee}$ in notation of Section 5.3. The latter gives a decomposition of the  bifunctor $\op{Hom}^\bullet$. The above family of bimodules define a homotopy which can be used for definition of an orientation data as in Section 5.2.

For any $z\in \C, \op{Im}z>0$ we have a stability condition $\sigma_z$ such that $E\in \CC^{ss}, \,Z(E):=Z(\op{cl}_\kk(E))=z,\, \op{Arg}(E)=\op{Arg}(z)\in (0,\pi)$.
For a strict sector $V$ such that $\op{Arg}(V)\subset (0,\pi)$ we have the category $\CC_V$ which is either trivial (if $z\notin V$) or consists of objects $0,E,E\oplus E,\dots$ (if $z\in V$).
Then $A_V^{mot}=1$ in the first case and
$$A_V^{mot}=\sum_{n\geqslant 0}\frac{{\mathbb L}^{n^2/2}}{[GL(n)]}\hat{e}_{\gamma_1}^n\,\,,$$
in the second case. Here  $\gamma_1:={\op{cl}_\kk(E)}$ is the generator of $\Gamma$.

Let us comment on the answer.
In this case $\op{Ext}^1(nE,nE)=0$, where we set $nE=E^{\oplus n}, n\geqslant 1$. Therefore $W_{nE}=0$ which implies that $MF(W_{nE})=0$.
The numerator is
$${\mathbb L}^{n^2/2}={\mathbb L}^{{1\over{2}}\,\dim \op{Ext}^0(nE,nE)}=
{\mathbb L}^{{1\over{2}}\sum_{i\le 1}(-1)^i \dim \op{Ext}^i(nE,nE)}\,\,,$$
since $\op{Ext}^{\ne 0}(nE,nE)=0$.

Let us consider the ``quantum dilogarithm" series
$${\bf E}(q^{1/2},x)=\sum_{n\geqslant 0}{{q}^{n^2/2}\over{(q^n-1)\dots(q^n-q^{n-1})}}x^n\in \Q(q^{1/2})[[x]]\,\,.$$

Since $[GL(n)]=({\mathbb L}^n-1)\dots({\mathbb L}^n-{\mathbb L}^{n-1})$, we conclude that
$$A_V^{mot}={\bf E}({\mathbb L}^{1/2},\hat{e}_{\gamma_1})\,\,.$$
In order to simplify the notation we will denote ${\bf E}(q^{1/2},x)$ simply by ${\bf E}(x)$.
 In Section 7.1 we will discuss the quasi-classical limit, and will associate numerical Donaldson-Thomas invariants $\Omega(\gamma)\in \Q$
 for any $\gamma\in \Gamma$ for given stability structure $\sigma\in Stab(\CC,\op{cl})$. In our basic example we have (for any $\sigma$)
$$\Omega(\pm \gamma_1)=1,\,\,\,\Omega(n \gamma_1)=0\mbox { for }n\ne \pm 1\,.$$

2) Assume that $\CC$ is generated by two spherical objects $E_1,E_2$ defined over $\k$ such that $\dim  \op{Ext}^i(E_2,E_1)=0$
if $i\ne 1$ and $\dim \op{Ext}^1(E_2,E_1)=1$. Notice that the unique (up to isomorphism) non-trivial extension $E_{12}$ appears in the exact triangle $E_1\to E_{12}\to E_2$ and it is a spherical object.

For any $z_1,z_2\in \C, \op{Im}z_i>0, i=1,2$ there is a unique stability condition $\sigma_{z_1,z_2}$ such that $Z(E_i):=Z(\op{cl}_{\kk}(E_i))=z_i, i=1,2$, and the category $\CC_V(\kk)$ in the case $z_1,z_2\in V, \op{Arg}(V)\subset (0,\pi)$ consists of subsequent extensions of the copies of $E_1$ and $E_2$.

If $\op{Arg}(z_1)>\op{Arg}(z_2)$ then the only $\sigma_{z_1,z_2}$-semistable objects are (up to shifts) $E_1, 2E_1,\dots,E_2,2E_2,\dots$, where we use the notation $nE$ for $E^{\oplus n}$, as before.
If $\op{Arg}(z_2)>\op{Arg}(z_1)$ then we have three groups of $\sigma_{z_1,z_2}$-semistable objects:
$nE_1,nE_2,nE_{12}, n\geqslant 1$.

\vspace{3mm}
\begin{picture}(350,150)(0,0)
\setlength{\unitlength}{.8pt}
\thicklines
\put(0,0){\line(1,0){200}}
\put(240,0){\line(1,0){200}}
\put(100,0){\circle{3}}
\multiput(70,30)(-30,30){3}{\circle*{4}}
\multiput(310,30)(-30,30){3}{\circle*{4}}
\multiput(130,30)(30,30){3}{\circle*{4}}
\multiput(370,30)(30,30){3}{\circle*{4}}
\multiput(340,60)(0,60){2}{\circle*{4}}
\thinlines
\put(100,0){\line(-1,1){100}}
\put(100,0){\line(1,1){100}}
\put(340,0){\line(-1,1){100}}
\put(340,0){\line(1,1){100}}
\put(340,0){\line(0,1){140}}
\put(53,20){$E_1$}
\put(15,50){$2E_1$}
\put(-15,80){$3E_1$}
\put(293,20){$E_2$}
\put(255,50){$2E_2$}
\put(225,80){$3E_2$}
\put(133,20){$E_2$}
\put(163,50){$2E_2$}
\put(193,80){$3E_2$}
\put(373,20){$E_1$}
\put(403,50){$2E_1$}
\put(433,80){$3E_1$}
\put(345,57){$E_{12}$}
\put(345,117){$2E_{12}$}

\put(340,0){\circle{3}}
\end{picture}
\vspace{3mm}

The wall-crossing formula implies the following well-known 
identity (see \cite{FK}) in the algebra $D_q\langle\langle x_1,x_2\rangle\rangle/(x_1x_2-q x_2x_1)$:
$${\bf E}(x_1){\bf E}(x_2)={\bf E}(x_2){\bf E}(x_{12}){\bf E}(x_1)\,\,,$$
where $x_{12}=q^{-1/2}x_1x_2=q^{1/2}x_2x_1$ and $x_i$ corresponds to $\hat{e}_{\op{cl}_\kk(E_i)},\, i=1,2,12$.
  Namely, both sides of the above identity are equal  to $A_{V_{big},q}$ for any sector $V_{big}$ in the upper half-plane containing $z_1,z_2$.
The LHS and the RHS of the identity come from the decompositions 
$$A_{V_{big},q}=A_{V_1,q} A_{V_2,q},\,\,\,A_{V_{big},q}=A_{V_2,q}A_{V_{12},q}A_{V_1,q}\,,$$
where $V_i,\,i=1,2,12$ are some narrow sectors containing $z_i$.

\begin{rmk} The function ${\bf E}(x)$ satisfies also the identity
$${\bf E}(x_2){\bf E}(x_1)={\bf E}(x_1+x_2)$$
for $x_1,x_2$ obeying the relations $x_1x_2=q x_2x_1$ as above. This follows from the formula
$${\bf E}(x)=\exp_q \left( \frac{q^{1/2}}{q-1} x \right)$$
where $\exp_q(x)$ is the usual $q$-exponent
$$\exp_q(x):=\sum_{n\geqslant 0}\frac{x^n}{[n]_q !}\,,\,\,\,[n]_q!:=\prod_{j=1}^n[j]_q,\,\,\,[j]_q:=1+q+\dots+q^{j-1}\,\,.$$
The exponential property of ${\bf E}(x)$ seems to play no role in our considerations.
\end{rmk}

If we denote by $\gamma_i \in \Gamma\simeq \Z^2$ the classes $\op{cl(E_i)}, \, i=1,2$, then the only 
non-trivial numerical Donaldson-Thomas invariants are
$$\Omega(\pm \gamma_1)=\Omega(\pm \gamma_2)=1$$
in the case $\op{Arg}(z_1)>\op{Arg}(z_2)$, and 
  $$\Omega(\pm \gamma_1)=\Omega(\pm \gamma_2)=\Omega(\pm (\gamma_1+\gamma_2))=1$$
in the case $\op{Arg}(z_1)<\op{Arg}(z_2)$.

\subsection{D0-D6 BPS bound states: an example related to the MacMahon function}

Let $X$ be a compact $3d$ Calabi-Yau manifold over $\k$, such that $H^1(X,{\cal O}_X)=0$.
We denote by $\CC^{(0,6)}$ the ind-constructible triangulated category generated by the structure sheaf ${\cal O}_X$ and torsion sheaves ${\cal O}_x, x\in X$.\footnote{This category is related to the counting of D0-D6 BPS bound states, compare with \cite{DM}, formula (6.1).}
 This category has a $t$-structure with the heart consisting of coherent sheaves on $X$ which are trivial vector bundles outside of a finite set.  Then
${\cal O}_X$ is the only spherical object in $\CC^{(0,6)}$. We choose $\Gamma=\Z\gamma_1\oplus \Z\gamma_2$, which is the image of $K_0(\CC^{(0,6)})$ under the Chern class in the quotient of the Chow group by the numerical equivalence, where $
 \gamma_1=\op{cl}_{\kk}({\cal O}_x))$  for any point $x\in X$, and    $ \gamma_2=\op{cl}_{\kk}({\cal O}_X)$.
We are going to consider a stability condition $\sigma=(Z,(\CC^{(0,6)})^{ss},\dots)$ on $\CC^{(0,6)}$ with the above $t$-structure and such that 
$$z_1:=Z(\gamma_1)=-1,\,\,\,z_2=Z(\gamma_2)=i=\sqrt{-1}\,\,.$$ Then $\sigma$-semistable objects in $\CC^{(0,6)}$ will be either pure torsion sheaves supported at finitely many points or torsion-free sheaves.

This corresponds to the following picture for $\Omega(\gamma)$.

\begin{picture}(250,200)(-40,20)
\thinlines
\put(0,20){\vector(1,0){250}}
\put(70,0){\vector(0,1){200}}
\multiput(20,20)(50,0){5}{\circle*{3}}
\multiput(20,70)(50,0){5}{\circle*{3}}
\multiput(20,120)(50,0){5}{\circle*{3}}
\multiput(20,170)(50,0){5}{\circle*{3}}
\put(10,28){$-\chi$}
\put(110,28){$-\chi$}
\put(110,78){$-\chi$}
\put(165,77){$\frac{\chi^2+5\chi}{2}$}
\put(215,77){$-\frac{\chi^3+15\chi^2+20\chi}{6}$}
\put(165,25){?}
\put(165,125){?}
\put(215,25){?}
\put(215,125){?}
\put(215,175){?}
\put(63,75){$1$}
\put(13,75){0}
\multiput(13,125)(50,0){3}{0}
\multiput(13,175)(50,0){4}{0}
\end{picture}

\vspace{6mm}
Let us comment on the last figure.

a) The vertical line corresponds to the subcategory generated by the spherical object ${\cal O}_X$, for which we know $\Omega(\gamma)$. Namely, $\Omega(\gamma_2)=1$ and $\Omega(n\gamma_2)=0, n\geqslant 2$.

b) Horizontal line $z_2-nz_1, n\geqslant 0$ corresponds to sheaves of ideals of $0$-dimensional subschemes. Then:

$$ \sum_{n\geqslant 0}\Omega(\gamma_2-n\gamma_1)t^n=M(-t)^{\chi(X)},$$
where $\chi(X)$ is the Euler characteristic of $X$ and $$M(x):=\prod_{n\geqslant 1}(1-x^n)^{-n}\in \Z[[x]] 
$$ is the MacMahon function (see \cite{MaNOP}, \cite{BF} about this identity).

c) The torsion sheaves ${\cal O}_x, x\in X$ are Schur objects in $\CC^{(0,6)}$. Their moduli space is canonically identified with $X$. By Behrend's formula (see \cite{B}) their contribution to the virtual fundamental number of objects is 
$$\Omega(\gamma_1)=(-1)^{\dim X}\chi(X)=-\chi(X)\,.$$

d) The numbers marked by $``?"$ correspond to (possibly non-Schur) objects. Notice that there are no semistable objects with the class $n\gamma_2-m\gamma_1$ with $0<m<n$. They correspond to the sector filled by $0$'s.

Let us now choose a path $\sigma_{z_1(\tau),z_2(\tau)}$ in the space of the above stability structures such that 
$$ z_1(\tau)=-\exp(i\tau),\,z_2(\tau)=i\,,$$ where $\tau\in [0,\pi/2+\varepsilon],\, i=\sqrt{-1}$ and $\varepsilon>0$ is sufficiently small. The heart of the $t$-structure for $\tau>0$ consists of complexes of sheaves $E$ such that there exists an exact triangle $E_1\to E\to E_2[-1]$, where $E_1$ is a torsion-free sheaf and $E_2$ is a torsion sheaf (indeed, the new $t$-structure is obtained from the initial one by the standard tilting procedure). 
This heart coincides with the category $\CC^{(0,6)}_V$ for any $\tau \in (0,\pi/2+\varepsilon]$ where 
$$V=\{z\in \C^\ast|\,0\le \op{Arg}(z)\le \pi/2+\varepsilon\}\,.$$
Object ${\cal O}_X\in \CC^{(0,6)}_V$ can not be represented as a non-trivial extension in $\CC^{(0,6)}_V$, hence it is semistable for
 any $\tau \in [0,\pi/2+\varepsilon]$. 

Let us now consider the case $\tau \in (\pi/2,\pi/2+\varepsilon]$. Then  object ${\cal O}_X$ has the minimal argument among all non-trivial
objects in $\CC^{(0,6)}_V$. Therefore, all other indecomposable semistable objects $E$ are strictly on the left of ${\cal O}_X$, and we have
 $$\op{Ext}^0(E, {\cal O}_X)=0\,.$$
Taking the long exact sequence of $Ext$-groups to the object ${\cal O}_X$ one easily shows that 
in the decomposition $E_1\to E\to E_2[-1]$ we have
$E_1=0$.  Hence in this new heart $\CC^{(0,6)}_V$  the left orthogonal to  ${\cal O}_X$ 
consists of objects $F[-1]$, where $F$ is a torsion sheaf.
We conclude that for the stability condition with $\tau \in (\pi/2,\pi/2+\varepsilon]$ the only semistable objects have classes which belong to
$\Z_{\ne 0}\gamma_2\sqcup \Z_{\ne 0}\gamma_1$. Therefore, all DT-invariants $\Omega_\tau (\gamma)$ for $\sigma_{z_1({\tau}),z_2(\tau)}$ with $\tau 
>\pi/2$ are completely determined by the numbers $a_n=\Omega_{\tau}(-n\gamma_1), n\geqslant 1$
(and known invariants $\Omega_{\tau}(m\gamma_2)=\delta_{m,1}, m\geqslant 1$). Then the wall-crossing formula determines all the invariants 
$\Omega(\gamma)=\Omega_0(\gamma)$ for the initial stability condition $\sigma_{z_1(0),z_2(0)}$ in terms of the numbers $a_n, n\geqslant 1$.

The wall-crossing formula  implies that the following identity:
$$\prod_{n\geqslant 1}T_{-n\gamma_1}^{a_n}T_{\gamma_2}=\prod_{m\geqslant 1,n\geqslant 0}^{\longrightarrow}T_{-n\gamma_1+m\gamma_2}^{\Omega(-n\gamma_1+m\gamma_2)}\prod_{n\geqslant 1}T_{-n\gamma_1}^{a_n}\,\,.$$
 Using the known result for special values $\Omega(\gamma_2-n \gamma_1),\,\,n\geqslant 1$ (in terms of the MacMahon function), one can deduce that 
all the numbers
$a_n=\Omega(n\gamma_1))$ for $n\geqslant 1$ are equal to $-\chi(X)$. We don't know how to prove this identity directly.
We see that invariants $\Omega_\tau(\gamma)$ for $\tau >\pi/2$ have a much simpler form than $\Omega(\gamma)=\Omega_0(\gamma)$.
 Moreover, it is now possible (in principle) to  work out a formula for $\Omega(-n\gamma_1+m\gamma_2)$ for any given $m\geqslant 2$.

\begin{rmk} One can try to generalize the above considerations to the case of D0-D2-D6 bound states. Mathematically this means that we consider the triangulated category generated by the sheaf ${\cal O}_X$ and sheaves with at most $1$-dimensional support (cf. \cite{MaNOP}). A problem arises here, since for the natural $t$-structure there is no central charge which gives a stability condition on the category. Presumably, in this case one can use the limit stability conditions (see \cite{Ba}, \cite{Tod1}).

\end{rmk}

\begin{rmk} Let $X$ be a $3d$ complex Calabi-Yau manifold, $C\simeq {\bf P}^1\subset X$ a rational curve with normal bundle isomorphic to ${\cal O}(-1)\oplus {\cal O}(-1)$ and $\CC$ be an ind-constructible $\A$-version of the category $Perf_C(X)$ of perfect complexes supported on $C$. Then $\Gamma:=K_0(\CC)\simeq \Z^2$ carries a trivial skew-symmetric (Euler) form. The lattice $\Gamma$ is generated by $\op{cl}({\cal O}_{pt})$ and $\op{cl}({\cal O}_{C})$.
It follows that there are no wall-crossings in this case, and hence our invariants $\Omega(\gamma)$ do not change under continuous deformations of a stability condition. In order to use this idea for computations one can choose two stability conditions by specifying the corresponding $t$-structures and central charges:

a) choose the $t$-structure with the heart consisting of coherent sheaves on $X$ supported on $C$ and the central charge $Z$ such that $$Z(\op{cl}({\cal O}_{pt}))\in \R_{<0},\,\, \op{Im} Z(\op{cl}({\cal O}_{C}))>0\,\,;$$

b) choose the $t$-structure given by the category of finite-dimensional representations of the quiver with two vertices and two double arrows in each direction and the potential $$W=a_1b_1a_2b_2-a_1b_2a_2b_1\,\,.$$

Then calculations from \cite{Sz} give the following formulas for the invariants $\Omega(\gamma)$:

$$\begin{array}{rcc}\Omega(n\op{cl}({\cal O}_{pt}))&=&-2, n\ne 0\,\,;\\
\Omega(n\op{cl}({\cal O}_{pt})\pm \op{cl}({\cal O}_{C}))&=&1, n\in \Z\,\,.
\end{array}$$
In all other cases $\Omega(\gamma)=0$.

For recent generalization of \cite{Sz} see \cite{MoRe2},\cite{NagNak},\cite{Nag}.

\end{rmk}

\section{Quasi-classical limit and integrality conjecture}

\subsection{Quasi-classical limit, numerical DT-invariants}

The elements $A_{V,q}\in {\cal R}_{V,q}$ corresponding to $A_V^{mot}$ are series in $\hat{e}_{\gamma},\gamma\in \Gamma$ with coefficients which are rational functions in $q^{1/2}$. They can have poles as $q^n=1$ for some $n\geqslant 1$. Hence it is not clear how to take the quasi-classical limit as $q^{1/2}\to -1$ (this corresponds to the taking of Euler characteristic of the corresponding motives).

Let us assume that the skew-symmetric form on $\Gamma$ is non-degenerate (otherwise we can replace $\Gamma$ by the symplectic lattice $\Gamma\oplus \Gamma^{\vee}$). The element $A_{V,q}$ defines an automorphism of an appropriate completion of
${\cal R}_{\Gamma,q}$. More precisely, it acts by the conjugation $x\mapsto A_{V,q}xA_{V,q}^{-1}$ on the subring
$$\prod_{\gamma\in C_0(V)\cap \Gamma}D_q\hat{e}_{\gamma}\,$$ where $C_0(V)=C_0(V,Z,Q)$ is the union of $0$ with the  
convex hull $C(V,Z,Q)$ of the set $Z^{-1}(V)\cap \{Q\geqslant 0\}$ (see Section 2).

Let us recall the  example of the category generated by two spherical objects from Section 6.4. We will use notation for sectors $V_1,V_2,V_{big}$ introduced there.
One has, for 
quantum variables $x_1 x_2= q x_2 x_1$ and $A_{V_1,q}={\bf E}(x_1)$:

$$\begin{array}{l}
x_1\mapsto {\bf E}(x_1)x_1{\bf E}(x_1)^{-1}=x_1\,\,;\\
x_2\mapsto {\bf E}(x_1)x_2{\bf E}(x_1)^{-1}=x_2(1+q^{1/2}x_1)\,\,.
\end{array}$$

This follows from the formula $f(x_1)x_2=x_2f(qx_1)$, where $f(x)$ is an arbitrary series as well from the formula
$${\bf E}(x)=\prod_{n\geqslant 0}(1+q^{(2n+1)/2}x)^{-1}\,\,,$$
which is valid for $0<q<1$. The latter formula implies the needed identity in $\Q(q^{1/2})[[x]]$:

$${\bf E}(qx)=(1+q^{1/2}x){\bf E}(x)\,\,.$$

A similar formula holds for the conjugation by $A_{V_2,q}$. 
We remark that in this example the conjugation by $A_{V,q}$ for $V_1, V_2$ or $V_{big} $ preserves the  subring
$\prod_{\gamma\in C_0(V)\cap \Gamma}\Z[q^{\pm 1/2}]\hat{e}_{\gamma}$. In particular,
 one can make a specialization at 
 $$q^{1/2}=-1\,.$$
 \begin{rmk} Recall that at the end of Section 6.3 we defined a homomorphism
 $${\cal R}_\Gamma\to {\cal R}_{\Gamma,q}$$
 as the composition of Serre polynomial with the involution $q^{1/2}\mapsto -q^{-1/2}$.
 In particular, the specialization $q^{1/2}=-1$ is well-defined on the subring of series
 in generators $\hat{e}_\gamma$ with coefficients in
 ${\cal M}^{\mu}(Spec(\k))[{\mathbb L}^{-1/2}]$ (see also Section 7.3), and it corresponds
  to the usual Euler characteristic. We  use the twisting $q^{1/2}\mapsto -q^{-1/2}$ in order to
  avoid a lot of minus signs in formulas.
 \end{rmk}
  The ``integer" quantum torus
$$\bigoplus_{\gamma\in C_0(V)\cap \Gamma} \Z[q^{\pm 1/2}]\hat{e}_{\gamma}\subset {\cal R}_{\Gamma,q}$$ has the quasi-classical limit\footnote{There is another quasi-classical limit $q^{1/2}\to +1$ which we do not consider here.} which is the Poisson algebra with basis ${e}_{\gamma}, \gamma\in C_0(V)\cap \Gamma$ with the product and Poisson bracket
given by
$${e}_{\gamma}{e}_{\mu}=(- 1)^{\langle \gamma,\mu\rangle}e_{\gamma+\mu},\,\,\,\{e_{\gamma},e_{\mu}\}=(- 1)^{\langle \gamma,\mu\rangle}\langle \gamma,\mu\rangle e_{\gamma+\mu}\,\,.$$ 
The Poisson bracket is the limit of a normalized bracket:
$$[\hat{e}_\gamma,\hat{e}_\mu]=\left(q^{1/2\langle \gamma, \mu\rangle}-q^{-1/2\langle \gamma, \mu\rangle}\right)\hat{e}_{\gamma+\mu}\,\,,$$
$$\lim_{q^{1/2}\to -1}
(q-1)^{-1}\cdot\left(q^{1/2\langle \gamma, \mu\rangle}-q^{-1/2\langle \gamma, \mu\rangle}\right)=(- 1)^{\langle \gamma,\mu\rangle}\langle \gamma,\mu\rangle\,.$$
 One can write informally
 $$e_\gamma=\lim_{q^{1/2}\to -1}\frac{\hat{e}_\gamma}{q-1}\,\,.$$
 
\begin{conj} For any $3d$ Calabi-Yau category with polarization and any strict sector $V$ the automorphism $x\mapsto A_{V,q}xA_{V,q}^{-1}$ preserves the subring $$\prod_{\gamma\in C_0(V)\cap \Gamma}D_q^+\hat{e}_{\gamma}\,\,,$$ where
$D_q^+ :=\Z[q^{\pm 1/2}]$.

\end{conj}

Later we will present arguments in favor of this conjecture as well as a stronger version.
Assuming the Conjecture  we can define ``numerical" DT-invariants of a $3d$ Calabi-Yau category with polarization in the following way. Consider the quasi-classical limit (i.e. specialization at $q^{1/2}=-1$) of the automorphism $x\mapsto A_{V,q}xA_{V,q}^{-1}$.
We will present (see Section 7.4. and Conjecture 10) an explicit conjectural formula for this ``quasi-classical limit" which does not depend on the orientation data.
The quasi-classical limit gives rise to a formal symplectomorphism of the torus ${\mathbb T}_{\Gamma}$ and therefore induces 
 the stability data on the graded Lie algebra $\g_{\Gamma}$ (see Section 2.5). Alternatively, we can define such data as
$$a(\gamma):=\lim_{q^{1/2}\to -1}(q-1)a(\gamma)_q$$
in the obvious notation.
For a generic central charge $Z$ the symplectomorphism can be written as
$$A_V=\prod_{Z(\gamma)\in V}^{\longrightarrow}T_{\gamma}^{\,\Omega(\gamma)},$$
where $$T_{\gamma}(e_{\mu})=(1-e_{\gamma})^{\langle \gamma,\mu\rangle}e_{\mu}$$ and $\Omega(\gamma)\in \Q$ (see Section 2.5). In the above example of the Calabi-Yau category generated by one spherical object $E$ we have $\Omega(n\,\op{cl}(E))=1$ if $n\ne 1$ and $\Omega(n\,\op{cl}(E))=0$ otherwise.

\begin{conj} For a generic central charge $Z$ all numbers $\Omega(\gamma),\gamma \in \Gamma\setminus \{0\}$ are integers.

\end{conj}

The collection $(\Omega(\gamma))_{\gamma \in \Gamma}$
seems to be the correct mathematical definition of the counting of BPS states in  String Theory.

Finally, we make a comment about the relationship with the work of Kai Behrend (see \cite{B}). Recall that he defined a $\Z$-valued invariant of a critical point $x$ of a function $f$ on $X$ which is equal to 
$$(-1)^{\dim X}(1-\chi(MF_x(f)))\,\,,$$ where $\chi$ denotes the Euler characteristic. By Thom-Sebastiani theorem this number does not change if we add to $f$ a function with a quadratic singularity at $x$ (stable equivalence).

Let $M$ be a scheme with perfect obstruction theory (see \cite{BF}). Thus $M$ is locally represented as a scheme of critical points of a function $f$ on a manifold $X$. Then the above invariant gives rise to a $\Z$-valued constructible function $B$ on $M$. The 
global invariant is
$$\int_MB\,d\chi:=\sum_{n\in \Z}n \chi(B^{-1}(n))\,,$$ where $\chi$ denotes the Euler characteristic. Behrend proved that for a {\it proper} $M$ the invariant $\int_M B\, d\chi$ coincides with the degree of the virtual fundamental class $[M]^{virt}\in H_0(M)$ given by $\int_{[M]^{virt}}1$.

Now let us assume that $M\subset \CC^{ss}$ consists of Schur objects $E$ (see Section 1.3), such that $\op{cl}(E)=\gamma\in \Gamma$ is a fixed primitive class. Let us look at the contribution of $M$ to the motivic DT-invariant $a(\gamma)_{mot}$.
By definition it is equal to
$$\int_M{{\mathbb L}^{{1\over{2}}(1-\dim \op{Ext}^1(E,E))}\over {\mathbb L}-1}(1-MF(E))(1-MF_0(Q_E)){\mathbb L}^{-{1\over{2}}\op{rk}Q_E}\hat{e}_{\gamma}\,\,.$$
Mapping it to the quantum torus and taking the quasi-classical limit $q^{1/2}\to -1$, and taking into account
the relation $-a(\gamma)=\Omega(\gamma)$ for primitive $\gamma\in \Gamma$ (see Section 2.5), we obtain 
 that  Behrend's formula implies that the contribution of $M$ to the value 
$\Omega(\gamma)$ is equal to $\int_{[M]^{virt}}1$.

\subsection{Deformation invariance and intermediate Jacobian}

We also expect the following (not very precise) conjecture to be true as well.

\begin{conj}

 The collection $(\Omega(\gamma))_{\gamma \in \Gamma}$ is invariant with respect to the ``polarization preserving" deformations of $\CC$, in the case when $\CC$ is homologically smooth in the sense of \cite{KoSo3}.

\end{conj}

The motivation for the last Conjecture is the deformation invariance of the virtual fundamental class in the ``classical" Donaldson-Thomas theory.
Recall that  homologically smooth $Ext$-finite categories can be thought as non-commutative analogs of smooth proper schemes. Hence, we can expect that the moduli stacks of semistable objects in such categories are also proper in some
 sense. Therefore, we can also expect that the degree of the virtual fundamental class
  is invariant under deformations.

Also, we expect the following generalization of our theory in the case when $\k=\C$ and the $3d$ Calabi-Yau category is homologically smooth (see \cite{KoSo3}).

1) First, we recall that even without imposing the Calabi-Yau condition one expects that a triangulated compact homologically smooth $\A$-category $\CC$ (possibly $\Z/2\Z$-graded) admits (conjecturally) a non-commutative pure Hodge structure (see \cite{Ko2}, \cite{KKP}, \cite{KoSo3} about motivations, definitions as well as some conjectures and applications of this notion).
In particular, periodic cyclic homology groups $HP_{even}(\CC)$ (resp. $HP_{odd}(\CC)$) carry descending Hodge filtrations  
$$\begin{array}{l}HP_{even}(\CC)\dots\supset F^i_{even}\supset F^{i-1}_{even}\supset\dots,\, i\in \Z\\
HP_{odd}(\CC)\dots\supset F^i_{odd}\supset F^{i-1}_{odd}\supset\dots,\, i\in \Z+{1\over{2}}\,\,.
\end{array}$$ 
In $3d$ Calabi-Yau case we assume that the smallest non-trivial term of the filtration $F^{\bullet}_{odd}$ is $F^{-3/2}, \dim F^{-3/2}=1$. Moreover, in general, it is expected that there are lattices $K_{top}^{even}(\CC)$ and $K_{top}^{odd}(\CC)$ which belong to the corresponding periodic cyclic homology groups (they represent the non-commutative version of the image of the topological $K$-theory in the de Rham cohomology).

2) If $\CC$ is homologically smooth Calabi-Yau category then it is easy to see that (assuming the degeneration of the Hodge-to-de Rham conjecture, see \cite{KoSo3}) the moduli space ${\cal M}$ of formal deformations of $\CC$ is smooth of dimension $\dim {\cal M}={1\over{2}}\dim HP_{odd}(\CC)$ (this is a corollary of the formality of the little disc operad as well as the fact that the action of the Connes differential is represented by the rotation of the circle, which is homotopically trivial under the assumption). It is expected that the global moduli space also exists.
Notice that the Calabi-Yau structure on $\CC$ induces a symplectic structure on the vector space $HP_{odd}(\CC)$ and in the $3d$ case the moduli space ${\cal M}$ is locally embedded into $HP_{odd}(\CC)$ as a Lagrangian cone.

3) We expect that for an arbitrary triangulated compact homologically smooth $\A$-category $\CC$ one has a non-commutative version of the Deligne
 cohomology $H_{D}(\CC)$ which  fits into a short exact sequence
$$0\to HP_{odd}(\CC)/(F^{1/2}_{odd}+ K^{top}_{odd}(\CC))\to H_D(\CC)\to F^0_{even}\cap K^{top}_{even}(\CC)\to 0\,\,.$$
Morally, $H_D(\CC)$ should be thought as zero cohomology group of the homotopy colimit of the following diagram of cohomology theories:
$$\begin{CD}  @. HC_\bullet^-(\CC)\\
  @. @VVV  \\
K^{top}_\bullet(\CC) @>>>  HP_\bullet(\CC)
\end{CD} $$
where $HC_\bullet^-(\CC)$ is the negative cyclic homology.

Any object of $\CC$ should have its characteristic class in $H_{D}(\CC)$. More precisely, there should be a homomorphism of groups $ch_D: K_0(\CC)\to H_{D}(\CC)$ (in the case of  Calabi-Yau manifold it is related to holomorphic Chern-Simons functional). The reason for this is that every object $E\in Ob(\CC)$
  has natural characteristic classes in $K^{top}_0(\CC)$ and in $HC_0^-(\CC)$ whose images in $HP_0(\CC)$ coincide with each other. 
The total space ${\cal M}^{tot}$ of the fibration ${\cal M}^{tot}\to {\cal M}$
with the fiber $H_{D}(\CC)$ over the point $[\CC]\in {\cal M}$ should be a holomorphic symplectic manifold (cf. \cite{DonM}). Moreover, any fiber of this fibration (i.e. the group $H_{D}(\CC)$ for given $[\CC]$) is a countable union of complex Lagrangian tori.
By analogy with the commutative case we expect that the locus ${\cal L}\subset {\cal M}^{tot}$ consisting of values of $ch_D$ is a countable union of Lagrangian subvarieties. Every such subvariety can be either a finite ramified covering of ${\cal M}$ or a fibration over a proper subvariety of ${\cal M}$ with the fibers which are abelian varieties.

4) For generic $[\CC]\in {\cal M} $ one can use the triple $(K_0(\CC),H_{D}(\CC),ch_D)$ instead of the triple $(K_0(\CC),\Gamma,\op{cl})$. Analogs of our motivic Donaldson-Thomas invariants $A_V^{mot}\in {\cal R}_V$ will be  formal countable sums of points in $H_{D}(\CC)$ with ``weights" which are elements of the motivic ring ${\overline{D^{\mu}}}$. The pushforward map from
$H_{D}(\CC)$ to $\Gamma=F^0_{even}\cap K^{top}_{0}(\CC)$ gives the numerical DT-invariants. The continuity of motivic DT-invariants means that after taking the quasi-classical limit the weights become integer-valued functions on the set of those irreducible components of ${\cal L}$ which are finite ramified coverings on ${\cal M}$.

These considerations lead to the following

\begin{que} Is there a natural extension of the numerical DT-invariants to those components of ${\cal L}$ which project to a proper subvariety of ${\cal M}$?

\end{que}

\begin{rmk} Let us notice the similarity of the above considerations with those in the theory of Gromov-Witten invariants. Suppose $X$ is a $3d$ complex compact Calabi-Yau manifold with $H^1(X,\Z)=0$. Then we have an exact sequence
$$0\to H^3_{DR}(X)/(F^2H^3_{DR}(X)+H^3(X,\Z))\to H^4_D(X)\to H^4(X,\Z)\to 0\,\,,$$
where $H^4_D(X)={\mathbb H}^4(X,\Z\to {\cal O}_X\to \Omega^1_X)$ is the Deligne cohomology.
Then any curve $C\subset X$ defines the class $[C]\in H^4_D(X)$. For a generic complex structure on $X$ the class is constant in any smooth connected family of curves. Moreover, a stable map to $X$ defines a class in $H^4_D(X)$. Then we have exactly the same picture with holomorphic symplectic fibration ${\cal M}^{tot}\to {\cal M}$ with the Lagrangian fibers, as we discussed above. Similarly to the case of DT-invariants the GW-invariants appear as infinite linear combinations of points in $H^4_D(X)$, but this time with rational coefficients. We expect that the well-known relationship ``GW=DT" (see \cite{MaNOP}) should be a statement about the equality of the above-discussed counting functions (assuming positive answer to the above question).

\end{rmk}

\subsection{Absence of poles in the series $A_V^{\op{Hall}}$}

Here we are going to discuss a stronger version of the Conjecture 5.

\begin{conj} Let $\overline{D^+}:=\overline{{\cal M}^{\mu}}(Spec(\k))[{\mathbb L}^{-1/2}]$ be the ring of equivalence classes of motivic functions. Then the automorphism of the motivic quantum torus given by $x\mapsto A_V^{mot}x(A_V^{mot})^{-1}$ preserves the subring $\prod_{\gamma\in C(V)\cap \Gamma}\overline{D^+}\hat{e}_{\gamma}$ for all strict sectors $V\subset \R^2$.

\end{conj}

It is enough to check the conjecture for all $x=\hat{e}_{\gamma}, \gamma\in \Gamma$. Moreover, because of Factorization Property it is enough to consider the case when $V=l$ is a ray. In the latter case we can split the infinite product into those corresponding to different arithmetic progression, hence reducing the conjecture to the case when $Z(\Gamma)\cap \,l=\Z_{>0}\cdot\gamma_0$ for some non-zero $\gamma_0\in \Gamma$. Then we have
$$A_l^{mot}=A_l^{mot}(\hat{e}_{\gamma_0})=1+\sum_{n\geqslant 1}c_n\hat{e}_{\gamma_0}^n\in \overline{D^{\mu}}[[\hat{e}_{\gamma_0}]]\,\,.$$
Using the commutation relations in the motivic quantum torus we have:
$$A_l^{mot}(\hat{e}_{\gamma_0})\,\hat{e}_{\gamma}(A_l^{mot}(\hat{e}_{\gamma_0}))^{-1}=
\hat{e}_{\gamma}\,A_l^{mot}({\mathbb L}^{\langle \gamma_0,\gamma\rangle} \hat{e}_{\gamma_0})A_l^{mot}(\hat{e}_{\gamma_0})^{-1}\,\,.$$

Since for any series $f(t)=1+\dots$ we have
$${f({\mathbb L}^nt)\over{f(t)}}={f({\mathbb L}^nt)\over{f({\mathbb L}^{n-1}t)}}\dots{f({\mathbb L}t)\over{f(t)}}\,\,,$$
in order to prove the conjecture it suffices to check that 
$$  A_l^{mot}({\mathbb L} \hat{e}_{\gamma_0})A_l^{mot}(\hat{e}_{\gamma_0})^{-1}\in \overline{D^+}[[\hat{e}_{\gamma_0}]]\,\,.$$
Since in that case we are dealing with objects whose central charges belong to the ray $l$, we can restrict ourselves to the subcategory $\CC_l$. The latter can be thought of as a heart of the $t$-structure of an ind-constructible $3d$ Calabi-Yau category with  {\it vanishing} Euler form. More precisely, $\CC_l(\kk)$ is an abelian artinian category with
$\op{Hom}_{\CC_l(\kk)}(E,F):=\op{Ext}^0_{\CC(\kk)}(E,F)$. Then $K_0(\CC_l(\kk))\simeq \oplus_{E\ne 0}\Z\cdot[E]$, where the sum runs over the set of non-zero simple objects of  $\CC_l(\kk)$.

Next, we can reduce the conjecture to a special case when $\op{cl}_{\k}(E)=\gamma_0$ is a fixed class for all simple objects $E$ of $\CC_l(\kk)$. Indeed, let us consider an ind-constructible homomorphism $\op{cl}_{\k}^{\prime}:K_0(\CC_l(\kk))\to \Gamma^{\prime}:=\Z\oplus \Z$ such that $\op{cl}_{\k}^{\prime}(E)=(1,0)$ if
$\op{cl}_{\k}(E)=\gamma_0$ and $\op{cl}_{\k}^{\prime}(E)=(0,1)$ if $\op{cl}_{\k}(E)\in \{2\gamma_0,3\gamma_0,\dots\}$ for a simple object $E$. Let choose two complex numbers $z_1,z_2$ in such a way that $0<\op{Arg}(z_1)<\op{Arg}(z_2)<\pi$ and define a central charge $Z^{\prime}:\Gamma^{\prime}\to \C$  by the formula $Z^{\prime}((1,0))=z_1, Z^{\prime}((0,1))=z_2$. In this way we obtain a new stability structure on the triangulated envelope of $\CC_l(\kk)$ with the same heart. In particular, the element $A_l^{mot}$ will be decomposed into an infinite product:
$$A_l^{mot}=\prod^{\longrightarrow}A_{l^{\prime}}^{mot}$$
of series $A_{l^\prime}^{mot}$ corresponding to abelian categories $\CC_{l^{\prime}}(\kk)$ for the new stability structure. One of these categories will be the subcategory generated by  simple objects $E$ such that $\op{cl}_{\k}^{\prime}(E)=\gamma_0$. Let us call such category {\it pure of class $\gamma_0$}. All other categories $\CC_{l^{\prime}}(\kk)$ do not contain objects with the class $\gamma_0$. Repeating the procedure we reduce the conjecture to the case of pure category of the class $m\gamma_0$ for some $m\geqslant 1$. Similarly to the above arguments we can reduce it further to the case $m=1$. In this case the conjecture follows from the one below which concerns Hall algebras of categories which are no longer required to carry a Calabi-Yau structure.

In order to formulate this new conjecture we are going to use the following set-up. Let $({\CC},{\cal A})$ be a pair consisting of an ind-constructible triangulated $\A$-category over a ground field $\k$ and ${\cal A}\subset Ob({\CC})$ be an ind-constructible subset such that ${\cal A}(\overline{\k})$ is the heart of a bounded $t$-structure in ${\CC}(\overline{\k})$. We assume that simple objects of the abelian category ${\cal A}(\overline{\k})$ form a constructible subset of $Ob({\CC})(\overline{\k})$ and every object in ${\cal A}(\overline{\k})$ is a finite extension of simple ones. These data are equivalent to a special kind of an ind-constructible category with a stability structure. Namely, let us take $\Gamma:=\Z$ and define $\op{cl}_{\overline{\k}}(E)=1$ for every simple object of ${\cal A}(\overline{\k})$. It follows that $\op{cl}_{\overline{\k}}(F)=length(F)$ for any object of ${\cal A}(\overline{\k})$. Furthermore, we choose a complex number $z_0$ in the upper-half plane and define a central charge $Z:\Gamma\to \C$ by the formula $Z(1)=z_0$. Then ${\cal A}=\CC_l$ for $l=\R_{>0}\cdot z_0$. Therefore the element $A_l^{\op{Hall}}$ defined for this stability structure can be thought of as a series in one variable:
$$A_l^{\op{Hall}}(t)=1+\sum_{n\geqslant 1}c_nt^n.$$

Let us define a subalgebra $H^+(\CC) \subset H(\CC)$ to be the set of linear combination of elements of the form
  $${\mathbb L}^n \cdot [Z\to Ob(\CC)]$$
where $n\in \Z$ and $Z\to Ob(\CC)$ is a 1-morphism of ind-constructible stacks (see Section 4.2) with $Z$ being an ordinary constructible set endowed with trivial action of the trivial group. The multiplication law in $H(\CC)$ preserves such class of elements.

\begin{conj} The element $F_l(t):=A_l^{\op{Hall}}({\mathbb L}t)A_l^{\op{Hall}}(t)^{-1}$ belongs to the completed Hall algebra
$\widehat{H^+}(\CC)$ (i.e. we do not invert motives $[GL(n)], n\geqslant 1$ of the general linear groups).

\end{conj}

Below we discuss two special cases in which the above conjecture holds. But first we present a similar motivating statement in the case of finite fields.
Let $R$ be finitely generated algebra over a finite field ${\bf F}_q$, and $R-mod^f$ denotes the category of finite-dimensional (over ${\bf F}_q$) left $R$-modules. We define the Hall algebra $H(R-mod^f)$ as a unital associative algebra over the ring $\Z[{1\over{q}}]$ generated by the isomorphism classes $[M]$ of objects of $R-mod^f$
with the multiplication
$$[E]\cdot [F]=q^{-\dim \op{Hom}(F,E)}\sum_{\alpha\in \op{Ext}^1(F,E)}[E_{\alpha}]\,\,,$$
where, as before, $ E_{\alpha}$ denotes an extension with the class $\alpha$.
\begin{prp}
Let
$$A(t):=\sum_{[M]\in Iso(R-mod^f)}{[M]\over{\#\op{Aut}(M)}}t^{\,\dim M}\,\,.$$
Then
$$F(t):=A(qt)A(t)^{-1}\in H(R-mod^f)[[t]]\,\,.$$
Moreover,
$$F(t)=\sum_{I\subset R, I=RI,\, \dim R/I<\infty}[R/I]t^{\,\dim R/I}\,\,.$$

\end{prp}

Hence the quotient $F(t)$ does not have denominator $(q^n-1), n\geqslant 1$ and can be represented in terms of the ``non-commutative Hilbert scheme" of left ideals in $R$ of finite codimension.

{\it Proof.}
Let us make use of the basis of ``renormalized" elements $$\widehat{[E]}:={[E]\over{\#\op{Aut}(E)}}$$ in the $\Q$-algebra $H(R-mod^f)\otimes \Q$. Then the product can be rewritten in a more familiar form:
$$\widehat{[E]}\cdot\widehat{[F]} =\sum_{[G]}c_{\widehat{[E]},\widehat{[F]}}^{\widehat{[G]}}\widehat{[G]},$$
where the structure constant $c_{\widehat{[E]},\widehat{[F]}}^{\widehat{[G]}}\in \Z$ denotes the number of subobjects in $G$ isomorphic to $E$ and such that the quotient
is isomorphic to $F$.
In these notation we have:
$$A(t)=\sum_{[M]}\widehat{[M]}t^{\dim M}\,\,.$$
Since first statement of the Proposition follows from the second one, we are going to show the latter. In the new notation it becomes:
$$\sum_{I\subset R, I=RI, \,\dim R/I<\infty}\widehat{[R/I]}\#\op{Aut}(R/I)\,t^{\dim R/I}\cdot\sum_{[M]}\widehat{[M]}\,t^{\dim M}=$$
$$=\sum_{[N]}\widehat{[N]}\,q^{\dim N}t^{\dim N}\,\,.$$
Let us fix an object $N$, and consider the coefficient of 
the term $\widehat{[N]}t^{\dim N}$. In the RHS it is equal to $q^{\dim N}$.
 It is easy to see that the corresponding coefficient in the LHS is of the form
$$\sum_{I\subset R, I=RI, \,\dim R/I<\infty}\,\,\sum_{N^{\prime}\subset N,\, N^{\prime}\simeq R/I}\#\op{Aut}(R/I)=$$
$$=\sum_{I\subset R, I=RI, \,\dim R/I<\infty}\,\,\sum_{R/I\hookrightarrow N}1=\#N=q^{\dim N}.$$
Notice that in the last sum we consider all possible embeddings of $R/I$ to $N$ and every summand corresponds to a choice of a  cyclic vector in a cyclic submodule in $N$. This proves the Proposition. $\blacksquare$

The above Proposition suggest to interpret our category as a category of modules and then apply  similar arguments which reduce the sum (or even the motivic integral) to the sum over all cyclic submodules. It is useful to keep this in mind when considering two examples in the  next subsection.

\begin{rmk}  The subalgebra $H^+(\CC)$ of the Hall algebra has the advantage that one can apply the Euler characteristic
 $\chi$ to its elements fiberwise over $Ob(\CC)$ , and  get a constructible $\Z$-valued function (with constructible support) on the ind-constructible
 set $Iso(\CC)$ of isomorphism classes of objects of $\CC(\kk)$. The multiplication in $H^+(\CC)$ descends to a multiplication
 on the abelian group of such functions. It is easy to see that this multiplication is commutative, and one has
 $$\nu_E \times \nu_F=\nu_{E\oplus F}$$
where $\nu_E$ etc. are delta-functions (see Section 6.1). This follows from the fact that for any non-zero $\alpha \in \op{Ext}^1(F,E)$
 all objects $E_{t \alpha}$ are isomorphic to each other for $t\in \kk^\times$, and the Euler characteristic of ${\mathbf G}_m$ is zero.
\end{rmk}

\subsection{Reduction to the case of category of modules}

Here we present two special cases when the conjecture holds.

1) Assume that the abelian category ${\cal A}(\kk)$ contains only one\footnote{The arguments below extend immediately to the case of finitely many such objects.} (up to an isomorphism) simple object $E\ne 0$, and this object is defined over the field $\k$.  Hence $\op{Ext}^0_{\CC(\k)}(E,E)\simeq \k$. We also assume that $\A$-algebra $\op{Hom}^\bullet(E,E)$ is minimal, i.e. $m_1=0$, and hence $\op{Hom}^\bullet(E,E)=\op
{Ext}^\bullet(E,E)$.
\begin{prp} The category ${\cal A}(\k)$ is equivalent to the category $B-mod^{f,cont}$ of continuous finite-dimensional representations of a finitely generated topological algebra $B$.

\end{prp}
{\it Proof.} There is a general way to construct the algebra $B$ from the $\A$-structure.
Let $x_1,\dots,x_m$ be a basis in the vector space $(\op{Ext}^1(E,E))^{\ast}$. Then the higher compostion maps
$m_n:\op{Ext}^1(E,E)^{\otimes n}\to \op{Ext}^2(E,E), n\geqslant 2$ define a  linear map
$$\sum_{n\geqslant 2}m_n: (\op{Ext}^2(E,E))^{\ast}\to \k\langle\langle x_1,\dots,x_m\rangle\rangle=\prod_{n\geqslant 0}((\op{Ext}^1(E,E))^{\ast})^{\otimes n}\,.$$
We define a topological algebra
$B_E:=B$ as the quotient of
$\k\langle\langle x_1,\dots,x_m\rangle\rangle$ by the closure of the $2$-sided ideal generated by the image of $\sum_{n\geqslant 2}m_n$.

Next we observe that any object $M$ of ${\cal A}(\k)$ is a finite extension of objects isomorphic to $E$. Hence, it can be thought of as deformation of an object $mE:=E\oplus E\oplus\dots\oplus E$ ($m$ summands) preserving the filtration 
$$E\subset E\oplus E\subset\dots\subset mE\,,$$ where $m=length(M)$. Every such a deformation is given by a solution to the Maurer-Cartan equation
$$\sum_{n\geqslant 2}m_n(\alpha,\dots,\alpha)=0\,\,,$$
where $\alpha=(\alpha_{ij})$ is an upper-triangular matrix with coefficients from $\op{Ext}^1(E,E)$. It is easy to see that such a solution gives rise to a representation of the algebra $B$ in the upper-triangular matrices of finite size. Furthermore one checks that this correspondence provides an equivalence of categories 
$$F:{\cal A}(\k)\simeq B-mod^{f,cont}\,.$$ This proves the Proposition. $\blacksquare$

Notice that $$length(M)=\dim F(M)$$ for any object $M$.

Using the framework of finite-dimensional continuous representations we can modify the arguments from the proof of the Proposition 13 to the case of motivic functions instead of finite fields and obtain the formula
$$A^{\op{Hall}}({\mathbb L}t)A^{\op{Hall}}(t)^{-1}=\sum_{n\geqslant 0}[\pi: Hilb_n(B)\to Ob(\CC)]t^n,$$
where $Hilb_n(B)$ is the scheme of closed left ideals in $B$ of codimension $n$ (non-commutative analog of Hilbert scheme) and $\pi(I)=B/I$ for any such ideal.

2) Let us assume that $\k={\bf F}_q$ and ${\cal A}$ is an abelian $\k$-linear category such that every object has finitely many subobjects. We define the map $\op{cl}:K_0({\cal A})\to \Z$ such that $\op{cl}([E])=n$ if $E$ is simple object and
$\op{End}(E)\simeq {\bf F}_{q^n}$.

\begin{prp} Assume that ${\cal A}$ is a heart of a $t$-structure of a triangulated $\op{Ext}$-finite ${\bf F}_q$-linear $\A$-category $\CC$.
Let us consider the series
$$A(t):=\sum_{[M]\in Iso({\cal A})}{[M]\over{\#\op{Aut}(M)}}t^{\op{cl}(M)}\,\,.$$
Then we claim that
$$F(t):=A(qt)A(t)^{-1}=\sum_{[M],\,\, M \rm{ \,\,\,is \,\,\,cyclic}}c_M[M]t^{\op{cl}(M)}\,\,,$$
where $c_M\in \Z[{1\over{q}}]$, and the notion of a cyclic object is introduced below.
\end{prp}

We are going to reduce the proof to the case of modules over an algebra. Moreover we will give an explicit formula for the coefficients $c_M$.
In order to do  that we need the following categorical definition of a cyclic object.
 
\begin{defn} We say that an object $N$  in an artinian abelian category is cyclic if there is no epimorphism
$N\to E\oplus E$ where $E\ne 0$ is simple.
\end{defn}
In the category of finite-dimensional modules over an associative algebra (over any field), cyclic objects are the same as cyclic modules.

Any object $M\in {\cal A}$ admits a decomposition $M=\oplus_{\alpha}M_{\alpha}$ into a direct sum of indecomposables.
For each indecomposable summand $M_{\alpha}$ we have  a decomposition $M_{\alpha}^{ss}=\oplus_iE_{\alpha,i}$
of its maximal semisimple factor $M_{\alpha}^{ss}$ (called the cosocle of $M_\alpha$)
into a direct sum of simple objects $E_{\alpha,i}$.

Let us assume that $M$ is a cyclic object. It is equivalent to the condition that
 all simple factors $E_{\alpha,i}$ are pairwise different.
 Notice
that 
$$\op{End}(M^{ss})=\oplus_{\alpha,i} \op{End}(E_{\alpha,i})\simeq \oplus_{\alpha,i}{\bf F}_{q^{m_{\alpha,i}}}\,$$ where $m_{\alpha,i}=\op{cl}(E_{\alpha,i})\in \Z_{> 0}$.
Also, it follows from the cyclicity of $M$  that 
$$\op{End}(M)^{ss}\simeq\oplus_{\alpha}{\bf F}_{q^{n_{\alpha}}}$$ for some
positive integers $n_{\alpha}$. It follows from the definition that $m_{\alpha,i}$ is divisible by  $n_{\alpha}$ for any pair $(\alpha,i)$. Observe that in the above notation 
$$\#\op{Aut}(M)=q^r\prod_{\alpha}(q^{n_{\alpha}}-1)\,,$$ where $r$ is the dimension over ${\bf F}_q$ of the radical of $\op{End}(M)$. Now we claim that in the above Proposition 15
$$c_M={q^{\op{cl}(M)}\cdot \prod_{\alpha,i}{q^{m_{\alpha,i}}-1\over{q^{m_{\alpha,i}}}}\over{q^r\cdot \prod_{\alpha}(q^{n_{\alpha}}-1)}
}\,\,.$$
The property $n_{\alpha}|m_{\alpha,i}$  implies that $c_M\in \Z[{1\over{q}}]$.
We are going to prove the Proposition together with the above formula for $c_M$.

{\it Proof.} We may assume that ${\cal A}$ is generated by finitely many simple objects (but they can be defined over different finite extensions of ${\bf F}_q$). First, we claim that ${\cal A}$ is equivalent to the category $B-mod^{f,cont}$ of finite-dimensional continuous representations of a topological algebra $B$, similarly to the previous example. More precisely, let $N=\oplus_iE_i$ be the direct sum of all simple objects $E_i$, and set $C:=\op{End}(N)$. Then $C$ is a semisimple associative unital ${\bf F}_q$-algebra, which is isomorphic to $\oplus_i{\bf F}_{q^{\op{cl}(E_i)}}$. Let us consider $\op{Ext}^1(N,N)$ as a $C$-bimodule and take
$$G:=\op{Hom}_{C\otimes C^{op}-\op{mod}}(\op{Ext}^1(N,N), C\otimes C^{op})$$
to be the dual bimodule. The topological free algebra 
$$\prod_{n\geqslant 0}
{G\otimes_C\otimes G\otimes_C\dots\otimes_CG}\,\,\,\,({n\mbox{ tensor factors}})$$ contains a closed two-sided ideal generated by the image of the map $\sum_{n\geqslant 2}m_n$ (here we use the ``$\A$-origin" of our abelian category). We denote by $B$ the quotient of the free algebra by this ideal. Then $B$ can be thought of as a completed path algebra of the quiver defined by simple objects $E_i$ with the arrows which correspond to a basis of 
$\left(\op{Ext}^1(E_i,E_j)\right)^\ast$.  Similarly to the previously considered example,
 we have an equivalence of categories $\Psi:{\cal A}\simeq B-mod^{f,cont}$. Under this equivalence  simple object $E_i$
maps to the direct summand 
${\bf F}_{q^{\op{cl}(E_i)}}$ of $C$, hence $\dim \Psi(E_i)=\op{cl}(E_i)$. It follows that for any object $M$ we have $\dim  \Psi(M)=\op{cl}(M)$.

\begin{lmm} Module $M\in B-mod^{f,cont}$ is cyclic iff
$M^{ss}$ is cyclic.
Moreover $v\in M$ is a generator iff its projection $\overline{v}\in M^{ss}$ is a generator.

\end{lmm}
{\it Proof of lemma.} The first statement follows directly from the definition of a cyclic object. 
 In order to prove the second statement assume that $\overline{v}\in M^{ss}$ is a generator. We want to prove that the quotient $M/Bv=0$. If this is not the case then we have an epimorphism $M/Bv\to E_{i_0}$ to a simple module $E_{i_0}$. It follows that we have an epimorphism
$M^{ss}\to E_{i_0}$ such that $\overline{v}\mapsto 0$. This contradicts to the assumption that $\overline{v}\in M^{ss}$ is a generator. The lemma is proved. $\blacksquare$

In order to finish the proof of the Proposition, it is enough to check that the coefficient $c_M$ given by a product formula on the previous page, 
 is equal to the number of isomorphism classes of generators $v\in M$ up to an automorphism of $M$. In order to do that we observe that
the product $\prod_{\alpha,i}(q^{m_{\alpha,i}}-1)$ from the formula for $c_M$ is in fact equal to the number of generators of $M^{ss}$. Furthermore, the factor ${q^{\op{cl}(M)}\over{\prod_{\alpha,i}q^{m_{\alpha,i}}}}$ is equal to the number of liftings of a generator of $M^{ss}$ to a generator of $M$ (this number is the number of elements in the kernel of the projection $M\to M^{ss}$). Finally, we recall that $q^r\prod_{\alpha}(q^{n_{\alpha}}-1)=\#\op{Aut}(M)$. Applying the above lemma we finish the proof of the Proposition. $\blacksquare$

\begin{rmk} It looks plausible that the Proposition holds without the assumption that ${\cal A}$ is a $t$-structure of an $\A$-category.

\end{rmk}

We do not know the ``motivic" analog of the above Proposition. In that case one should replace ${\cal A}$ by an ind-constructible 
abelian category over  any field. There is a notion of  semisimple and cyclic modules, it is preserved under field extensions\footnote{Notice that notions of simple or indecomposable objects are {\it not} preserved under the field extension.}.
It looks natural to expect that an analog of the quotient ${\prod_{\alpha,i}(q^{m_{\alpha,i}}-1)\over{\prod_{\alpha}(q^{n_{\alpha}}-1)}}$ is the motive
 $\underline{\op{Aut}}(M^{ss})/\underline{\op{Aut}}'(M)$ where $\underline{\op{Aut}}(M^{ss})$ is the affine group scheme of automorphisms of $M^{ss}$
 and $\underline{\op{Aut}}'(M)$ is the image of the scheme of automorphisms of $M$ in $\underline{\op{Aut}}(M^{ss})$. Both groups schemes are algebraic tori.
 Although the motivic version is not absolutely clear, we can write down the ``numerical" version, which is the result of the quasi-classical limit $q^{1/2}\to -1$ (equivalently, this is the result of taking the Euler characteristic of the corresponding motives).

It follows from the Proposition that in the quasi-classical limit only those terms in the formula for $c_M$ are non-zero for which $\underline{\op{Aut}}(M^{ss})=\underline{\op{Aut}}'(M)$. Let us call such objects {\it special cyclic}.  A cyclic object is special cyclic iff under the extension of scalars
to $\kk$ the cococles of all indecomposable summands (i.e. objects $M^{ss}_\alpha$ in our notation) are simple.

In the case of finite-dimensional modules over an associative algebra $A$, a cyclic object (or module) $M$ is special iff the scheme of left ideals $I\subset A$ such that $M\simeq A/I$ has Euler characteristic $1$. For non-special cyclic objects the corresponding Euler characteristic vanishes.

Let us return to our considerations in the case of $3d$ ind-constructible Calabi-Yau category over a field $\k$ of characteristic zero.  We reduced the main conjecture to the case of a single ray, hence ${\cal A}(\kk)$ is the heart of a $t$-structure of $\CC_l$. In this case isomorphism classes of special cyclic objects $M$ with the fixed class $\op{cl}(M)$ form a constructible set ${\cal SC}_n$. Thus, we arrive to the following formula
$$\begin{array}{rcl}\chi_\Phi(F_l(t))&=&\chi_\Phi(A_l^{mot}({\mathbb L}t)A_l^{mot}(t)^{-1})=\\
&=&\sum_{n\geqslant 0}t^n\int_{{\cal SC}_n}(-1)^{(M,M)_{\leqslant 1}}(1-\chi(MF(M)))d\chi\,\,,
\end{array}$$
where 
$\int_V f d\chi=\sum_{n\in\Z}n \chi(f^{-1}(n))$ denotes the ``integral over Euler characteristic" $\chi$ of the map $f:V\to \Z$, and $\chi_\Phi$
 is the composition of the homomorphism $\Phi$ from the motivic Hall algebra to the motivic quantum torus (restricted to subalgebra $H^+(\CC)\subset H(\CC)$), and of the Euler characteristic morphism acting on coefficients as $\chi:\overline{D^+}\to \Z$.

We remark that the RHS does not depend on the orientation data.

\begin{conj} In case if the category $\CC$ is not endowed with  orientation data the above procedure gives rise to  well-defined stability data on the graded Lie algebra $\g_{\Gamma}$ of Poisson automorphisms of the algebraic Poisson torus $\op{Hom}(\Gamma,{\bf G}_m)$ as well as a continuous local homeomorphism $Stab(\CC,\op{cl})\to Stab(\g_{\Gamma})$.

\end{conj}

\subsection{Evidence for the integrality conjecture}

In this section we present arguments in favor of the integrality of the ``numerical" DT-invariants $\Omega(\gamma)$. Recall that if $E$ is an object of a $\k$-linear triangulated category, then we say that $E$ is a Schur object if 
$$\op{Ext}^{<0}(E,E)=0,\,\,\,\op{Ext}^0(E,E)\simeq \k\cdot id_E\,\,.$$

Let us assume now that $\CC$ is an ind-constructible $3d$ Calabi-Yau category generated by a Schur object $E\in \CC(\k)$ in the sense that the category $\CC(\overline{\k})$ consists of finite extensions of the shifts $E[i], i\in \Z$. In this case $K_0(\CC(\overline{\k}))\simeq \Z\cdot \op{cl}_{\kk}(E)$. We take
$\Gamma=K_0(\CC(\overline{\k}))$ and the trivial skew-symmetric form on $\Gamma$.
For any $z\in \C, \op{Im}z>0$ our category carries an obvious stability condition $\sigma_z$ such that $Z(E):=Z(\op{cl}_{\kk}(E))=z,\,\, \op{Arg}(E)=\op{Arg}(z)\in (0,\pi)$. All objects $F\in \CC^{ss}(\overline{\k})$ with $\op{Arg}(F)=\op{Arg}(E)$ are $n$-fold extensions of copies of $E$ for some $n\geqslant 1$. We denote by $l$ the ray $\R_{>0}\cdot z$.

In the previous section we obtained a formula for ${\chi}_{\Phi}(F_l(t))$ in terms of the integral over Euler characteristic over the moduli space of special cyclic objects of $\CC_l(\kk)$. We are going to make it more explicit further, by using the potential of $E$. Let us recall (see Section 3.3) that with the object $E$ we associate a collection of cyclically invariant polylinear maps 
$$W_N: (\op{Ext}^1(E,E)^{\otimes N})^{\Z/N\Z}\to \k, N\geqslant 3\,\,,$$
$$ a_1\otimes\dots\otimes a_N\mapsto W_N(a_1,\dots,a_N)=(m_{N-1}(a_1,\dots,a_{N_1}),a_N)\,\,.$$
 Let us choose a basis $x^1,\dots,x^m$ in $\op{Ext}^1(E,E)$. Then to any $n\geqslant 0$ and collection of matrices $X_1,\dots,X_m\in Mat(n\times n,\kk)$ we associate the number
$$W_N^{(n)}(X_1,\dots,X_m)={1\over{N}}\sum_{1\leqslant i_1,\dots,i_N\leqslant m}W_N(x^{i_1},\dots,x^{i_N})\op{Tr}(X_{i_1}\dots
X_{i_N})\,\,.$$
Thus we have a polynomial on $\kk^{mn^2}$. The series
$$W^{(n)}=\sum_{N\geqslant 3}W_N^{(n)}$$
is a formal function on the formal neighborhood of the reduced closed subscheme $Nilp_{m,n}\subset {\bf A}^{mn^2}$ whose $\kk$-points are collections of matrices 
$$(X_1,\dots,X_m)\in Mat(n\times n,\kk)$$
 which satisfy the property that there exists a basis in which all $(X_i)_{i=1,\dots, m}$
 are strictly upper triangular. Equivalently, $\op{Tr}(X_{i_1}\dots X_{i_r})=0$ for any sequence of indices $i_\bullet \in \{1,\dots,m\}$ with $r\ge (n+1)$.
 This property ensures that $W_N^{(n)}$ is well-defined in a formal neighborhood of $Nilp_{m,n}$.
Then
$$
\chi_\Phi(F_l(t))=1+\sum_{n\geqslant 1}t^n\int_{ Nilp_{m,n}^{SC}/PGL(n)}(-1)^{n^2(1-m)}(1-\chi(MF_{(X_1,\dots,X_m)}(W^{(n)})))d\chi\,\,,
$$
where $Nilp_{m,n}^{SC}$ for $n\geqslant 1$ is a subscheme of $Nilp_{m,n}$ whose $\kk$-points consists  of those collections $(X_1,\dots,X_m)$ for which 
$$\op{codim}\left(\sum_i \op{Im}(X_i)\right)=1\,.$$

Let us comment on this formula. First we remark that it is sufficient to integrate over the set $Crit(W^{(n)})$ of critical points of $W^{(n)}$, since for all non-critical points $\chi(MF_{(X_1,\dots,X_m)}(W^{(n)}))=1$. Repeating the arguments of the previous section we obtain that $\CC(\kk)$ is equivalent to the category $B_W-mod^{f,cont}$ of continuous finite-dimensional representations over $\kk$ of the topological $\k$-algebra
$B_W=\k\langle\langle x_1,\dots,x_m\rangle\rangle/\overline{(\partial_{x_i}W)}, 1\leqslant i\leqslant m$, where $x_i,1\leqslant i\leqslant m$ are the coordinates corresponding to the chosen basis $x^i, 1\leqslant i\leqslant m$, and  $\overline{(\partial_{x_i}W)}$ denotes the closure of the $2$-sided ideal generated by the cyclic derivatives of the cyclic potential
$W=\sum_{N\geqslant 3}N ^{-1} W_N$. Indeed, it is straightforward to see that a point $(X_1,\dots,X_m)\in Nilp_{m,n}(\kk)$ gives rise to a continuous $n$-dimensional representation of $W$ if an only if it belongs to $Crit(W^{(n)})$. In terms of the category
$\CC(\kk)$ these points correspond to  $n$-fold extensions of the Schur object $E$ by itself. Special cyclic modules correspond to critical points belonging to $Nilp_{m,n}^{SC}\subset Nilp_{m,n}$.

Considering an object $M$ of length $n$ as an upper-triangular deformation of the ``free" object $nE=E\oplus\dots\oplus E$ ($n$-times) we see that
$$(M,M)_{\leqslant 1}=(nE,nE)_{\leqslant 1}+r,\,\,\,r:=\dim \op{Im}(W^{(n)})^{\prime\prime}_{(X_1,\dots,X_m)}\,\,.$$
Then
$$1-MF_{(X_1,\dots,X_m)}(W^{(n)})=(1-MF(E))(1-MF_0(Q_E))\,\,,$$
where $Q_E$ is a quadratic form and $\op{rk}Q_E=r$.
Indeed, $W^{(n)}$ coincides with the potential $W_{nE}$ of the object $nE$ under the isomorphism $\op{Ext}^1(nE,nE)\simeq \kk^{mn^2}$. Thus we see that

$$(-1)^{(M,M)_{\leqslant 1}}=(-1)^{n^2(1-m)+r}\,.$$ Since $\chi(1-MF_0(Q_E))=(-1)^{\op{rk}Q_E}$ we obtain the desired formula for $\chi_\Phi(F_l(t))$.
Alternatively, in the integral one can replace the quotient $Nilp_{m,n}^{SC}/PGL(n)$ by $Nilp_{m,n}^{cycl}/GL^{(1)}(n)$. Here $Nilp_{m,n}^{cycl}\subset Nilp_{m,n}$
consists of collection of matrices such that $$\kk\langle\langle X_1,\dots,X_m\rangle\rangle v_1=\kk^n$$ where $v_1:=(1,0,\dots,0)$ is the first base vector, and the group $GL^{(1)}(n)\subset GL(n)$ is the stabilizer of $v_1$. Notice that $GL^{(1)}(n)$ acts freely on $Nilp_{m,n}^{cycl}$.
   The reason is that the contribution of non-special cyclic objects vanishes as follows from the vanishing of the Euler characteristic
 of corresponding schemes of modules with chosen cyclic generators.

\begin{conj} We have:
$$\chi_\Phi(F_l(t))=\prod_{n\geqslant 1}(1-t^n)^{n\,\Omega(n)}\,\,,$$
where all $\Omega(n)=\Omega(n\op{cl}_{\kk}(E))$ are integer numbers (see Section 1.4).

\end{conj}

Let us illustrate the conjecture in  few examples.

1) Let $m=0$ (i.e. the case of just one spherical object). Then $W=0$ and $B_W=\k$. The only non-trivial cyclic representation have dimension one, hence 
$\chi(F_l(t))=1-t$. Then we have 
$$\Omega(1)=1,\,\,\,\Omega(n)=0\mbox{ for }n>1\,.$$

2) Let $m=1, W^{(n)}(X_1)=\op{Tr}(X_1^d)$ for $d=3,4,\dots$.
Then $B_W=\k[x_1]/(x_1^{d-1})$. There is a unique isomorphism class of cyclic $B_W$-modules in an dimension $n=0,1,2,\dots,d-1$. One can show directly that 
$$\chi(F_l(t))=(1-t)^{d-1},\,\,\, 
\Omega(1)=d-1,\,\,\,\Omega(n)=0\mbox{ for }n>1\,.$$

3) Let $m\ge 1 $ be arbitrary and $W=0$. In this case
$$\chi(F_l(t))=1+\sum_{n\geqslant 1}(-1)^{n^2(1-m)}\chi(Nilp_{m,n}^{cycl}/GL^{(1)}(n))t^n\,\,.$$

Euler characteristic $\chi(Nilp_{m,n}^{cycl}/GL^{(1)}(n))$ coincides with the Euler characteristic of the non-commutative Hilbert scheme $H_{n,1}^{(m)}$ from \cite{Re2}. The latter parametrizes left ideals of codimension $n$ in the free algebra
$\kk\langle x_1,\dots,x_m\rangle$. The reason why we can disregard all non-nilpotent collections $(X_1,\dots,X_m)$ of  matrices is that the latter carries a free action of the group ${\bf G}_m$, such that $X_i\mapsto \lambda X_i, 1\leqslant i\leqslant m$ where $\lambda\in {\bf G}_m(\kk)$. 
Hence the corresponding Euler characteristic is trivial. Then using explicit formulas from \cite{Re2} we obtain
$$G_{(m)}(t):=\chi(F_l(t))=\sum_{n\geqslant 0}{(-1)^{n(1-m)}\over{(m-1)n+1}}\binom {mn}{n}t^n.$$
Notice that this series can be written
as
$$\exp\left(\sum_{n\geqslant 1}{(-1)^{n(1-m)}\over{mn}}\binom {mn}{n}t^n\right)\,\,.$$
For $m=1$ we have $\Omega(1)=-1, \Omega(n)=0, n\geqslant 2$.
In general
$$\Omega(n)={1\over{mn^2}}\sum_{d|n}\mu({n}/{d})\binom {md}{d}(-1)^{(m-1)d+1}\,\,,$$
where $\mu(k)$ is the M\"obius function (for $m=2$ see the entry $A131868$ in the online Encyclopedia of integer sequences).

\begin{rmk}
One can check that the generating function $G=G_{(m)}$ is algebraic: it satisfies the equation\footnote{Compare with the algebraic series in the Introduction, section 1.4.}
$$G(t)+t(-1)^m(G(t))^m-1=0\,\,.$$
An interesting question arises: which algebraic functions admit multiplicative factorization of the form $\prod_{n\geqslant 1}(1-t^n)^{n\,\Omega(n)}$,
where all $\Omega(n)$ are integer numbers?

\end{rmk}

\section{Donaldson-Thomas invariants and cluster transformations}

\subsection{Spherical collections and mutations}

Let ${\cal C}$ be a $3$-dimensional  ind-constructible Calabi-Yau category over a field $\k$ of characteristic zero. Assume that it is endowed with a finite collection of spherical generators
${\cal E}=\{E_i\}_{i\in I}$ of ${\cal C}$ defined over $\k$.  Then
$\op{Ext}_{\cal C(\k)}^\bullet(E_i,E_i)$ is
isomorphic to $H^{\bullet}(S^3,\k ), \,\,i\in I$.
The  matrix of the Euler form (taken with the {\it minus} sign)
$$a_{ij}:=-\chi\left(\op{Ext}_{\cal C(\k)}^{\bullet}(E_i,E_j)\right)$$
is integer and skew-symmetric.
 In fact, the ind-constructible category   ${\cal C}$ can be canonically reconstructed from the (plain, i.e. not ind-constructible)
 $\k$-linear Calabi-Yau $\A$-category $\CC(\K)$, or even from its full subcategory consisting of the collection ${\cal E}$.
 In what follows we will omit the subscript $\CC(\k)$ in the notation for $\op{Ext}^\bullet$-spaces.

\begin{defn} The collection ${\cal E}$ is called cluster if for any $i\ne j$ the graded space $\oplus_{m\in \Z} \op{Ext}^m(E_i,E_j)$ is either zero, or it is concentrated
 in one of two degrees $m=1$ or $m=2$ only.

\end{defn}

We will assume that our collection is cluster. In that case $K_0(\CC(\kk))\simeq \Z^I$ with the basis formed by the isomorphism classes $[E_i], i\in I$.

With the cluster collection ${\cal E}$ we  associate a quiver $Q$ such that
$Q$ does not have oriented cycles of lengths $1$ and $2$, and $a_{ij}\geqslant 0$ is the number of arrows from $i$ to $j$ (we identify the set of arrows from $i$ to $j$ with a basis in $\op{Ext}^1(E_i,E_j)$).
Then the potential for the object $E=\oplus_{i\in I}E_i$ gives rise to the potential $W=W_Q$
of the quiver $Q$, i.e. the restriction of the potential to $\oplus_{i,j\in I} \op{Ext}^1(E_i,E_j)$. 
The latter is an infinite linear combination of cyclic words (see \cite{DWZe}, \cite{Ze} where the potential appears abstractly without the relation with Calabi-Yau categories). Any such linear combination is called a potential of $Q$. In our case the potential is automatically minimal, i.e. all words have length at least $3$. The group of continuous automorphisms of the completed path algebra of $Q$ preserving the projectors $pr_i,\,\,i\in I$, acts on the set of potentials of $Q$.
 We call it the gauge action.
Let us state the following general result.

\begin{thm} Let $\CC$ be a 3-dimensional $\k$-linear Calabi-Yau category generated by a finite  collection ${\cal E}=\{E_i\}_{i\in I}$ of generators satisfying the condition that
\begin{itemize}
\item
 $\op{Ext}^0(E_i,E_i)=\k \,id_{E_i}$, 
\item $\op{Ext}^0(E_i,E_j)=0$ for any $i\ne j$,
\item $\op{Ext}^{<0}(E_i,E_j)=0$, for any $i,j$.
\end{itemize}
The equivalence classes
of such categories with respect to $\A$-transformations preserving the Calabi-Yau structure and the collection ${\cal E}$,
 are in one-to-one correspondence with the gauge equivalence classes of pairs
$(Q,W)$ where $Q$ is a finite oriented quiver (possibly with cycles of length $1$ or $2$) and $W$
is a minimal potential of $Q$ (i.e. its Taylor decomposition starts with terms of degree at least $3$).

\end{thm}

The case of cluster collections corresponds to quivers without oriented cycles of length $1$ and $2$.

 {\it Proof.} We will present  the proof of the Theorem in the case of the category with single object $E$ (i.e. $A=\op{Hom}(E,E)$ is a $3d$ Calabi-Yau algebra). The general case  can be proved in a similar way.

Let $Q$ be a quiver with one vertex and $|J|$ loops, where $J$ is a finite set. We assume that $Q$ is endowed with the potential $W_0$. We would like to construct a $3d$ Calabi-Yau category with a single object $E$ such that the number of loops in $Q$ is equal to $\op{Ext}^1(E,E)$ and the restriction of the potential of the category to $\op{Ext}^1(E,E)$
coincides with the given $W_0$. Our considerations proceed such as follows. Assuming that such a category exists we will find an explicit formula for the potential on $A=\op{Hom}(E,E)$. Then we simply take this explicit formula as the definition.
If the desired category is constructed then we can consider the graded vector space $\op{Ext}^{\bullet}(E,E)[1]$ which decomposes as 
$$\op{Ext}^{0}(E,E)[1]\oplus \op{Ext}^1(E,E)\oplus \op{Ext}^2(E,E)[-1]\oplus \op{Ext}^3(E,E)[-2]\,\,.$$
The first and the last summand are isomorphic to $\k[1]$ and $\k[-2]$ respectively, and the middle two summands are dual two each other. We introduce graded coordinates on $\op{Ext}^{\bullet}(E,E)[1]$ and denote them such as follows:

a) the coordinate $\alpha$ of degree $+1$ on $\op{Ext}^{0}(E,E)[1]$;

b) the coordinate $a$ of degree $-2$ on $\op{Ext}^3(E,E)[-2]$;

c) the coordinates $x_i, i\in J$ of degree $0$ on $\op{Ext}^1(E,E)$;

d) the coordinates $\xi_i, i\in J$ of degree $-1$ on $\op{Ext}^2(E,E)[-1]$.

The Calabi-Yau structure on $A$ gives rise to the minimal potential $W=W(\alpha, x_i,\xi_i,a)$, which is a series in cyclic words on the space  $\op{Ext}^{\bullet}(E,E)[1]$. If it arises from the pair $(Q,W_0)$, then the restriction $W(0,x_i,0,0)$ must coincide with $W_0=W_0(x_i)$. Furthermore, $A$ defines a non-commutative formal pointed graded manifold endowed with a symplectic structure (see \cite{KoSo3}). The potential $W$ satisfies the ``classical BV equation" $\{W,W\}=0$, where $\{\bullet,\bullet \}$ denotes the corresponding Poisson bracket.

With these preliminary considerations we see what problem should be solved. We need to construct an extension of $W_0$ to the formal series  $W$ of degree $0$ in cyclic words on the graded vector space $\k[1]\oplus \k^J\oplus \k^J[-1]\oplus \k[-2]$,  satisfying the classical BV-equation with respect to the Poisson bracket
$$\{W,W\}=\sum_i\partial W/\partial x_i\,\partial W/\partial \xi_i+\partial W/\partial \alpha \,\partial W/\partial a \,.$$

Here is the construction. Let us start with the potential
$$W_{can}=\alpha^2 a+\sum_{i\in J}(\alpha x_i \xi_i-\alpha \xi_i x_i)\,\,.$$
This potential makes the above graded vector space into a $3d$ Calabi-Yau algebra with associative product and the unit. The multiplication vanishes on the graded components
$$\op{Ext}^1(E,E)\otimes \op{Ext}^1(E,E)\to \op{Ext}^2(E,E)$$ 
and 
is a non-degenerate bilinear form on components 
 $$\op{Ext}^1(E,E)\otimes \op{Ext}^2(E,E)\to \op{Ext}^3(E,E)\simeq \k\,\,.$$ Now we see that starting with an arbitrary minimal potential $W_0$ on $\op{Ext}^1(E,E)$ we can lift it to the minimal potential on $\op{Ext}^{\bullet}(E,E)$ by setting
$$\widehat{W}:=W_{can}+W_0\,.$$ We claim that $\{\widehat{W},\widehat{W}\}=0$. Indeed, we have
$\{W_{can},W_{can}\}=\{W_0,W_0\}=0$. Moreover,
$$\{W_{can},W_0\}=\alpha\sum_{i\in J}[x_i,\partial W_0/\partial x_i]=0$$ (we use here the  well-known identity  $\sum_{i\in J}[x_i,\partial W_0/\partial x_i]=0$).

Next we need to check compatibility of the above construction with the gauge group action. Let $G_0$ be the subgroup of the grading preserving automorphisms of the group of continuous automorphisms of the algebra of formal series $\k\langle\langle \alpha,x_i,\xi_i,a\rangle\rangle, i\in J$.
Let ${\cal J}\subset  \k\langle\langle \alpha,x_i,\xi_i,a\rangle\rangle$ be a closed $2$-sided ideal generated by $\alpha,a$ and $\xi_i$ for $ i\in J$. Since every generator of ${\cal J}$ has non-zero degree we conclude that the group $G_0$ preserves ${\cal J}$ (it can be deduced from the fact that it transforms generators into series of non-zero degrees).  Therefore we obtain a homomorphism of groups $G_0\to \op{Aut}(\k\langle\langle x_i\rangle\rangle), i\in J$. The restriction of the potential $W$ to $\op{Ext}^1(E,E)$ defines a surjection from the set of $\A$-equivalence classes of $3d$ Calabi-Yau algebras to the gauge equivalence classes of $(Q,W_0)$, where $Q$ is a quiver with one vertex endowed with the minimal potential $W_0$. Such algebras can be thought of as deformations of the ``ansatz", which is a $3d$ Calabi-Yau algebra $A_{can}$ corresponding to the potential $W_{can}$.

Finally we are going to show that the above surjection is in fact a bijection. The latter will follow from the equivalence of the corresponding deformation theories.
The deformation theory of the Calabi-Yau algebra $A_{can}$ is controlled by a  DGLA  $\g_{A_{can}}=\oplus_{n\in \Z}\g_{A_{can}}^n$, which is a DG Lie subalgebra of the DG Lie algebra 
$$\widehat{\g_{A_{can}}}=\left(\prod_{n\geqslant 1}Cycl^n((A[1])^{\ast})\right)[-1] $$
 of all cyclic series in the variables $\alpha,x_i,\xi_i,a, i\in J$ (the Lie bracket is given by the Poisson bracket and the differential is given by $\{W_{can},\bullet\}$). Namely, the component of $\g_{A_{can}}$ of degree $N$ consists of all cyclic series which contain at least $2+N$ letters $\alpha,x_i,\xi_i,a, i\in J$. We will call the degree defined in terms of these letter a {\it cyclic degree} in order to distinguish it from the {\it cohomological degree} of complexes.
Notice that the set of $\A$-equivalence classes of minimal $3d$ Calabi-Yau agebras can be identified with the set of gauge equivalence classes of solutions $\gamma\in \g_{A_{can}}^1$ to the Maurer-Cartan equation
$$d\gamma+{1\over{2}}[\gamma,\gamma]=0\,\,.$$
Similarly, the set of gauge equivalence classes of minimal potentials on $A^1=\op{Ext}^1(E,E)$ can be identified with the set of gauge equivalence classes of solutions to the Maurer-Cartan equation in the DGLA ${\h}={\h}^0\oplus {\h}^1$, where
$${\h}^0=\prod_{n\geqslant 1}((A^1)^{\ast})^{\otimes n}\otimes A^1\,\,,\,\,\,\,
{\h}^1=\prod_{n\geqslant 3}Cycl^n((A^1)^{\ast})\,\,.$$
Here we identify ${\h}^0$ with the Lie algebra of continuous derivations of the topological algebra $\k\langle\langle x_1,\dots,x_n\rangle\rangle$ preserving the augmentation ideal $(x_1,\dots,x_n)$, and we identify ${\h}^1$ the ${\h}^0$-module of minimal cyclic potentials on $A^1$.

The above construction of the ``lifting" $\widehat{W}=W_0+W_{can}$ can be interpreted as a homomorphism of DGLAs $\psi:{\h}\to \g_{A_{can}}$. Namely, ${\h}^0$ is identified (after the shift $[1]$) with the space of such cyclic series in $x_i,\xi_i, i\in J$ which contain exactly one of the variables $\xi_i$ and at least one of the variables $x_j$ for some $i,j\in J$. Similarly ${\h}^1$ is identified with the space of cyclic series in $x_i,i\in J$ which has terms of degree at least $3$.

We claim that $\psi$ induces an epimorphism (previous considerations ensure that it is a monomorphism) of cohomology groups in degree $1$, and for both DGLAs ${\h}$ and  $\g_{A_{can}}$ there is no cohomology in degree greater or equal than $2$. This would imply the desired surjectivity of $\psi$.

Notice that the differential $\{W_{can},\bullet\}$ preserves the difference between cyclic and cohomological degree. It follows that  the complex $\g_{A_{can}}$ is a direct summand of the complex $\widehat{\g_{A_{can}}}$. The latter is dual to the cyclic complex $CC_{\bullet}(A_{can})$. Let $A_{can}^+\subset A_{can}$ be a non-unital $\A$-subalgebra consisting of terms of positive cohomological degree. Then, one has for the cyclic homology: $HC_{\bullet}(A_{can})\simeq HC_{\bullet}(A_{can}^+)\oplus HC_{\bullet}(\k)$. In terms of the dual complex this isomorphism means the decomposition into a direct sum of the space of cyclic series in variables $x_i,\xi_i,a, i\in J$ (corresponds to $(HC_{\bullet}(A_{can}^+))^{\ast}$) and the space of cyclic series in the variable $\alpha$ of odd cyclic degree (corresponds to $(HC_{\bullet}(\k)^{\ast}$). It is easy to see that the series in the variable $\alpha$ do not contribute to the cohomology of $\g_{A_{can}}\subset\widehat{\g_{A_{can}}} $.

The cohomological degree of series in variables $x_i,\xi_i,a, i\in I$ is non-positive. Recall that we shifted the grading in Lie algebras by $1$ with respect to the cohomological grading.  Hence $H^{\geqslant 2}(\g_{A_{can}})=0$. Also, it is immediate that $H^1(\g_{A_{can}})$ is isomorphic to the space of cyclic series in the variables $x_i,i\in I$  with terms of degree at least $3$. Hence $H^1(\g_{A_{can}})\simeq {\h}^1\simeq H^1({\h})$ (the latter holds since the differential on ${\h}$ is trivial). This concludes  the proof. $\blacksquare$

Next, we will introduce the notion of a {\it mutation} on the set of cluster collections in a given category $\CC$. Let us choose an element of $I$ which we will denote by $0$. We are going to write $i<0$ if $a_{i0}>0$, and $i>0$ if $i\ne 0$ and $a_{i0}\leqslant 0$. The mutation of ${\cal E}$ at the object $E_0$ is defined as a new spherical collection ${\cal E}^{\prime}=(E_i^{\prime})_{i\in I}$ such that:

$$\begin{array}{ll}
E_i^{\prime}=E_i, & i<0,\\
 E_0^{\prime}=E_0[-1], &\\
  E_i^{\prime}=R_{E_0}(E_i), &i>0\,\,.
  \end{array}$$
where $R_{E_0}(E_i):=
Cone(E_0\otimes \op{Ext}^{\bullet}(E_0,E_i)\to E_i)$ is the reflection functor given by 
 the cone of the natural evaluation map. Explicitly, 
the object $E_i^{\prime}$ for  $i>0$ fits in the exact triangle $$E_i\to E_i^{\prime}\to E_0\otimes \op{Ext}^1(E_0,E_i)\,\,.$$
 Notice that all objects $E_i^{\prime},i\ne 0$ belong to the abelian category generated by $E_i,i\in I$.
We remark that the spherical collection  ${\cal E}^{\prime}$ is not necessarily a cluster one.

At the level of the lattice $\Gamma:=\Z^I$ the change of the spherical collections
${\cal E}\to {\cal E}^{\prime}$ corresponds to the following
relation between the  basis $v_i:= \op{cl}_{\kk}(E_i),\,\,i\in I$ and the mutated basis $ v_i^{\prime}=\op{cl}_{\kk}(E_i^{\prime}),\,\,i\in I$:

$$ \begin{array}{l}v_i^{\prime}=v_i, i< 0, \\
v_0^{\prime}=-v_0,\\
v_i^{\prime}=v_i-\langle v_0,v_i \rangle v_0=v_i+a_{0i}v_0, i>0\,\,.
\end{array}$$
We recall that $a_{ij}=-\langle v_i,v_j\rangle$.
The mutated matrix $(a_{ij}^{\prime})$ is given by

$$\begin{array}{l}a_{ij}^{\prime}=a_{ij}+a_{i0}a_{0j}\,\mbox{ if } i<0<j,\\
a_{i0}^{\prime}=-a_{i0},\\
a_{0i}^{\prime}=-a_{0i},\\
a_{ij}^{\prime}=a_{ij},\, \mbox{otherwise}.\end{array}$$
Thus we see that the mutation at $E_0$ gives rise to the mutation
of the matrix $(a_{ij})$ in the sense of cluster algebras (see \cite{Ze}).
Notice that at the categorical level the mutation is not an involution. The composition of the mutation at $E_0$ and of the mutation at
 $E'_0=E_0[-1]$ is the reflection functor  $R_{E_0}$ 
applied to all elements of the cluster collection.

Identifying Calabi-Yau categories endowed with cluster collections with quivers with potentials  we obtain the well-known notion of mutation of a quiver with potential (see \cite{Ze}). Then we have the following result.

\begin{thm} In the scheme (an infinite-dimensional affine space) of potentials ${\cal PT}$ there is a countable set of algebraic hypersurfaces
$X_i, i\geqslant 1$ invariant under the gauge group action, such that
for any potential belonging to the  set  ${\cal PT}\setminus \cup_{i\geqslant 1}X_i$ one can make mutations indefinitely, obtaining each time a potential
from ${\cal PT}\setminus \cup_{i\geqslant 1}X_i$. In particular, all corresponding quivers do not have oriented cycles of length one or two.

\end{thm}

{\it Sketch of the proof.} The mutated spherical collection fails to be cluster if for some $i\ne j$ we have simultaneously $\op{Ext}^1(E_i^{\prime},E_j^{\prime})\ne 0$ and $\op{Ext}^2(E_i^{\prime},E_j^{\prime})\ne 0$.
This property is not stable under deformations of 3-dimensional Calabi-Yau $\A$-category, since we can add a quadratic term to the potential
$W_{E_i^{\prime}\oplus E_j^{\prime}}$ reducing the dimension of $\op{Ext}^1(E_i^{\prime},E_j^{\prime})$ and $\op{Ext}^2(E_i^{\prime},E_j^{\prime})$. Therefore, the property that the mutated collection is also a cluster one holds on a Zariski open non-empty subset of the space of all potentials. Moreover,
the mutation induces a birational identification between varieties (maybe infinite-dimensional) of gauge equivalence classes of generic potentials
 for quivers corresponding to skew-symmetric matrices $(a_{ij})$ and $(a'_{ij})$. $\blacksquare$

Any cluster collection ${\cal E}=\{E_i\}_{i\in I}$ defines an open domain $U_{{\cal E}}\subset Stab(\CC,\op{cl})$, where $\Gamma=K_0(\CC(\kk)), \op{cl}=id$. Namely, for any collection
$z_i\in \C, \op{Im}z_i>0, i\in I$ we have a stability condition $\sigma_{(z_i)}:=\sigma_{(z_i)_{i\in I}}$ with the $t$-structure defined by $(E_i)_{i\in I}$ and the central charge $Z$ such that $Z(E_i):=Z(\op{cl}(E_i))=z_i, i\in I$. The heart of the $t$-structure is an abelian category ${\cal A}_{\cal E}$ generated by $(E_i)_{i\in I}$, which is artinian with simple objects $E_i, i\in I$.  This abelian category is equivalent to the category of continuous finite-dimensional representations of the algebra $B_W$ where $W$ is the potential of the path algebra of the quiver $Q$.
If ${\cal E}^{\prime}$ is a cluster collection obtained from ${\cal E}$ by the mutation at $E_0$ then the domains
$U_{{\cal E}}$ and $U_{{\cal E}^{\prime}}$ do not have common interior points, but have a common part of the boundary which is the wall of second kind. The common boundary corresponds to the stability structure with $Z(E_0)\in \R_{<0}$.

\vspace{15mm}
\begin{picture}(200,200)(-80,0)
\thicklines
\put(0,100){\line(1,0){200}}
\put(100,100){\circle*{3}}
\thinlines
\put(100,100){\vector(-4,1){90}}
\put(100,100){\vector(4,-1){90}}
\put(100,100){\vector(-1,4){20}}
\put(100,100){\vector(1,3){25}}
\put(100,100){\vector(0,1){80}}
\put(0,130){$Z(E_0)$}
\put(165,60){$Z(E_0[-1])$}
\put(80,190){$Z(E_i),\,\,i\ne 0 $}
\end{picture}

Category ${\cal A}_{{\cal E}'}$ is obtained from ${\cal A}_{\cal E}$ by tilting. Namely, any object $M$ of ${\cal A}_{\cal E}$ admits a unique presentation as an extension 
$$0\to nE_0\to M\to N\to 0$$
where $N\in {\cal B}:=\{E\in {\cal A}_{\cal E}\,|\,\,\op{Hom}(E,E_0)=0\}$. Similarly, any object $M'$ of ${\cal A}_{{\cal E}'}$ admits a unique presentation as an extension 
$$0\to N\to M'\to nE_0[-1]\to 0$$
with $N\in {\cal B}$.

\subsection{Orientation data for cluster collections}

Let ${\cal E}=(E_i)_{i\in I}$ be a cluster collection. We set $R:=R_{\cal E}= \op{Ext}^\bullet(E,E)$, where $E=\oplus_{i\in I}E_i$. Then $R$ is an $\A$-algebra. We denote by $M:=M_{\cal E}$ the algebra $R$ considered as  $R$-bimodule. Using the truncation functors $\tau_{\leqslant i}$ and $\tau_{\geqslant i}$ we define a sub-bimodule $M_{\geqslant 2}=\tau_{\geqslant 2}M$ as well as a quotient bimodule $M/M_{\geqslant 2}$, which is isomorphic to 
$M_{\leqslant 1}=\tau_{\leqslant 1}M$.  Then we can deform the extension $M_{\geqslant 2}\to M\to M_{\leqslant 1}$ into the direct sum of bimodules $M_{\geqslant 2}\oplus M_{\leqslant 1}$. Moreover, one can check 
 that there exists a deformation which consists of self-dual bimodules (i.e. they give rise to self-dual functors in the sense of Section 5.3). Thus we would like to define an orientation data using the splitting given by the bifunctor $F$ which corresponds to the bimodule $M_{\leqslant 1}$, i.e.
$(E_i,E_j)\mapsto \tau_{\leqslant 1} \op{Ext}^\bullet(E_i,E_j)$.
Let ${\cal E}^{\prime}=(E_i^{\prime})_{i\in I}$ be the cluster collection obtained by a mutation at $i=0$.
One can check directly that $\Z/2\Z$-valued quadratic form defined on $K_0(\CC(\kk))$ by
$$[E]=\sum_{I\in I}n_i[E_i]\mapsto \sum_{i\in I}n_i^2-\sum_{i,j\in I, a_{ij}>0}a_{ij}n_in_j \op{mod} 2\Z$$ is invariant under mutations. This means that the parity of the super line bundle $\sqrt{\cal D}_F=\op{sdet}(\tau_{\leqslant 1} \op{Ext}^\bullet(F,F))$ is preserved under mutations. This makes plausible the following conjecture.

\begin{conj} Bifunctors $M_{\leqslant 1}^{\cal E}$ and $M_{\leqslant 1}^{{\cal E}^{\prime}}$ define isomorphic orientation data on $\CC$.

\end{conj}

In order to check the conjecture one needs to find a self-dual ${\bf A}^1$-deformation of $M_{\leqslant 1}^{\cal E}\oplus  (M_{\leqslant 1}^{{\cal E}^{\prime}})^{\vee}$ to a bifunctor of the type $N\oplus N^{\vee}$ (we identify bifunctors with bimodules here).

\subsection{Quantum DT-invariants for quivers}

For any $\sigma\in U_{{\cal E}}$ (recall domain $U_{{\cal E}}$ introduced at the end of Section 8.1) we have the corresponding element $A_V^{\op{Hall}}$, where $V$ is any strict sector containing all $Z(E_i), i\in I$. The element $A_V^{\op{Hall}}$ does not depend on $\sigma$. Moreover,  this element depends only on the gauge equivalence class of the corresponding potential. The associated element $A_{V,q}:=A_{{\cal E},q}$ of the quantum torus ${\cal R}_{V,q}$ depends (for a {\it generic} potential) on the matrix $(a_{ij})$ only.

Let us associate with our quiver $Q$ the quantum torus ${\cal R}_{Q,q}$. By definition it is an associative unital algebra over the field $\Q(q^{1/2})$ of rational functions, with invertible generators $\hat{e}_i^{\pm 1},\,\,i\in I$ subject to the relations
$$\hat{e}_i\hat{e}_j=q^{a_{ji}}\hat{e}_j\hat{e}_i\,\,.$$
We are going to use its double $D({\cal R}_{Q,q})$, which is generated by ${\cal R}_{Q,q}$, new set of generators $\hat{e}_i^{\vee}, i\in I$ subject to the additional set of relations:
$$\hat{e}_i^{\vee}\hat{e}_j^{\vee}=\hat{e}_j^{\vee}\hat{e}_i^{\vee}, \,\,\hat{e}_i^{\vee}\hat{e}_j=q^{-\delta_{ij}}\hat{e}_j \hat{e}_i^{\vee}, \,\,i,j\in I\,\,.$$
The corresponding quasi-classical limits are Poisson tori which we will denote by ${\mathbb T}_Q$ and
$D({\mathbb T}_Q)$ respectively.

Identifying ${\cal R}_{\Gamma,q}$ with ${\cal R}_{Q,q}$ in the obvious way we obtain an element
$${\bf E}_Q=1+\dots\in \widehat{{\cal R}_{Q,q}}$$
corresponding to $A_{{\cal E},q}$. We observe that ${\bf E}_Q$ is a series in non-negative powers of $\hat{e}_i,i\in I$.

Conjugation with ${\bf E}_Q$ gives rise to an automorphism of the quantum torus $D({\cal R}_{Q,q})$. By the ``absence of poles" conjecture it does not have poles at $q^n=1, n\geqslant 1$. In particular it defines a formal symplectomorphism of the double torus
$D({\mathbb T}_Q)$ (see  Section 2.6, with the notation $b_{ij}:=-a_{ij}$).

\subsection{Quivers and cluster transformations}

The formal power series  ${\bf E}_Q$ in $\hat{e}_i,\,\,i\in I$ defined in the previous section satisfy a number of remarkable properties.

1) If $|I|=1$ then $Q$ is a quiver with one vertex $i$. We have
$${\bf E}_Q={\bf E}(\hat{e}_i)\,\,,$$
where ${\bf E}$ is the  quantum dilogarithm function.

2) Let $I=I_1\sqcup I_2$ , and we assume that $a_{i_1,i_2}<0$ for any $i_1\in I_1, i_2\in I_2$. Then we have two subquivers $Q_1$ and $Q_2$ of $Q$ with the sets of vertices $I_1$ and $I_2$ correspondingly, and all the arrows connecting $Q_1$ and $Q_2$ go only in the
 direction from $Q_2$ to $Q_1$ (i.e. there is no arrows from $Q_1$ to $Q_2$).

\begin{prp} One has:
$${\bf E}_Q={\bf E}_{Q_1}{\bf E}_{Q_2}$$
where we embed ${\cal R}_{Q_j,q},\,j=1,2$ into ${\cal R}_{Q,q}$ in the obvious way: $\hat{e}_i \mapsto \hat{e}_i$ for $i\in I_1$ or $i\in I_2$.

\end{prp}
{\it Proof.} Consider the stability condition $\sigma\in U_{{\cal E}}$ on the Calabi-Yau category $\CC_Q$ associated with $Q$ and a generic potential. Let ${\cal E}=\{E_i\}_{i\in I}$ be the corresponding cluster collection. We choose a stability condition $\sigma\in U_{{\cal E}}$  in such a way that $\op{Arg}(E_{i_1})>\op{Arg}(E_{i_2})$ for $i_1\in I_1, i_2\in I_2$. In this case $\CC^{ss}_Q=\CC^{ss}_{Q_1}\sqcup \CC^{ss}_{Q_2}$. This implies the desired identity. $\blacksquare$

\begin{rmk} It follows from the Properties 1) and 2) that for any acyclic quiver $Q$ the element
${\bf E}_Q$ can be expressed as the product of ${\bf E}(\hat{e}_i), i\in I$. In particular, the conjugation by ${\bf E}_Q$ has a well-defined quasi-classical limit as $q^{1/2}\to -1$, which is a birational symplectomorphism of the torus $D({\mathbb T}_Q)$.

\end{rmk}

3) Let $Q^{\prime}$ be the quiver obtained from $Q$ by the mutation at $0\in I$. We denote the standard generators of the corresponding quantum tori by $(\hat{e}_i^{\prime})_{i\in I}, \hat{e}_i^{\prime}=\hat{e}_{\op{cl}_{\kk}(E_i^{\prime})}$ and $(\hat{e}_i)_{i\in I},\hat{e}_i=\hat{e}_{\op{cl}_{\kk}(E_i)}$ respectively. Let us introduce the elements
$$R_Q={\bf E}(\hat{e}_0)^{-1}\cdot{\bf E}_Q,\,\,R_{Q^{\prime}}={\bf E}_{Q^{\prime}}\cdot{\bf E}(\hat{e}_0^{\prime})^{-1}\,\,.$$ Here
$R_Q$ is a series in variables $\hat{e}_i$ for $ i<0$, and in (dependent) variables  $\hat{e}_i,\hat{e}_i\hat{e}_0,\dots,\hat{e}_i\hat{e}_0^{a_{0i}}$ for $ j>0$.
 Similarly, 
$R_{Q^{\prime}}$ is a series in variables $\hat{e}_i^{\prime}$ for $ i<0$ and
$\hat{e}_i^{\prime},\hat{e}_i^{\prime}\hat{e}_0^{\prime}\dots,\hat{e}_i^{\prime}(\hat{e}_0^{\prime})^{ a_{0i}}$ for $i>0$.

Then $R_Q=R_{Q^{\prime}}$ under the identification 
$$\begin{array}{l}
\hat{e}_i^{\prime}=\hat{e}_i, \,\,\,i<0\,,\\
\hat{e}_0^{\prime}=\hat{e}_0^{-1}\,,\,\,\\
\hat{e}_i^{\prime}=q^{-{1\over{2}}a_{0i}^2}\hat{e}_i\hat{e}_0^{a_{0i}}, \,\,\,i>0\,\,.
\end{array}$$ This follows from the above-discussed picture of tilting via the wall-crossing, more precisely, from the formula
$${\bf E}(\hat{e}_{\op{cl}_{\kk}(E_0)})^{-1}A_{{\cal E},q}=A_{{\cal E}^{\prime},q}{\bf E}(\hat{e}_{-\op{cl}_{\kk}(E_0)})^{-1}.$$
Element $R_Q=R_{Q^{\prime}}$ corresponds to the integral over the space of objects of category ${\cal B}$ in notation at the end of Section 8.1.

For the convenience of the reader we give also the formulas comparing
dual coordinates on the double quantum torus:
$$  \begin{array}{l}
{\hat{e}^\vee_i}{}^{\prime}=\hat{e}^\vee_i,\,\,\,\forall i\ne 0 \,\,,\\
{\hat{e}^\vee_0}{}^{\prime}=(\hat{e}^\vee_0)^{-1}\cdot \prod_{i>0}\left(\hat{e}^\vee_i\right)^{a_{0i}}\,\,.
\end{array}  $$

Let us now consider the minimal class ${\cal P}$ of oriented finite quivers which satisfies the following properties:

a) the trivial quiver (one vertex no arrows) belongs to ${\cal P}$;

b) class ${\cal P}$ is closed under mutations;

c) if $Q_1,Q_2\in {\cal P}$ then
a  quiver $Q$ obtained from the disjoint union of $Q_1$ and $Q_2$ by inserting a finite number of arrows from $Q_2$ to  $Q_1$ (without changing anything else for $Q_1$ and $Q_2$) also belongs to ${\cal P}$.
We will say in this case that $Q$ is an extension of $Q_1$ by $Q_2$.
 At the level of categories this means that any object $J$
of the category ${\cal A}({\cal E})$ generated by $E_i\in {\cal E},i\in I$ is an extension $F_1\to J\to F_2$ where $F_1$ (resp. $F_2$) is an object of the abelian category generated by $E_i, i\in I_1$ (resp. $E_i,i\in I_2$).

This class ${\cal P}$ enjoys the property that the gauge group associated with $Q\in {\cal P}$ when acting on the space of potentials on $Q$ has one open orbit (this can be shown by induction), hence the corresponding 3-dimensional Calabi-Yau category is {\it rigid}.
Moreover for any $Q\in {\cal P}$ the element ${\bf E}_Q$ is a finite product of the elements
${\bf E}(f)$, where $f=\hat{e}_\gamma$ is a monomial.  In particular, the conjugation with ${\bf E}_Q$ has a quasi-classical limit as $q^{1/2}\to -1$, which is a birational transformation.

One of the first nontrivial examples of a quiver $Q$ which is not in the class ${\cal P}$ is
the quiver $Q_3$ which has three vertices and two parallel arrows between any two vertices (see the Figure). This quiver is stable under mutations. The element ${\bf E}_{Q_3}$ satisfies an overdetermined system of equations. The computer check shows that the conjugation with ${\bf E}_{Q_3}$ has the quasi-classical limit which is not rational. It is not clear whether it admits an analytic continuation.

\vspace{3mm} 

\begin{picture}(100,100)(-120,0)
\thicklines
\put(-1,11){\circle*{3}}
\put(100,11){\circle*{3}}
\put(90,13.1){\vector(-1,0){80}}
\put(90,7.9){\vector(-1,0){80}}
\put(8,20){\vector(2,3){38}}
\put(3.5,23.5){\vector(2,3){38}}
\put(59, 79.5){\vector(2,-3){38}}
\put(54.5,76){\vector(2,-3){38}}
\put(50,87){\circle*{3}}

\end{picture}

\vspace{3mm}

The mutation property of $Q_3$ has the following explicit corollary. Namely, there exist collections 
$$c_{i,j,k},\,\,b_{m_1,m_2,n}\in \Q(q^{1/2}), \,\,i,j,k\in \Z_{\geqslant 0},\,\, m_1,n\geqslant 0,\,\, -m_1\leqslant m_2\leqslant m_1$$ such that the following system of equations is satisfied:

$$c_{0,0,0}=b_{0,0,0}=1, c_{i,j,k}=c_{j,k,i}=c_{k,i,j}\,\,,$$
$$c_{n_0,n_1,n_2}=\sum_{l\geqslant 0}\varepsilon_lq^{l(n_2-n_1)}b_{n_1,n_0-l-n_1,n_2}\,\,,$$
$$c_{n_0,n_1,n_2}=\sum_{l\geqslant 0}\varepsilon_lq^{l(n_2-n_0)}b_{n_0,n_0+l-n_1,n_2}\,\,,$$
where
$$\varepsilon_l={q^{l^2/2}\over{(q^l-1)\dots(q^l-q^{l-1})}}$$
are coefficients of the series $\mathbf{E}$.
To have a solution of this system of equations is the same as to write the element
$${\bf E}_Q=\sum_{i,j,k}c_{i,j,k}\hat{e}_{(i,j,k)}\,\,,$$
where  we identified $\Gamma$ with $\Z^3$.
The above system of equations follows from the identity $R_{Q_3}=R_{Q_3^{\prime}}$ since $Q_3=Q_3^{\prime}$ after the mutation. The elements $b_{m_1,m_2,n}$ are derived from $c_{i,j,k}$.

Notice that the above system of equations has a solution which is not unique. Therefore the element ${\bf E}_Q$ is determined  non-uniquely, but only up to a multiplication by a series of the type $$1+\sum_{n\geqslant 1}a_n\,\hat{e}_{1,1,1}^n,\,\,\, a_n\in \Q(q^{1/2})$$ which belongs to the center of the quantum torus ${\cal R}_{Q_3,q}$.

Let as before ${\cal E}=(E_i)_{i\in I}$ be a cluster collection in $\CC$ such that the corresponding potential is generic.  We make an additional assumption that the conjugation $\op{Ad}_{A_{{\cal E},q}}:x\mapsto A_{{\cal E},q}xA_{{\cal E},q}^{-1}$ is a birational transformation of the double quantum torus ${\cal R}_{\Gamma\oplus \Gamma^{\vee},q}\simeq D({\cal R}_{Q,q})$.
This means that it is an automorphism of the (well-defined) skew field $K_{\Gamma\oplus \Gamma^{\vee},q}$ of fractions of this quantum torus. In the equivalent language of quivers it suffices to require that $Q\in {\cal P}$.

Let us denote by $\Phi_{\cal E}$ the automorphism of $K_{\Gamma\oplus \Gamma^{\vee},q}$ given by
$$\Phi_{\cal E}(x)=(\op{Ad}_{A_{{\cal E},q}}^{-1}\circ \tau)(x)\,\,,$$
where $\tau$ is the involution induced by the antipodal involution $\gamma\mapsto -\gamma$ of $\Gamma\oplus \Gamma^{\vee}$.

\begin{prp} If ${\cal E}^{\prime}=(E_i^{\prime})_{i\in I}$ is the cluster collection obtained by the mutation at $E_0$ then
$$\op{Ad}_{{\bf E}(\hat{e}_{\op{cl}_{\kk}(E_0)})}^{-1}\circ\Phi_{\cal E}\circ \op{Ad}_{{\bf E}(\hat{e}_{\op{cl}_{\kk}(E_0)})}=\Phi_{{\cal E}^{\prime}}\,\,.$$

\end{prp}

{\it Proof.} From the known identity
$$\op{Ad}_{{\bf E}(\hat{e}_{\op{cl}_{\kk}(E_0)})}^{-1}\circ \op{Ad}_{A_{{\cal E},q}}=\op{Ad}_{A_{{\cal E}^{\prime},q}}
\circ \op{Ad}_{{\bf E}(\hat{e}_{-\op{cl}_{\kk}(E_0)})}^{-1}$$
we obtain the desired one by multiplying it from the right by $\tau\circ \op{Ad}_{{\bf E}(\hat{e}_{\op{cl}_{\kk}(E_0)})}$. $\blacksquare$

Now we can state a similar result for a quiver $Q$ which satisfies the condition that $\op{Ad}_{{\bf E}_Q}$ is a birational transformation of the skew field $K_Q$ of fractions of the double quantum torus $D({\cal R}_{Q,q})$.
Let us define $\Phi_Q:=\op{Ad}_{{\bf E}_Q}^{-1}\circ \tau$ where $\tau$ is the obvious
 involution: $$\tau(\hat{e}_i)=\hat{e}_i^{-1},\,\,\tau(\hat{e}^\vee_i)=(\hat{e}^\vee_i)^{-1}\,\,.$$ Let $Q^{\prime}$ be the quiver obtained as a mutation of $Q$ at the vertex $0\in I$.
Then we have the following corollary of the above Proposition.
\begin{cor} Let us define the map $C_{Q,0}: K_{Q,q}\to K_{Q^{\prime},q}$ as the composition
$$K_{Q,q}\to K_{\Gamma\oplus \Gamma^{\vee},q}\to K_{\Gamma\oplus \Gamma^{\vee},q}\to K_{Q^{\prime},q}\,\,,$$ where
the middle arrow is the automorphism $\op{Ad}^{-1}_{{\bf E}(\hat{e}_{\op{cl}_{\kk}(E_0)})}$ while the other maps are obvious isomorphisms of skew fields. Then
$$C_{Q,0}\circ \Phi_Q=\Phi_{Q^{\prime}}\circ C_{Q,0}\,\,.$$

\end{cor}

{\it Proof.} It is just a reformulation of the previous Proposition in the language of quivers. $\blacksquare$

Let us compute $C_{Q,0}(\hat{e}_i)$, where $\hat{e}_i=\hat{e}_{\op{cl}_{\kk}(E_i)}, i\in I$, as well as $C_{Q,0}(\hat{e}_i^{\vee}), i\in I$.
We have to compute the action of $\op{Ad}^{-1}_{{\bf E}(\hat{e}_{\op{cl}_{\kk}(E_0)})}$ on these generators. Thus we obtain
$$\begin{array}{l}\hat{e}_0\mapsto (\hat{e}_0^{\prime})^{-1}\,\,,\\
\hat{e}_i\mapsto \hat{e}_{i}^{\prime}\cdot \prod_{0\leqslant n\leqslant a_{i0}-1}(1+q^{n+1/2}(\hat{e}_0^{\prime})^{-1})^{-1},\, i<0\,\,,\\
\hat{e}_i\mapsto \hat{e}_{i}^{\prime}\cdot \prod_{0\leqslant n\leqslant a_{i0}-1}(1+q^{n+1/2}\hat{e}_0^{\prime}), \,i>0\,\,.
\end{array}$$
Similarly we obtain that
$$\begin{array}{l}\hat{e}_i^{\vee}\mapsto {\hat{e}_i^{\vee}}{}^{\prime},\, i\ne 0\,\,,\\
\hat{e}_0^{\vee}\mapsto ({\hat{e}_0^{\vee}}{}^{\prime})^{-1}\cdot \prod_{i>0}({\hat{e}_i^{\vee}}{}^{\prime})^{a_{0i}}\cdot (1+q^{1/2}(\hat{e}_0^{\prime})^{-1})^{-1}\,\,.\end{array}$$

Under quasi-classical limit the generators $\hat{e}_i, i\in I$ go to the coordinates $y_i, i\in I$ and $\hat{e}_i^{\vee}$
go to the coordinates $x_i, i\in I$ of the symplectic double torus (see Section 2.6).
Then in those coordinates we obtain
$$\begin{array}{l}y_i\mapsto {y_i^{\prime}\over{(1-1/y_0^{\prime})^{a_{i0}}}},\,i<0\,,\\
y_0\mapsto (y_0^{\prime})^{-1}\,,\\
y_i\mapsto y_i^{\prime}(1-y_0^{\prime})^{a_{0i}},\, i>0\,\,.
\end{array}$$
For the dual coordinates we have:
$$\begin{array}{l}x_i\mapsto x_i^{\prime}, \,i\ne 0\,,\\
x_0\mapsto (x_0^{\prime})^{-1}\cdot \prod_{i>0}(x_i^{\prime})^{a_{0i}}\cdot (1-1/y_0^{\prime})^{-1}
\,\,.\end{array}$$

Up to a change of sign these are cluster transformations. Namely, if we set
$X_i=-y_i, X_i^{\prime}=-y_i^{\prime}, A_i=1/x_i, i\in I$ then our formulas become formulas (17) and (18) from \cite{FG}
(in the notation from loc. cit).

\begin{rmk} Let us recall the variety $N$ from Section 2.6 defined by the equations
$N=\{y_i=-\prod_{j\in I}x_j^{a_{ij}}, i\in I\}$,
and let $N^{\prime}$ be a similar variety defined for the transformed coordinates $x_i^{\prime}, y_i^{\prime}, i\in I$. One can check that the quasi-classical limit of $C_{Q,0}$ transforms $N$ into $N^{\prime}$. Furthermore, the quasi-classical limit of the automorphism $\Phi_Q$ preserves $N$.

\end{rmk}

\begin{rmk} 1) Let us assume  that $ \op{Ad}_{A_{{\cal E},q}}$ is birational (e.g. for $Q\in {\cal P}$).The above considerations show that the conjugacy class of the element ${\Phi}_Q$ is an invariant of the quiver $Q$ under mutations. Passing to quasi-classical limit we obtain an invariant of a quiver (under mutations) which is a conjugacy class in the group of birational transformations of the classical double torus.
2) The categorical version of the above remark holds in a greater generality. Namely, let us assume that $\CC$ has a $t$-structure generated by finitely many objects. Then we can define the motivic DT-invariant $A_{\CC}^{mot}:=A_{V}^{mot}$ (and its quantum and semi-classical relatives) for every stability condition such that all the generators of the $t$-structure are stable. Here $V$ can be any strict sector containing their central charges, so we can replace it by the upper-half plane.   Then  the conjugacy class of the automorphism $\Phi_{\CC}:= \op{Ad}_{A_\CC}^{-1}\circ \tau$ (if it makes sense) will be 
independent (under appropriate conditions) of the choice of stability condition.

\end{rmk}

\vspace{5mm}

Addresses:

M.K.: IHES, 35 route de Chartres, F-91440, France

{maxim@ihes.fr}\\

Y.S.: Department of Mathematics, KSU, Manhattan, KS 66506, USA

{soibel@math.ksu.edu}


\addcontentsline{toc}{section}{References}


\begin{thebibliography}{XX}

\bibitem {Ba} A.~Bayer, {\it Polynomial Bridgeland stability conditions and the large volume limit}, arXiv:0712.1083.


\bibitem
{B} K.~Behrend, {\it  Donaldson-Thomas invariants via microlocal geometry}, math/0507523.

\bibitem
{BDh} K.~Behrend, A.~Dhillon, {\it On the Motive of the Stack of Bundles}, math/0512640.

\bibitem
{BF} K.~Behrend, B.~Fantechi, {\it  Symmetric obstruction theories and Hilbert schemes of points on threefolds},
math/0512556.

\bibitem
{BiMi} E.~Bierstone, P.~Milman, {\it  Functoriality in resolution of singularities}, ArXiv:math/0702375.


\bibitem
{Bit} F.~Bittner, {\it The universal Euler characteristic for varieties of characteristic zero}, math. AG/0111062.


\bibitem
 {BVdB} A.~Bondal, M.~Van den Bergh, {\it Generators and representability of functors in commutative and non-commutative geometry}, math.AG/0204218.


\bibitem
 {Bor} R.~Borcherds, {\it Automorphic forms with singularities on Grassmannians}, arXiv:alg-geom/9609022.


\bibitem
{Br1} T.~Bridgeland, {\it Stability conditions on triangulated categories}, math.AG/0212237.


\bibitem
{BrTL} T.~Bridgeland, V.~Toledano-Laredo, {\it Stability conditions and Stokes factors}, arXiv:0801.3974.


\bibitem
 {CeVa} S.~Cecotti, C.~Vafa, {\it On Classification of N=2 Supersymmetric Theories},  arXiv:hep-th/9211097.


\bibitem
 {ChVer} M.~Cheng, E.~Verlinde, {\it  Wall Crossing, Discrete Attractor Flow, and Borcherds Algebra},
arXiv:0806.2337.

\bibitem
 {DM} F.~Denef, G.~Moore, {\it Split States, Entropy Enigmas, Holes and Halos},
hep-th/0702146.


\bibitem
 {DeLo1} J.~Denef, F.~Loeser, {\it Geometry on arc spaces of algebraic varieties},
arXiv:math/0006050.


\bibitem
 {DeLo2} J.~Denef, F.~Loeser, {\it Motivic exponential integrals and a motivic Thom-Sebastiani Theorem},
arXiv:math/9803048.


\bibitem
 {DeLo3} J.~Denef, F.~Loeser, {\it Germs of arcs on singular algebraic varieties and motivic integration},
arXiv:math/9803039.



\bibitem
 {DWZe} H.~Derksen, J.~Weyman, A.~Zelevinsky, {\it Quivers with potentials and their representations I: Mutations},
arXiv:0704.0649.



\bibitem
 {DiM} E.~Diaconescu, G.~Moore, {\it Crossing the Wall: Branes vs. Bundles}, arXiv:0706.3193.

\bibitem
 {DoT} S.~Donaldson, R.~Thomas, {\it Gauge theory in higher dimensions},  ``The geometric universe:
science, geometry and the work of Roger Penrose", Oxford University Press, 1998.

\bibitem
 {DonM} R.~Donagi, E.~Markman, {\it Cubics, integrable systems, and Calabi-Yau threefolds}, arXiv:alg-geom/9408004.


\bibitem
{Dou} M.~Douglas, {\it Dirichlet branes, homological mirror symmetry, and stability}, arXiv:math/0207021.


\bibitem{FK} L.~Faddeev, R.~Kashaev, {\it Quantum Dilogarithm}, Mod.Phys.Lett. A9 (1994) 427-434, arXiv hep-th/9310070.

\bibitem
 {FG} V.~Fock, A.~Goncharov, {\it The quantum dilogarithm and representations quantum cluster varieties } , arXiv:math/0702397.



\bibitem
 {GaMNeit} D.~Gaiotto, G.~Moore, A.~Neitzke, {\it Four-dimensional wall-crossing via three-dimensional field theory},
arXiv:0807.4723.


\bibitem
{Gi} V.~Ginzburg, {\it Calabi-Yau algebras}, arXiv:math/0612139.


\bibitem
  {GSV} M.~Gekhtman, M.~Shapiro, A.~Vainshtein, {\it Cluster algebras and Poisson geometry},
math.QA/0208033.


\bibitem
 {GS1} M.~Gross, B.~Siebert, {\it Mirror Symmetry via Logarithmic Degeneration Data I},
arXiv:math/0309070.


\bibitem
 {GS2} M.~Gross, B.~Siebert, {\it Mirror Symmetry via Logarithmic Degeneration Data II},
 arXiv:0709.2290.


\bibitem
 {HM} J.~Harvey, G.~Moore, {\it On the algebra of BPS states}, hep-th/9609017.


\bibitem
 {HIVa} K.~Hori, A.~Iqbal, C.~Vafa, {\it D-branes and mirror symmetry}, arXiv:hep-th/0005247.


\bibitem
 {In} M.~Inaba, {\it  Moduli of stable objects in a triangulated category}, arXiv:math/0612078.



\bibitem
 {Jo1} D.~Joyce, {\it Holomorphic generating functions for invariants counting
coherent sheaves on Calabi-Yau $3$-folds}, hep-th/0607039.



\bibitem
 {Jo2} D.~Joyce, {\it Configurations in abelian categories. III. Stability conditions and identities},
math.AG/0410267.



\bibitem
 {Jo3} D.~Joyce, {\it Configurations in abelian categories. IV. Invariants and changing stability
conditions}, math.AG/0410268.



\bibitem
 {Jo4} D.~Joyce, {\it Motivic invariants of Artin stacks and ``stack functions"}, Arxiv:math/0509722.

\bibitem
 {KKP} L.~Katzarkov, M.~Kontsevich, T.~Pantev, {\it Hodge theoretic aspects of mirror symmetry}, arXiv:0806.0107.

\bibitem
 {Kel} B.~Keller, {\it A-infinity algebras, modules and functor categories}, arXiv:math/0510508.


\bibitem
 {Ko1} M.~Kontsevich, {\it Formal non-commutative symplectic geometry}, in: Gelfand Mathematical Seminars, 1990-1992, Birkhauser 1993, 173-187.


\bibitem
 {Ko2} M.~Kontsevich,  {\it XI Solomon Lefschetz Memorial Lecture Series: Hodge structures in non-commutative geometry}, arXiv:0801.4760.


\bibitem
 {KoSo1} M.~Kontsevich, Y.~Soibelman, {\it Affine structures and non-archimedean analytic spaces}, math.AG/0406564.


\bibitem
 {KoSo2} M.~Kontsevich, Y.~Soibelman, {\it Deformation theory} (book in preparation, preliminary draft is available at
{www.math.ksu.edu/$\sim$soibel}).


\bibitem
 {KoSo3} M.~Kontsevich, Y.~Soibelman, {\it Notes on A-infinity algebras, A-infinity categories and non-commutative geometry. I}, math.RA/0606241.





\bibitem
 {Lau} G.~Laumon, {\it Transformation de Fourier, constantes d'equations fonctionelle et conjecture de Weil}. Publ. Math. IHES {\bf 65}, 1987,131-210.



\bibitem
 {Laz} C.~Lazaroiu, {\it  Generating the superpotential on a D-brane category: I}, arXiv:hep-th/0610120.


\bibitem
{Loo} E.~Looijenga,
{\it Motivic measures.} S\'eminaire Bourbaki, Asterisque {\bf 42} (1999-2000), Exposé No. 874, 31 p.

\bibitem
 {LyuMan}
 V.~Lyubashenko, O.~Manzyuk, {\it Unital ${A}_{\infty}$-categories}, arXiv:0802.2885.



\bibitem
 {MaNOP} D.~Maulik, N.~Nekrasov, A.~Okounkov, R.~Pandharipande, {\it Gromov-Witten theory and Donaldson-Thomas theory, I}, arXiv:math/0312059.


\bibitem {MoRe1} S.~Mozgovoy, M.~Reineke, {\it On the number of stable quiver representations over finite fields},
arXiv:0708.1259.

\bibitem
 {MoRe2} S.~Mozgovoy, M.~Reineke, {\it  On the noncommutative Donaldson-Thomas invariants arising from brane tilings},arXiv:0809.0117.

\bibitem {NagNak} K.~Nagao, J.~Nakajima, {\it Counting invariant of perverse coherent sheaves and its wall-crossing}, arXiv:0809.2992. 

\bibitem {Nag} K.~Nagao, {\it Derived categories of small toric Calabi-Yau 3-folds and counting invariants}, arXiv:0809.2994.


\bibitem
 {NaYo} 
H.~Nakajima, K.~Yoshioka, {\it Perverse coherent sheaves on blow-up, II. Wall-crossing and Betti numbers formula},
arXiv:0806.0463.

\bibitem  {MikNekSet} A.~Mikhailov, N.~Nekrasov, S.~Sethi,
{\it Geometric realization of BPS states in N=2 theories},
arXiv:hep-th/9803142


\bibitem {NicS} J.~Nicaise, J.~Sebag, {\it Motivic Serre invariants, ramification, and the analytic Milnor fiber}, arXiv:math/0703217.

\bibitem
 {O} So Okada, {\it Topologies on a triangulated category}, math.AG/0701507.


\bibitem
{PR1}  R.~Pandharipande, R.~P.~Thomas,  {\it Curve counting via stable pairs in the derived category},
arXiv:0707.2348.


\bibitem
 {PR2} R.~Pandharipande, R.~P.~Thomas, {\it The 3-fold vertex via stable pairs}, arXiv:0709.3823.


\bibitem
 {Re1} M.~Reineke, {\it Poisson automorphisms and quiver moduli}, arXiv:0804.3214.


\bibitem
 {Re2} M.~Reineke, {\it Cohomology of non-commutative Hilbert schemes},
arXiv:math/0306185.



\bibitem
 {Re3} M.~Reineke, {\it The Harder-Narasimhan system in quantum groups and cohomology of quiver moduli}, arXiv:math/0204059.

\bibitem
 {Sa} M.~Saito, {\it  Mixed Hodge Modules}, Publ. RIMS, Kyoto Univ.
{\bf 26}, 1990, 221-333.


\bibitem
 {Seg} E.~Segal,  {\it The A-infinity Deformation Theory of a Point and the Derived Categories of Local Calabi-Yaus},
arXiv:math/0702539.

\bibitem
 {SeibW} N.~Seiberg, E.~Witten, {\it Electric-magnetic duality, monopole condensation, and confinement in $N=2$ supersymmetric Yang-Mills theory}, hep-th/9407087.




\bibitem
 {SGA} SGA 7, vol.2, Lecture Notes in Mathematics, 340, 1973.



\bibitem
 {So1} Y.~Soibelman, {\it  Mirror symmetry and non-commutative geometry}, J. Math. Phys., 45:10, 2004.





\bibitem
 {Sz} B.~Szendr\"oi, {\it Non-commutative Donaldson-Thomas theory and the conifold},
arXiv:0705.3419.


\bibitem
 {T1} R.~Thomas, {\it Moment maps, monodromy and mirror manifolds},
arXiv:math/0104196.


\bibitem
 {T2} R.~Thomas, {\it A holomorphic Casson invariant for Calabi-Yau 3-folds, and bundles on K3 fibrations},
 arXiv:math/9806111.


\bibitem
 {T3} R.~Thomas, {\it  Stability conditions and the braid group}  arXiv:math/0212214.



\bibitem
 {To} B.~To\"en, {\it Derived Hall Algebras}, arXiv:math/0501343.

\bibitem {ToVa} B.~To\"en, M.~Vaquie, {\it Moduli of objects in dg-categories}, arXiv:math/0503269.

\bibitem
 {ToVe} B.~To\"en, G.~Vezzosi, {\it Homotopical Algebraic Geometry II: geometric stacks and applications}, arXiv:math/0404373.

\bibitem
 {Tod1} Y.~Toda, {\it Limit stable objects on Calabi-Yau 3-folds}, arXiv:0803.2356.


\bibitem
 {Tod2} Y.~Toda, {\it Moduli stacks and invariants of semistable objects on K3 surfaces}, arXiv:math/0703590.



\bibitem
 {Tod3} Y.~Toda, {\it  Birational Calabi-Yau 3-folds and BPS state counting}, arXiv:0707.1643.


\bibitem
 {XiXu} Jie Xiao, Fan Xu, {\it  Hall algebras associated to triangulated categories},
arXiv:math/0608144.

\bibitem
 {Ze} A.~Zelevinsky, {\it Mutations for quivers with potentials: Oberwolfach talk, April 2007}, arXiv:0706.0822.

\bibitem
 {Zo} A.~Zorich, {\it Flat Surfaces}, arXiv:math/0609392.

\end{thebibliography}
\end{document}